\documentclass[preprint,11pt,number,sort]{elsarticle}
\usepackage[a4paper,margin=2.5cm]{geometry}
\usepackage{graphicx}
\usepackage{subcaption}
\usepackage{float}
\usepackage{bm}
\usepackage{
    amsfonts,
    amsmath,
    amssymb,
    amsthm,
    dsfont,
    tabularx,
    makecell,
    multirow
}

\usepackage{
    booktabs,
    siunitx
}

\usepackage{xcolor}
\usepackage{hyperref}

\hypersetup{
    colorlinks=true,
    linkcolor=blue,
    citecolor=blue,
    urlcolor=blue
}

\usepackage{orcidlink}

\newcolumntype{Y}{>{\centering\arraybackslash}X}

\DeclareSIUnit{\corehour}{core-hour}

\newcommand{\verti}[1]{\ensuremath{\left\lvert #1 \right\rvert}}
\newcommand{\vertii}[1]{\ensuremath{\left\lVert #1 \right\rVert}}
\newcommand{\indicator}[1]{\ensuremath{\mathds{1}_{#1}}}
\newcommand{\R}{\ensuremath{\mathbb{R}}}

\newcommand{\dy}{\ensuremath{\textup{d}y}}

\begin{document}

\begin{frontmatter}

\title{ADEx-FNO: A Unified Ambient-Domain Framework for Fourier Neural
Operators on Varying Geometries}

\author[kaust-cemse]{Roberto Nuca\orcidlink{0000-0002-9031-5668}\corref{cor1}}
\ead{roberto.nuca@kaust.edu.sa}

\author[kaust-cemse]{Giovanni Testa\orcidlink{0009-0008-5906-7433}}
\ead{giovanni.testa@kaust.edu.sa}

\author[kaust-cemse,polimi-daer]{Luca Galimberti\orcidlink{0009-0001-1145-0539}}
\ead{luca.galimberti@kaust.edu.sa}

\author[kaust-cemse,kaust-pse]{Matteo Parsani\orcidlink{0000-0001-7300-1280}}
\ead{matteo.parsani@kaust.edu.sa}

\cortext[cor1]{Corresponding author}

\address[kaust-cemse]{
Computer, Electrical and Mathematical Sciences and Engineering Division,
King Abdullah University of Science and Technology,
Thuwal 23955-6900, Saudi Arabia
}

\address[polimi-daer]{
Dipartimento di Scienze e Tecnologie Aerospaziali,
Politecnico di Milano,
Milan, Italy
}

\address[kaust-pse]{
Physical Science and Engineering Division,
King Abdullah University of Science and Technology,
Thuwal 23955-6900, Saudi Arabia
}

\begin{abstract}
Fourier neural operators (FNOs) combine nonlocal spectral learning with an efficient operator architecture, but their application to varying geometries and independently chosen numerical discretizations remains challenging. We introduce the ambient-domain extension Fourier neural operator (ADEx-FNO), a unified and deterministic framework for incorporating geometry into an FNO pipeline without replacing its defining Fourier-operator layers. Each admissible physical domain is embedded in a fixed ambient hypercube and represented by its signed distance function. The problem data and solution fields are extended deterministically to the ambient domain, transferred to a common, potentially nonuniform, rectilinear latent grid, and processed by the underlying Fourier operator. After inference, the predicted field is interpolated to an independently chosen target discretization and restricted to the physical domain. All geometry-transfer operations remain outside the optimization loop and introduce no trainable graph, point-cloud, deformation, or geometry-decoding modules.

ADEx-FNO attains relative $\ell^2$ errors of $0.32\%$--$0.77\%$ on held-out smooth-domain tests for nonlinear Poisson and advection-reaction-diffusion problems in 2D and 3D, and is further assessed on unseen nonsmooth geometries. We then use one ADEx-FNO inference to initialize conventional computational fluid dynamics solvers. Across all $29$ converged 2D and 3D RANS cases, the pseudo-time iterations decrease in every pair, with mean reductions of $44.17\%$ and $43.03\%$, respectively; comparable gains persist across three mesh resolutions. The URANS cases reduce post-window physical-time advances by $18.52\%$--$27.51\%$. Finally, in a stringent transfer from two-dimensional URANS training data to a DNS at different Mach and Reynolds numbers, the bootstrap interval is reduced by $23.47\%$--$48.21\%$, depending on the statistic required. In every CFD test, ADEx-FNO supplies only the initial field; the governing-equation solver retains complete control of the subsequent solution process, and physical or statistical consistency is assessed separately from the reduction in numerical work.
\end{abstract}

\begin{keyword}
Fourier Neural operator (FNO)
\sep Ambient-domain extension Fourier neural operator (ADEx-FNO)
\sep Functional learning
\sep Scientific machine learning (SciML)
\sep Computational fluid dynamics (CFD)
\sep Reynolds-averaged Navier--Stokes equations (RANS)
\sep Unsteady Reynolds-averaged Navier--Stokes equations (URANS)
\sep Direct numerical simulation (DNS)
\sep 
\end{keyword}

\end{frontmatter}

\section{Introduction}
\label{sec:introduction}

Partial differential equations (PDEs) provide the mathematical foundation for
a broad range of models in science and engineering. Reliable numerical methods
can approximate their solutions with high fidelity, but repeated simulations
remain expensive in many-query settings such as design optimization,
uncertainty quantification, inverse problems, control, and digital-twin
construction. Scientific machine learning offers complementary strategies for
reducing this cost. Physics-informed methods incorporate governing equations
into the learning process, whereas operator-learning methods seek a reusable
map from functional problem data to the corresponding solution field
\cite{RaissiEtAl2019,LuEtAl2021DeepONet,KovachkiEtAl2023,
LiEtAl2024PINO}. Once trained, such a map can amortize the cost of generating
high-fidelity data over many subsequent evaluations.

Among operator-learning architectures, DeepONet
\cite{LuEtAl2021DeepONet} and the Fourier Neural Operator (FNO)
\cite{LiEtAl2021FNO} have become widely used representatives. The FNO
parameterizes a nonlocal integral operator through Fourier multipliers and
evaluates the resulting spectral convolutions efficiently with Fourier
transforms. Its standard implementation is therefore particularly effective
when the input and output fields are represented on structured rectangular
grids. Many practical PDE problems, however, are solved on unstructured meshes
whose geometry, connectivity, point count, and local resolution vary among
samples. Geometric variability then changes not only the discrete data but
also the physical domains on which the input and output function spaces are
defined. Extending Fourier operator learning to this setting requires a
principled interface between varying physical domains, irregular native
discretizations, and a common spectral representation.

Several approaches address this interface through architectures that operate
directly on irregular data. Graph neural operators use learned integral
kernels on point sets \cite{LiEtAl2020MGNO}, while GNOT and Transolver employ
attention-based representations for irregular meshes and query-dependent
fields \cite{HaoEtAl2023GNOT,WuEtAl2024Transolver}. CORAL represents PDE
states through coordinate-based neural fields
\cite{SerranoEtAl2023CORAL}, and point-cloud approaches such as
Geom-DeepONet, the Point Cloud Neural Operator, and the Geometry-Informed
Neural Operator Transformer encode geometry from sampled coordinates or
surface point clouds
\cite{HeEtAl2024GeomDeepONet,ZengEtAl2025PCNO,LiuEtAl2026GINOT,LiEtAl2023GINO}.
These methods provide flexible mechanisms for processing irregular domains and
variable query locations, with geometry incorporated through learned graph,
attention, point-cloud, or neural-field representations.

A second family of methods constructs a common computational representation
before applying the operator model. Geo-FNO learns a deformation from the
physical geometry to a regular latent domain
\cite{LiEtAl2023GeoFNO}, whereas NUNO partitions nonuniform data and
interpolates the resulting subsets onto uniform subgrids
\cite{LiuEtAl2023NUNO}. \cite{LiEtAl2023GINO} combines graph neural operators with an FNO,
using learned input and output operators to transfer information between
irregular point clouds and a regular latent grid
\cite{LiEtAl2023GINO}. The Domain Agnostic Fourier Neural Operator of
Liu et al.\ incorporates a smoothed characteristic function directly into the
Fourier integral layers~\cite{LiuJafarzadehYu2023DAFNO}. Reference-domain
approaches such as DIMON, the Diffeomorphism Neural Operator, and the
Harmonic-Mapping Operator instead transport geometry, problem data, and
solutions to a template domain through diffeomorphic or harmonic maps
\cite{YinEtAl2024DIMON,ZhaoEtAl2025DNO,YanEtAl2026HMO}. These methods
provide important routes to geometry-dependent operator learning, but differ
in whether the physical-to-latent interface is learned, mask-based,
decomposition-based, or constructed through a geometry-dependent coordinate
map.

We introduce the Ambient-Domain Extension Fourier Neural Operator
(ADEx-FNO), which follows a different construction. Its novelty is not a claim
that the learned operator is independent of geometry: geometry enters
explicitly as a functional input. Instead, ADEx-FNO provides a unified,
deterministic way to incorporate varying geometries into an FNO pipeline while
preserving the original architecture. Each admissible physical
domain $\Omega\subset\mathbb{R}^d$ is embedded in a fixed ambient hypercube,
$\mathrm{Q}$, and represented by its signed distance function. The problem data
and solution fields are extended from $\Omega$ to $\mathrm{Q}$ through a
prescribed procedure selected according to the regularity of the boundary.
This construction defines an auxiliary problem on common function
spaces over $\mathrm{Q}$ whose restriction to $\Omega$ recovers the original
physical solution. The extended fields are then transferred deterministically
to a common rectilinear latent grid, where the Fourier operator is applied.

The geometry-transfer operations remain outside the optimization loop.
ADEx-FNO introduces no trainable graph, point-cloud, deformation, or
geometry-decoding modules beyond the underlying Fourier operator, and it does
not require a one-to-one map from each physical domain to a reference
geometry. For external-flow applications, the latent grid can be nonuniform,
allowing degrees of freedom to be concentrated near solid boundaries and other
regions of large spatial variation while retaining the spectral structure of
the FNO. During inference, the predicted extended field is interpolated to
independently chosen target points and restricted to the physical domain. The
same learned representation can therefore interface with data generated on
different unstructured meshes and provide initial or surrogate fields on
target discretizations that need not coincide with those used during training.

We first evaluate this construction on nonlinear Poisson and
advection-reaction-diffusion problems in two and three dimensions (2D and 3D). The tests
assess accuracy on held-out smooth domains, transfer to unseen nonsmooth
geometries, and preservation of the same data-generation and learning pipeline
across equations and spatial dimensions. We then examine a distinct practical
use of ADEx-FNO: one-shot initialization of established computational fluid
dynamics (CFD) solvers. Reynolds-averaged Navier--Stokes (RANS) and unsteady
Reynolds-averaged Navier--Stokes (URANS) models are used to assess the method in
engineering workflows involving routine analysis, design, optimization,
database generation, and repeated changes in geometry or operating conditions.
Large-eddy simulation (LES), wall-modeled LES, and hybrid RANS-LES increasingly
complement these approaches when important turbulent structures must be
resolved, whereas direct numerical simulation (DNS) remains primarily a
research and reference tool because of its substantially greater resolution
and sampling requirements
\cite{CFDVision2030Update2024,AltmannEtAl2023,GoinisEtAl2026}.

The validation hierarchy comprises 2D and 3D RANS cases,
transfer across CFD mesh resolutions, 2D and 3D URANS cases,
and a stringent transfer from 2D RANS training data to a DNS at substantially different Mach and Reynolds numbers.
The RANS studies examine reductions in nonlinear and cumulative linear-solver
work while verifying consistency of the converged aerodynamic coefficients.
The URANS studies distinguish reduction of the initial transient from agreement
of the resulting time-averaged quantities. The DNS study asks whether the
inferred large-scale field shortens the flow-development interval required
before statistics can be accumulated, while the 3D
instabilities and turbulent scales are still generated entirely by the DNS
solver. This final test also provides evidence relevant to future LES and
hybrid RANS-LES initialization, without claiming that those applications are
demonstrated here. In every CFD case,
ADEx-FNO is evaluated only once to supply the initial field; the
governing-equation solver then retains complete control of the numerical
solution and its convergence or statistical development.

\section{Mathematical setup and model definition}
\label{sec:architecture}

\subsection{Mathematical setting of the Fourier NO}
Let $\Omega\subset\R^d$ be a bounded open set, and let
\begin{equation*}
\mathcal{A}=\mathcal{A}(\Omega;\R^{d_a})
\qquad\text{and}\qquad
\mathcal{U}=\mathcal{U}(\Omega;\R^{d_u})
\end{equation*}
denote Banach spaces of functions defined on $\Omega$ and taking values in $\R^{d_a}$ and $\R^{d_u}$, respectively. Furthermore, let
\begin{equation*}
\mathcal{G}^\dagger:\mathcal{A}\to\mathcal{U}
\end{equation*}
be an arbitrary nonlinear operator. In the present setting, $\mathcal{G}^\dagger$ represents the solution operator associated with a parametric partial differential equation, mapping the problem data to the corresponding solution field. NOs constitute a recent class of machine-learning models specifically designed to approximate mappings between function spaces.
Given a collection of observations 
$$\left\{(a_j,u_j)\in\mathcal{A}\times\mathcal{U}\right\}_{j=1}^n,$$
the objective is to construct an approximation $\mathcal{G}^\theta$ of $\mathcal{G}^\dagger$ by optimizing the model parameters $\theta\in\Theta$, for some finite-dimensional parameter space $\Theta$, with respect to a prescribed loss functional. Unlike conventional neural networks, which typically learn mappings between finite-dimensional vectors, NOs aim to approximate operators between function spaces. A key property of these models is discretization invariance. Suppose that the training data consist of input-output function pairs $(a_j,u_j)$ observed at a finite set of sampling locations $x_i\in\Omega_j=\{x_1,\dots,x_{n_j}\}\subset\Omega$. Once trained, the learned operator can be evaluated at arbitrary locations $x\in\Omega$, potentially including sampling locations that were not present in the original discretization datasets.

The NO architecture introduced in \cite{KovachkiEtAl2023} is formulated in terms of a sequence of latent functions
\begin{equation*}
v_0,\ldots,v_L,
\qquad
v_\ell:\Omega\to\R^{d_v}.
\end{equation*}
The input field is first lifted to a latent feature space through the pointwise transformation
\begin{equation*}
v_0(x)=P\bigl(a(x)\bigr),
\end{equation*}
where $P:\R^{d_a}\to\R^{d_v}$ is a learnable map, typically parameterized by a shallow neural network. The latent representation is subsequently updated through a sequence of operator layers according to
\begin{equation*}
v_{t+1}(x)=\sigma\left(Wv_t(x)+(\mathcal{K}(a;\phi)v_t)(x)\right),\quad\forall x\in\Omega,
\end{equation*}
where $\sigma:\R\to\R$ denotes a component-wise nonlinear activation function, $W:\R^{d_v}\to\R^{d_v}$ is a learnable linear transformation, and $\mathcal{K}:\mathcal{A}\times\Theta_{\mathcal{K}}\to\mathcal{L}(\mathcal{U}(\Omega;\R^{d_v}),\mathcal{U}(\Omega;\R^{d_v}))$ is a parametrized integral operator acting on the latent representation. Specifically,
\begin{equation}
\label{eq:kernel}
(\mathcal{K}(a;\phi)v_t)(x):=\int_\Omega k_{\phi}(x,y,a(x),a(y)) \, v_t(y) \,\dy,\quad\forall x\in \Omega.
\end{equation}
After the final operator layer, a projection map $Q:\R^{d_v}\to\R^{d_u}$ transforms pointwise the latent representation back into the physical space, yielding the predicted solution field $u(x)=Q\bigl(v_L(x)\bigr)$.

Among the NO architectures proposed in the literature, the Fourier NO (FNO) \cite{LiEtAl2021FNO} has received considerable attention because of its computational efficiency and strong performance across a broad range of scientific-computing applications, see for example \cite{LiEtAl2023GeoFNO, FNO-reaction, FNO-enhanced}.

A central idea of the FNO is to restrict the generic integral kernel in \eqref{eq:kernel} to a translation-invariant form by setting
$k_{\phi}\bigl(x,y,a(x),a(y)\bigr)=k_{\phi}(x-y)$. Under this assumption, the integral operator becomes a convolution operator and can therefore be evaluated in Fourier space using the convolution theorem:
\begin{equation}
\bigl(\mathcal{K}_{\phi}v_t\bigr)(x)=
\mathcal{F}^{-1}
\left(
\mathcal{F}(k_{\phi})\cdot
\mathcal{F}(v_t)
\right)(x),\quad\forall x\in\Omega.
\label{eq:fno_convolution}
\end{equation}
Here, the product in Fourier space is understood modewise. For vector-valued latent functions, it corresponds to multiplication by a matrix-valued Fourier multiplier.

This spectral representation constitutes the principal computational building block of the FNO architecture. At each Fourier layer, the latent representation is first transformed to the frequency domain using the Fast Fourier Transform (FFT). For each retained Fourier mode, the corresponding coefficients are multiplied by a learnable complex-valued matrix acting across the latent channels, while modes outside the prescribed truncated set are discarded. The inverse FFT then maps the transformed representation back to the spatial domain. By evaluating the nonlocal component of the operator in Fourier space, the FNO can efficiently represent long-range interactions and global dependencies. For a fixed latent dimension and a prescribed number of retained modes, the dominant computational cost associated with the Fourier transforms scales as $\mathcal{O}(N\log N)$, where $N$ denotes the number of spatial degrees of freedom.

The FNO has demonstrated strong performance in approximating solution operators for a broad class of parametric partial differential equations, with applications including fluid dynamics, porous-media flow, elasticity, and wave propagation. Such operator surrogates are particularly valuable in settings that require repeated evaluations of computationally expensive numerical models. 

\paragraph{Nonuniform Fourier NO}
Despite its excellent approximation capabilities, the standard Fourier NO presents a significant limitation in the context of CFD applications. The latent structured grid on which the FFT operates is assumed to be uniformly spaced along each coordinate direction. While this assumption is well suited for problems posed on regular Cartesian domains, it becomes highly inefficient when applied to external aerodynamic simulations. In such problems, the computational domain typically extends far beyond the body of interest in order to accurately capture wakes, shocks, and other flow structures while avoiding spurious boundary effects. Consequently, the underlying unstructured CFD mesh is strongly nonuniform, featuring a high concentration of cells in the vicinity of the body and a progressively coarser discretization in the far-field region.

Constructing a uniformly spaced latent grid over the entire computational domain would therefore allocate a substantial fraction of its degrees of freedom to regions where the solution varies smoothly, while providing insufficient resolution in the vicinity of the body, where the relevant physical phenomena occur. Since the computational cost of the Fourier NO scales directly with the number of latent degrees of freedom, such an approach would result in an inefficient use of both memory and computational resources. For this reason, we employ a rectilinear grid whose spacing is allowed to vary independently along each coordinate direction. The resulting latent discretization concentrates grid points where the solution exhibits the largest spatial gradients, while maintaining a considerably coarser resolution in regions of limited physical interest.

This modification, however, prevents the direct application of the classical FNO architecture. Indeed, the FFT assumes that the underlying samples are uniformly spaced, and its mathematical derivation is no longer valid on nonuniform grids. To accommodate the proposed latent discretization, we therefore replace the FFT with a Nonuniform Discrete Fourier Transform (NDFT), following the framework introduced in \cite{NUFNO}. 

\begin{figure}[h]
    \centering
    \includegraphics[width=0.85\linewidth]{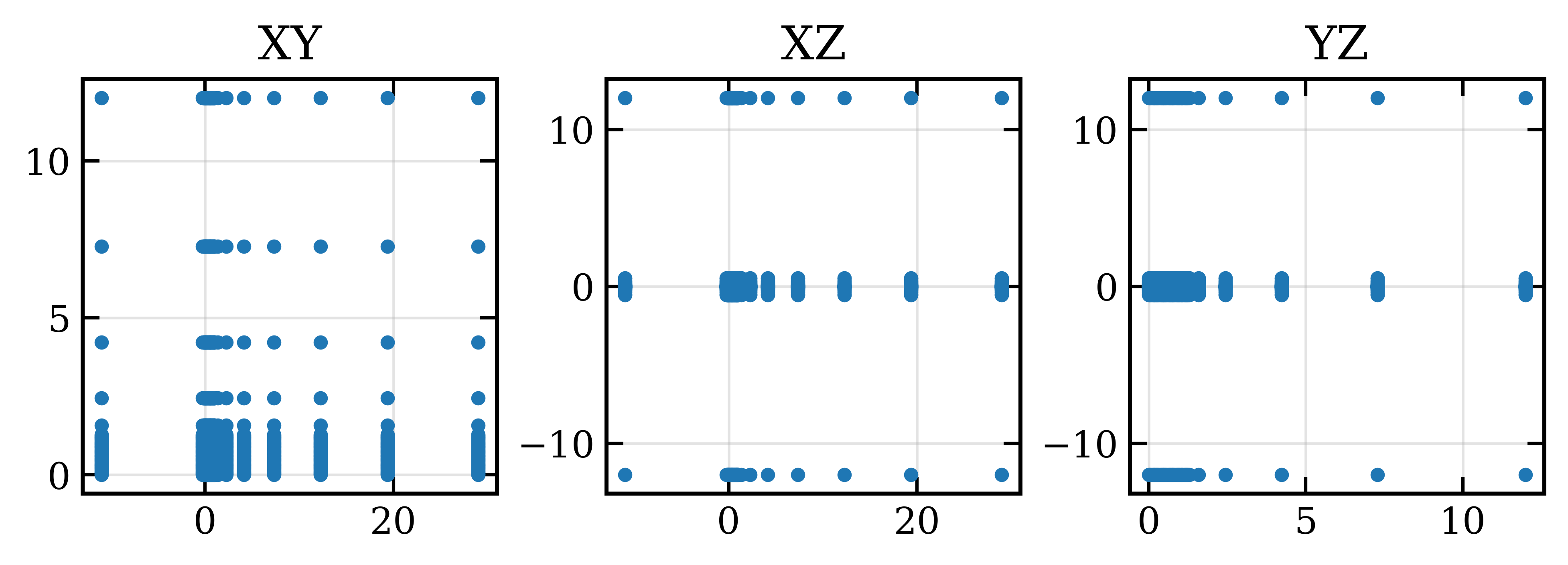}
    \caption{Nonuniform grids along coordinate axes used for our 3D model.}
    \label{fig:nonuniform-grid-sample}
\end{figure}

The spectral convolution layer is then implemented exactly as in the classical FNO: the latent representation is transformed to the frequency domain, the learnable weights act on a truncated set of Fourier modes, and the updated representation is finally mapped back to the physical domain by means of the inverse transform. Consequently, the overall architecture remains unchanged, with the only modification consisting of replacing the FFT-based spectral transform with its nonuniform counterpart NDFT. This allows the NO to exploit the computational advantages of spectral learning while efficiently operating on a latent grid whose resolution is tailored to the underlying physics of the problem. Figure~\ref{fig:nonuniform-grid-sample} shows an example of the nonuniform grid spacing used along each of the three coordinate directions for a representative 3D problem.

\subsection{Geometry-Informed NO}
The Geometry-Informed NO \cite{LiEtAl2023GINO} was introduced to extend the applicability of Fourier-based NOs to partial differential equations posed on complex geometries. While the classical FNO assumes that the input and output fields are defined on structured Cartesian grids, Li \textit{et al.} \cite{LiEtAl2023GINO} addresses the more general setting in which the computational domain is represented by an irregular point cloud or an unstructured mesh. To bridge this discrepancy, the architecture augments the standard FNO with two learnable graph-based operators. An input graph operator transfers information from the irregular physical discretization to a structured latent grid, where operator learning is performed by the Fourier NO, while an output graph operator maps the latent representation back to arbitrary query points in the physical domain. The geometry is encoded through a combination of the signed distance function and a point-cloud representation of the boundary, enabling the model to handle complex 3D domains while retaining the computational advantages of spectral convolutions. Owing to this hybrid graph-spectral formulation, \cite{LiEtAl2023GINO} has demonstrated remarkable performance in learning solution operators for large-scale PDEs on irregular geometries, substantially broadening the range of applications accessible to Fourier-based NOs.

Unlike \cite{LiEtAl2023GINO}, the proposed framework does not introduce any additional trainable components outside the Fourier NO itself. The geometry encoding and decoding stages are entirely deterministic and are performed offline through interpolation and evaluation of the signed distance function. Consequently, the number of learnable parameters, the optimization problem, and the computational cost of training remain unchanged relative to the underlying FNO architecture.

The proposed framework shares \cite{LiEtAl2023GINO}'s objective of extending Fourier-based NOs to problems posed on irregular computational domains. The two approaches, however, pursue this goal through fundamentally different strategies. \cite{LiEtAl2023GINO} introduces two learnable graph-based operators that bridge the irregular physical discretization and the structured latent grid required by the FNO. An input graph operator transfers information from the computational mesh to the latent representation, while a symmetric output graph operator maps the latent prediction back to arbitrary query locations. These graph operators are optimized jointly with the FNO, allowing the entire architecture to learn both the geometric representation and the solution operator in an end-to-end fashion.

The philosophy adopted in the present work is instead entirely deterministic. Rather than learning the transfer between the physical and latent discretizations, we explicitly construct it through a sequence of preprocessing and postprocessing operations based on the signed distance function and high-fidelity interpolation. Consequently, the geometry encoding and decoding stages are performed offline and do not introduce any additional trainable parameters beyond those of the underlying Fourier NO. The optimization problem therefore remains identical to that of a standard FNO, while the geometry-processing pipeline is completely decoupled from the learning stage.

A second important distinction concerns the construction of the latent computational grid. \cite{LiEtAl2023GINO} relies on the standard Fourier NO and therefore inherits the requirement for a uniformly sampled latent discretization. In contrast, the proposed framework employs a nonuniform rectilinear grid, allowing the latent degrees of freedom to be concentrated in regions where the solution exhibits the largest spatial variations while maintaining a coarser resolution elsewhere. This design is particularly advantageous for CFD applications, where the computational domain typically contains extensive far-field regions with relatively smooth solution fields and highly localized flow features in the vicinity of solid boundaries. By replacing the classical FFT with the corresponding nonuniform spectral transform, the proposed architecture preserves the Fourier-based formulation while enabling a substantially more efficient allocation of the latent degrees of freedom.

The two methodologies should therefore be regarded as complementary rather than competing. \cite{LiEtAl2023GINO} provides a highly flexible, end-to-end trainable framework capable of learning the transfer between arbitrary discretizations and the latent Fourier representation. The present approach, on the other hand, leverages deterministic geometric preprocessing to preserve the simplicity of the original FNO optimization problem while exploiting prior knowledge of the computational domain. Since all geometry-related operations are performed offline, the training cost and the number of learnable parameters remain unchanged with respect to the underlying FNO architecture. Moreover, the mesh-independent nature of the resulting NO enables efficient transfer of predictions to arbitrary computational meshes, including discretizations substantially finer than those employed during training. This opens the possibility of generating relatively inexpensive datasets on coarse meshes while deploying the trained model on high-resolution simulations, thereby reducing the cost of data generation without sacrificing the accuracy required by practical CFD applications.

\subsection{Ambient-Domain Extension FNO (ADEx-FNO)}
In this section, we introduce the principal concept proposed in this work. In the previous subsection, we described the Fourier NO, which is capable of approximate maps between function spaces defined on hypercubes in $\mathbb{R}^n$. From a computational perspective, the Fourier layers defining the FNO architecture rely on a Fourier basis defined on a uniform grid and use the fast Fourier transform to perform the spectral convolution. The nonuniform variant was then introduced to accommodate nonuniform rectilinear grids, allowing the grid points to be concentrated in regions of the domain where high accuracy is required. Finally, the GNO mapping was combined with the FNO to approximate operators acting on function spaces defined over sufficiently regular, bounded, open, connected, and bounded subsets of $\mathbb{R}^n$, leading to the \cite{LiEtAl2023GINO} architecture.

A fundamental limitation of FNO architectures is their inability to directly operate on data defined over unstructured discretizations. Indeed, the spectral convolution underlying the FNO relies on the FFT, which inherently assumes that the input and output fields are sampled on a structured grid. In contrast, numerical simulations of partial differential equations on complex geometries are typically performed on unstructured meshes, whose flexibility allows for an accurate representation of arbitrary domains and localized solution features. Bridging these two fundamentally different data representations, therefore, constitutes a key challenge in applying FNOs to practical scientific computing problems. To overcome this limitation, the numerical solution must first be mapped from the original unstructured mesh onto a structured latent discretization. This is achieved by projecting the raw data onto a rectilinear grid through interpolation. Although this preprocessing step is essential to obtain data compatible with the FNO architecture, it inevitably introduces an approximation error associated with the projection of the original solution onto the latent grid. Once the interpolation is performed, the physical quantities of interest are represented as separate channels of the structured input tensor provided to the FNO. In the proposed framework, the signed distance function is systematically included as a primary input feature, allowing the network to encode the geometry of the computational domain together with the physical fields. This structured representation forms the dataset used to train the ADEx-FNO. 

To illustrate the methodology, let us consider the following prototype boundary-value problem:
\begin{equation}
    \label{eq:protoproblem}
    \begin{cases}
        \mathbb{A}(u)=f\quad\text{in}\,\Omega,\\
        \mathcal{B}(u)=0\quad\text{on}\,\partial\Omega,
    \end{cases}
\end{equation}
where $\mathbb{A}(\cdot)$ is a linear or nonlinear differential operator, and $\mathcal{B}(\cdot)$ is a compact representation to denote boundary conditions, e.g. Dirichlet, Neumann, Robin. In general
\begin{equation*}
    \mathcal{B}(u)
    :=
    (u-g)\indicator{\Gamma_D}
    +
    (\partial_{\mathbf{n}}u-h)\indicator{\Gamma_N}
    +
    (\alpha u+\beta\partial_{\mathbf{n}}u-w)\indicator{\Gamma_R}.
\end{equation*}
Here,
\begin{equation*}
    \partial\Omega = \Gamma_D\mathbin{\dot{\cup}}\Gamma_N\mathbin{\dot{\cup}}\Gamma_R.
\end{equation*}
We consider classical hypotheses for $\Omega$, namely, that it is open, bounded, connected, and has a Lipschitz boundary. The well-posedness of \eqref{eq:protoproblem} is intended in the sense of Hadamard \cite{Evans2010}.

We aim to generalize the architecture of $\mathcal{G}^\theta$ by including $\Omega$ among the problem inputs directly within the mathematical setting of the function spaces on which the NO acts. The well-posedness of a partial differential equation depends on the differential operator, the problem data, the boundary conditions, and the geometry of the physical domain on which the solution is defined. From a practical perspective, a NO capable of predicting solutions over arbitrary domains broadens the range of potential applications and increases the value of the computational investment required for training.

Our approach to include $\Omega$ as an input is based on its representation through the signed distance function and on a suitable prolongation of functions in $\mathcal{A}$ and $\mathcal{U}$ outside of $\Omega$. For a given $\Omega\subset\mathbb{R}^d$, the \textit{Signed Distance Function} (SDF) is defined by
\begin{equation}
    \label{eq:sdf}
    d(x):=
    \begin{cases}
        - \underset{y\in\partial\Omega}{\inf}\vertii{x-y}_2 & \text{for}\,x\in\Omega,\\
        \underset{y\in\partial\Omega}{\inf}\vertii{x-y}_2 & \text{for}\,x\in\Omega^c.
    \end{cases}
\end{equation}
Now, consider the hypercube $\mathrm{Q}:=[-M,M]^d$, for some $M>0$, such that $\Omega \subsetneq \mathrm{Q}$ and $\partial\Omega$ is sufficiently far from $\partial\mathrm{Q}$, in a sense that will be specified later. As a consequence, we can identify
\begin{align*}
    \Omega&=\left\{x\in\mathrm{Q}\,\big|\,d(x)<0\right\},\\
    \partial\Omega&=\left\{x\in\mathrm{Q}\,\big|\,d(x)=0\right\},\\
    \Omega^c&=\left\{x\in\mathrm{Q}\,\big|\,d(x)>0\right\}.
\end{align*}
We denote by $\mathcal{D}\subset\mathcal{C}^0(\mathrm{Q};\R)$ the set of signed distance functions associated with all admissible computational domains contained within the hypercube $\mathrm{Q}$.

Next, we address the problem of prolonging functions outside of their domain of definition. Let $v\in\mathcal{V}$, where $\mathcal{V}$ is any space of sufficiently regular functions defined on $\Omega$. We seek a prolongation $\tilde{v}$ defined on $\mathrm{Q}$ such that
\begin{equation}
    \label{eq:pbmext1}
    \tilde{v}(x)=
    \begin{cases}
        v(x)  &\text{for}\,x\in\overline{\Omega},\\
        v^c(x)&\text{for}\,x\in\mathrm{Q}\setminus\overline{\Omega}.
    \end{cases}
\end{equation}
We propose two approaches for constructing such prolongations, depending on the regularity of $\partial\Omega$.

\paragraph{Smooth boundary}
If $\partial\Omega$ is of class $\mathcal{C}^k$, with $k\geq2$, then the signed distance function $d$ has the same regularity (at least) in a sufficiently small tubular neighborhood of $\partial\Omega$; see \cite{Foote1984}. This allows us to prolong the function smoothly outside of its original domain of definition, without introducing sharp gradients that would degrade the Fourier approximation. Assuming that $v$ possesses normal derivatives up to order $N$ on
$\partial\Omega$, with $N\leq k$, we define the exterior prolongation by
\begin{equation}
    \label{eq:smooth-extension}
    v^c(x):= \varphi\left(\frac{d(x)}{\varepsilon}\right)\left[\sum_{m=0}^N\frac{d(x)^m}{m!}\,\partial_{\mathbf{n}}^{(m)}v(q(x))\right]\quad\text{for}\,x\in\mathrm{Q}\setminus\overline{\Omega},
\end{equation}
where $\varphi$ is the bump function compactly supported in $[0,1]$, and $q(x)$ denotes the unique nearest-point projection of $x$ onto $\partial\Omega$ within the tubular neighborhood.

\begin{figure}
    \centering
    \includegraphics[scale=1]{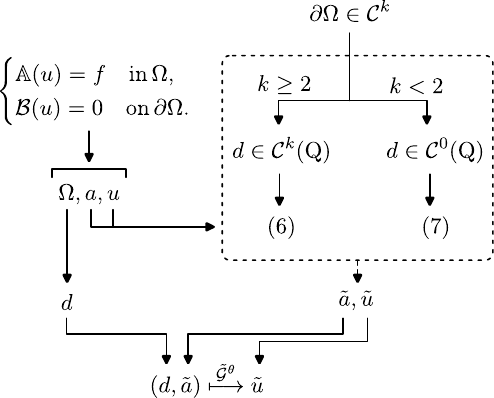}
    \caption{Flowchart of the proposed-function extension procedure used to construct the dataset for training the Domain Agnostic Fourier Neural Operator.}
    \label{fig:extension-pipeline}
\end{figure}

\paragraph{Lipschitz Boundary}
When $\partial\Omega$ is only Lipschitz and may contain singular points, such as sharp corners, we define the exterior prolongation by solving the auxiliary problem:
\begin{equation}
    \label{eq:auxproblem}
    \begin{cases}
        \Delta v^c=0&\text{in}\,\mathrm{Q}\setminus\overline{\Omega},\\
        v^c=0&\text{on}\,\partial\mathrm{Q},\\
        v^c=v|_{\partial\Omega}&\text{on}\,\partial\Omega.
    \end{cases}
\end{equation}
The boundary conditions in \eqref{eq:auxproblem} are understood in the trace sense. For functions defined only over the boundary $\partial\Omega$, e.g. boundary conditions, we first prolong inside by solving
\[
\begin{cases}
    \Delta v = 0 & \text{in}\,\Omega\\
    v = g& \text{on}\,\Gamma_D\\
    v = h& \text{on}\,\Gamma_N\\
    v = w& \text{on}\,\Gamma_R\\
\end{cases}
\]
then, prolong outside accordingly to the regularity of $\partial\Omega$. An overview of the proposed function extension procedure, which constitutes a fundamental component of the dataset generation pipeline, is provided in Figure~\ref{fig:extension-pipeline}.

The two constructions above provide a unique prolongation of $v$ over $\mathrm{Q}\setminus\Omega$, in the sense that for every pair $(\Omega,v)$, where $v$ is defined on $\Omega$, the proposed construction associates a uniquely determined extension $\tilde{v}$ defined on the ambient hypercube $\mathrm{Q}$. The extension is given explicitly by the prescribed procedure, so that no additional choices or ambiguities are involved in its definition. This property allows to uniquely identify each sample in 
\begin{equation*}
    \left\{(a_j,u_j)_\Omega\in\mathcal{A}\times\mathcal{U}\right\}_{j=1}^n,
\end{equation*}
with their prolongations as an element of
\begin{equation*}
    \left\{(d,\tilde{a}_j,\tilde{u}_j)\in\tilde{\mathcal{A}}\times\tilde{\mathcal{U}}\right\}_{j=1}^n,
\end{equation*}
where $\widetilde{\mathcal{A}}$ and $\widetilde{\mathcal{U}}$ denote the corresponding spaces of prolonged input and output functions defined on the entire hypercube $\mathrm{Q}$. We remark that the elements of $\tilde{\mathcal{A}}$ and $\tilde{\mathcal{U}}$ are functions defined over the entire $\mathrm{Q}$.
Finally, this extended mathematical framework allows us to use the FNO to
construct an approximation $\widetilde{\mathcal{G}}^{\theta}$ of the extended
solution operator
\begin{equation*}
\widetilde{\mathcal{G}}^{\dagger}:
\mathcal{D}\times\widetilde{\mathcal{A}}
\longrightarrow
\widetilde{\mathcal{U}}.
\end{equation*}

\paragraph{Inference}
During inference, the model predicts the solution directly on the same latent structured grid employed during training. To recover the solution on an arbitrary set of query points, an interpolant is constructed from the predicted field using a procedure analogous to that adopted during preprocessing and is subsequently evaluated at the desired spatial locations. This final step enables the prediction to be transferred back to an arbitrary mesh-based discretization while preserving the computational efficiency afforded by the structured latent representation.

An important advantage of the proposed ADEx-FNO framework is that the geometry-transfer operations are deterministic and remain outside the optimization loop of the neural operator. Geometry encoding is performed during preprocessing, whereas decoding is performed after inference. Consequently, the framework introduces no additional trainable geometry-processing modules beyond the underlying Fourier-based NO. These deterministic operations do not require an auxiliary optimization problem, although they incur preprocessing and postprocessing costs. Their role is to transform the original mesh-based data into a structured latent representation compatible with the Fourier-based architecture and, after inference, to map the predicted solution back to the desired computational discretization. Thus, the optimization procedure remains that of the underlying NO, while the geometry-transfer overhead is confined to deterministic preprocessing and postprocessing operations.

The proposed architecture relies fundamentally on the SDF, which encodes a sharp geometric representation of the domain through a scalar function. In particular, the model is expected to exploit the fact that regions satisfying $d>0$ contain limited information, since the function rapidly decays to zero outside $\Omega$. Conversely, in the region where $d<0$, the network is required to approximate the PDE solution operator. Hence, the SDF naturally separates the computational domain into an informative region, where the solution must be learned, and an exterior region carrying substantially less relevant information.

\paragraph{Computation of the SDF}
In practice, the SDF associated with the continuous boundary is approximated using a piecewise-linear polynomial representation of the boundary. Given a finite element approximation of $\partial\Omega$, the boundary is represented by a collection of linear segments in 2D and by a collection of planar triangular facets in 3D. In the 2D case, the distance from a point to each boundary segment is computed using the closed-form expression for the distance to the supporting line. The point is first projected onto the line defined by the segment endpoints, after which the projection is clipped to the segment whenever it falls outside its endpoints. In the 3D case, an analogous procedure is employed: the point is projected onto the plane supporting each triangular facet, and the projection is retained whenever it lies inside the triangle; otherwise, the distance is computed to the closest point on the triangle boundary. The unsigned distance is then obtained by taking the minimum over all boundary elements. Then, its sign is assigned by determining whether the point lies inside or outside $\Omega$, consistently with the convention in \eqref{eq:sdf}.


We emphasize that our implementation evaluates the distance to the discrete boundary using the exact geometric formulas for line segments in 2D and triangular facets in 3D. Consequently, up to floating-point roundoff, correct inside-outside classification, and the requirement that the spatial partitioning preserve the exact nearest-element search, no additional geometric approximation is introduced relative to the piecewise-linear boundary representation. The discrepancy between the computed distance and the true SDF associated with the continuous boundary $\partial\Omega$ is therefore attributable to the piecewise-linear approximation of the boundary, and consequently, the error is determined by the mesh size.

An alternative approach is to compute the unsigned distance separately inside and outside $\Omega$ by solving the Eikonal equation $\vertii{\nabla \delta}_2=1$ with $\delta=0$ on $\partial\Omega$. However, such an approach would further complexity in the pipeline, which is the reason why we opted for the direct computation of the distance for each point.

The signed distance is then obtained by assigning a negative sign in $\Omega$ and a positive sign outside $\Omega$. When the exterior problem is solved on the bounded region $\mathrm{Q}\setminus\overline{\Omega}$, an appropriate treatment of the outer boundary $\partial\mathrm{Q}$ must also be specified. Under the standard assumptions for the corresponding Hamilton--Jacobi problem, the distance-to-boundary function is characterized as its viscosity solution.

\section{Numerical validation}
To assess the effectiveness of the algorithm presented in the previous sections, we selected a set of challenging benchmark problems designed to evaluate its accuracy, robustness, and generalization capabilities. The test cases span a range of representative applications, including nonlinear elliptic PDEs and CFD simulations, thereby covering both academic benchmark problems and more realistic engineering scenarios. Although our numerical investigation focuses on these examples, the proposed methodology is not restricted to this particular class of problems. With suitably generated training datasets, we expect the same pipeline to be readily applicable to a broad variety of PDE models, making it a flexible and effective framework for learning solution operators across diverse and increasingly complex physical systems. We adopt the following two metrics to evaluate the accuracy of the proposed method:
\begin{itemize}
    \item discrete relative $\ell^2$ error:
        \begin{equation*}
            \text{err}_{\text{rel},2} (\bm{y}^{\textup{true}},\bm{y}^{\textup{pred}}):=\frac{\vertii{\bm{y}^{\textup{true}}-\bm{y}^{\textup{pred}}}_2}{\vertii{\bm{y}^{\textup{true}}}_2},
        \end{equation*}
    \item range-normalized discrete $\ell^\infty$ error:
        \begin{equation*}
            \text{err}_{\text{rel},\infty}(\bm{y}^{\textup{true}},\bm{y}^{\textup{pred}}):=\frac{\vertii{\bm{y}^{\textup{true}}-\bm{y}^{\textup{pred}}}_\infty}{\texttt{max}(\bm{y}^{\textup{true}})-\texttt{min}(\bm{y}^{\textup{true}})},
        \end{equation*}
\end{itemize}
where $\bm{y}^{\textup{true}}$ denotes the numerical reference solution and $\bm{y}^{\textup{pred}}$ denotes the corresponding model prediction. In the second metric, the maximum and minimum are taken over all entries of $\bm{y}^{\textup{true}}$. The ADEx-FNO models are trained using an $L^2(\mathrm{Q})$-based loss. In our experiments, this choice yielded lower test errors than the conventional mean-squared-error (MSE) loss. The AdamW is used with initial learning rate of $1e-3$ and weight decay $1e-4$. A \texttt{CosineAnnealingLR} learning-rate scheduler is also used. The code implementing the neural model is developed starting from the original modules provided in the \texttt{neuralop} library \cite{NO-library}. Each model was trained on a single NVIDIA A100 GPU hosted on the Ibex cluster at KAUST. Table~\ref{tab-1} summarizes the training configuration adopted for each benchmark problem, together with the main characteristics of the resulting ADEx-FNO models. In addition to the dataset size and training hyperparameters, the table reports the model complexity, final training and test errors, and the computational cost associated with both training and inference.

\begin{table}[H]
    \centering
    \caption{ADEx-FNO training configurations and model characteristics for all benchmark cases. The table reports the dataset sizes, batch sizes, numbers of training epochs, input resolutions, model parameter counts, final training and test losses, training time per epoch, and inference time per solution. All timings were measured on a single NVIDIA A100 GPU.}
    \renewcommand{\arraystretch}{1.15}
    \begin{tabularx}{\textwidth}{lccccc}
        \toprule
        & \makecell{2D nonlinear\\Poisson} & \makecell{2D advect.-\\diff.} & \makecell{3D nonlinear\\Poisson} & \makecell{2D\\airfoil} & \makecell{3D\\wing} \\ \midrule
        Training samples & 1,700 & 1,700 & 819 & 497 & 382 \\ \midrule
        Test samples & 300 & 300 & 205 & 88 & 68 \\ \midrule
        Batch size & 6 & 6 & 1 & 5 & 2 \\ \midrule
        Epochs & 1,500 & 1,000 & 1,500 & 1,500 & 1,200 \\ \midrule
        Resolution & $1,024^2$ & $1,024^2$ & $128^3$ & $1,024^2$ & $152^3$ \\ \midrule
        \makecell[l]{Number of\\parameters} & \num{1.74e6} & \num{1.74e6} & \num{2.80e7} & \num{4.94e7} & \num{2.66e8} \\ \midrule
        Training loss & \num{8.62e-5} & \num{9.13e-5} & \num{1.77e-5} & \num{7.40e-4} & \num{5.42e-3} \\ \midrule
        Test loss & \num{8.92e-4} & \num{4.84e-4} & \num{3.85e-5} & \num{1.89e-3} & \num{3.43e-2} \\ \midrule
        \makecell[l]{Training time\\(s/epoch)} & \num{2.69e2} & \num{2.69e2} & \num{2.59e2} & \num{1.11e2} & \num{1.92e2} \\ \midrule
        \makecell[l]{Inference time\\(s/solution)} & \num{7.27e-3} & \num{7.44e-3} & \num{8.75e-3} & \num{2.44e-2} & \num{1.28e-3} \\
        \bottomrule
    \end{tabularx}
    \label{tab-1}
\end{table}

\subsection{Nonlinear Poisson problem}
\begin{figure}[htbp]
    \centering
    \includegraphics[width=0.75\linewidth]{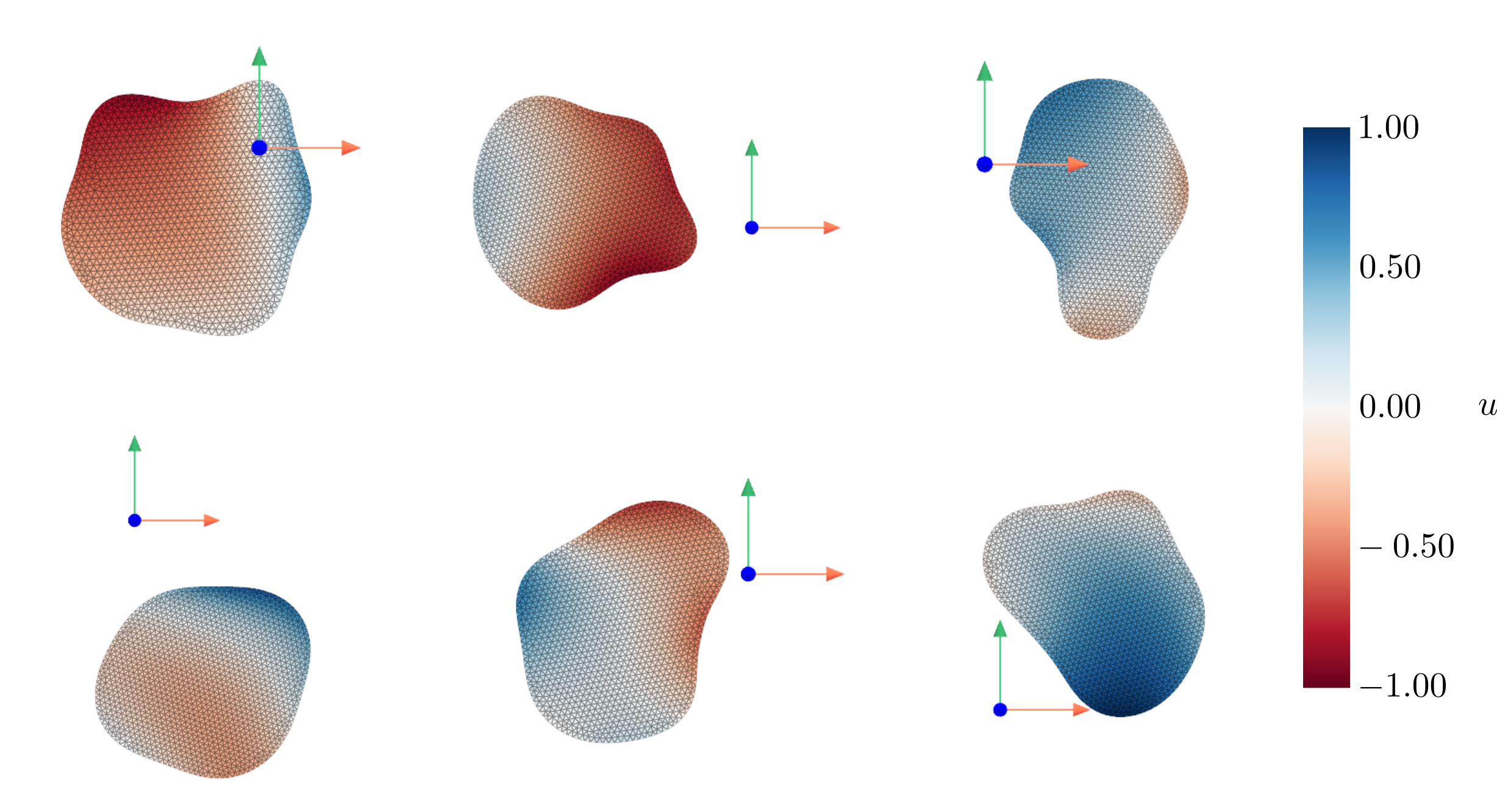}
    \caption{Examples of the 2D dataset used to train the ADEx-FNO model. Each sample represents a domain with smooth boundary.}
    \label{fig:smooth_domains_2D}
\end{figure}

As first test, we consider the nonlinear Poisson equation
\begin{equation}
\label{poisson_eq}
    \begin{cases}
        -\operatorname{div}\left((0.05+\verti{u})\nabla u\right)=f, & \text{in } \Omega,\\
        u=g, & \text{on } \partial\Omega.
    \end{cases}
\end{equation}
We constructed the dataset from point clouds representing the boundary of the PDE domain, following the \textit{Elasticity} test case described in \cite[Section 3.2]{LiEtAl2023GINO}. In the reference setting, the physical domain is embedded in the unit box $[0,1]^2$. In contrast, our construction introduces additional geometric variability: each sample is generated by first applying a spatial dilation, followed by a planar rotation and a rigid translation. The transformations parameters are sampled independently for each instance in the dataset, thereby enriching the diversity of the geometries encountered during training. Figure~\ref{fig:smooth_domains_2D} shows six instances of such geometries. The boundary is obtained by fitting a B-Spline through the points defining $\partial\Omega$. The spline construction is chosen so that $\partial\Omega$ is at least of class $\mathcal{C}^2$.

We consider the source terms $f(x,y)=y/2+e^{-x^2}$ and $g(x,y)=y/2\tanh(x/2)$. The corresponding reference solutions are computed with FEniCSx using a finite-element discretization based on $P_1$ elements, yielding an unstructured computational mesh defined on the physical domain $\Omega\subset\R^2$. To obtain a representation suitable for operator learning, the numerical solution is subsequently extended from $\Omega$ to the larger domain $\mathrm{Q}=[-4,4]^2$ by means of the smooth extension procedure described in the previous section.

The extended solution field is then interpolated onto a finer structured grid, which serves as the latent discretization for the Fourier NO. Throughout this work, linear interpolation is employed for this projection. For the datasets considered here, the measured transfer error is negligible relative to the overall prediction error of the proposed pipeline. Consequently, this preprocessing stage can be regarded as effectively lossless while providing a structured representation compatible with Fourier-based NOs.

The resulting structured dataset is used to train the ADEx-FNO models. The target output consists of the extended solution field represented on the latent grid as a single channel. The input, on the other hand, is composed of multiple channels encoding both the geometry of the computational domain and the governing physical problem. The signed distance function associated with $\Omega$ is provided as the primary geometric descriptor, complemented by two additional channels containing the Cartesian coordinates $(x,y)$ of the latent grid. Furthermore, the extended source term $f$ is included as an additional input channel, ensuring that the forcing term defining the PDE is explicitly available to the model. Finally, the harmonic lifting obtained by solving the Laplace equation with the prescribed boundary data is provided as an additional input feature. This auxiliary field conveys information about the influence of the prescribed boundary conditions independently of the source term, thereby supplying the ADEx-FNO with a richer physical representation of the underlying problem and facilitating the learning of the solution operator.

\begin{figure}[H]
    \centering
    \includegraphics[width=0.75\linewidth]{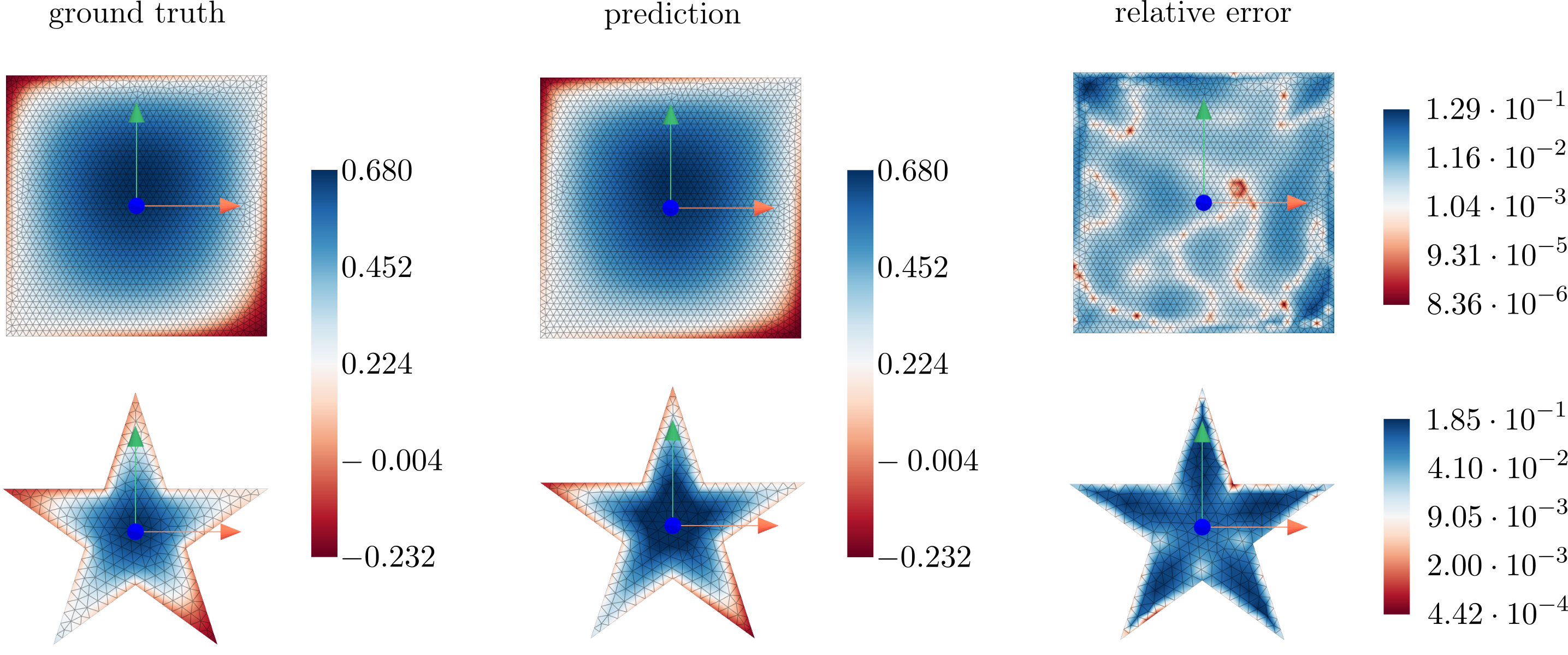}
    \caption{ADEx-FNO predictions of the 2D ADEx-FNO model for OOD square and star-shaped geometries governed by the nonlinear Poisson problem \eqref{poisson_eq}.}
    \label{fig:poisson_2D}
\end{figure}

To further assess the robustness of the proposed framework, we evaluated the trained model on geometries that were deliberately excluded from the training distribution (or out-of-distribution (OOD)). As shown in Figure~\ref{fig:poisson_2D}, although it is trained exclusively on smooth domains, the ADEx-FNO produces very good predictions for domains with sharp corners and highly nonconvex boundaries, including square and star-shaped geometries. Nevertheless, the quantitative errors reported below show that the generalization accuracy depends on the geometric departure from the training distribution.

\begin{table}[H]
\centering
\caption{Summary of the ADEx-FNO inference performance of the proposed machine learning model for 2D geometries on PDE problem \eqref{poisson_eq}. The reported metrics characterize both prediction accuracy (relative and maximum errors) and computational efficiency (average inference time).}
\begin{tabularx}{\textwidth}{lYYY}
    \toprule
    & Test samples & Square & Star \\
    \midrule
    $\text{err}_{\text{rel},2}$       & \makecell{0.77 \%\\(std: \num{3.73e-3})}   & \makecell{2.81 \%\\\,}      & \makecell{26.28 \%\\\,}   \\
    \midrule
    $\text{err}_{\text{rel},\infty}$       & \makecell{2.41 \%\\(std: \num{6.74e-3})}   & \makecell{12.88 \%\\\,}      & \makecell{18.54 \%\\\,}   \\
    \midrule
    Inference time (s)                 & \makecell{\num{7.27e-3} \\(std: \num{7.08e-5})}    & \makecell{\num{7.03e-3} \\\,}    & \makecell{\num{7.34e-3} \\\,} \\
    \midrule
    \makecell[l]{Rectilinear interpolation\\time (s)}  & \makecell{\num{7.82e-3} \\(std: \num{4.53e-4})}    & \makecell{\num{6.85e-3} \\\,}    & \makecell{\num{6.52e-3} \\\,} \\
    \bottomrule
\end{tabularx}
\label{Poisson-inference}
\end{table}

\begin{figure}[H]
    \centering
    \includegraphics[width=0.75\linewidth]{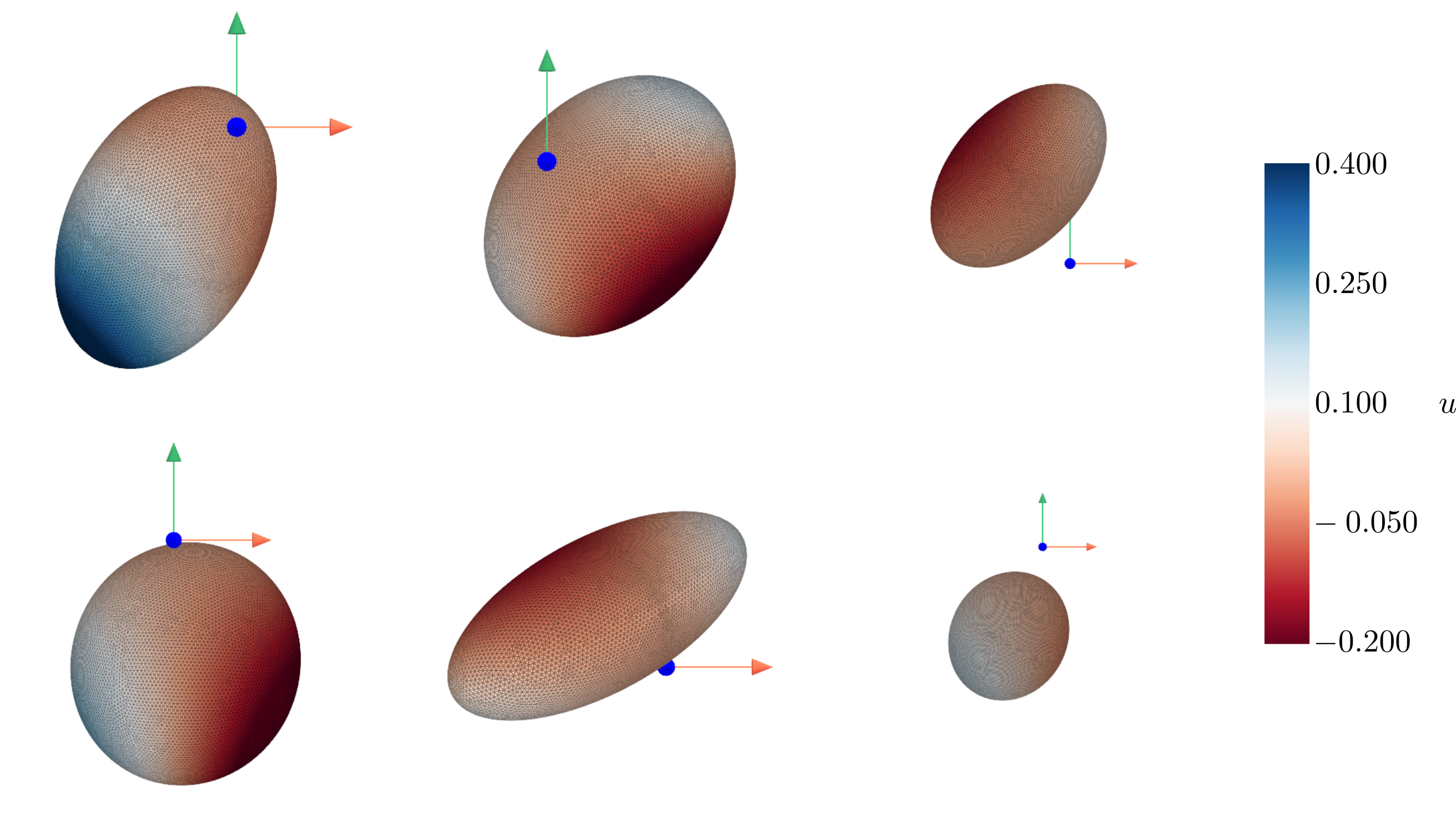}
    \caption{Examples of the 3D dataset used to train the ADEx-FNO model. Each sample represents a domain with smooth boundary.}
    \label{fig:smooth_domains_3D}
\end{figure}

To assess the dimensional robustness of the proposed framework, we considered the same partial differential equation posed in 3D. An independent dataset was generated following exactly the same pipeline adopted in the 2D case, allowing the ADEx-FNO to be trained and evaluated under analogous modeling assumptions. In this setting, the computational domains are ellipsoids parametrized by seven degrees of freedom. Specifically, we vary the lengths of the three semi-axes, a restricted 3D rotation parameterized by two angles, and finally translate the geometry using a fixed-magnitude offset vector whose direction is parameterized by two additional angles. A selection of six representative geometries, together with their corresponding solution fields, is shown in Figure \ref{fig:smooth_domains_3D}.

\begin{figure}[htbp]
    \centering
    \includegraphics[width=0.85\linewidth]{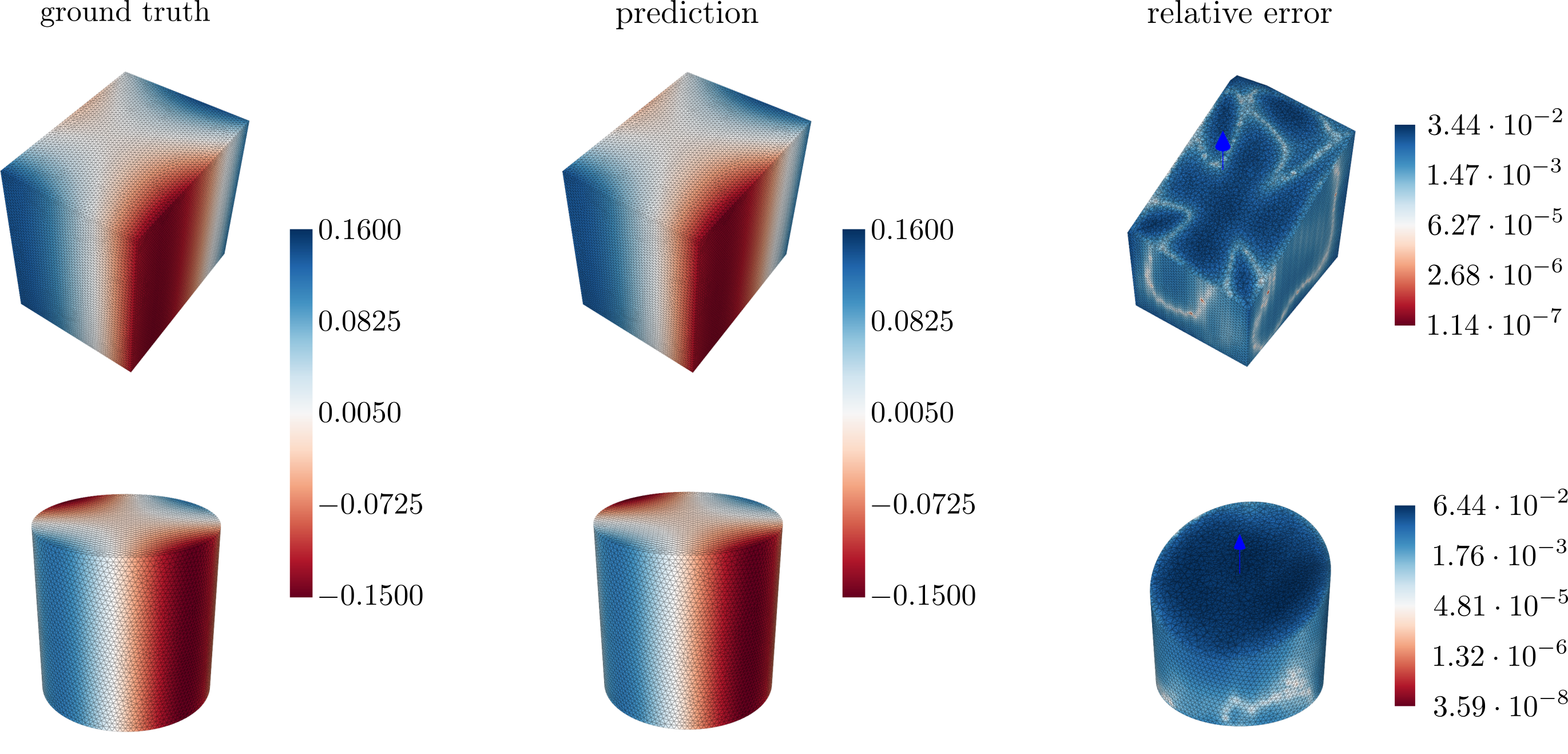}
    \caption{ADEx-FNO inference of the 3D NO model on \textit{out of distribution} geometries for PDE problem \eqref{poisson_eq}.}
    \label{fig:inference_3D}
\end{figure}

As in the 2D experiments, the training dataset consists exclusively of smooth geometries. To evaluate the model's ability to generalize beyond the training distribution, we performed inference on geometries with substantially different characteristics, namely a cylinder and a parallelepiped. The corresponding predictions are reported in Figure \ref{fig:inference_3D}, demonstrating that the proposed architecture successfully generalizes to previously unseen 3D domains. A quantitative summary of the ADEx-FNO inference accuracy is provided in Table \ref{tab-6}. Taken together with the 2D results, these experiments show that the proposed framework maintains its predictive capabilities across different spatial dimensions while preserving the same data-generation and learning pipeline.


\begin{table}[H]
    \centering
    \caption{Summary of the ADEx-FNO inference performance of the proposed machine learning model for 3D geometries on PDE model \eqref{poisson_eq}. The reported metrics characterize both prediction accuracy (relative and maximum errors) and computational efficiency (average inference time).}
    \renewcommand{\arraystretch}{1.15}
    \begin{tabularx}{\textwidth}{lYYY}
        \toprule
        & Test samples & Parallelepiped & Cylinder \\
        \midrule
        $\mathrm{err}_{\mathrm{rel},2}$ & \makecell{0.32\,\% \\ (std: \num{2.48e-3})} & \makecell{2.54\,\% \\ \,} & \makecell{9.19\,\% \\ \,} \\

        $\mathrm{err}_{\mathrm{rel},\infty}$ & \makecell{0.36\,\% \\ (std: \num{2.88e-3})} & \makecell{3.36\,\% \\ \,} & \makecell{6.01\,\% \\ \,} \\

        Inference time (s) & \makecell{\num{8.35e-3} \\ (std: \num{1.22e-3})} & \makecell{\num{8.81e-3} \\ \,} & \makecell{\num{8.68e-3} \\ \,} \\

        \makecell[l]{Rectilinear interpolation\\time (s)} & \makecell{\num{3.93e-2} \\ (std: \num{4.95e-3})} & \makecell{\num{4.23e-2} \\ \,} & \makecell{\num{2.88e-2} \\ \,} \\

        \bottomrule
    \end{tabularx}
    \label{tab-6}
\end{table}

\subsection{Nonlinear advection-reaction-diffusion problem}
To further assess the generality of the proposed framework, we next consider a different governing equation, namely a nonlinear advection-reaction-diffusion problem posed on 2D domains. In contrast to the previous experiment, which considered the nonlinear Poisson problem in 2D and 3D, the objective here is to assess whether the same framework can accommodate a different class of partial differential equations without modifying the overall data-generation and learning pipeline.

Consider the following problem:
\begin{equation}
\label{adv_eq}
    \begin{cases}
        -\eta\,\Delta u
        +(1+u)\,\mathbf{v}\cdot\nabla u
        +\gamma u
        =f,
        & \text{in } \Omega,\\
        u=g,
        & \text{on } \partial\Omega.
    \end{cases}
\end{equation}
where $\eta=0.5$, $\gamma=0.01$, and
$\mathbf{v}=\frac{1}{\sqrt{2}}(1,1)^\top$. The Dirichlet datum and source term are given by
$g(x,y)=\sin\left(\frac{\pi x}{2}\right)\cos(y)$ and
$f(x,y)=y e^{-x^2}$, respectively. 

To ensure a controlled comparison, the training dataset was generated using the same family of smooth geometries adopted for the nonlinear Poisson problem, while changing the governing PDE and its associated source and boundary data. This allows the influence of the physical model to be assessed while keeping the geometric distribution unchanged. As before, the NO was trained exclusively on smooth domains represented through their signed distance functions.

\begin{figure}[htbp]
    \centering
    \includegraphics[width=0.85\textwidth]{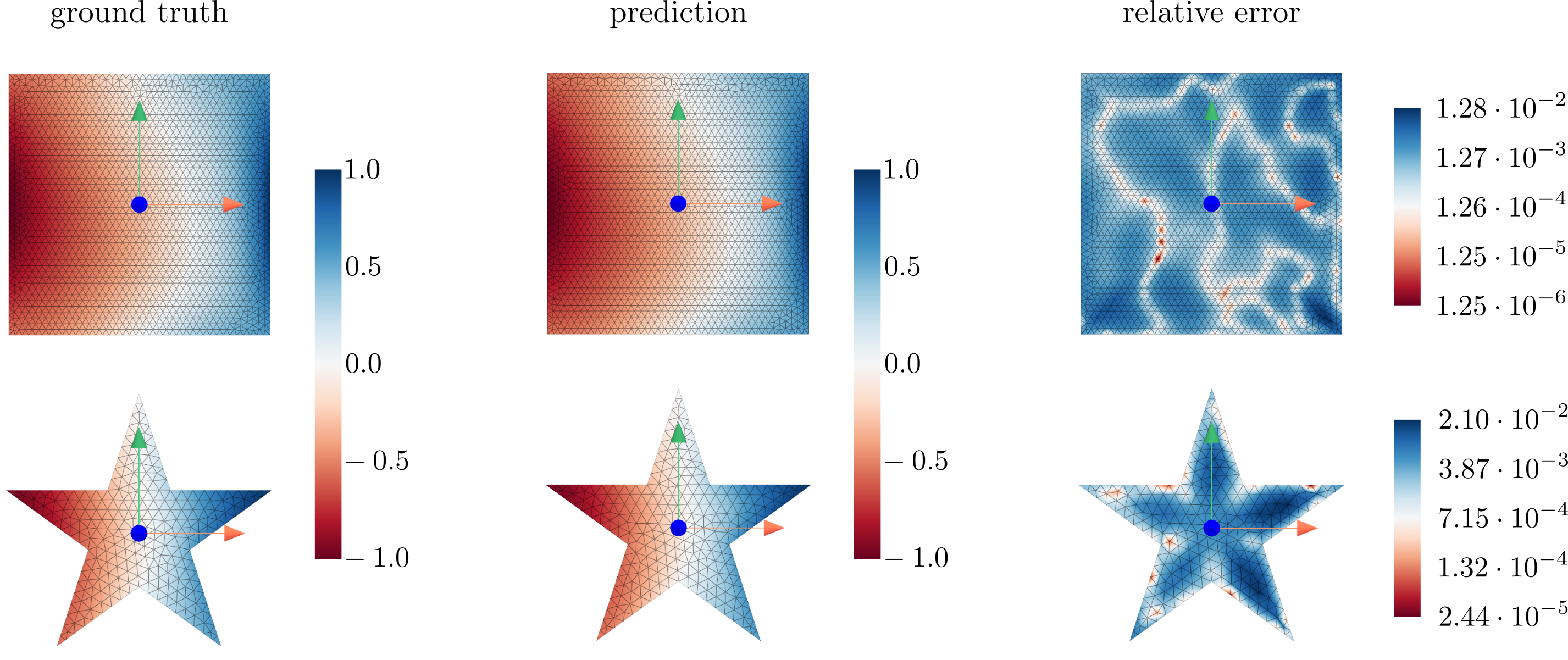}
    \caption{ADEx-FNO predictions of the 2D NO model for OOD square and star-shaped geometries governed by the nonlinear advection-reaction-diffusion problem \eqref{adv_eq}.}
    \label{fig:inference_2D_advection}
\end{figure}

The geometric generalization capability of the trained model was then evaluated on geometries lying outside the training distribution. In particular, inference was performed on square and star-shaped domains, whose sharp corners and nonsmooth boundaries were not represented during training. The corresponding predictions are presented in Figure~\ref{fig:inference_2D_advection}, while the quantitative errors are reported in Table~\ref{tab:adv-inference}. Together, these results assess the model's ability to approximate the solution on geometries whose characteristics differ substantially from those represented in the training dataset. The same evaluation metrics adopted for the nonlinear Poisson problem are used here. The results demonstrate that the same preprocessing and learning framework can be applied to another class of PDEs and that the independently trained model can generalize to previously unseen geometries.

\begin{table}[H]
\centering
\caption{ADEx-FNO inference performance for the 2D nonlinear advection-reaction-diffusion problem \eqref{adv_eq}. The test-set column reports aggregate statistics over the held-out smooth-domain test set, whereas the square and star columns correspond to individual OOD geometries. The reported quantities characterize the prediction accuracy and the computational costs of inference and rectilinear interpolation.}
\label{tab:adv-inference}
\begin{tabularx}{\textwidth}{lYYY}
    \toprule
    & Test set & Square & Star \\
    \midrule
    $\text{err}_{\text{rel},2}$ & \makecell{0.49 \%\\(std: \num{7.10e-3})} & \makecell{0.80 \%\\\,} & \makecell{2.84 \%\\\,} \\
    \midrule
    $\text{err}_{\text{rel},\infty}$ & \makecell{0.44 \%\\(std: \num{5.80e-3})} & \makecell{1.28 \%\\\,} & \makecell{2.21 \%\\\,} \\
    \midrule
    Inference time (s) & \makecell{\num{7.44e-3} \\(std: \num{4.73e-5})} & \makecell{\num{7.37e-3} \\\,} & \makecell{\num{7.25e-3} \\\,} \\
    \midrule
    \makecell[l]{Rectilinear interpolation\\ time (s)}  & \makecell{\num{5.82e-3} \\(std: \num{5.33e-4})} & \makecell{\num{7.15e-3} \\\,} & \makecell{\num{6.39e-3} \\\,} \\
    \bottomrule
\end{tabularx}
\end{table}

\subsection{RANS and URANS simulations}
\label{sec:rans_urans}

Reynolds-averaged Navier--Stokes (RANS) and unsteady Reynolds-averaged
Navier--Stokes (URANS) provide complementary levels of turbulence modeling.
RANS targets statistically stationary mean flows and remains widely used for
routine analysis, design exploration, optimization, uncertainty
quantification, and aerodynamic-database generation. URANS advances the
averaged equations in physical time and is appropriate when coherent
large-scale unsteadiness, separation, vortex shedding, moving geometry, or
time-dependent loading must be represented while the unresolved turbulence is
modeled.

These formulations remain central to engineering CFD despite advances in
hardware and artificial intelligence. Large-eddy simulation (LES) and direct
numerical simulation (DNS) provide greater scale resolution, but their mesh,
time-integration, and sampling requirements remain substantial for realistic
Reynolds numbers, wall-bounded flows, complex geometries, and repeated
many-query studies
\cite{CFDVision2030Update2024,AltmannEtAl2023,GoinisEtAl2026}. RANS and URANS
therefore retain a practical balance among physical fidelity, robustness,
turnaround time, and compatibility with established verification and
validation procedures
\cite{ManiDorgan2023,MaueryEtAl2021,SteinerEtAl2022,VisonneauEtAl2022,CroneEtAl2024,Ricci2025}.

Artificial intelligence and scientific machine learning are most naturally
positioned as complementary technologies within the CFD multi-fidelity hierarchy,
rather than as immediate replacements for established solvers \cite{VinuesaBrunton2022,CFDVision2030Update2024}.
Even a
moderate reduction for one calculation can accumulate into a substantial
saving when simulations are repeated hundreds or thousands of times for design
optimization, aerodynamic-database generation, uncertainty quantification,
digital-twin construction, control studies, or certification-oriented analyses
\cite{DuensingKenway2020,CFDVision2030Update2024,MaueryEtAl2021,DuensingKenway2020}. 
Scientific machine learning is also used here as a complementary technology: the ADEx-FNO
is evaluated once to provide an initial field, after which the conventional
solver retains complete control.

\subsubsection{RANS test cases}
\label{sec:rans_cases}

We consider 2D low-Mach-number airfoil flows and 3D transonic wing flows. In
both settings, SU2~\cite{SU2} solves the compressible RANS equations closed
with a Spalart--Allmaras model. The principal physical and numerical settings
are reported in Table~\ref{tab:rans_su2_configurations} in
\ref{app:rans_configurations}. The ADEx-FNO supplies only the initial flow
field; SU2 performs all subsequent pseudo-time and linear iterations and
applies the prescribed convergence conditions.

Since external-aerodynamic meshes concentrate resolution near the body and
large gradients, we use the NUFNO architecture described in
Section~\ref{sec:architecture}. Converged CFD fields are projected to the
latent grid by linear interpolation, and the inferred field is transferred to
the target CFD nodes with \texttt{RegularGridInterpolator}. The target mesh
therefore need not coincide with either the latent discretization or a mesh
used to generate the training data.

Each uniform-flow calculation is paired with an ADEx-FNO-based calculation using the
same mesh, physical parameters, numerical methods, linear-solver settings, and
convergence criterion. Let $N$ be the number of executed pseudo-time
iterations, and let $L_{\mathrm{flow}}$ and $L_{\mathrm{turb}}$ be the
cumulative mean-flow and turbulence linear iterations. For any work measure
$Q$, we define
\begin{equation}
    G_Q
    =
    100\left(1-\frac{Q_{\mathrm{ADEx-FNO}}}{Q_{\mathrm{UF}}}\right),
    \qquad
    Q\in\left\{N,L_{\mathrm{flow}},L_{\mathrm{turb}},L_\Sigma\right\},
    \qquad
    L_\Sigma=L_{\mathrm{flow}}+L_{\mathrm{turb}},
    \label{eq:rans_work_reduction}
\end{equation}
where the subscripts $\mathrm{UF}$ and $\mathrm{ADEx-FNO}$ denote uniform-flow and ADEx-FNO-based
initialization. The unweighted sum $L_\Sigma$ is an algebraic-work diagnostic,
not an exact cost-equivalent sum.

\begin{figure}[H]
\centering
\includegraphics[width=0.99\linewidth]{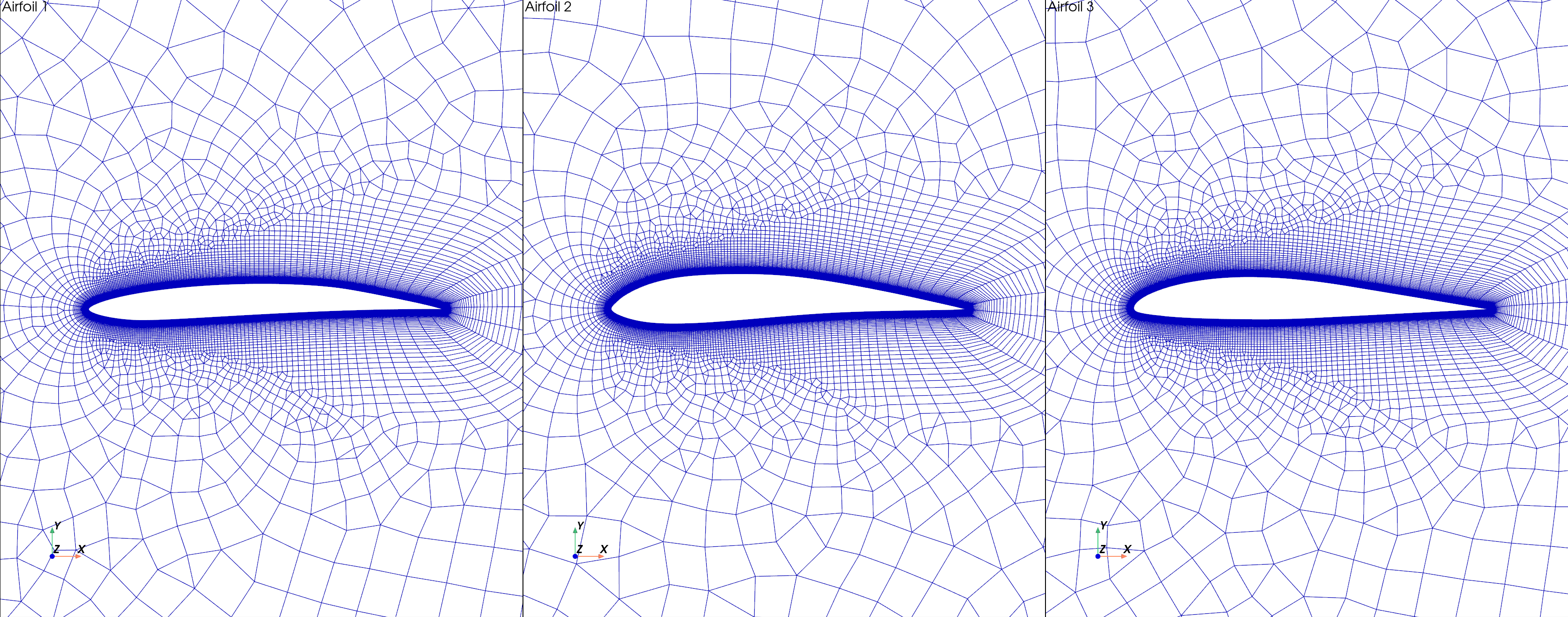}
\caption{Held-out airfoil geometries used to evaluate ADEx-FNO-based initialization
of the 2D RANS solver.}
\label{fig:rans_2d_airfoils}
\end{figure}

\paragraph{2D airfoil cases}
The ADEx-FNO was trained on $45$ airfoil geometries, each evaluated at $13$ angles
of attack in $[-5^\circ,5^\circ]$. Incidence was encoded by rotating the
airfoil during mesh generation, so the signed distance function represents
both shape and orientation. The model is evaluated on the three unseen
airfoils (shown in Figure~\ref{fig:rans_2d_airfoils}) at
\begin{equation*}
    \alpha\in
    \left\{-7^\circ,-2.5^\circ,0^\circ,2.5^\circ,7^\circ\right\},
\end{equation*}
for $15$ geometry-operating-condition pairs. The three angles inside the
training interval are in-distribution (ID), whereas $\alpha=\pm7^\circ$ are
OOD. A run is converged when
$r_\rho=\log_{10}(\mathrm{RMS}_\rho)$ first falls below $-12$ after monitoring
begins. The stored pseudo-time index is zero based, so the executed count is
the final index plus one.

\begin{table}[H]
\centering
\caption{Executed pseudo-time iterations for the 2D RANS airfoil cases.
Bold entries identify OOD angles. NC denotes the Airfoil~2,
$\alpha=0^\circ$ pair, for which neither initialization satisfies
$r_\rho<-12$ within $29{,}999$ iterations.}
\label{tab:rans_2d_iterations}
\begin{tabular}{@{}lrrrr@{}}
\toprule
& $\alpha$ ($^\circ$) & \makecell{Uniform-flow\\iterations}
& \makecell{ADEx-FNO-based\\iterations} & $G_{\mathrm{outer}}$ (\%) \\
\midrule
Airfoil~1 & \textbf{-7.0} & \textbf{8,402} & \textbf{4,969} & \textbf{40.86} \\
& -2.5 & 5,999 & 3,990 & 33.49 \\
& 0.0 & 5,855 & 3,980 & 32.02 \\
& 2.5 & 6,072 & 3,421 & 43.66 \\
& \textbf{7.0} & \textbf{7,704} & \textbf{4,590} & \textbf{40.42} \\
\midrule
Airfoil~2 & \textbf{-7.0} & \textbf{8,167} & \textbf{3,964} & \textbf{51.46} \\
& -2.5 & 6,538 & 3,652 & 44.14 \\
& 0.0 & 29,999~(NC) & 29,999~(NC) & -- \\
& 2.5 & 6,260 & 3,284 & 47.54 \\
& \textbf{7.0} & \textbf{7,325} & \textbf{3,652} & \textbf{50.14} \\
\midrule
Airfoil~3 & \textbf{-7.0} & \textbf{7,745} & \textbf{3,640} & \textbf{53.00} \\
& -2.5 & 6,186 & 3,487 & 43.63 \\
& 0.0 & 6,057 & 3,490 & 42.38 \\
& 2.5 & 6,358 & 3,477 & 45.31 \\
& \textbf{7.0} & \textbf{7,324} & \textbf{3,644} & \textbf{50.25} \\
\bottomrule
\end{tabular}
\end{table}

\begin{figure}[H]
\centering
\includegraphics[width=0.75\linewidth]{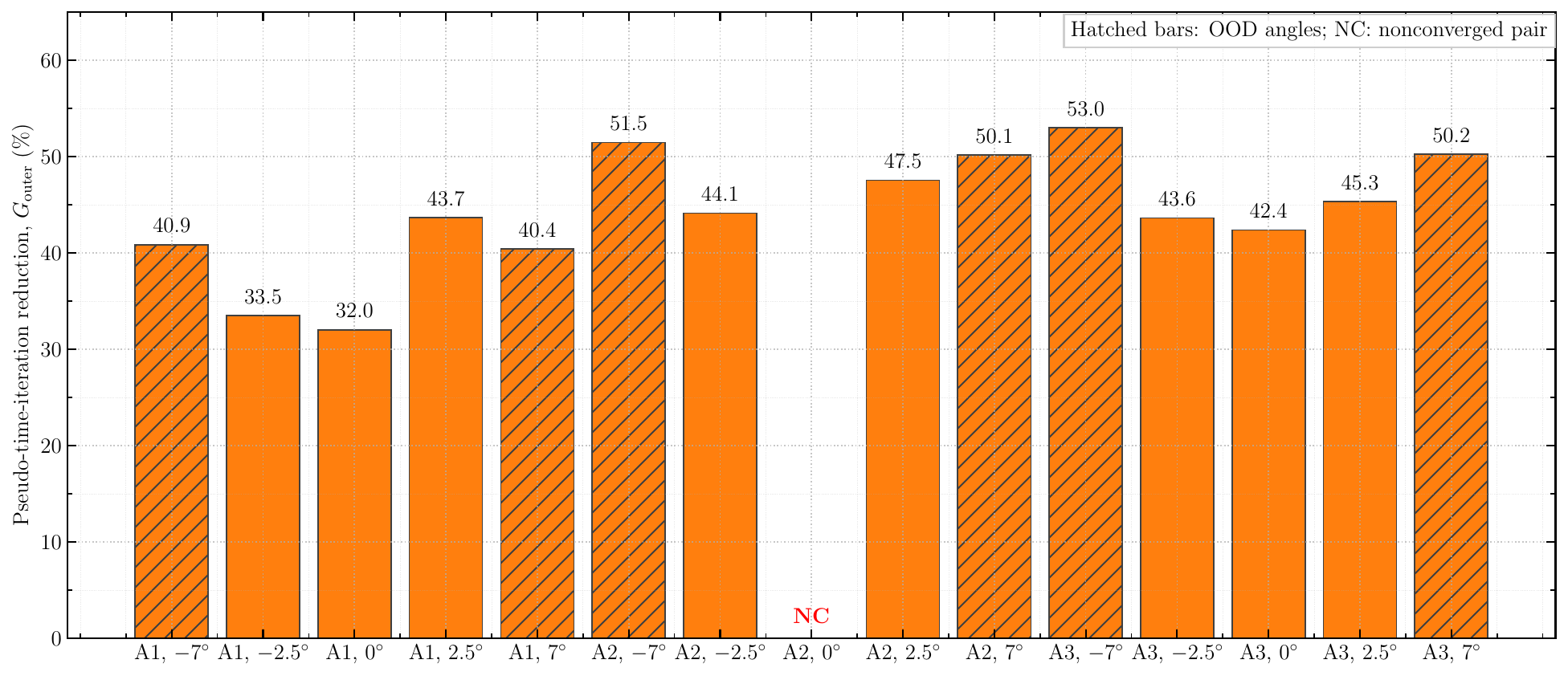}
\caption{Reduction in executed pseudo-time iterations for the complete 2D RANS
test matrix. Hatched bars denote OOD angles; NC identifies the Airfoil~2,
$\alpha=0^\circ$ pair for which neither initialization reaches the prescribed
threshold.}
\label{fig:rans_2d_outer_reduction}
\end{figure}

Both initializations converge in $14$ of the $15$ cases, and the ADEx-FNO-based run
requires fewer pseudo-time iterations in every converged pair. The pairwise
reduction $G_{\mathrm{outer}}$ has a mean of $44.17\%$ and ranges from
$32.02\%$ to $53.00\%$. Figure~\ref{fig:rans_2d_outer_reduction} shows the
complete matrix and retains the nonconverged pair.

The mean reduction is $41.52\%$ for the eight converged ID pairs and $47.69\%$
for the six OOD pairs. This difference is descriptive: the OOD cases also have
the largest $|\alpha|$ and generally larger uniform-flow iteration counts. The
geometry-specific means are $38.09\%$, $48.32\%$, and $46.91\%$ for
Airfoils~1--3, respectively; three held-out geometries are insufficient to
establish a general geometric trend. Detailed grouped statistics are given in
\ref{app:rans_2d_grouped}.

For Airfoil~2 at $\alpha=0^\circ$, both calculations reach the limit of
$29{,}999$ iterations without satisfying $r_\rho<-12$. Their minimum residuals
are $-10.8615$ and $-11.0630$ for uniform-flow and ADEx-FNO-based initialization,
respectively. This pair is excluded from statistics requiring a convergence
iteration.

Aggregated over the $14$ converged pairs, the pseudo-time count decreases from
$95{,}992$ to $53{,}240$ ($44.54\%$), while the cumulative mean-flow and
turbulence linear counts decrease by $44.37\%$ and $43.84\%$. The unweighted
combined count decreases by $44.23\%$. Pairwise, $G_\Sigma$ has a mean of
$43.84\%$ and ranges from $31.29\%$ to $52.69\%$.

\begin{table}[H]
\centering
\caption{Aggregate pseudo-time and linear-solver work over the $14$
converged 2D RANS pairs. $L_\Sigma$ is the unweighted sum of the mean-flow
and turbulence linear iterations.}
\label{tab:rans_2d_aggregate_work}
\begin{tabular}{@{}lrrr@{}}
\toprule
Quantity & \makecell{Uniform-flow\\initialization}
& \makecell{ADEx-FNO-based\\initialization} & Reduction (\%) \\
\midrule
Executed pseudo-time iterations & 95,992 & 53,240 & 44.54 \\
Mean-flow linear iterations & 7,170,150 & 3,988,806 & 44.37 \\
Turbulence linear iterations & 2,562,832 & 1,439,249 & 43.84 \\
Unweighted sum $L_\Sigma$ & 9,732,982 & 5,428,055 & 44.23 \\
\bottomrule
\end{tabular}
\end{table}

The aggregate mean combined linear count per pseudo-time iteration changes only
from $101.39$ to $101.95$. Therefore, the cumulative reduction arises from
fewer pseudo-time iterations rather than cheaper individual iterations. The
interval from the first $r_\rho<-10$ crossing to $r_\rho<-12$ is also shorter
in all $14$ converged pairs, with a mean reduction of $47.33\%$. Finally, the
maximum absolute differences between the converged aerodynamic coefficients
are $2.9176\times10^{-5}$ for $C_L$ and $4.7816\times10^{-7}$ for $C_D$; the
corresponding maximum relative differences are $0.0474\%$ and $0.00263\%$.
The complete residual, linear-solver, and coefficient evidence is given in
\ref{app:rans_2d_appendix}.

\paragraph{Transfer across CFD mesh resolutions}
We next apply the same trained ADEx-FNO, without retraining, to Airfoil~2 at
$\alpha=-7^\circ$ on the base mesh \texttt{L1} and the refined meshes
\texttt{L2} and \texttt{L4}. Their volume-element counts are $24{,}403$,
$44{,}968$, and $93{,}171$, respectively. Each pair uses the same target mesh
and the 2D RANS configuration listed in
Table~\ref{tab:rans_su2_configurations} in~\ref{app:rans_configurations}.

\begin{figure}[H]
\centering
\includegraphics[width=0.99\linewidth]
{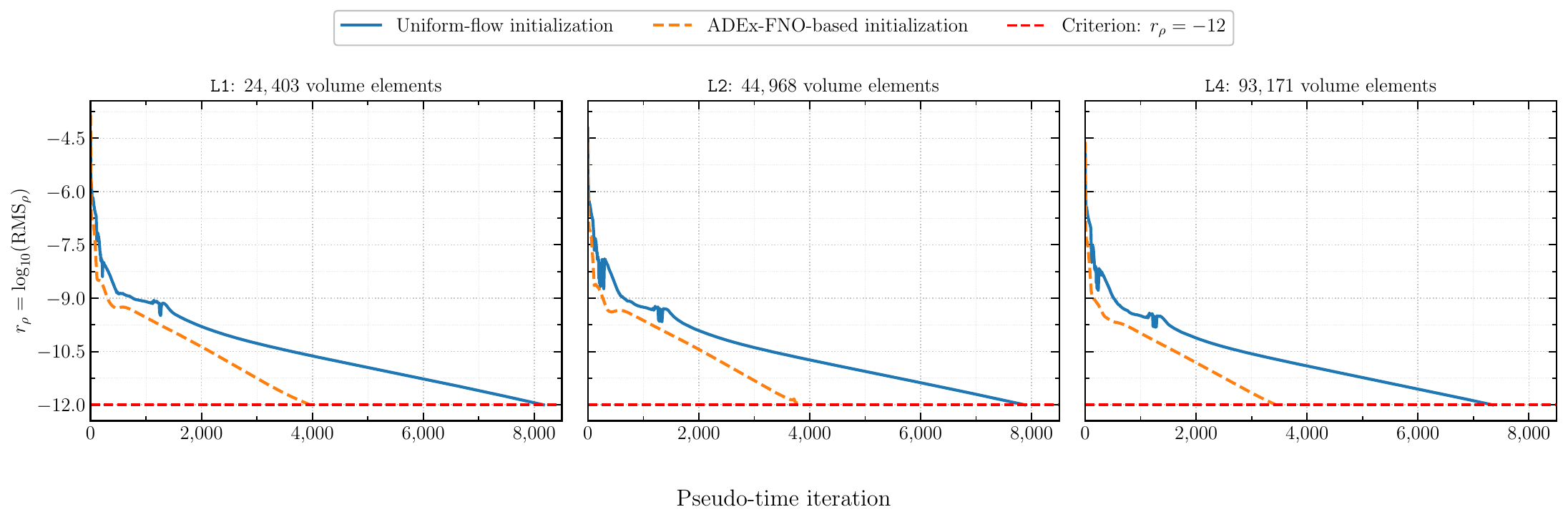}
\caption{Paired RMS-density residual histories for the mesh-refined 2D RANS
Airfoil~2 case at $\alpha=-7^\circ$. Blue solid curves denote uniform-flow
initialization, orange dashed curves denote ADEx-FNO-based initialization, and the
red dashed line marks the prescribed $r_\rho=-12$ threshold.}
\label{fig:rans_2d_mesh_refinement}
\end{figure}

\begin{table}[H]
\centering
\caption{Mesh size and solver-work reductions obtained by applying the same
trained NO to the mesh-refined 2D RANS case.}
\label{tab:rans_2d_mesh_refinement}
\begin{tabular}{@{}lrrrrr@{}}
\toprule
Level & \makecell{Volume\\elements} & $N_{\mathrm{UF}}$ & $N_{\mathrm{ADEx-FNO}}$
& $G_{\mathrm{outer}}$ (\%) & $G_\Sigma$ (\%) \\
\midrule
\texttt{L1} & 24,403 & 8,167 & 3,964 & 51.46 & 50.96 \\
\texttt{L2} & 44,968 & 7,875 & 3,851 & 51.10 & 50.13 \\
\texttt{L4} & 93,171 & 7,340 & 3,440 & 53.13 & 52.31 \\
\bottomrule
\end{tabular}
\end{table}

All six calculations satisfy $r_\rho<-12$. Across the three levels, the
pseudo-time reduction remains between $51.10\%$ and $53.13\%$, and the
reduction in $L_\Sigma$ remains between $50.13\%$ and $52.31\%$. The mean
combined linear count per pseudo-time iteration is similar, and slightly
larger, for the ADEx-FNO-based runs, so the cumulative benefit again follows from
fewer pseudo-time iterations. The interval from $r_\rho=-10$ to final
convergence is reduced by $56.27\%$--$58.59\%$. On each fixed mesh, the two
initializations converge to consistent coefficients, with maximum relative
differences of $0.00471\%$ for $C_L$ and $0.00867\%$ for $C_D$. The detailed
results are reported in~\ref{app:rans_2d_mesh_refinement}. This study
demonstrates transfer over the tested resolutions, not universal mesh
independence or formal grid convergence.

\begin{figure}[H]
    \centering
    \begin{subfigure}[t]{0.49\textwidth}
        \centering
        \includegraphics[width=\linewidth]{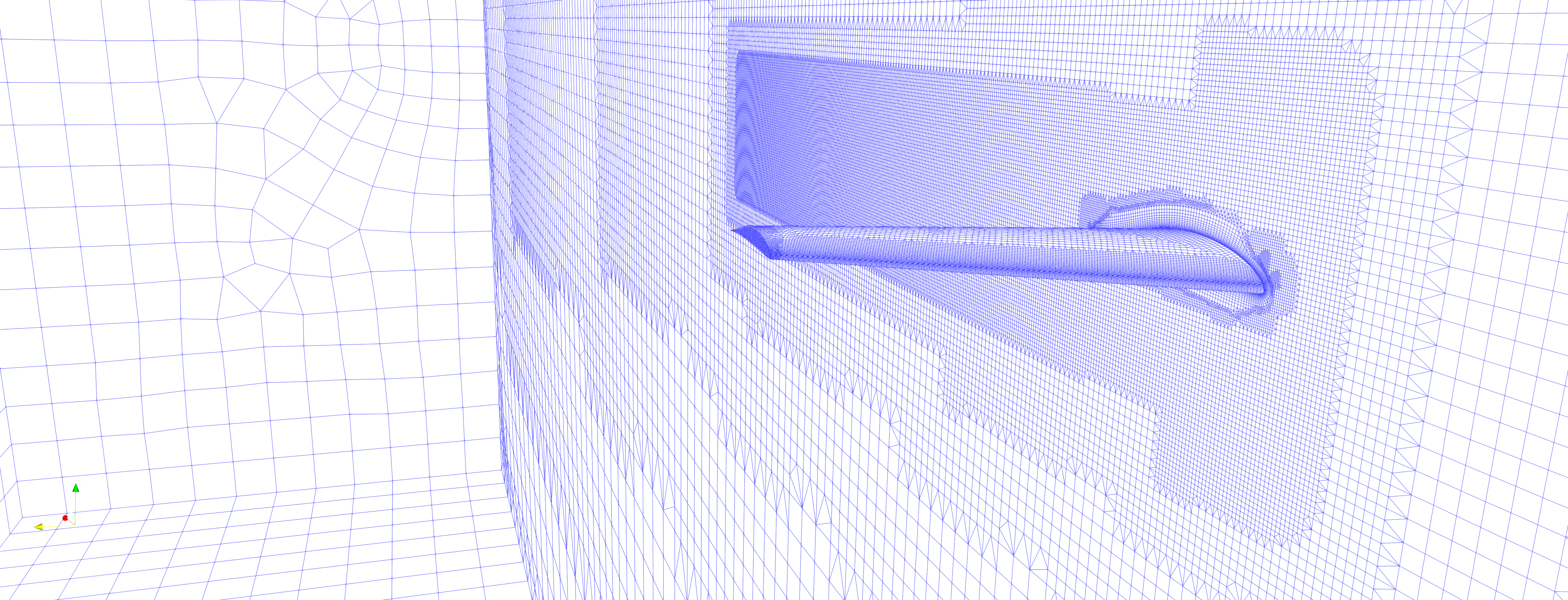}
        \caption{Front view of the mesh on the domain boundaries.}
        \label{fig:rans_3d_mesh_front}
    \end{subfigure}
    \begin{subfigure}[t]{0.49\textwidth}
        \centering
        \includegraphics[width=\linewidth]{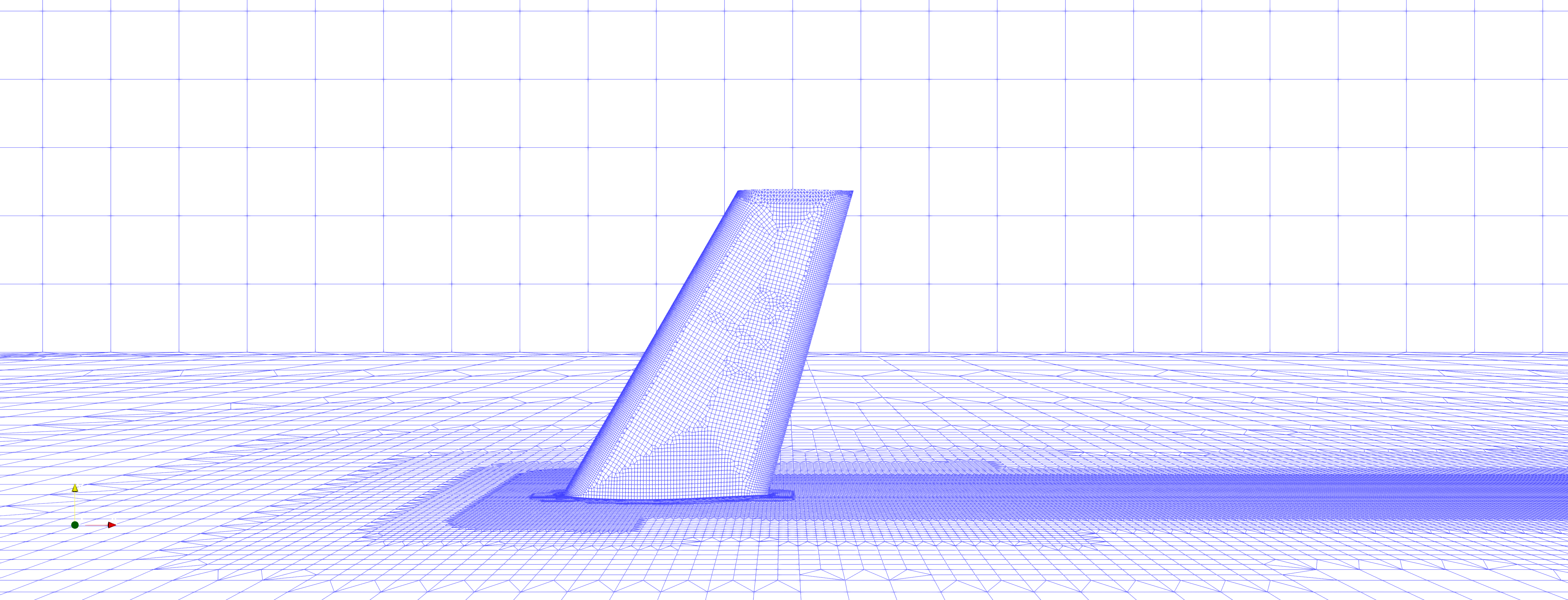}
        \caption{Top view of the mesh on the domain boundaries.}
        \label{fig:rans_3d_mesh_top}
    \end{subfigure}
    \par\medskip
    \begin{subfigure}[t]{0.60\textwidth}
        \centering
        \includegraphics[width=\linewidth]{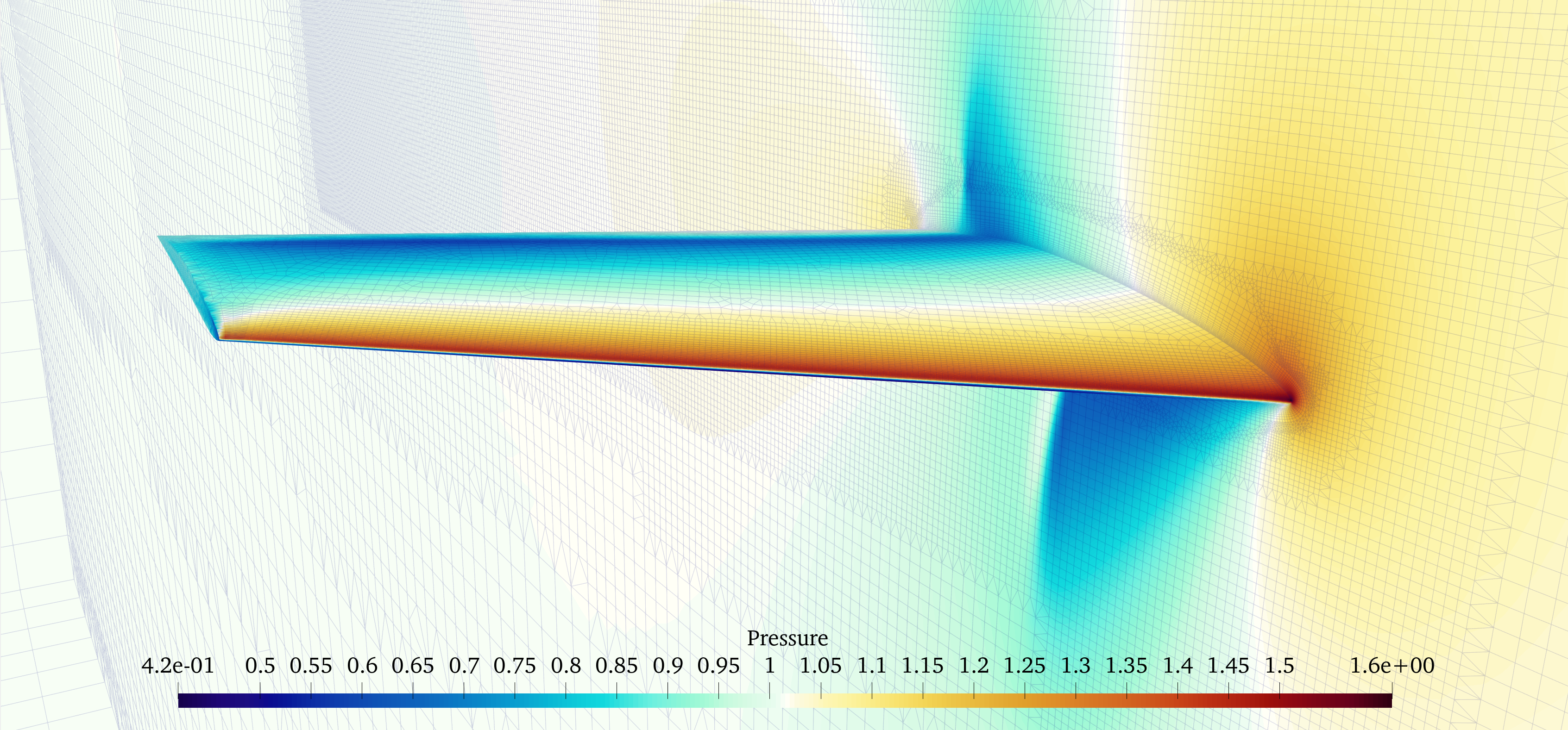}
        \caption{Pressure field over the computational domain boundaries.}
        \label{fig:rans_3d_pressure}
    \end{subfigure}
    \caption{Representative boundary-mesh views and pressure field for the 3D
transonic-wing RANS cases.}
    \label{fig:rans_3d_mesh}
\end{figure}

\paragraph{3D transonic wing cases}
The second RANS test matrix comprises three held-out ONERA~M6-based wing
geometries~\cite{schmitt1979oneraM6}, each evaluated at the same five incidence
conditions as the 2D airfoils. Incidence is encoded through the orientation of
the case-specific mesh, while the configured freestream angle remains zero.
The meshes contain $4{,}047{,}165$--$4{,}971{,}712$ grid points and
$5{,}975{,}312$--$7{,}056{,}614$ volume elements. Figure~\ref{fig:rans_3d_mesh}
shows a representative configuration and RANS solution, and
Table~\ref{tab:rans_3d_meshes} in~\ref{app:rans_3d_appendix} lists all mesh
sizes. A run terminates only when
$r_\rho<-12$ and both the drag and lift Cauchy diagnostics are below $10^{-7}$.

\begin{table}[H]
\centering
\caption{Executed pseudo-time iterations for the 3D RANS wing cases. Bold
entries identify OOD incidence conditions.}
\label{tab:rans_3d_iterations}
\begin{tabular}{@{}ccrrr@{}}
\toprule
Wing & $\alpha$ ($^\circ$) & \makecell{Uniform-flow\\iterations}
& \makecell{ADEx-FNO-based\\iterations} & $G_{\mathrm{outer}}$ (\%) \\
\midrule
1 & \textbf{-7.0} & \textbf{4,369} & \textbf{2,513} & \textbf{42.48} \\
& -2.5 & 4,429 & 3,217 & 27.37 \\
& 0.0 & 4,688 & 3,186 & 32.04 \\
& 2.5 & 4,577 & 3,060 & 33.14 \\
& \textbf{7.0} & \textbf{3,717} & \textbf{2,910} & \textbf{21.71} \\
\midrule
2 & \textbf{-7.0} & \textbf{4,383} & \textbf{1,659} & \textbf{62.15} \\
& -2.5 & 4,394 & 2,293 & 47.82 \\
& 0.0 & 4,391 & 1,605 & 63.45 \\
& 2.5 & 4,604 & 2,434 & 47.13 \\
& \textbf{7.0} & \textbf{3,781} & \textbf{2,969} & \textbf{21.48} \\
\midrule
3 & \textbf{-7.0} & \textbf{4,163} & \textbf{1,430} & \textbf{65.65} \\
& -2.5 & 4,558 & 1,652 & 63.76 \\
& 0.0 & 4,954 & 2,751 & 44.47 \\
& 2.5 & 5,305 & 2,724 & 48.65 \\
& \textbf{7.0} & \textbf{3,895} & \textbf{2,956} & \textbf{24.11} \\
\bottomrule
\end{tabular}
\end{table}

All $15$ pairs satisfy the complete stopping condition, and the ADEx-FNO-based run
requires fewer pseudo-time iterations in every case. The mean pairwise
reduction is $43.03\%$, with a range of $21.48\%$--$65.65\%$.

\begin{figure}[H]
\centering
\includegraphics[width=0.75\linewidth]{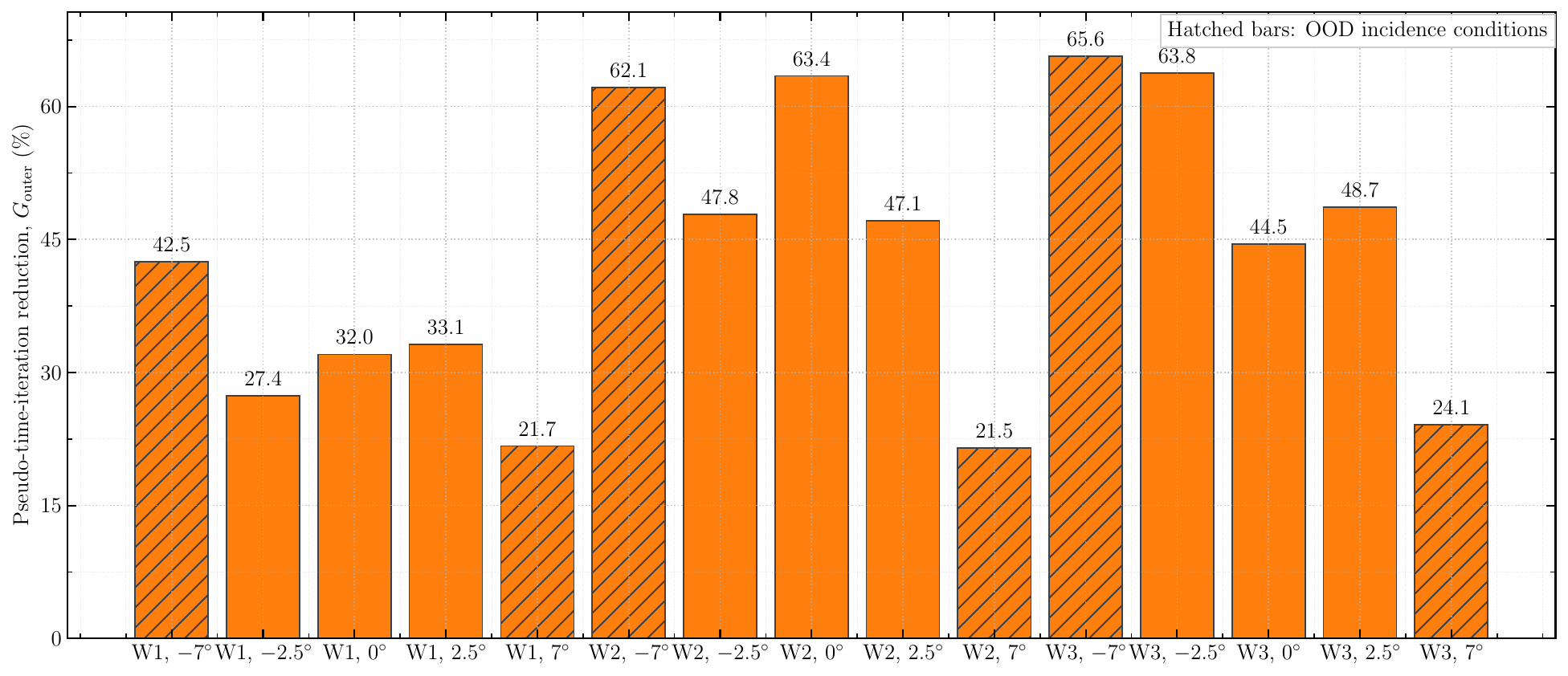}
\caption{Reduction in executed pseudo-time iterations for the complete 3D RANS
test matrix. Hatched bars denote OOD incidence conditions.}
\label{fig:rans_3d_outer_reduction}
\end{figure}

The mean reduction is $45.31\%$ for the nine ID pairs and $39.60\%$ for the
six OOD pairs; every OOD case retains a positive reduction. The wing-specific
means are $31.35\%$, $48.40\%$, and $49.33\%$ for Wings~1--3. These values are
descriptive because geometry, incidence, and baseline convergence behavior are
not independent. Grouped statistics are given in
\ref{app:rans_3d_grouped}.

The aggregate pseudo-time count decreases by $43.57\%$, and $L_\Sigma$
decreases by $43.50\%$. The mean combined linear count per pseudo-time
iteration changes only from $36.09$ to $36.14$, again showing that the
cumulative reduction is driven by fewer pseudo-time iterations.

\begin{table}[H]
\centering
\caption{Aggregate pseudo-time and linear-solver work over all $15$
converged 3D RANS pairs. $L_\Sigma$ is the unweighted sum of the mean-flow
and turbulence linear iterations.}
\label{tab:rans_3d_aggregate_work}
\begin{tabular}{@{}lrrr@{}}
\toprule
Quantity & \makecell{Uniform-flow\\initialization}
& \makecell{ADEx-FNO-based\\initialization} & Reduction (\%) \\
\midrule
Executed pseudo-time iterations & 66,208 & 37,359 & 43.57 \\
Mean-flow linear iterations & 1,324,160 & 747,180 & 43.57 \\
Turbulence linear iterations & 1,065,087 & 602,857 & 43.40 \\
Unweighted sum $L_\Sigma$ & 2,389,247 & 1,350,037 & 43.50 \\
\bottomrule
\end{tabular}
\end{table}

The gain is concentrated in the initial and intermediate convergence phases.
Although every ADEx-FNO-based run begins with a higher $r_\rho$, it reaches
$r_\rho=-10$ using $34.99\%$--$84.39\%$ fewer iterations. Beyond that crossing,
the remaining interval is shorter in ten pairs and longer in five; the data do
not support a universal improvement in the late-stage decay rate. Across all
pairs, the maximum absolute coefficient differences are
$4.1349\times10^{-5}$ for $C_L$ and $4.2101\times10^{-6}$ for $C_D$, with
maximum relative differences of $0.0382\%$ and $0.00577\%$. Detailed
convergence and coefficient data are provided in
\ref{app:rans_3d_convergence} and
\ref{app:rans_3d_coefficients}.

\subsubsection{URANS test cases}
\label{sec:urans_cases}

URANS adds a three-level work hierarchy: physical-time advancement, inner
pseudo-time iterations at each physical step, and linear iterations within
each inner solve. Let $N_{\mathrm{phys}}^{\mathrm{win}}$ be the physical-time
advances from averaging-window activation to termination,
$N_{\mathrm{inner}}^{\mathrm{win}}$ the corresponding inner iterations, and
$L_{\mathrm{flow}}^{\mathrm{win}}$ and
$L_{\mathrm{turb}}^{\mathrm{win}}$ the cumulative linear iterations. We use
\begin{equation}
    G_Q=100\left(1-\frac{Q_{\mathrm{ADEx-FNO}}}{Q_{\mathrm{UF}}}\right),
    \qquad
    L_\Sigma=L_{\mathrm{flow}}+L_{\mathrm{turb}},
    \label{eq:urans_work_reduction}
\end{equation}
with the initialization subscripts defined above. The principal 2D and 3D
settings are reported in Table~\ref{tab:urans_su2_configurations} in
\ref{app:urans_configurations}. In both cases, the ADEx-FNO is evaluated
only once, and SU2 performs all subsequent physical-time, pseudo-time, and linear
iterations.

\begin{figure}[H]
\centering
\includegraphics[width=0.75\linewidth]{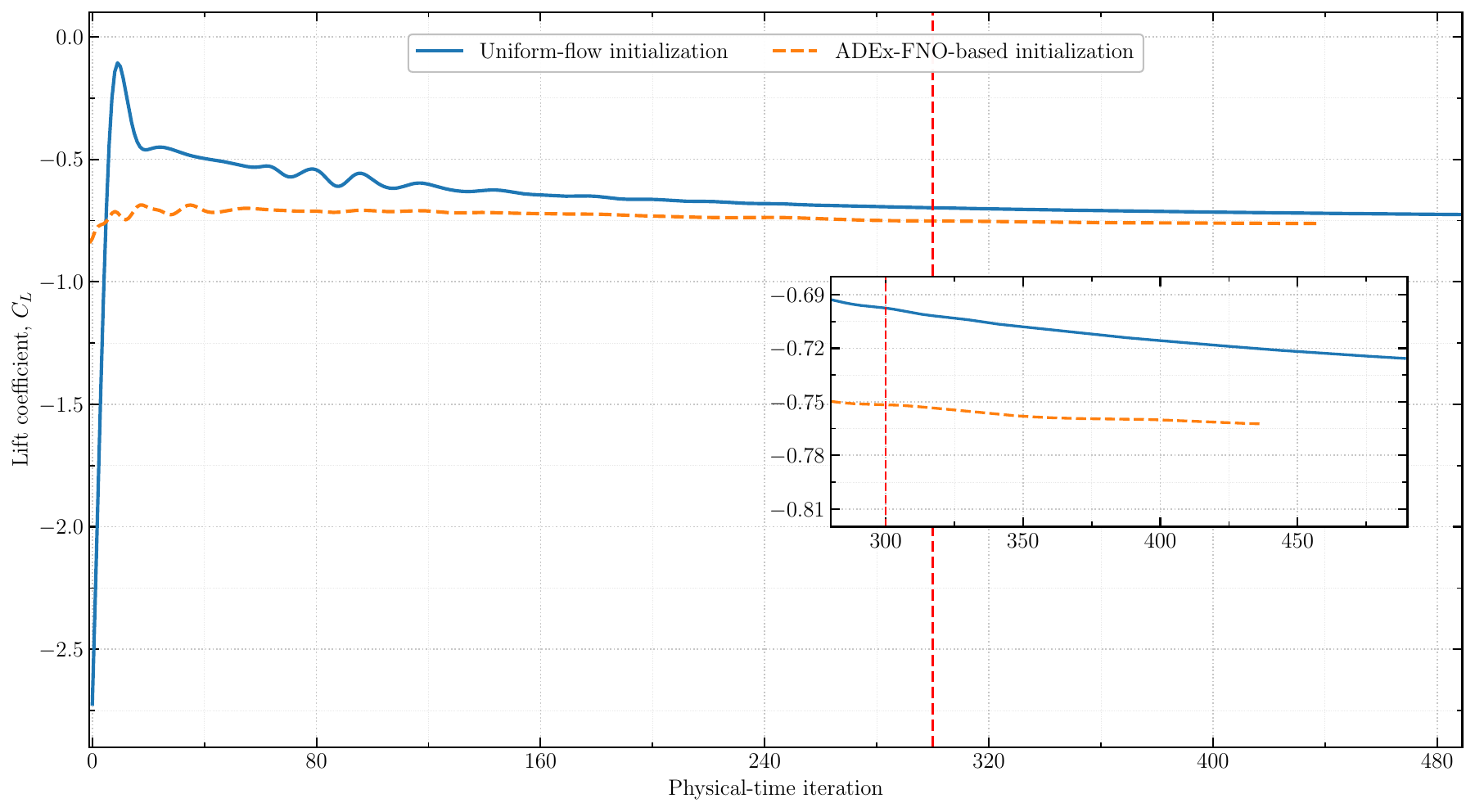}
\caption{Instantaneous lift-coefficient history for the 2D NACA~0012 URANS
case. Blue solid and orange dashed curves denote uniform-flow and ADEx-FNO-based
initialization, respectively. The vertical red dashed line marks activation of
the averaging window after $300$ physical-time iterations. Each curve ends at
its SU2 termination point. The inset zooms in on the late-time evolution.}
\label{fig:urans_2d_cl}
\end{figure}

\begin{figure}[H]
\centering
\includegraphics[width=0.75\linewidth]{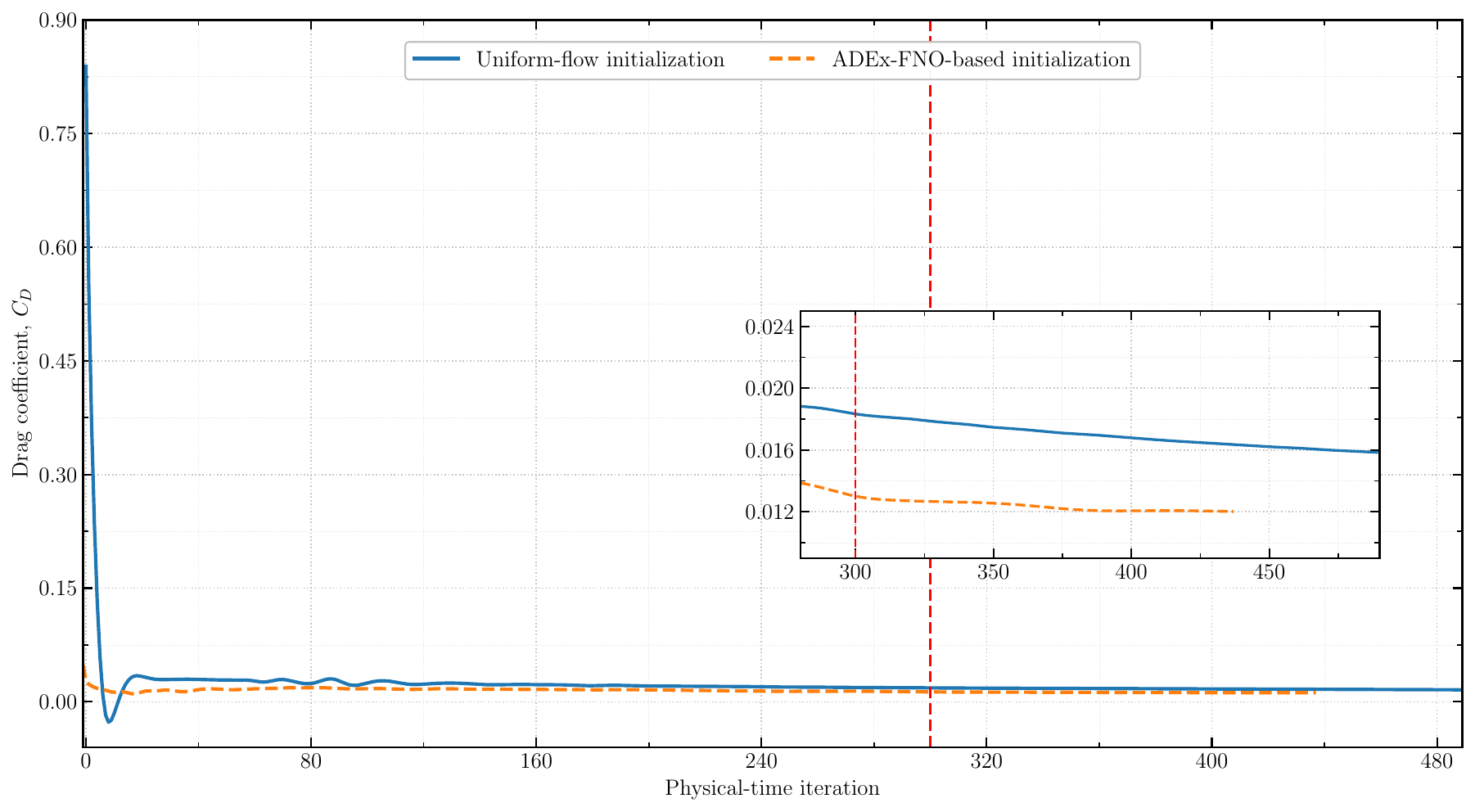}
\caption{Instantaneous drag-coefficient history for the 2D NACA~0012 URANS
case. The line styles, averaging-window marker, and curve termination have the
same meaning as in Figure~\ref{fig:urans_2d_cl}. The inset zooms in on the late-time evolution.}
\label{fig:urans_2d_cd}
\end{figure}

\paragraph{2D NACA~0012 airfoil}
The 2D case uses a mesh with $44{,}346$ grid points and $44{,}835$ volume
elements. Time convergence is assessed from squared-Hann-windowed averages of
$C_D$ and $C_L$. Averaging starts after $300$ physical-time iterations, and
both $50$-element Cauchy diagnostics must fall below $10^{-4}$. This criterion
measures stabilization of each running average, and it neither estimates error
relative to an exact mean nor enforces equality between the two simulations.
Figures~\ref{fig:urans_2d_cl} and~\ref{fig:urans_2d_cd} show the coefficient
histories, the averaging-window start, and each solver termination point.

After window activation, the uniform-flow and ADEx-FNO-based calculations require
$189$ and $137$ physical-time advances, respectively. Thus,
\begin{equation*}
    G_{N_{\mathrm{phys}}^{\mathrm{win}}}
    =100\left(1-\frac{137}{189}\right)=27.51\%.
\end{equation*}
The corresponding simulated intervals are $0.0945\,\mathrm{s}$ and
$0.0685\,\mathrm{s}$. The ADEx-FNO-based trajectory avoids $52$ advances, or
$0.0260\,\mathrm{s}$ of simulated physical time.

\begin{table}[H]
\centering
\caption{Physical-time, pseudo-time, and cumulative linear-solver work for the
2D URANS case.}
\label{tab:urans_2d_effort}
\begin{tabular}{@{}lrrr@{}}
\toprule
Quantity & \makecell{Uniform-flow\\initialization}
& \makecell{ADEx-FNO-based\\initialization} & $G_Q$ (\%) \\
\midrule
Post-window physical-time advances & 189 & 137 & 27.51 \\
Post-window simulated time ($\mathrm{s}$) & 0.0945 & 0.0685 & 27.51 \\
Post-window inner iterations & 1,649 & 1,223 & 25.83 \\
Post-window unweighted linear count $L_\Sigma$ & 56,451 & 29,703 & 47.38 \\
Complete inner iterations & 10,055 & 7,158 & 28.81 \\
Complete unweighted linear count $L_\Sigma$ & 339,275 & 180,008 & 46.94 \\
\bottomrule
\end{tabular}
\end{table}

The post-window inner count decreases by $25.83\%$, while the mean inner count
per remaining physical-time advance changes only from $8.72$ to $8.93$. The
larger $47.38\%$ reduction in $L_\Sigma$ reflects a decrease in mean combined
Krylov work per inner iteration from $34.23$ to $24.29$, driven primarily by
the turbulence system. The complete decomposition is given in
\ref{app:urans_2d_effort}.

Both calculations satisfy the configured Cauchy tolerance, with the lift
diagnostic controlling termination. Their terminal squared-Hann averages are
\begin{equation*}
\begin{aligned}
    \overline{C_L}^{\mathrm{UF}}&=-0.711,
    &\qquad \overline{C_D}^{\mathrm{UF}}&=0.0168,\\
    \overline{C_L}^{\mathrm{ADEx-FNO}}&=-0.753,
    &\overline{C_D}^{\mathrm{ADEx-FNO}}&=0.0122.
\end{aligned}
\end{equation*}
The relative differences are $5.97\%$ for lift and $27.12\%$ for drag.
Accordingly, this case demonstrates less SU2 work to reach the configured
self-convergence condition, but not convergence to a verified common
long-time-averaged state. 
The stopping diagnostics and
terminal statistics are reported in~\ref{app:urans_2d_window} and
\ref{app:urans_2d_statistics}.

\begin{figure}[H]
\centering
\includegraphics[width=0.75\linewidth]{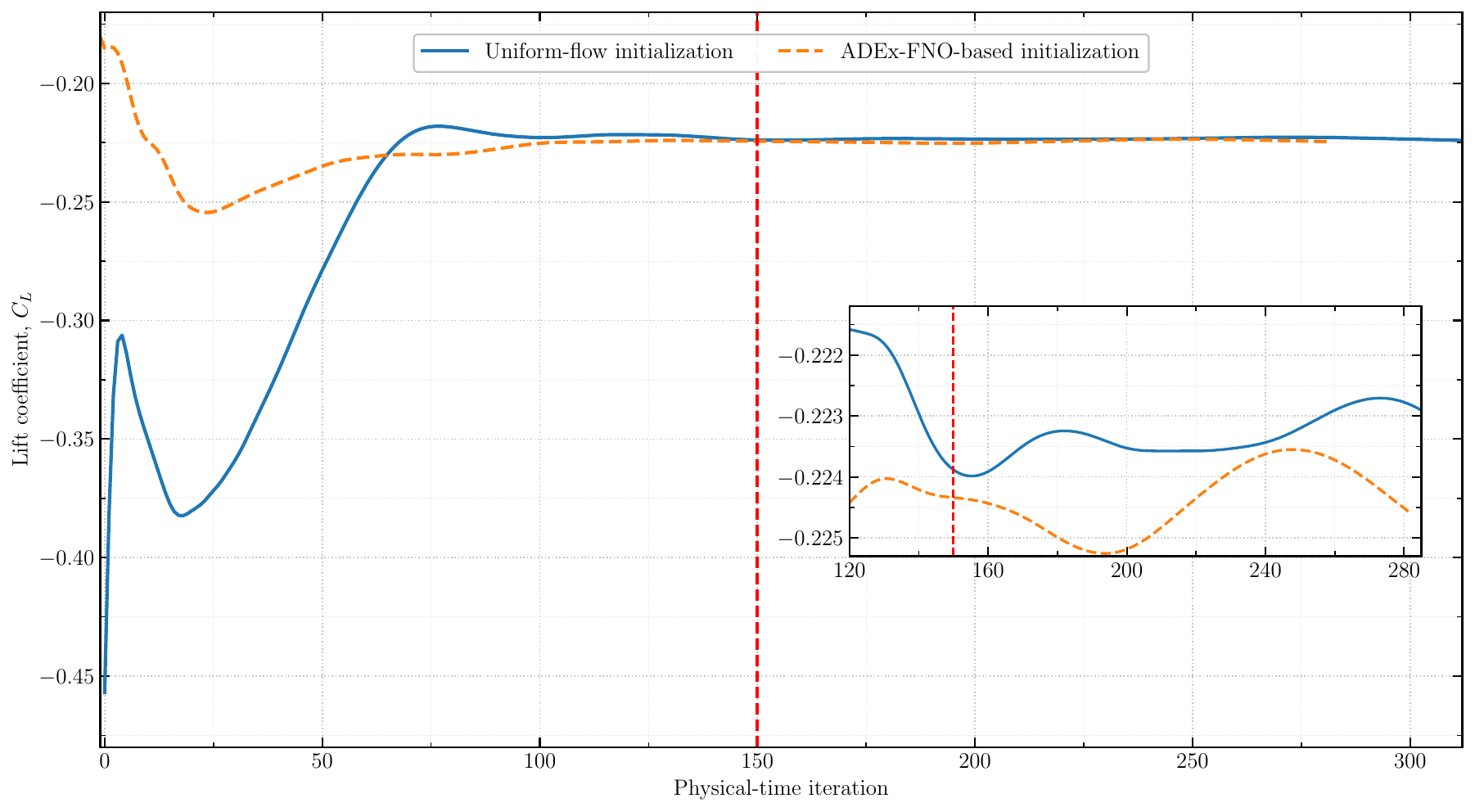}
\caption{Instantaneous lift-coefficient history for the 3D transonic-wing
URANS case. Blue solid and orange dashed curves denote uniform-flow and
ADEx-FNO-based initialization, respectively. The vertical red dashed line marks
activation of the averaging window after $150$ physical-time iterations. Each
curve ends at its SU2 termination point. The inset zooms in on the late-time evolution.}
\label{fig:urans_3d_cl}
\end{figure}

\begin{figure}[H]
\centering
\includegraphics[width=0.75\linewidth]{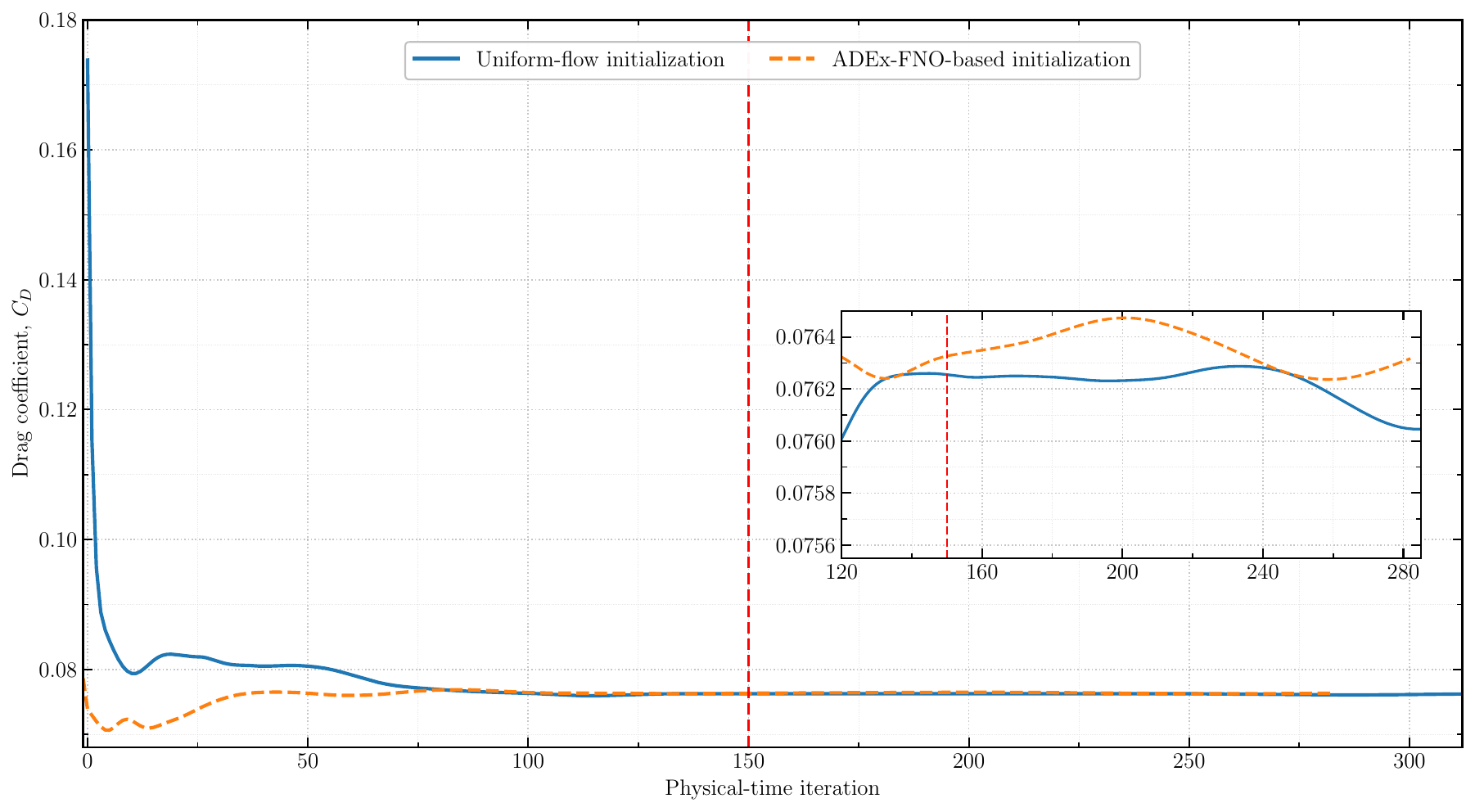}
\caption{Instantaneous drag-coefficient history for the 3D transonic-wing
URANS case. The line styles, averaging-window marker, and curve termination
have the same meaning as in Figure~\ref{fig:urans_3d_cl}. The inset resolves
the late-time evolution.}
\label{fig:urans_3d_cd}
\end{figure}

\paragraph{3D transonic wing}
The 3D case uses the held-out Wing~3 geometry at the OOD incidence
$\alpha=-7^\circ$ on a mesh with $4{,}052{,}234$ grid points and
$5{,}975{,}820$ volume elements. Averaging starts after $150$ physical-time
iterations, and both $50$-element Cauchy diagnostics must fall below $10^{-5}$.
Figures~\ref{fig:urans_3d_cl} and~\ref{fig:urans_3d_cd} show the coefficient
histories, window activation, and termination.

After window activation, the uniform-flow and ADEx-FNO-based calculations require
$162$ and $132$ physical-time advances, giving
\begin{equation*}
    G_{N_{\mathrm{phys}}^{\mathrm{win}}}
    =100\left(1-\frac{132}{162}\right)=18.52\%.
\end{equation*}
The corresponding simulated intervals are $0.0810\,\mathrm{s}$ and
$0.0660\,\mathrm{s}$. The ADEx-FNO-based trajectory avoids $30$ advances, or
$0.0150\,\mathrm{s}$.

\begin{table}[H]
\centering
\caption{Physical-time, pseudo-time, and cumulative linear-solver work for the
3D URANS case.}
\label{tab:urans_3d_effort}
\begin{tabular}{@{}lrrr@{}}
\toprule
Quantity & \makecell{Uniform-flow\\initialization}
& \makecell{ADEx-FNO-based\\initialization} & $G_Q$ (\%) \\
\midrule
Post-window physical-time advances & 162 & 132 & 18.52 \\
Post-window simulated time ($\mathrm{s}$) & 0.0810 & 0.0660 & 18.52 \\
Post-window inner iterations & 2,382 & 1,710 & 28.21 \\
Post-window mean-flow linear iterations & 46,169 & 32,245 & 30.16 \\
Post-window turbulence linear iterations & 46,771 & 28,231 & 39.64 \\
Post-window unweighted linear count $L_\Sigma$ & 92,940 & 60,476 & 34.93 \\
Complete inner iterations & 7,860 & 5,409 & 31.18 \\
Complete unweighted linear count $L_\Sigma$ & 311,407 & 199,182 & 36.04 \\
\bottomrule
\end{tabular}
\end{table}

The post-window inner count decreases by $28.21\%$, while the mean inner count
per remaining physical-time advance decreases from $14.70$ to $12.95$. The
mean combined linear count per inner iteration also decreases from $39.02$ to
$35.37$, producing a $34.93\%$ reduction in post-window $L_\Sigma$. The
complete decomposition is reported in~\ref{app:urans_3d_effort}.

Both calculations satisfy the same coefficient-based stopping condition, with
the lift diagnostic controlling termination. Their terminal squared-Hann
averages are
\begin{equation*}
\begin{aligned}
    \overline{C_L}^{\mathrm{UF}}&=-0.222,
    &\qquad \overline{C_D}^{\mathrm{UF}}&=0.0758,\\
    \overline{C_L}^{\mathrm{ADEx-FNO}}&=-0.223,
    &\overline{C_D}^{\mathrm{ADEx-FNO}}&=0.0758.
\end{aligned}
\end{equation*}
The relative differences are $0.360\%$ for lift and $0.072\%$ for drag. Thus,
this case combines reductions in post-window physical-time advances
($18.52\%$), inner iterations ($28.21\%$), and combined linear work
($34.93\%$) with sub-percent agreement of the terminal windowed statistics.
Detailed stopping and coefficient data are given in
\ref{app:urans_3d_window} and~\ref{app:urans_3d_statistics}.

\subsection{Direct numerical simulation bootstrap test}
\label{sec:dns_bootstrap}

The RANS and URANS studies quantify reductions in iterative solver work. For DNS, the relevant question is whether the ADEx-FNO field shortens the initial flow-development, or bootstrap, interval that must be discarded before useful statistical sampling. The paired calculations consider the compressible flow around a 3D NACA~0012 wing at
\begin{equation*}
    M_\infty=0.4,
    \qquad
    Re_\infty=\num{5e4},
    \qquad
    \alpha=\SI{5}{\degree}.
\end{equation*}
The calculations use the entropy-stable discontinuous collocated Galerkin (SSDC) solver presented in~\cite{ParsaniEtAl2021SSDC}. The physical and numerical configuration and reference flow data are based on
\cite{karp2026EffectsOfLowerFloatingPointPrecision,jones2008naca0012,sandberg2011naca0012}.
Both simulations use the same geometry, $504{,}896$-cell mesh, degree-seven discontinuous Galerkin discretization, boundary conditions, time integration, and physical parameters, and are advanced by the DNS solver to $t^*=10$. The only intended difference is the initial field.

This is a demanding transfer test: the ADEx-FNO was trained on 2D URANS data at $M_\infty=0.15$, $Re_\infty=\num{6e6}$, and $\alpha\in[-5^\circ,5^\circ]$, whereas the target is a 3D DNS at substantially different Mach and Reynolds numbers. The ADEx-FNO cannot provide the resolved 3D turbulent fluctuations. It can only supply a large-scale pressure and velocity field from which the DNS must develop the physical instabilities and smaller scales. 

Time is expressed in convective time units (CTU),
\begin{equation}
t^*=\frac{tU_\infty}{c},
\label{eq}
\end{equation}
where $c$ denotes the chord length of the NACA 0012 section used to construct the wing. We assess three progressively broader notions of readiness for statistical sampling: stabilization of the integral forces, development of the local skin-friction distribution, and removal of initialization-generated structures from the computational domain. These criteria correspond to different statistical objectives and are therefore reported separately. This distinction is important because DNS startup and statistical sampling can require comparable computational resources. The reference calculation reported in \cite{karp2026EffectsOfLowerFloatingPointPrecision} required nearly $4\times10^6$ core-hours on Shaheen~III CPU at KAUST.

\begin{figure}[H]
    \centering
    \includegraphics[width=0.75\textwidth]{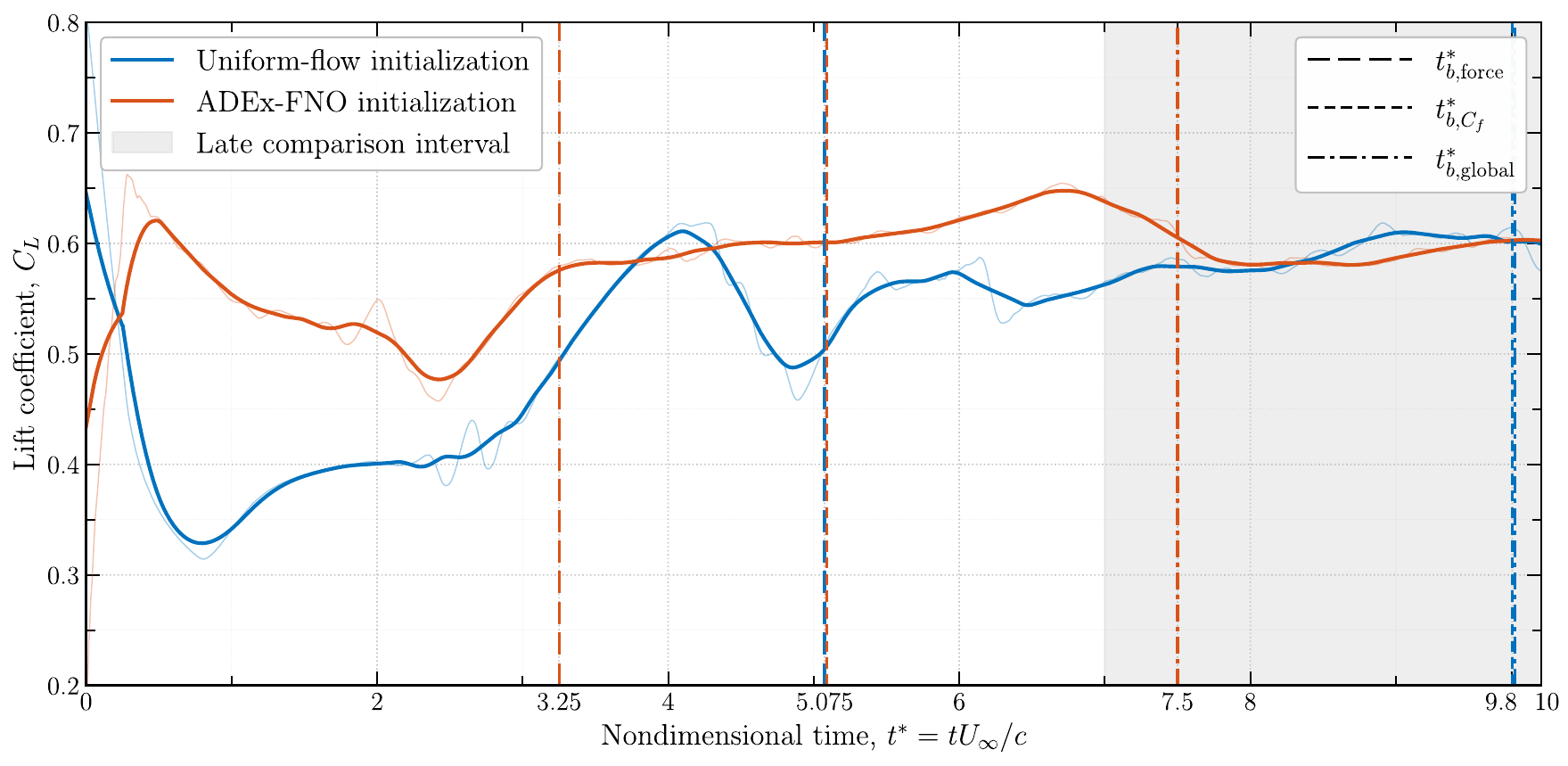}
    \caption{DNS lift-coefficient histories for uniform-flow and ADEx-FNO initialization. Thin curves show the instantaneous values, and thick blue solid and orange dashed curves show moving means computed over $0.5$ CTU. The dotted vertical lines mark the central MSER bootstrap estimates. The gray band denotes the common late-time comparison interval $7\leq t^*\leq10$.}
    \label{fig:dns_cl_history}
\end{figure}

\begin{figure}[H]
    \centering
    \includegraphics[width=0.75\linewidth]{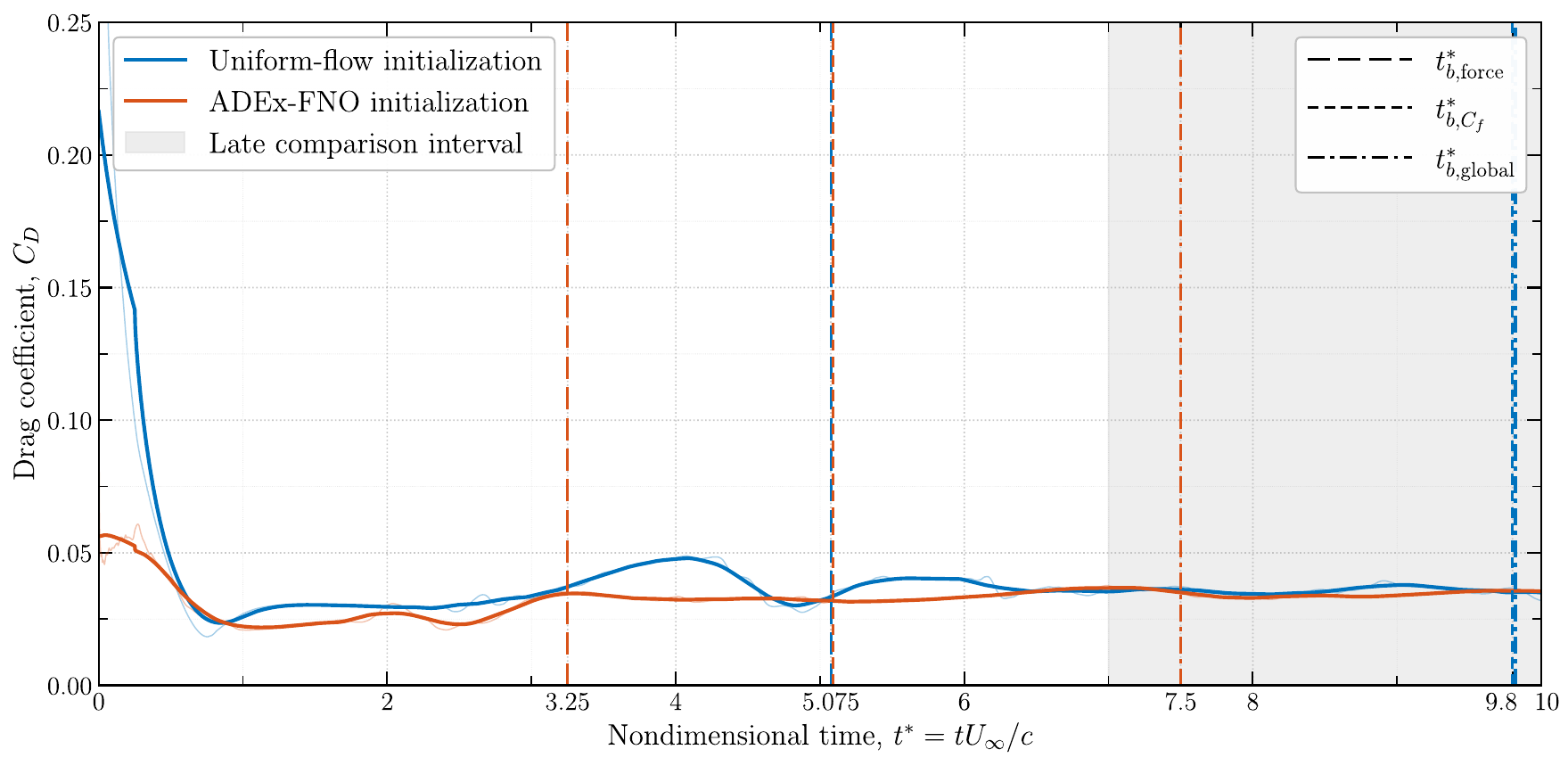}
    \caption{DNS drag-coefficient histories for uniform-flow and ADEx-FNO initialization. Colors, line styles, MSER markers, and the late-time comparison interval have the same meanings as in Figure~\ref{fig:dns_cl_history}.}
    \label{fig:dns_cd_history}
\end{figure}

\paragraph{Integral-force criterion}
The lift and drag histories are analyzed with the marginal standard error rule (MSER) applied to common block-averaged signals \cite{MockettEtAl2010,BergmannEtAl2022}. For each block duration, MSER gives a truncation time for $C_L$ and $C_D$. The later value defines the bootstrap time for that calculation. Identical block durations, candidate intervals, and objectives are used for both initializations. Figures~\ref{fig:dns_cl_history} and~\ref{fig:dns_cd_history} show the instantaneous histories, $0.5$-CTU moving means, and central bootstrap estimates.

For a quantity of interest $q$, we define
\begin{equation}
    G_{b,q}
    =100\left(
        1-\frac{t_{b,q,\mathrm{ADEx-FNO}}^*}
                 {t_{b,q,\mathrm{uniform}}^*}
    \right).
    \label{eq:dns_bootstrap_reduction}
\end{equation}
Six common block durations, $0.02\leq B\leq0.50$, are used, and the central estimate is the median over the six analyses. The resulting force-based times are $5.075$ CTU for uniform-flow initialization and $3.250$ CTU for ADEx-FNO initialization. The saving is $1.825$ CTU, corresponding to
\begin{equation*}
    G_{b,\mathrm{force}}
    =100\left(1-\frac{3.250}{5.075}\right)
    =35.96\%.
\end{equation*}
The reduction remains between $34.15\%$ and $38.10\%$ across the tested block durations.

\begin{figure}[H]
    \centering
    \includegraphics[width=0.70\linewidth]{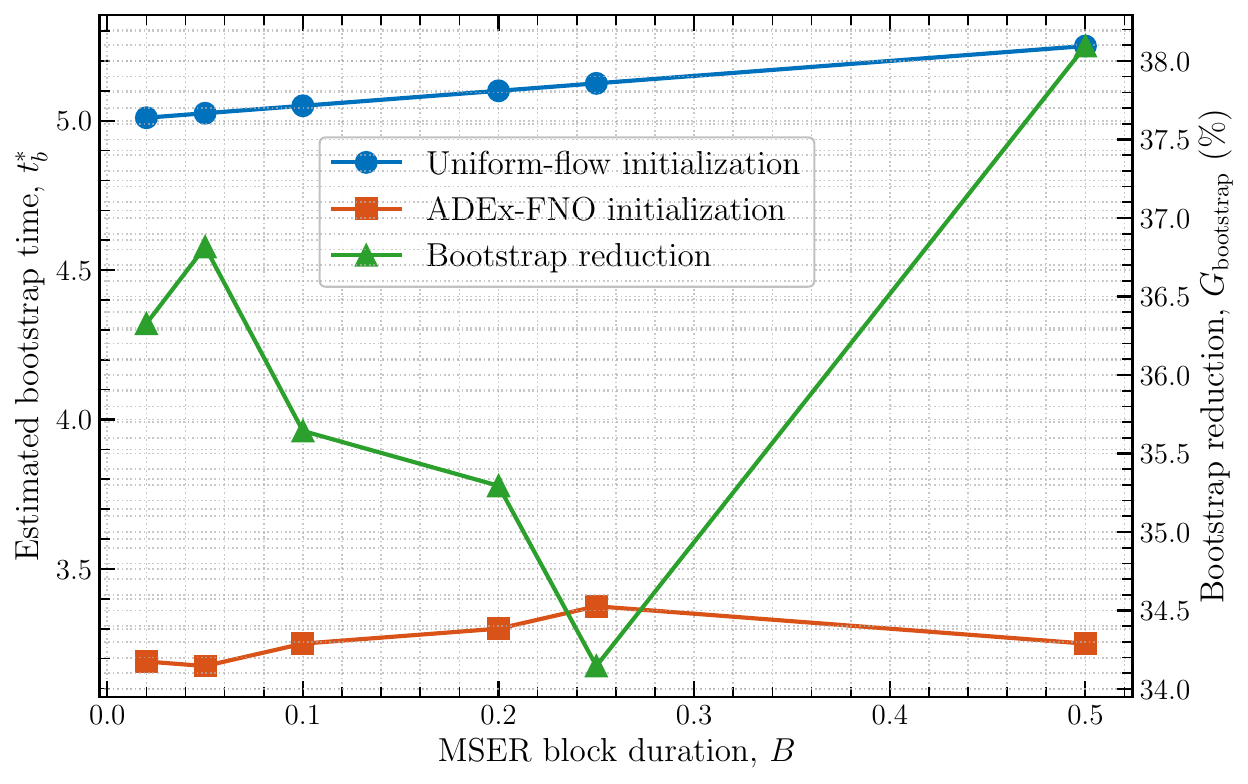}
    \caption{Sensitivity of the force-based DNS bootstrap estimate to the common MSER block duration $B$. Blue circles and orange squares denote the uniform-flow and ADEx-FNO bootstrap times, respectively, and green triangles denote the reduction defined by Eq.~\eqref{eq:dns_bootstrap_reduction}. The ADEx-FNO estimate remains below the uniform-flow estimate for every tested block duration.}
    \label{fig:dns_mser_sensitivity}
\end{figure}

The corresponding adaptive explicit-step indices decrease from $1{,}243{,}957$ to $829{,}132$, avoiding $414{,}825$ steps, or $33.35\%$. The difference between the CTU- and step-based reductions reflects the distinct adaptive time-step histories.

\begin{figure}[H]
    \centering
    \includegraphics[width=\linewidth]{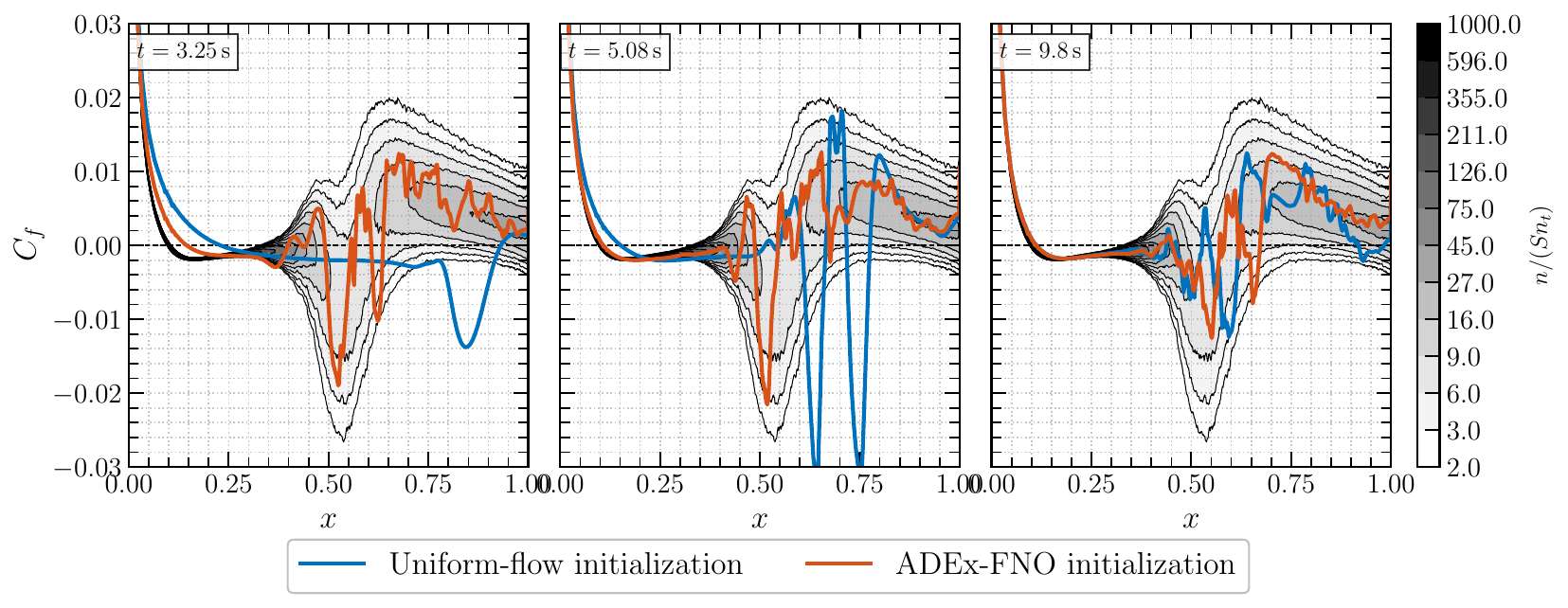}
    \caption{DNS skin-friction coefficient at $t^*=3.25$, $5.08$, and $9.80$ for uniform-flow initialization (blue) and ADEx-FNO initialization (orange). The grayscale background and black contours show the reference probability-density distribution of $C_f$; the dashed horizontal line marks $C_f=0$. The comparison assesses the development of the laminar separation region, transition, and turbulent reattachment.}
    \label{fig:Cf_DNS}
\end{figure}

\paragraph{Local skin-friction criterion}
The integral loads do not establish that the separated boundary layer is ready for surface statistics. Figure~\ref{fig:Cf_DNS} therefore compares instantaneous skin-friction profiles with the reference probability-density distribution of $C_f$. The $C_f=0$ crossings identify the instantaneous separation and reattachment locations. At $t^*=3.25$, the ADEx-FNO profile already exhibits developed trailing-edge fluctuations, while the uniform-flow profile remains largely laminar. By $t^*=5.075$ (shown as $5.08$), the ADEx-FNO profile is consistent with the developed distribution over the chord; the first displayed uniform-flow profile with comparable behavior occurs at $t^*=9.80$. Therefore, the snapshot-resolved estimates are
\begin{equation*}
    t_{b,C_f,\mathrm{UF}}^*=9.80,
    \qquad
    t_{b,C_f,\mathrm{ADEx-FNO}}^*=5.075,
\end{equation*}
which gives a saving of $4.725$ CTU and
\begin{equation*}
    G_{b,C_f}
    =100\left(1-\frac{5.075}{9.80}\right)
    =48.21\%.
\end{equation*}

\paragraph{Full-domain flow-field criterion}
Near-body development is not sufficient when statistics are required throughout the domain. The $Q$-criterion sequence in~\ref{app:dns_bootstrap} shows earlier formation of 3D wake structures under ADEx-FNO initialization, while the spanwise-vorticity sequence tracks the downstream convection of the initial structures. At $t^*=5.075$, those structures remain in the ADEx-FNO far field; they have left the displayed domain by $t^*=7.50$. Comparable flushing for uniform-flow initialization occurs at approximately $t^*=9.80$. Thus,
\begin{equation*}
    t_{b,\mathrm{global},\mathrm{UF}}^*=9.80,
    \qquad
    t_{b,\mathrm{global},\mathrm{ADEx-FNO}}^*=7.50,
\end{equation*}
with a saving of $2.30$ CTU and
\begin{equation*}
    G_{b,\mathrm{global}}
    =100\left(1-\frac{7.50}{9.80}\right)
    =23.47\%.
\end{equation*}
This flow-field estimate is more restrictive and remains based on the displayed snapshots.

\begin{table}[H]
\centering
\caption{DNS bootstrap estimates for the three quantities of interest. Each reduction is normalized by the corresponding uniform-flow bootstrap time.}
\label{tab:dns_bootstrap_criteria}
\begin{tabular}{@{}lrrrr@{}}
\toprule
Criterion
& \makecell{Uniform-flow\\$t_b^*$}
& \makecell{ADEx-FNO\\$t_b^*$}
& \makecell{Saving\\(CTU)}
& Reduction (\%) \\
\midrule
Integral forces ($C_L$, $C_D$) & 5.075 & 3.250 & 1.825 & 35.96 \\
Local skin friction ($C_f$)    & 9.80  & 5.075 & 4.725 & 48.21 \\
Full-domain wake flushing      & 9.80  & 7.50  & 2.30  & 23.47 \\
\bottomrule
\end{tabular}
\end{table}

The criteria should not be collapsed into one percentage because they address different statistical objectives. Nevertheless, each criterion yields a shorter bootstrap interval for the ADEx-FNO-initialized simulation. Over the common late-time interval $7\leq t^*\leq10$, the mean lift coefficients are $0.590544$ and $0.594988$, respectively, for the uniform-flow and ADEx-FNO initializations, while the corresponding mean drag coefficients are $0.035878$ and $0.034582$. The standard deviations differ by $4.19\%$ for $C_L$ and $5.29\%$ for $C_D$. Correlation-aware $95\%$ confidence intervals for the mean differences include zero~\cite{OliverEtAl2014}; however, only approximately $3$--$7$ effectively independent force samples are available, and the final block means continue to exhibit drift. Correlation-aware $95\%$ confidence intervals for the mean differences include zero~\cite{OliverEtAl2014}. However, only approximately $3$--$7$ effectively independent force samples are available, and the final block means continue to exhibit drift.

Overall, the supplied ADEx-FNO field reduces the bootstrap interval under all three definitions: $35.96\%$ for integral-force readiness, $48.21\%$ for the snapshot-based local $C_f$ criterion, and $23.47\%$ for the removal of initialization-generated wake structures from the full computational domain. These results indicate that the transferred field provides a more developed large-scale initial state despite the mismatch between the URANS training regime and the target DNS. These results indicate that the transferred field provides a more developed large-scale initial state despite the mismatch between the URANS training regime and the target DNS. The supporting right-hand-side and flow-field diagnostics are reported in \ref{app:dns_bootstrap}.

\section{Conclusions}
\label{sec:conclusions}

We introduced the Ambient-Domain Extension Fourier Neural Operator
(ADEx-FNO), a unified and deterministic framework for incorporating varying
geometries into an FNO pipeline without replacing its defining
Fourier-operator architecture. Each admissible physical domain is embedded in
a fixed ambient hypercube and represented through its signed distance
function. The problem data and solution fields are extended to the ambient
domain through prescribed operations and transferred to a common, potentially
nonuniform, rectilinear latent grid. After inference, the predicted field can
be interpolated to an independently selected target discretization and
restricted to the physical domain. These geometry-transfer operations remain
outside the optimization loop and introduce no trainable graph, point-cloud,
deformation, or geometry-decoding modules beyond the underlying FNO.

The nonlinear Poisson and advection-reaction-diffusion experiments
demonstrate that the same construction can be applied across different
equations and spatial dimensions. On the held-out smooth-domain test sets,
ADEx-FNO attains relative $\ell^2$ errors between $0.32\%$ and $0.77\%$.
Models trained exclusively on smooth geometries were also evaluated on unseen
domains containing corners, nonconvex boundaries, and shapes substantially
different from those used during training. These tests demonstrate the
flexibility of the common ambient-domain pipeline, while also showing that
accuracy under geometric extrapolation depends on the departure from the
training distribution.

The CFD studies examine a complementary use of ADEx-FNO: one-shot
initialization of established governing-equation solvers. ADEX-FNO supplies
only the initial flow field and does not modify the discrete equations,
convergence criteria, time integration, or linear solvers. Across all $29$
converged 2D and 3D RANS cases, the ADEx-FNO-based
initialization reduces the pseudo-time count in every pair, with mean
reductions of $44.17\%$ and $43.03\%$, respectively. The corresponding
reductions in cumulative linear-solver work are of similar magnitude, showing
that the benefit arises primarily from placing the nonlinear solver on a more
favorable convergence trajectory. Without retraining, the same model also
retains pseudo-time reductions of $51.10\%$--$53.13\%$ across three target
mesh resolutions. On each fixed mesh, the paired calculations converge to
consistent aerodynamic coefficients.

The URANS results demonstrate both the potential of the initialization and the
need to distinguish self-convergence from agreement of the resulting
time-averaged state. In the 2D case, ADEx-FNO reduces the
post-window physical-time advances by $27.51\%$ and the cumulative
post-window linear work by $47.38\%$. However, the terminal windowed
coefficients differ appreciably, so this result establishes a reduction in
the work required to satisfy the configured stopping condition, not
acceleration toward a verified common long-time-averaged state. In the
3D case, the post-window physical-time advances, inner
iterations, and cumulative linear work decrease by $18.52\%$, $28.21\%$,
and $34.93\%$, respectively, while the terminal windowed lift and drag
coefficients differ by only $0.360\%$ and $0.072\%$. This comparison therefore
combines a reduction in solver work with sub-percent agreement of the monitored
late-time statistics.

The DNS experiment provides the most stringent transfer test. An ADEx-FNO
trained using 2D RANS data is used to initialize a DNS at substantially different Mach and Reynolds numbers.
The inferred field cannot contain the resolved 3D turbulence;
it can only provide a more developed large-scale state from which the physical
instabilities and smaller scales must emerge under the DNS solver. Relative to
uniform-flow initialization, the estimated bootstrap interval is reduced by
$35.96\%$ according to the integral-force criterion, by $48.21\%$ according
to the local skin-friction criterion, and by $23.47\%$ according to the more
restrictive full-domain wake-flushing criterion. These values are reported
separately because they correspond to different statistical objectives. The
late-time force levels are compatible within the uncertainty of the available
records, but the integration interval remains too short to establish tightly
converged equality of all mean fields, Reynolds stresses, probability
distributions, spectra, or acoustic statistics.

Overall, ADEx-FNO provides a deterministic interface between Fourier operator
learning, varying physical geometries, and independently selected numerical
discretizations. Its role in the CFD studies is not to replace the
governing-equation solver, but to provide a physically informed starting field
from which that solver retains complete control. The reported gains are
therefore expressed in pseudo-time iterations, cumulative linear-solver work,
simulated physical-time advancement, and bootstrap duration, rather than as
unqualified end-to-end wall-clock speedups. 

\section*{CRediT authorship contribution statement}\noindent
\textbf{Roberto Nuca:} Conceptualization, Formal analysis, Investigation, Methodology, Supervision, Coding, Visualization, Writing - original draft, Writing - review \& editing.\\
\textbf{Giovanni Testa:} Formal analysis, Investigation, Methodology, Coding, Software, Data curation, Validation, Visualization, Writing - original draft, Writing - review \& editing.\\
\textbf{Luca Galimberti:} Data curation, Validation, Visualization, Writing - review \& editing.\\
\textbf{Matteo Parsani:} Methodology, Investigation, Supervision, Visualization, Writing - original draft, Writing - review \& editing, Funding acquisition.

\section*{Acknowledgments}
This work was supported by King Abdullah University of Science and Technology (KAUST) through grant no. BAS/1/1663-01-01. The authors also gratefully acknowledge the KAUST Supercomputing Laboratory for providing computational resources on the CPU partition of the Shaheen III supercomputer and the Ibex cluster. Roberto Nuca is a member of the GNCS/INdAM research group.

\section*{Conflict of interest}
The authors declare that they have no conflict of interest.

\newpage

\bibliographystyle{elsarticle-num}
\bibliography{references}

\newpage

\appendix

\section{Physical and numerical parameters of the RANS test cases}
\label{app:rans_configurations}

Table~\ref{tab:rans_su2_configurations} collects the physical, spatial,
implicit-solution, and convergence settings used in the 2D and 3D RANS
studies. Here, JST denotes the Jameson--Schmidt--Turkel scheme,
MUSCL denotes the monotonic upstream-centered scheme for conservation laws,
ILU denotes incomplete lower-upper factorization, FGMRES denotes the flexible
generalized minimal residual method, and CFL denotes the
Courant--Friedrichs--Lewy number.

\begin{table}[H]
\centering
\caption{Principal physical and numerical parameters of the 2D and 3D RANS
cases.}
\label{tab:rans_su2_configurations}
\renewcommand{\arraystretch}{1.06}
\setlength{\tabcolsep}{4pt}
\begin{tabularx}{\linewidth}{@{}
    >{\raggedright\arraybackslash}p{0.26\linewidth}
    >{\raggedright\arraybackslash}X
    >{\raggedright\arraybackslash}X
@{}}
\toprule
Parameter & \textbf{2D RANS} & \textbf{3D RANS} \\
\midrule

Reference configuration
& NACA~0012-based airfoil family
& ONERA~M6-based wing family \\

Flow regime
& Low-Mach-number external flow
& Transonic external flow \\

Mach number
& $M_\infty=0.15$
& $M_\infty=0.8395$ \\

Reynolds number
& $Re=6.0\times10^{6}$
& $Re=1.5\times10^{5}$ \\

Reference length
& $L_{\mathrm{ref}}=1$
& $L_{\mathrm{ref}}=0.64607$ \\

Freestream temperature
& $T_\infty=300\,\mathrm{K}$
& $T_\infty=288.15\,\mathrm{K}$ \\

Boundary conditions
& Adiabatic wall and far field
& Adiabatic wall, far field, and symmetry plane \\

\addlinespace[2pt]

Turbulence closure
& Spalart--Allmaras with \texttt{NEGATIVE} and
  \texttt{EXPERIMENTAL}
& Spalart--Allmaras-\texttt{noft2} \\

Mean-flow discretization
& Roe flux with MUSCL reconstruction
& JST scheme with coefficients $(0.5,0.02)$ \\

Turbulence discretization
& Scalar-upwind flux with MUSCL reconstruction
& Scalar-upwind flux without MUSCL reconstruction \\

Gradient evaluation
& Weighted least squares
& Green--Gauss \\

Pseudo-time integration
& Implicit Euler
& Implicit Euler \\

Linear solver and preconditioner
& ILU-preconditioned FGMRES
& ILU(0)-preconditioned FGMRES \\

Linear-solver tolerance
& $10^{-10}$
& $10^{-10}$ \\

Maximum linear iterations
& $75$
& $20$ \\

CFL strategy
& Adaptive, with initial value $100$ and bounds $[0.1,100]$
& Adaptive, with initial value $100$ and bounds $[10,200]$ \\

\addlinespace[2pt]

Monitored convergence quantities
& RMS density
& RMS density, drag, and lift \\

Residual threshold
& $r_\rho=\log_{10}(\mathrm{RMS}_\rho)<-12$
& $r_\rho<-12$ \\

Coefficient-Cauchy criterion
& Not used
& Drag and lift below $10^{-7}$ over $100$ entries \\

Monitoring start
& Pseudo-time iteration $10$
& Pseudo-time iteration $10$ \\

Maximum pseudo-time iterations
& $29{,}999$
& $29{,}999$ \\

\bottomrule
\end{tabularx}
\end{table}

\section{Details of the 2D RANS initialization study}
\label{app:rans_2d_appendix}

\subsection{Counting convention and convergence definition}

The executed pseudo-time count is the number of nonlinear iterations performed
by SU2. Because the solver index is zero based, a terminal index $k$
corresponds to $k+1$ executed iterations. The cumulative mean-flow and
turbulence linear counts are denoted by $L_{\mathrm{flow}}$ and
$L_{\mathrm{turb}}$, and their unweighted sum is $L_\Sigma$.

Convergence monitoring begins after pseudo-time iteration $10$, and a 2D RANS
run is converged when $r_\rho=\log_{10}(\mathrm{RMS}_\rho)$ first falls below
$-12$. Fourteen paired cases satisfy this criterion under both
initializations. Airfoil~2 at $\alpha=0^\circ$ reaches the iteration limit under
both initializations and is retained explicitly as a nonconverged case.

\subsection{Detailed linear-solver effort}
\label{app:rans_2d_linear}

Table~\ref{tab:rans_2d_linear_totals} reports the cumulative linear iterations
for each pair. For a converged pair, the reductions are computed with
Eq.~\eqref{eq:rans_work_reduction}. For Airfoil~2 at $\alpha=0^\circ$, the
counts are the work accumulated up to the imposed limit; no reduction is
assigned because neither run converged. The combined count $L_\Sigma$ is an
unweighted diagnostic and should not be interpreted as an exact cost-equivalent
sum.

\begin{table}[H]
\centering
\caption{Per-case cumulative linear-solver iterations for the 2D RANS
airfoil cases. Superscripts UF and ADEx-FNO denote uniform-flow and ADEx-FNO-based
initialization. The nonconverged pair is retained without an assigned
reduction.}
\label{tab:rans_2d_linear_totals}
\resizebox{\textwidth}{!}{%
\begin{tabular}{@{}ccrrcrrcrrc@{}}
    \toprule
    & & \multicolumn{3}{c}{Mean-flow system}
      & \multicolumn{3}{c}{Turbulence system}
      & \multicolumn{3}{c}{Unweighted combined count} \\
    \cmidrule(lr){3-5}\cmidrule(lr){6-8}\cmidrule(lr){9-11}
    Airfoil & $\alpha$ & $L_{\mathrm{flow}}^{\mathrm{UF}}$
    & $L_{\mathrm{flow}}^{\mathrm{ADEx-FNO}}$ & $G_{\mathrm{flow}}$ (\%)
    & $L_{\mathrm{turb}}^{\mathrm{UF}}$
    & $L_{\mathrm{turb}}^{\mathrm{ADEx-FNO}}$ & $G_{\mathrm{turb}}$ (\%)
    & $L_\Sigma^{\mathrm{UF}}$ & $L_\Sigma^{\mathrm{ADEx-FNO}}$
    & $G_\Sigma$ (\%) \\
    \midrule
    1 & -7.0 & 624,939 & 372,675 & 40.37 & 235,204 & 139,898 & 40.52 & 860,143 & 512,573 & 40.41 \\
    1 & -2.5 & 449,080 & 299,250 & 33.36 & 154,977 & 104,131 & 32.81 & 604,057 & 403,381 & 33.22 \\
    1 & 0.0 & 438,269 & 298,500 & 31.89 & 152,890 & 107,700 & 29.56 & 591,159 & 406,200 & 31.29 \\
    1 & 2.5 & 454,548 & 253,960 & 44.13 & 156,135 & 88,393 & 43.39 & 610,683 & 342,353 & 43.94 \\
    1 & 7.0 & 576,951 & 344,250 & 40.33 & 207,971 & 123,361 & 40.68 & 784,922 & 467,611 & 40.43 \\
    2 & -7.0 & 609,897 & 297,300 & 51.25 & 222,470 & 110,893 & 50.15 & 832,367 & 408,193 & 50.96 \\
    2 & -2.5 & 478,292 & 273,900 & 42.73 & 165,134 & 96,084 & 41.81 & 643,426 & 369,984 & 42.50 \\
    2 & 0.0 & 2,249,066 & 2,249,925 & -- & 764,116 & 766,152 & -- & 3,013,182 & 3,016,077 & -- \\
    2 & 2.5 & 468,645 & 246,300 & 47.44 & 167,202 & 88,865 & 46.85 & 635,847 & 335,165 & 47.29 \\
    2 & 7.0 & 548,522 & 272,321 & 50.35 & 209,861 & 103,923 & 50.48 & 758,383 & 376,244 & 50.39 \\
    3 & -7.0 & 580,023 & 273,000 & 52.93 & 194,491 & 93,427 & 51.96 & 774,514 & 366,427 & 52.69 \\
    3 & -2.5 & 463,098 & 261,525 & 43.53 & 165,740 & 94,483 & 42.99 & 628,838 & 356,008 & 43.39 \\
    3 & 0.0 & 453,426 & 261,750 & 42.27 & 157,786 & 93,196 & 40.94 & 611,212 & 354,946 & 41.93 \\
    3 & 2.5 & 476,010 & 260,775 & 45.22 & 168,163 & 92,723 & 44.86 & 644,173 & 353,498 & 45.12 \\
    3 & 7.0 & 548,450 & 273,300 & 50.17 & 204,808 & 102,172 & 50.11 & 753,258 & 375,472 & 50.15 \\
    \bottomrule
\end{tabular}}
\end{table}

Table~\ref{tab:rans_2d_linear_per_outer} reports the mean and maximum linear
iterations per pseudo-time step. The flow solve reaches its configured maximum
of $75$ iterations in every run, whereas the turbulence solve remains below
that limit. Aggregated over the converged pairs, the mean combined count per
pseudo-time step changes only from $101.39$ to $101.95$. Hence, the cumulative
linear-work reduction follows mainly from the reduction in pseudo-time
iterations.

\begin{table}[H]
\centering
\caption{Mean and maximum linear iterations per executed pseudo-time step for the 2D RANS airfoil cases.}
\label{tab:rans_2d_linear_per_outer}
\resizebox{\textwidth}{!}{%
\begin{tabular}{@{}ccrrrrrrrr@{}}
    \toprule
    & & \multicolumn{4}{c}{Mean-flow system}
      & \multicolumn{4}{c}{Turbulence system} \\
    \cmidrule(lr){3-6}\cmidrule(lr){7-10}
    Airfoil & $\alpha$
    & $\overline{l}_{\mathrm{flow}}^{\mathrm{UF}}$
    & $\overline{l}_{\mathrm{flow}}^{\mathrm{ADEx-FNO}}$
    & $l_{\mathrm{flow,max}}^{\mathrm{UF}}$
    & $l_{\mathrm{flow,max}}^{\mathrm{ADEx-FNO}}$
    & $\overline{l}_{\mathrm{turb}}^{\mathrm{UF}}$
    & $\overline{l}_{\mathrm{turb}}^{\mathrm{ADEx-FNO}}$
    & $l_{\mathrm{turb,max}}^{\mathrm{UF}}$
    & $l_{\mathrm{turb,max}}^{\mathrm{ADEx-FNO}}$ \\
    \midrule
    1 & -7.0 & 74.38 & 75.00 & 75 & 75 & 27.99 & 28.15 & 31 & 31 \\
    1 & -2.5 & 74.86 & 75.00 & 75 & 75 & 25.83 & 26.10 & 27 & 30 \\
    1 & 0.0 & 74.85 & 75.00 & 75 & 75 & 26.11 & 27.06 & 27 & 31 \\
    1 & 2.5 & 74.86 & 74.24 & 75 & 75 & 25.71 & 25.84 & 27 & 30 \\
    1 & 7.0 & 74.89 & 75.00 & 75 & 75 & 27.00 & 26.88 & 29 & 30 \\
    2 & -7.0 & 74.68 & 75.00 & 75 & 75 & 27.24 & 27.98 & 29 & 30 \\
    2 & -2.5 & 73.16 & 75.00 & 75 & 75 & 25.26 & 26.31 & 27 & 31 \\
    2 & 0.0 & 74.97 & 75.00 & 75 & 75 & 25.47 & 25.54 & 26 & 29 \\
    2 & 2.5 & 74.86 & 75.00 & 75 & 75 & 26.71 & 27.06 & 27 & 29 \\
    2 & 7.0 & 74.88 & 74.57 & 75 & 75 & 28.65 & 28.46 & 30 & 31 \\
    3 & -7.0 & 74.89 & 75.00 & 75 & 75 & 25.11 & 25.67 & 26 & 30 \\
    3 & -2.5 & 74.86 & 75.00 & 75 & 75 & 26.79 & 27.10 & 28 & 30 \\
    3 & 0.0 & 74.86 & 75.00 & 75 & 75 & 26.05 & 26.70 & 27 & 30 \\
    3 & 2.5 & 74.87 & 75.00 & 75 & 75 & 26.45 & 26.67 & 27 & 32 \\
    3 & 7.0 & 74.88 & 75.00 & 75 & 75 & 27.96 & 28.04 & 30 & 31 \\
    \bottomrule
\end{tabular}}
\end{table}

\begin{figure}[h]
    \centering
    \includegraphics[width=0.75\linewidth]
    {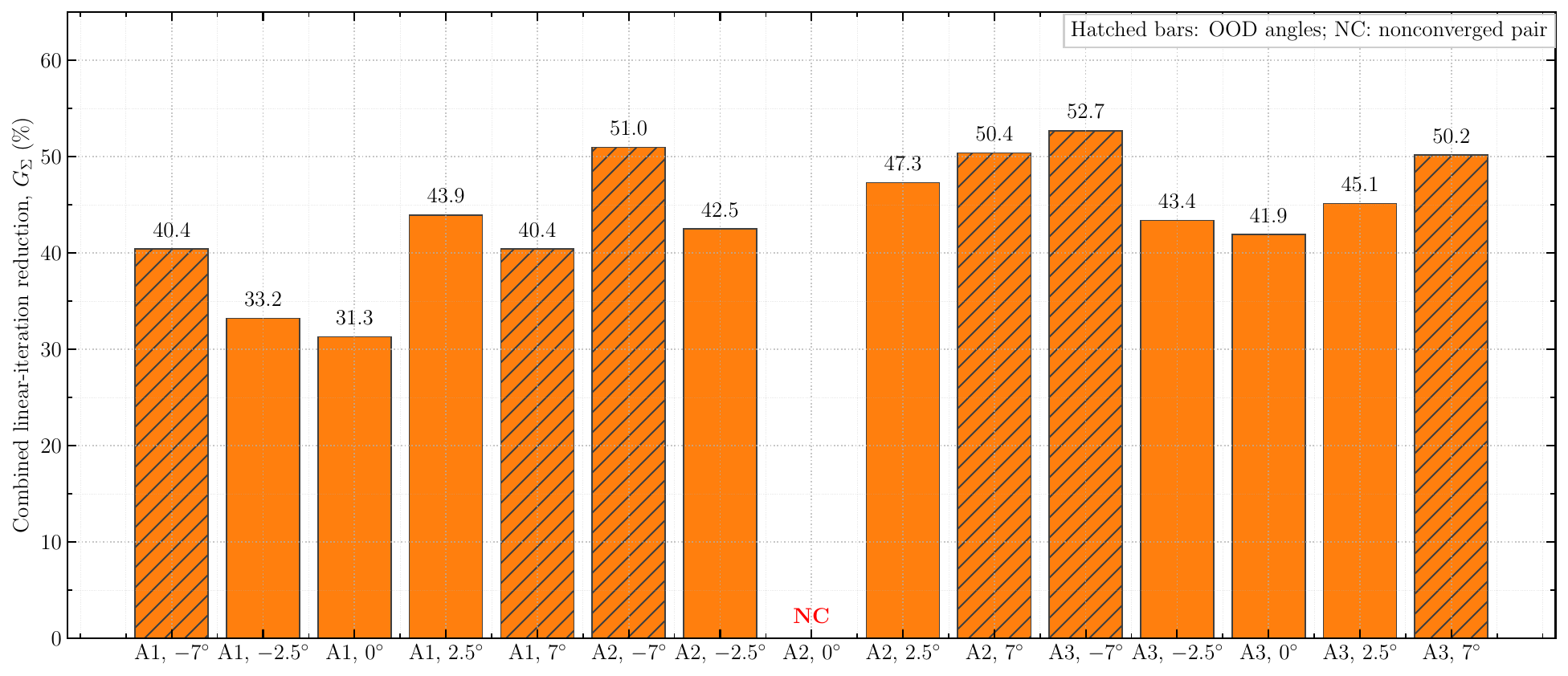}
    \caption{Per-case reduction in the unweighted cumulative linear-solver count
    $L_\Sigma$ for the 2D RANS airfoil cases. Hatched bars denote OOD
    angles, and NC denotes the nonconverged Airfoil~2, $\alpha=0^\circ$
    pair.}
    \label{fig:rans_2d_linear_reduction}
\end{figure}

Figure~\ref{fig:rans_2d_linear_reduction} summarizes the pairwise
reductions in $L_\Sigma$ and retains the nonconverged pair explicitly.

\subsection{Grouped statistics and dependence on the test parameters}
\label{app:rans_2d_grouped}

Table~\ref{tab:rans_2d_group_statistics} reports descriptive statistics for
$G_{\mathrm{outer}}$ and $G_\Sigma$. The OOD angles have a
larger (descriptive) mean reduction than
the ID angles, but the distinction is confounded with angle
magnitude and baseline convergence difficulty. It does not establish that extrapolation intrinsically improves the initialization. The sign of $\alpha$
also does not yield a uniform trend: for $|\alpha|=2.5^\circ$, the positive-angle
case has the larger reduction for all three geometries, whereas for
$|\alpha|=7^\circ$ the negative-angle case has the larger reduction for all
three. Geometry-specific means differ, but three held-out geometries are
insufficient to infer a general geometric dependence.

\begin{table}[H]
\centering
\caption{Descriptive statistics for the converged 2D RANS airfoil pairs.}
\label{tab:rans_2d_group_statistics}
\begin{tabular}{@{}lcrrrr@{}}
    \toprule
    Group & Metric & $n$ & Mean & Minimum & Maximum \\
    \midrule
    All converged pairs
    & $G_{\mathrm{outer}}$ & 14 & 44.17 & 32.02 & 53.00 \\
    All converged pairs
    & $G_\Sigma$ & 14 & 43.84 & 31.29 & 52.69 \\
    \addlinespace[2pt]

    ID angles
    & $G_{\mathrm{outer}}$ & 8 & 41.52 & 32.02 & 47.54 \\
    ID angles
    & $G_\Sigma$ & 8 & 41.08 & 31.29 & 47.29 \\
    \addlinespace[2pt]

    OOD angles
    & $G_{\mathrm{outer}}$ & 6 & 47.69 & 40.42 & 53.00 \\
    OOD angles
    & $G_\Sigma$ & 6 & 47.50 & 40.41 & 52.69 \\
    \addlinespace[2pt]

    Airfoil 1
    & $G_{\mathrm{outer}}$ & 5 & 38.09 & 32.02 & 43.66 \\
    Airfoil 1
    & $G_\Sigma$ & 5 & 37.86 & 31.29 & 43.94 \\
    \addlinespace[2pt]

    Airfoil 2
    & $G_{\mathrm{outer}}$ & 4 & 48.32 & 44.14 & 51.46 \\
    Airfoil 2
    & $G_\Sigma$ & 4 & 47.78 & 42.50 & 50.96 \\
    \addlinespace[2pt]

    Airfoil 3
    & $G_{\mathrm{outer}}$ & 5 & 46.91 & 42.38 & 53.00 \\
    Airfoil 3
    & $G_\Sigma$ & 5 & 46.66 & 41.93 & 52.69 \\
    \bottomrule
\end{tabular}
\end{table}

\begin{figure}[H]
    \centering
    \includegraphics[width=0.65\linewidth]
    {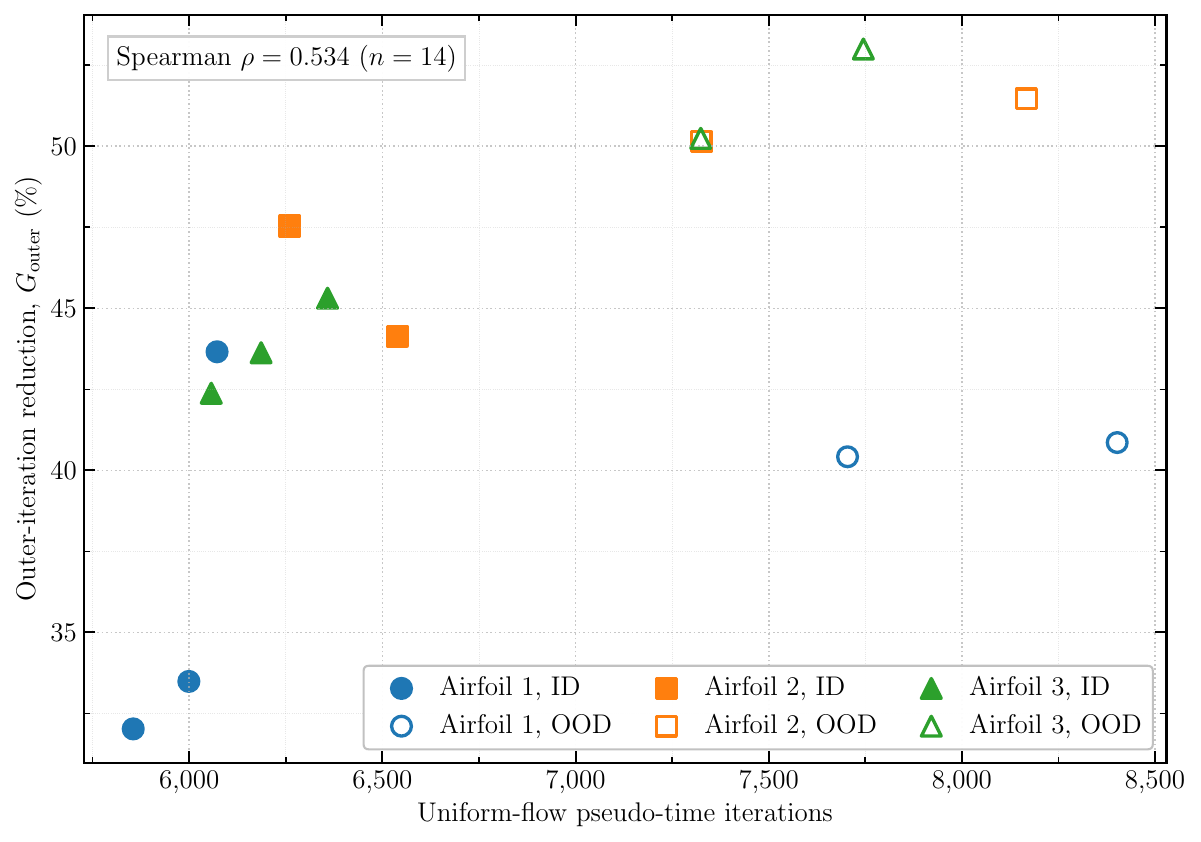}
    \caption{Outer-iteration reduction versus uniform-flow convergence effort for the
    $14$ converged 2D RANS airfoil pairs. Filled symbols denote ID angles
    within the training interval $[-5^\circ,5^\circ]$, open symbols denote
    OOD angles $\alpha=\pm7^\circ$, and color and marker shape identify
    the airfoil geometry.}
    \label{fig:rans_2d_reduction_vs_difficulty}
\end{figure}

Figure~\ref{fig:rans_2d_reduction_vs_difficulty} plots the outer-iteration
reduction against the uniform-flow convergence effort and separates the cases by
airfoil geometry and distribution status. Filled symbols denote ID angles
$\alpha\in\{-2.5^\circ,0^\circ,2.5^\circ\}$ within the training interval
$[-5^\circ,5^\circ]$; open symbols denote OOD angles $\alpha=\pm7^\circ$. Across the $14$ converged pairs, the Spearman correlation between the
uniform-flow pseudo-time count and $G_{\mathrm{outer}}$ is $0.534$
($p=0.049$), while the Pearson correlation is $0.496$ ($p=0.071$). The moderate
positive association indicates that cases requiring more pseudo-time iterations
from the uniform-flow state tend, within this dataset, to exhibit larger
percentage reductions. However, the point clouds overlap and the OOD cases also
correspond to the largest angle magnitude; geometry, angle magnitude,
distribution status, and baseline difficulty are therefore confounded. The
relation is reported as descriptive rather than as a general predictive law.

\subsection{Residual histories and late-stage convergence}

\begin{figure}[H]
    \centering
    \begin{subfigure}[t]{0.485\textwidth}
        \centering
        \includegraphics[width=\linewidth]{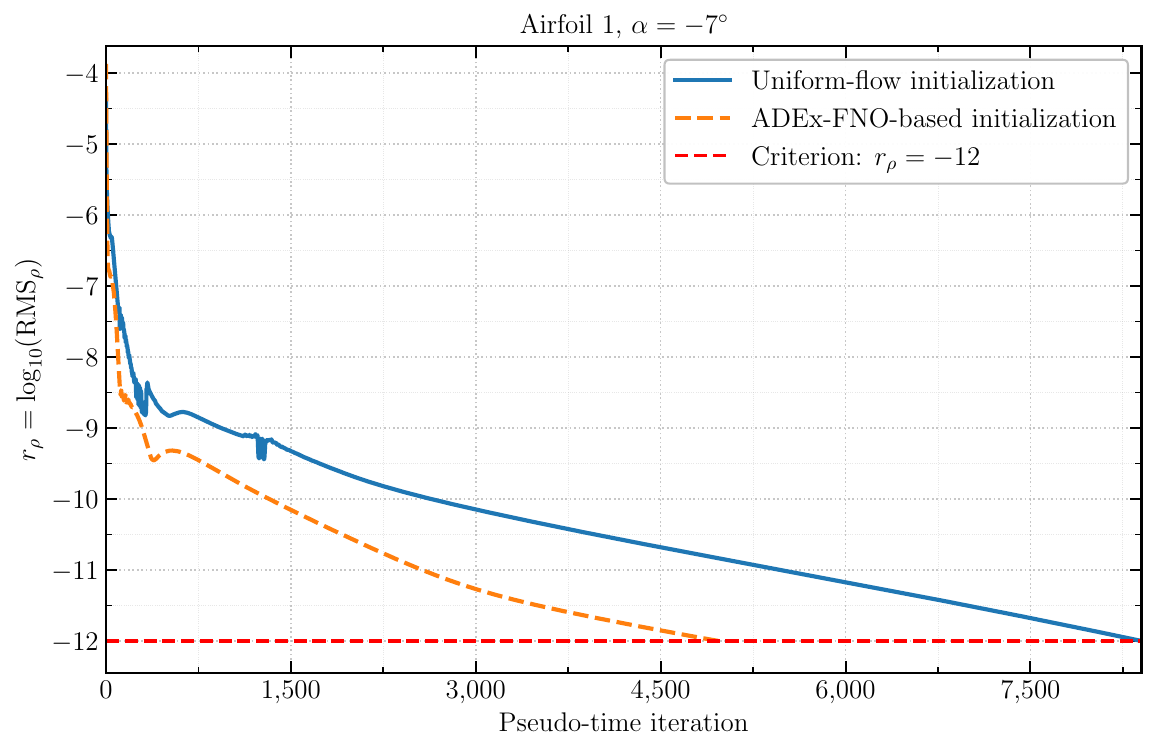}
    \end{subfigure}
    \hfill
    \begin{subfigure}[t]{0.485\textwidth}
        \centering
        \includegraphics[width=\linewidth]{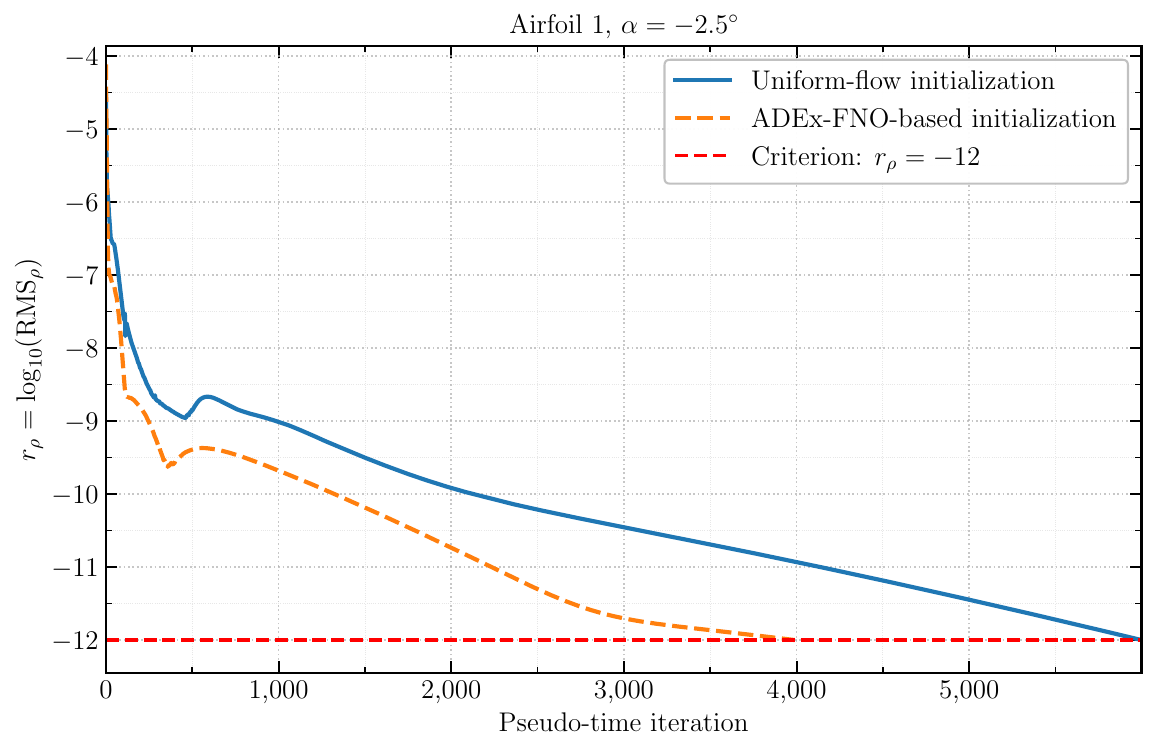}
    \end{subfigure}
    \par\medskip
    \begin{subfigure}[t]{0.485\textwidth}
        \centering
        \includegraphics[width=\linewidth]{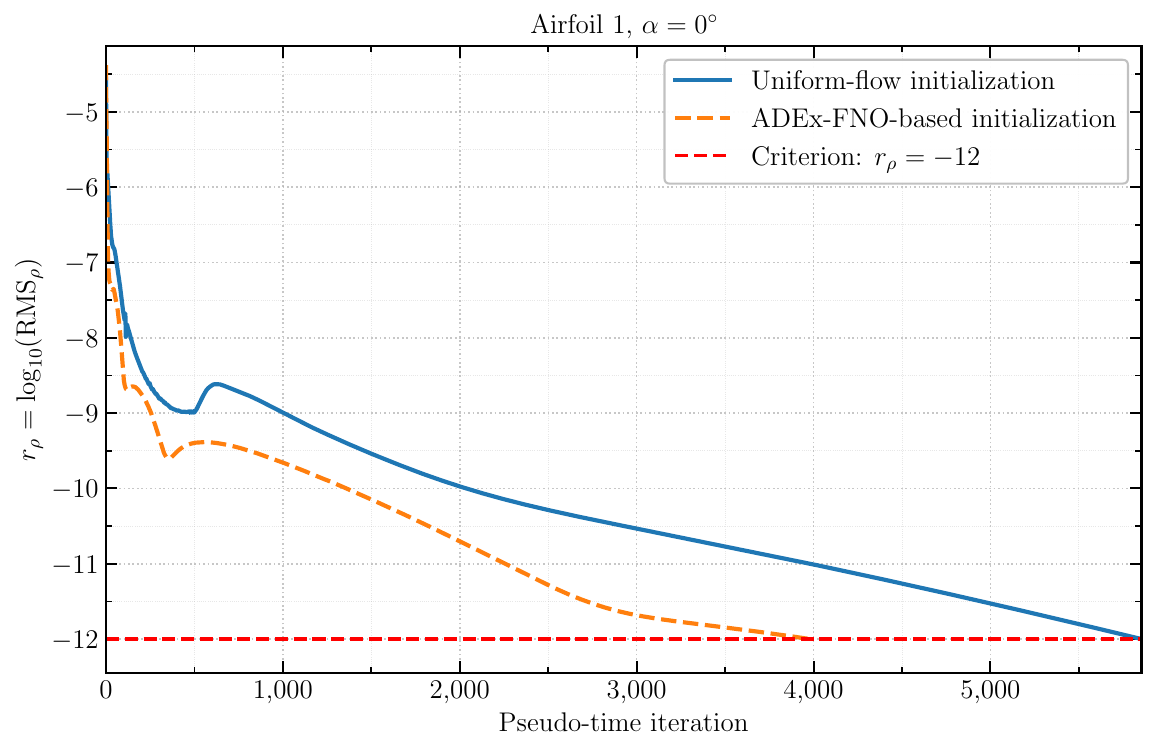}
    \end{subfigure}
    \hfill
    \begin{subfigure}[t]{0.485\textwidth}
        \centering
        \includegraphics[width=\linewidth]{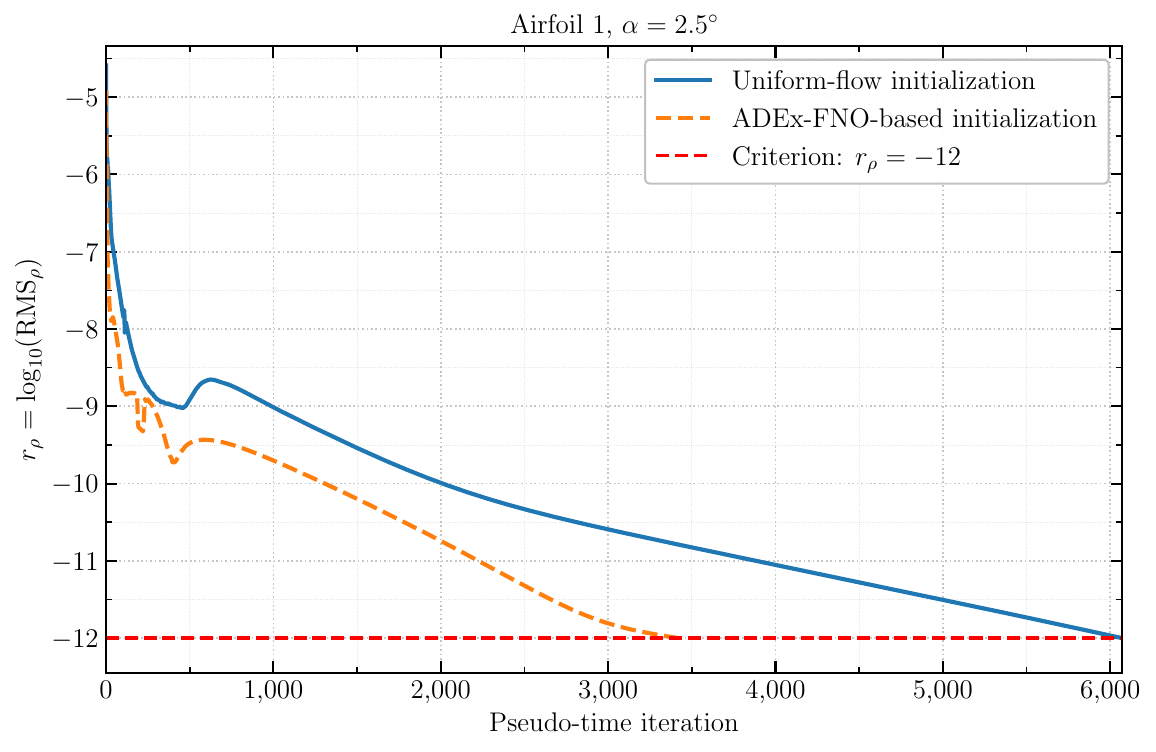}
    \end{subfigure}
    \par\medskip
    \begin{minipage}{\textwidth}
        \centering
        \begin{subfigure}[t]{0.485\textwidth}
        \centering
        \includegraphics[width=\linewidth]{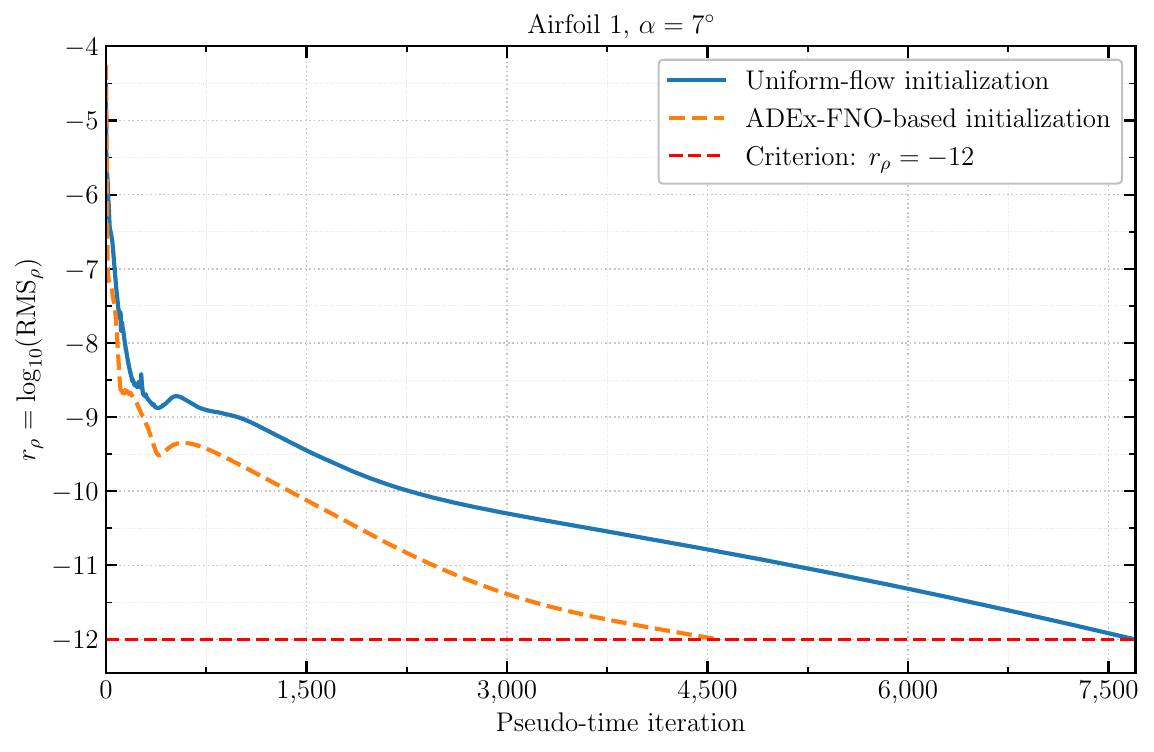}
    \end{subfigure}
    \end{minipage}
    \caption{Paired RMS-density residual histories for the 2D RANS Airfoil~1 cases. Blue solid curves denote uniform-flow initialization, orange dashed curves denote ADEx-FNO-based initialization, and the red dashed line is the prescribed $r_\rho=-12$ threshold.}
    \label{fig:rans_2d_residual_airfoil1}
\end{figure}

\begin{figure}[H]
    \centering
    \begin{subfigure}[t]{0.485\textwidth}
        \centering
        \includegraphics[width=\linewidth]{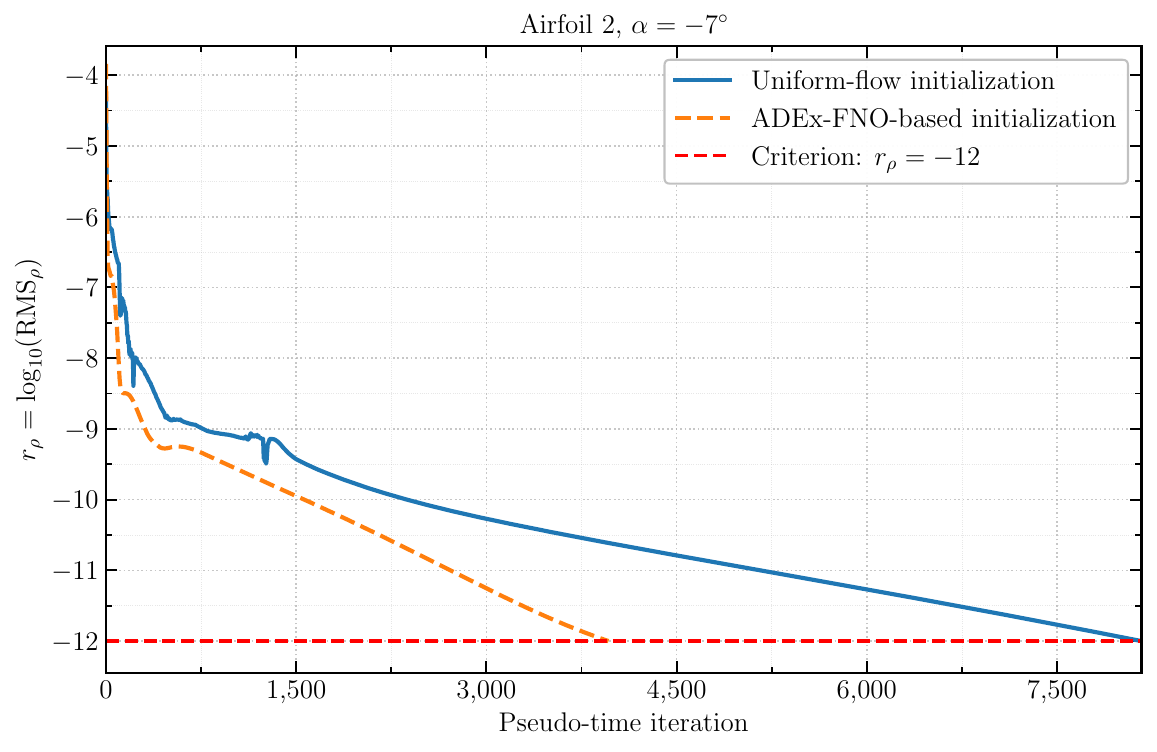}
    \end{subfigure}
    \hfill
    \begin{subfigure}[t]{0.485\textwidth}
        \centering
        \includegraphics[width=\linewidth]{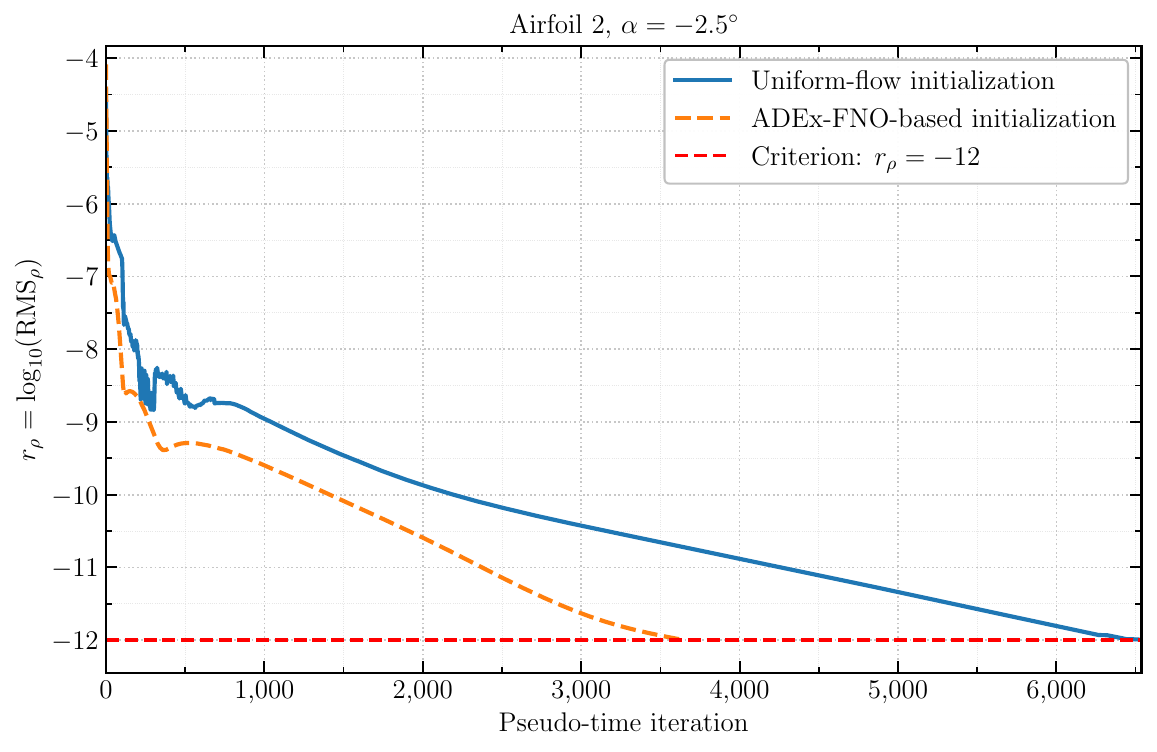}
    \end{subfigure}
    \par\medskip
    \begin{subfigure}[t]{0.485\textwidth}
        \centering
        \includegraphics[width=\linewidth]{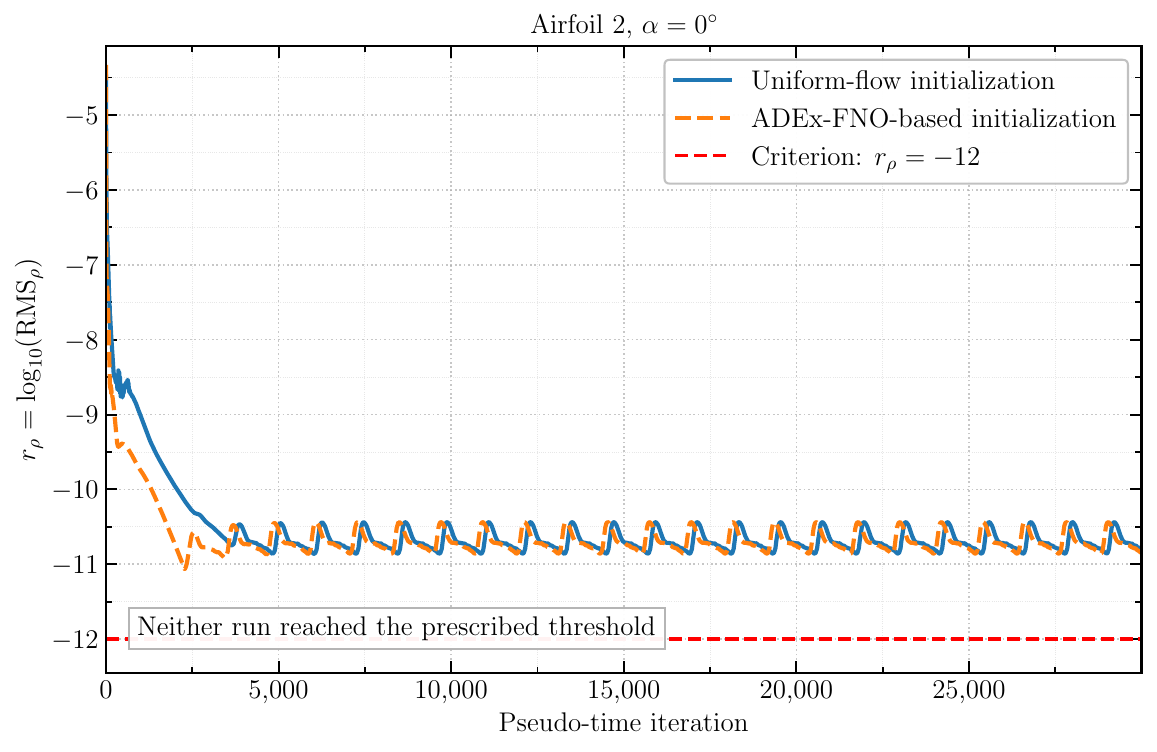}
    \end{subfigure}
    \hfill
    \begin{subfigure}[t]{0.485\textwidth}
        \centering
        \includegraphics[width=\linewidth]{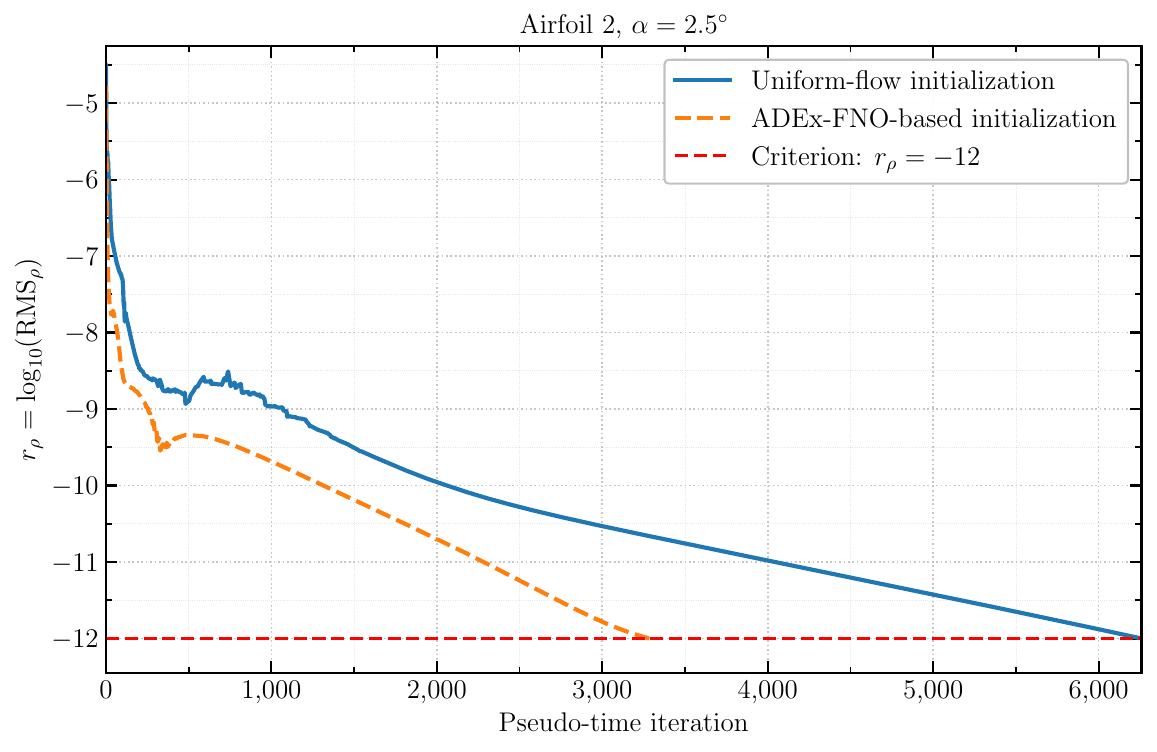}
    \end{subfigure}
    \par\medskip
    \begin{minipage}{\textwidth}
        \centering
        \begin{subfigure}[t]{0.485\textwidth}
        \centering
        \includegraphics[width=\linewidth]{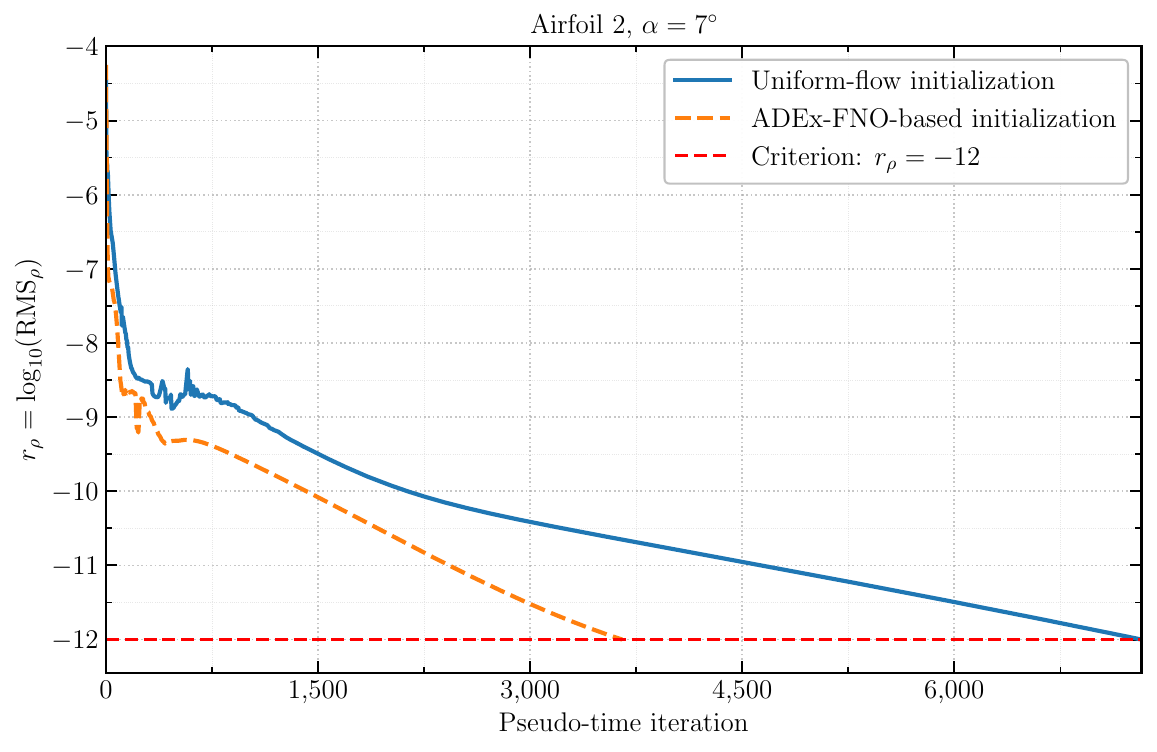}
    \end{subfigure}
    \end{minipage}
    \caption{Paired RMS-density residual histories for the 2D RANS Airfoil~2 cases. Blue solid curves denote uniform-flow initialization, orange dashed curves denote ADEx-FNO-based initialization, and the red dashed line is the prescribed $r_\rho=-12$ threshold.}
    \label{fig:rans_2d_residual_airfoil2}
\end{figure}

\begin{figure}[H]
    \centering
    \begin{subfigure}[t]{0.485\textwidth}
        \centering
        \includegraphics[width=\linewidth]{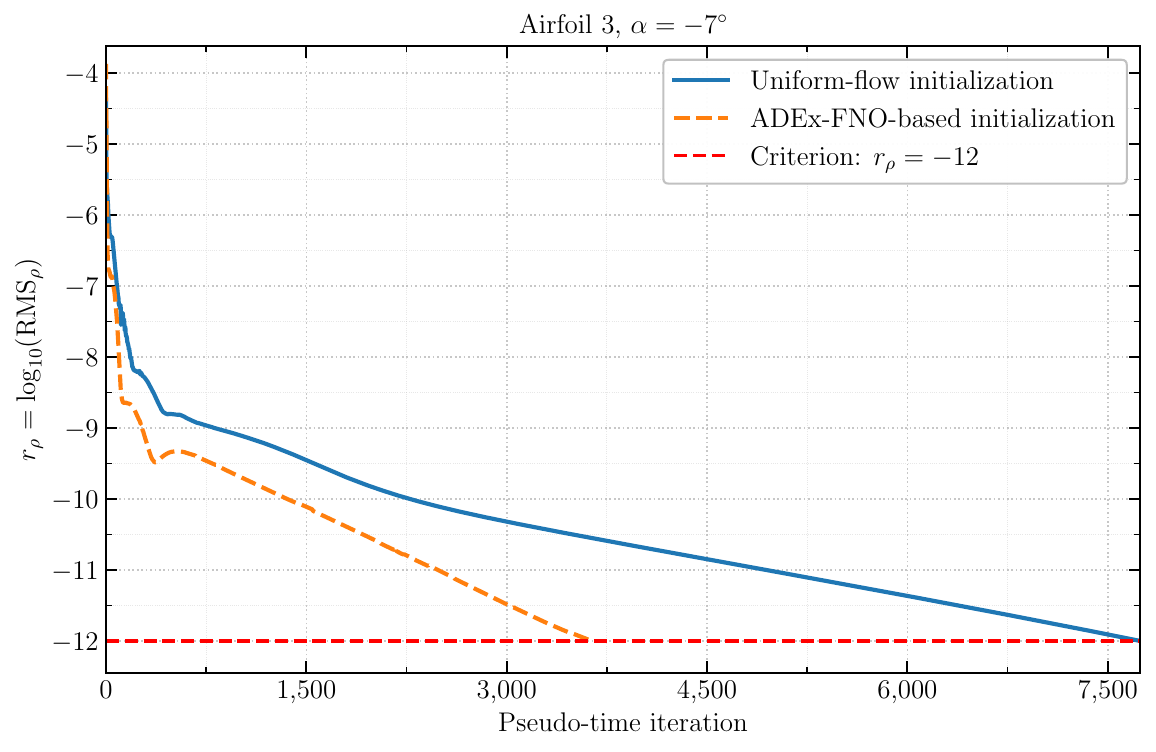}
    \end{subfigure}
    \hfill
    \begin{subfigure}[t]{0.485\textwidth}
        \centering
        \includegraphics[width=\linewidth]{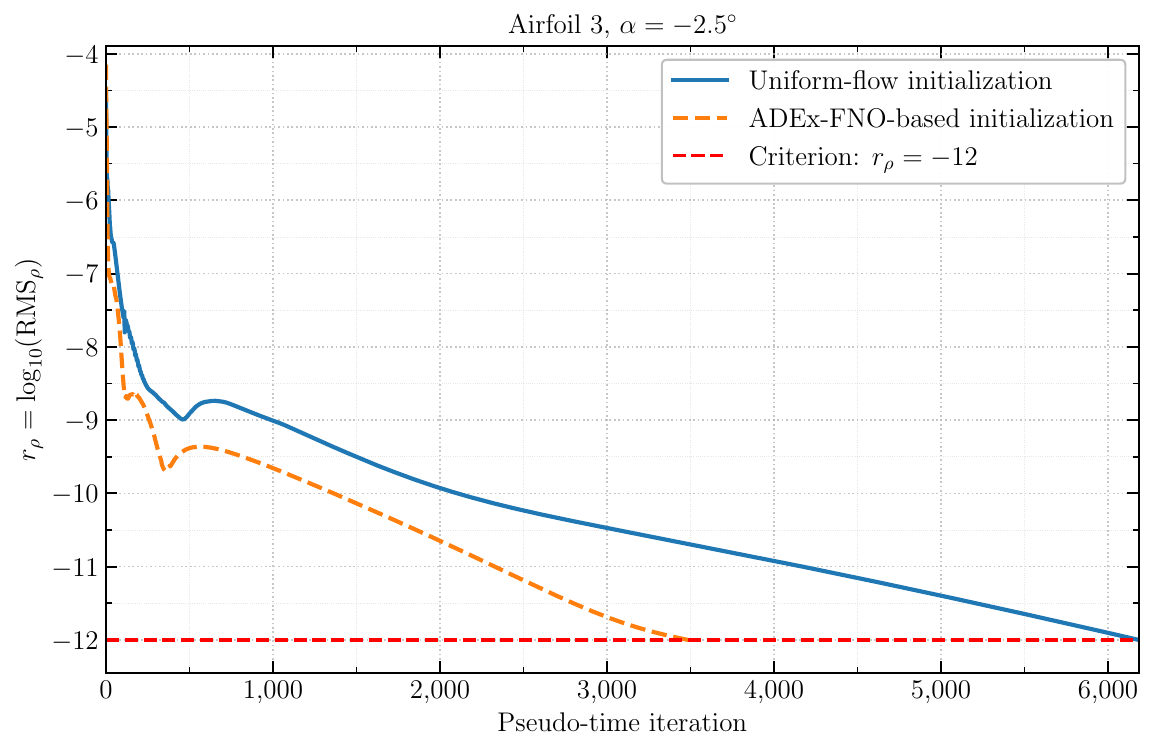}
    \end{subfigure}
    \par\medskip
    \begin{subfigure}[t]{0.485\textwidth}
        \centering
        \includegraphics[width=\linewidth]{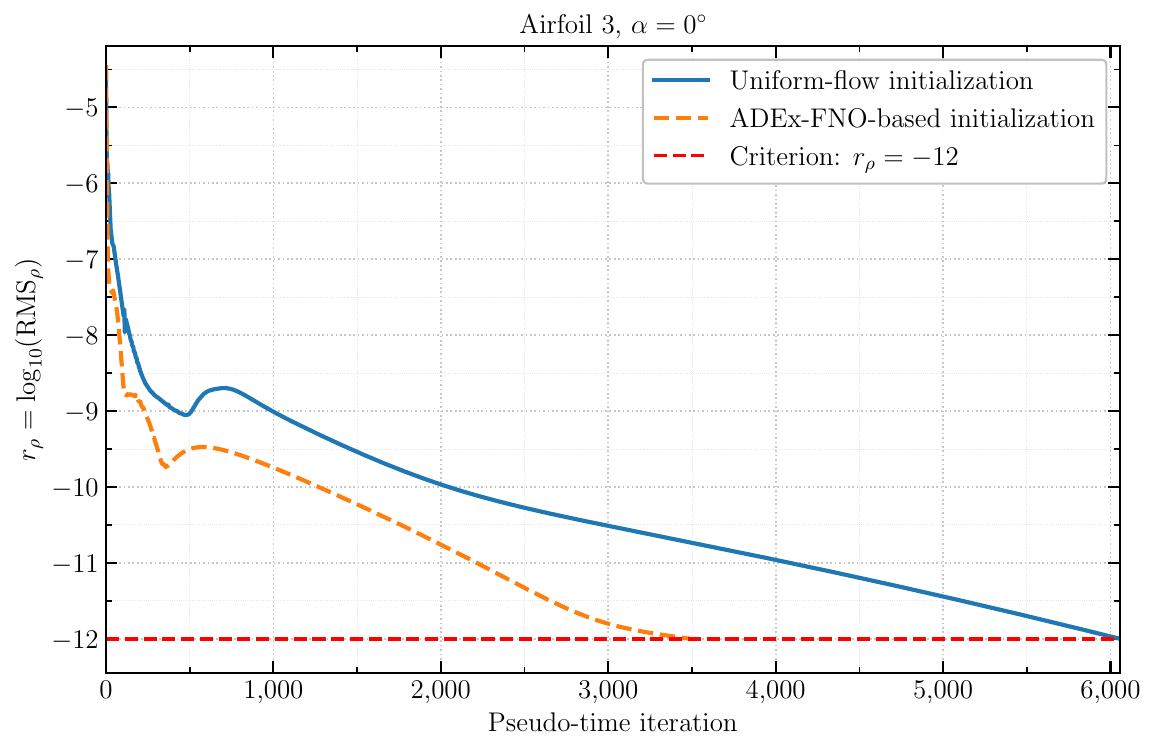}
    \end{subfigure}
    \hfill
    \begin{subfigure}[t]{0.485\textwidth}
        \centering
        \includegraphics[width=\linewidth]{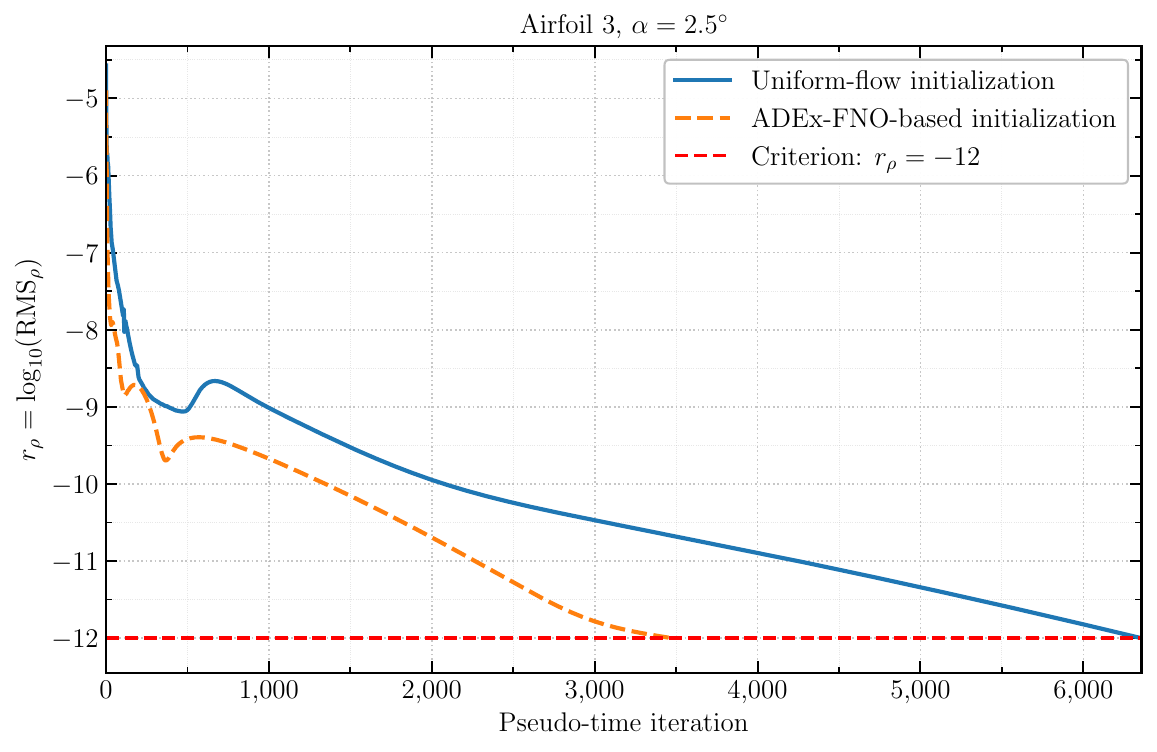}
    \end{subfigure}
    \par\medskip
    \begin{minipage}{\textwidth}
        \centering
        \begin{subfigure}[t]{0.485\textwidth}
        \centering
        \includegraphics[width=\linewidth]{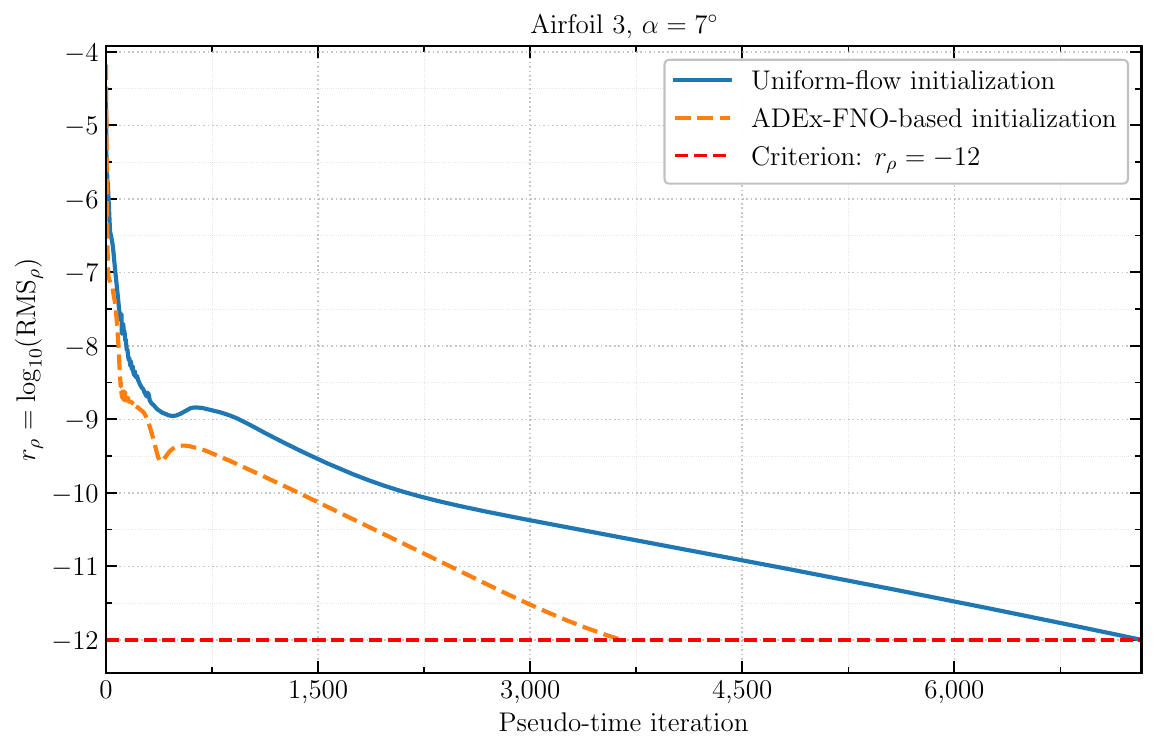}
    \end{subfigure}
    \end{minipage}
    \caption{Paired RMS-density residual histories for the 2D RANS Airfoil~3 cases. Blue solid curves denote uniform-flow initialization, orange dashed curves denote ADEx-FNO-based initialization, and the red dashed line is the prescribed $r_\rho=-12$ threshold.}
    \label{fig:rans_2d_residual_airfoil3}
\end{figure}

Figures~\ref{fig:rans_2d_residual_airfoil1},
\ref{fig:rans_2d_residual_airfoil2}, and
\ref{fig:rans_2d_residual_airfoil3} show the complete paired residual histories. Table~\ref{tab:rans_2d_late_stage}
quantifies the portion after each run first reaches $r_\rho=-10$. In every
converged pair, the ADEx-FNO-based initialization also shortens the subsequent
interval needed to reach $r_\rho=-12$. The reduction in this late-stage interval
ranges from $31.16\%$ to $58.59\%$, with a mean of $47.33\%$. Thus, the observed
benefit is not confined to the earliest residual drop. Because adaptive CFL
histories depend on the evolving state, this observation should not be
interpreted as an initialization-independent asymptotic convergence rate.

\begin{table}[H]
\centering
\caption{Pseudo-time iterations required after first reaching
$r_\rho=-10$ in the 2D RANS airfoil cases. The reduction is undefined for
Airfoil~2 at $\alpha=0^\circ$ because neither initialization converges.}
\label{tab:rans_2d_late_stage}
\begin{tabular}{@{}ccrrrrr@{}}
    \toprule
    Airfoil & $\alpha$
    & \makecell{UF index at\\$r_\rho=-10$}
    & \makecell{ADEx-FNO \\ index at\\$r_\rho=-10$}
    & \makecell{UF iterations\\$-10\rightarrow-12$}
    & \makecell{ADEx-FNO \\ iterations\\$-10\rightarrow-12$}
    & Reduction (\%) \\
    \midrule
    1 & -7.0 & 2,635 & 1,322 & 5,766 & 3,646 & 36.77 \\
    1 & -2.5 & 2,125 & 1,323 & 3,873 & 2,666 & 31.16 \\
    1 & 0.0 & 2,036 & 1,359 & 3,818 & 2,620 & 31.38 \\
    1 & 2.5 & 2,011 & 1,308 & 4,060 & 2,112 & 47.98 \\
    1 & 7.0 & 2,271 & 1,374 & 5,432 & 3,215 & 40.81 \\
    2 & -7.0 & 2,373 & 1,564 & 5,793 & 2,399 & 58.59 \\
    2 & -2.5 & 2,193 & 1,414 & 4,344 & 2,237 & 48.50 \\
    2 & 0.0 & 2,057 & 1,317 & -- & -- & -- \\
    2 & 2.5 & 2,060 & 1,315 & 4,199 & 1,968 & 53.13 \\
    2 & 7.0 & 2,133 & 1,418 & 5,191 & 2,233 & 56.98 \\
    3 & -7.0 & 2,274 & 1,356 & 5,470 & 2,283 & 58.26 \\
    3 & -2.5 & 2,104 & 1,364 & 4,081 & 2,122 & 48.00 \\
    3 & 0.0 & 2,045 & 1,271 & 4,011 & 2,218 & 44.70 \\
    3 & 2.5 & 2,078 & 1,348 & 4,279 & 2,128 & 50.27 \\
    3 & 7.0 & 2,130 & 1,364 & 5,193 & 2,279 & 56.11 \\
    \bottomrule
\end{tabular}
\end{table}

\subsection{Consistency of the terminal aerodynamic coefficients}

\begin{figure}[H]
    \centering
    \begin{subfigure}[t]{0.48\textwidth}
        \centering
        \includegraphics[width=\linewidth]
        {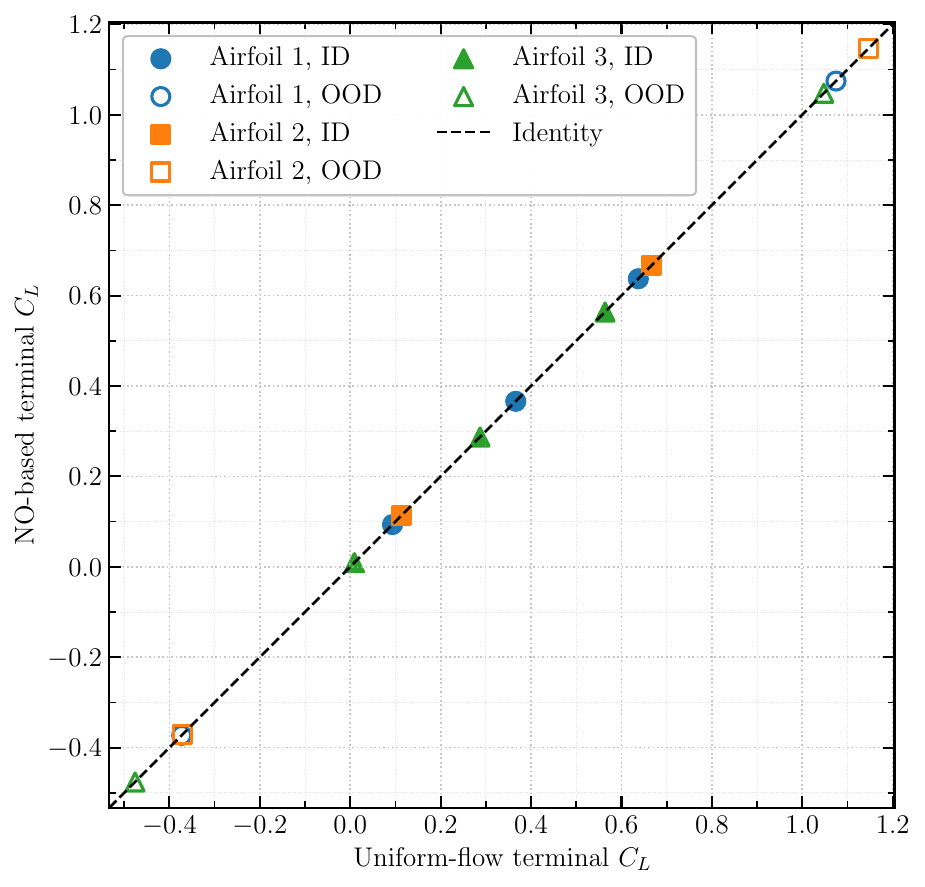}
        \caption{Lift coefficient.}
    \end{subfigure}
    \hfill
    \begin{subfigure}[t]{0.48\textwidth}
        \centering
        \includegraphics[width=\linewidth]
        {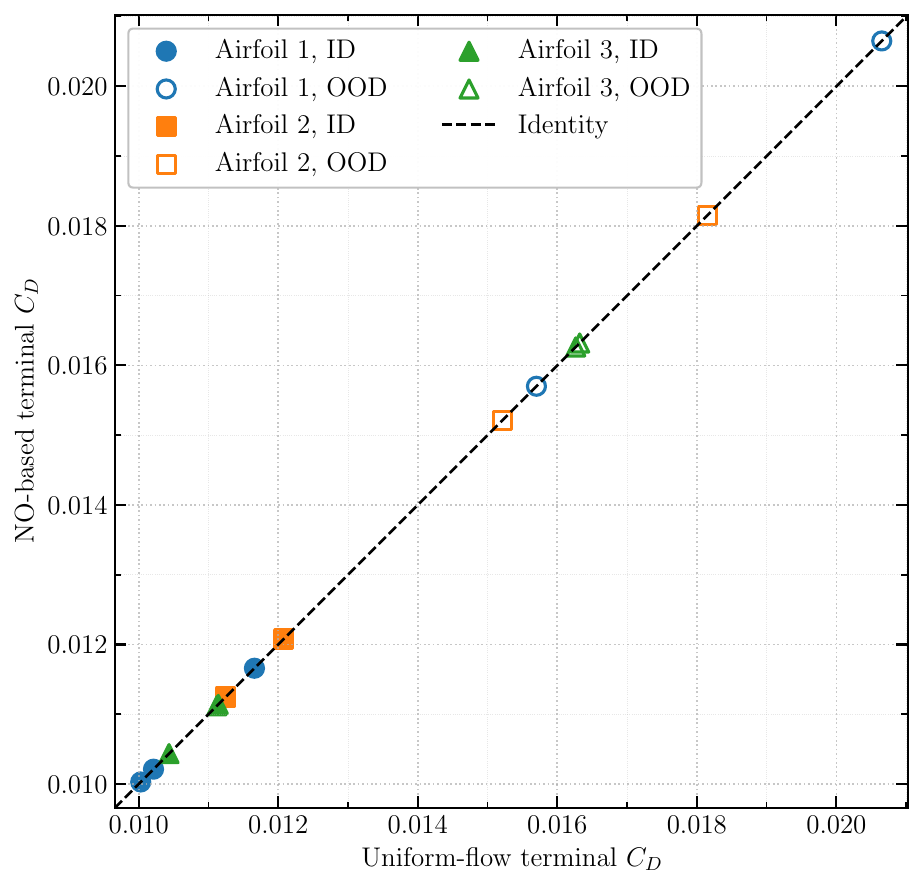}
        \caption{Drag coefficient.}
    \end{subfigure}
    \caption{Consistency of the terminal aerodynamic coefficients for the converged 2D RANS airfoil pairs.}
    \label{fig:rans_2d_coefficient_parity}
\end{figure}

Figure~\ref{fig:rans_2d_coefficient_parity} compares the terminal lift and drag
coefficients for the converged pairs. The points lie very close to the identity
line. Table~\ref{tab:rans_2d_terminal_quantities} gives the complete numerical
values. Among the converged pairs, the maximum absolute differences are
$2.9176\times10^{-5}$ for $C_L$ and $4.7816\times10^{-7}$ for $C_D$. The
maximum relative differences are $0.0474\%$ and $0.00263\%$, respectively. The
terminal quantities for Airfoil~2 at $\alpha=0^\circ$ are retained in the table
but are marked as nonconverged and are not used as evidence of
solution agreement.

\begin{table}[H]
\centering
\caption{Terminal residuals and aerodynamic coefficients for the 2D RANS airfoil
cases. The symbol $\dagger$ marks the nonconverged Airfoil~2,
$\alpha=0^\circ$ pair, which is excluded from convergence and
solution-consistency statistics.}
\label{tab:rans_2d_terminal_quantities}
\resizebox{\textwidth}{!}{%
\begin{tabular}{@{}ccrrrrrrrr@{}}
    \toprule
    Airfoil & $\alpha$
    & $r_\rho^{\mathrm{UF}}$ & $r_\rho^{\mathrm{ADEx-FNO}}$
    & $C_L^{\mathrm{UF}}$ & $C_L^{\mathrm{ADEx-FNO}}$ & $|\Delta C_L|$
    & $C_D^{\mathrm{UF}}$ & $C_D^{\mathrm{ADEx-FNO}}$ & $|\Delta C_D|$ \\
    \midrule
    1 & -7.0 & $-12.0$           & $-12.0$           & $ 1.075$           & $ 1.075$           & $\num{1.341e-6}$             & $0.0206$           & $0.0206$           & $\num{3.541e-8}$             \\
    1 & -2.5 & $-12.0$           & $-12.0$           & $ 0.637$           & $ 0.637$           & $\num{2.118e-6}$             & $0.0117$           & $0.0117$           & $\num{3.843e-8}$             \\
    1 &  0.0 & $-12.0$           & $-12.0$           & $ 0.366$           & $ 0.366$           & $\num{1.564e-6}$             & $0.0100$           & $0.0100$           & $\num{7.560e-9}$             \\
    1 &  2.5 & $-12.0$           & $-12.0$           & $ 0.093$           & $ 0.093$           & $\num{5.036e-6}$             & $0.0102$           & $0.0102$           & $\num{1.766e-7}$             \\
    1 &  7.0 & $-12.0$           & $-12.0$           & $-0.373$           & $-0.373$           & $\num{5.260e-6}$             & $0.0157$           & $0.0157$           & $\num{3.738e-8}$             \\
    2 & -7.0 & $-12.0$           & $-12.0$           & $ 1.147$           & $ 1.147$           & $\num{1.558e-5}$             & $0.0181$           & $0.0181$           & $\num{4.782e-7}$             \\
    2 & -2.5 & $-12.0$           & $-12.0$           & $ 0.667$           & $ 0.667$           & $\num{8.189e-6}$             & $0.0121$           & $0.0121$           & $\num{8.578e-8}$             \\
    2 &  0.0 & $-10.8^{\dagger}$ & $-10.9^{\dagger}$ & $ 0.390^{\dagger}$ & $ 0.390^{\dagger}$ & $(\num{1.378e-4})^{\dagger}$ & $0.0110^{\dagger}$ & $0.0110^{\dagger}$ & $(\num{8.366e-7})^{\dagger}$ \\
    2 &  2.5 & $-12.0$           & $-12.0$           & $ 0.113$           & $ 0.113$           & $\num{4.077e-6}$             & $0.0112$           & $0.0112$           & $\num{2.414e-7}$             \\
    2 &  7.0 & $-12.0$           & $-12.0$           & $-0.371$           & $-0.371$           & $\num{2.918e-5}$             & $0.0152$           & $0.0152$           & $\num{1.322e-7}$             \\
    3 & -7.0 & $-12.0$           & $-12.0$           & $ 1.048$           & $ 1.048$           & $\num{1.892e-5}$             & $0.0163$           & $0.0163$           & $\num{3.639e-7}$             \\
    3 & -2.5 & $-12.0$           & $-12.0$           & $ 0.564$           & $ 0.564$           & $\num{1.557e-5}$             & $0.0111$           & $0.0111$           & $\num{4.722e-8}$             \\
    3 &  0.0 & $-12.0$           & $-12.0$           & $ 0.287$           & $ 0.287$           & $\num{5.125e-6}$             & $0.0104$           & $0.0104$           & $\num{5.983e-8}$             \\
    3 &  2.5 & $-12.0$           & $-12.0$           & $ 0.009$           & $ 0.009$           & $\num{4.375e-6}$             & $0.0111$           & $0.0111$           & $\num{1.173e-7}$             \\
    3 &  7.0 & $-12.0$           & $-12.0$           & $-0.476$           & $-0.476$           & $\num{2.479e-5}$             & $0.0163$           & $0.0163$           & $\num{1.389e-8}$             \\
    \bottomrule
\end{tabular}}
\end{table}

\paragraph{Nonconverged case}
For Airfoil~2 at $\alpha=0^\circ$, both calculations stop after $29{,}999$
executed pseudo-time iterations. The uniform-flow run has terminal and minimum
residuals $-10.7830$ and $-10.8615$, respectively. The ADEx-FNO-based run has terminal
and minimum residuals $-10.8552$ and $-11.0630$. The accumulated unweighted
linear counts are $3{,}013{,}182$ and $3{,}016{,}077$. Because these calculations
do not reach the prescribed threshold (i.e., $r_\rho=-12$), they are not considered to assess
the convergence of $C_L$ and $C_D$.

\subsection{Transfer across CFD mesh resolutions}
\label{app:rans_2d_mesh_refinement}

The mesh-resolution study uses the Airfoil~2,
$\alpha=-7^\circ$ case at three nominal refinement levels. The supplied
\texttt{L1}, \texttt{L2}, and \texttt{L4} meshes contain, respectively,
$24,113$, $44,634$, and $92,730$ grid points and $24,403$, $44,968$, and
$93,171$ volume elements. Thus, the realized volume-element counts are
$1.00$, $1.84$, and $3.82$ times the \texttt{L1} count. The same
zero-based iteration convention is used at all three levels, and all six runs
terminate when $r_\rho$ first falls below $-12$.

Table~\ref{tab:rans_2d_mesh_refinement_linear_totals} reports the cumulative
linear-solver work. At every level, the ADEx-FNO-based initialization reduces the
mean-flow, turbulence, and unweighted combined linear counts. The mean-flow
reductions range from $50.36\%$ to $52.54\%$, the turbulence reductions from
$49.46\%$ to $51.57\%$, and the combined reductions from $50.13\%$ to
$52.31\%$.

\begin{table}[H]
\centering
\caption{Cumulative linear-solver iterations for the mesh-refined 2D RANS case. The
superscripts UF and ADEx-FNO identify uniform-flow and ADEx-FNO-based initialization.}
\label{tab:rans_2d_mesh_refinement_linear_totals}
\resizebox{\textwidth}{!}{%
\begin{tabular}{@{}lrrcrrcrrc@{}}
    \toprule
    & \multicolumn{3}{c}{Mean-flow system}
    & \multicolumn{3}{c}{Turbulence system}
    & \multicolumn{3}{c}{Unweighted combined count} \\
    \cmidrule(lr){2-4}\cmidrule(lr){5-7}\cmidrule(lr){8-10}
    Level
    & $L_{\mathrm{flow}}^{\mathrm{UF}}$
    & $L_{\mathrm{flow}}^{\mathrm{ADEx-FNO}}$
    & $G_{\mathrm{flow}}$ (\%)
    & $L_{\mathrm{turb}}^{\mathrm{UF}}$
    & $L_{\mathrm{turb}}^{\mathrm{ADEx-FNO}}$
    & $G_{\mathrm{turb}}$ (\%)
    & $L_\Sigma^{\mathrm{UF}}$
    & $L_\Sigma^{\mathrm{ADEx-FNO}}$
    & $G_\Sigma$ (\%) \\
    \midrule
    \texttt{L1} & 609,897 & 297,300 & 51.25
                 & 222,470 & 110,893 & 50.15
                 & 832,367 & 408,193 & 50.96 \\
    \texttt{L2} & 581,825 & 288,825 & 50.36
                 & 203,550 & 102,875 & 49.46
                 & 785,375 & 391,700 & 50.13 \\
    \texttt{L4} & 543,670 & 258,000 & 52.54
                 & 174,202 & 84,360 & 51.57
                 & 717,872 & 342,360 & 52.31 \\
    \bottomrule
\end{tabular}}
\end{table}

Table~\ref{tab:rans_2d_mesh_refinement_linear_per_outer} gives the mean and
maximum linear iterations per pseudo-time step. The mean-flow solve reaches
its configured maximum of $75$ iterations in every run. The mean combined
count per pseudo-time step changes from $101.92$ to $102.98$ at
\texttt{L1}, from $99.73$ to $101.71$ at \texttt{L2}, and from $97.80$
to $99.52$ at \texttt{L4}. Therefore, the ADEx-FNO-based initialization does not
reduce the average linear work of an individual pseudo-time step. Its
cumulative benefit follows from the smaller number of executed steps.

\begin{table}[H]
\centering
\caption{Mean and maximum linear iterations per executed pseudo-time step for the
mesh-refined 2D RANS case.}
\label{tab:rans_2d_mesh_refinement_linear_per_outer}
\resizebox{\textwidth}{!}{%
\begin{tabular}{@{}lrrrrrrrr@{}}
    \toprule
    & \multicolumn{4}{c}{Mean-flow system}
    & \multicolumn{4}{c}{Turbulence system} \\
    \cmidrule(lr){2-5}\cmidrule(lr){6-9}
    Level
    & $\overline{l}_{\mathrm{flow}}^{\mathrm{UF}}$
    & $\overline{l}_{\mathrm{flow}}^{\mathrm{ADEx-FNO}}$
    & $l_{\mathrm{flow,max}}^{\mathrm{UF}}$
    & $l_{\mathrm{flow,max}}^{\mathrm{ADEx-FNO}}$
    & $\overline{l}_{\mathrm{turb}}^{\mathrm{UF}}$
    & $\overline{l}_{\mathrm{turb}}^{\mathrm{ADEx-FNO}}$
    & $l_{\mathrm{turb,max}}^{\mathrm{UF}}$
    & $l_{\mathrm{turb,max}}^{\mathrm{ADEx-FNO}}$ \\
    \midrule
    \texttt{L1} & 74.68 & 75.00 & 75 & 75
                 & 27.24 & 27.98 & 29 & 30 \\
    \texttt{L2} & 73.88 & 75.00 & 75 & 75
                 & 25.85 & 26.71 & 27 & 30 \\
    \texttt{L4} & 74.07 & 75.00 & 75 & 75
                 & 23.73 & 24.52 & 26 & 27 \\
    \bottomrule
\end{tabular}}
\end{table}

The initial ADEx-FNO-based density residual is less negative than the corresponding
uniform-flow value at every level. In the order
$(r_\rho^{\mathrm{UF}},r_\rho^{\mathrm{ADEx-FNO}})$, the initial pairs are
\begin{equation*}
\begin{aligned}
    \texttt{L1}:&\quad(-4.3312,-3.8401),\\
    \texttt{L2}:&\quad(-4.5001,-4.1793),\\
    \texttt{L4}:&\quad(-4.6721,-4.6015).
\end{aligned}
\end{equation*}
The convergence reduction therefore cannot be attributed simply to a smaller
initial residual. Figure~\ref{fig:rans_2d_mesh_refinement} shows that the
ADEx-FNO-based histories instead enter the later residual-decay regime earlier.

Table~\ref{tab:rans_2d_mesh_refinement_late_stage} quantifies the interval
from the first occurrence of $r_\rho<-10$ to the final $r_\rho<-12$
crossing. The reduction remains between $56.27\%$ and $58.59\%$, showing
that the benefit persists well beyond the earliest residual transient. As in
the broader 2D study, this observation is descriptive and should not be
interpreted as an initialization-independent asymptotic convergence rate.

\begin{table}[H]
\centering
\caption{Pseudo-time iterations required after first reaching $r_\rho=-10$ in the
mesh-refined 2D RANS case.}
\label{tab:rans_2d_mesh_refinement_late_stage}
\begin{tabular}{@{}lrrrrr@{}}
    \toprule
    Level
    & \makecell{UF index at\\$r_\rho=-10$}
    & \makecell{ADEx-FNO \\ index at\\$r_\rho=-10$}
    & \makecell{UF iterations\\$-10\rightarrow-12$}
    & \makecell{ADEx-FNO \\ iterations\\$-10\rightarrow-12$}
    & Reduction (\%) \\
    \midrule
    \texttt{L1} & 2,373 & 1,564 & 5,793 & 2,399 & 58.59 \\
    \texttt{L2} & 2,152 & 1,467 & 5,722 & 2,383 & 58.35 \\
    \texttt{L4} & 1,803 & 1,018 & 5,536 & 2,421 & 56.27 \\
    \bottomrule
\end{tabular}
\end{table}

Table~\ref{tab:rans_2d_mesh_refinement_terminal} compares the terminal
residuals and aerodynamic coefficients within each mesh level. All values are
taken from runs satisfying the same $r_\rho<-12$ criterion. The maximum
absolute differences over the three levels are $5.4055\times10^{-5}$ for
$C_L$ and $1.4840\times10^{-6}$ for $C_D$; the corresponding maximum
relative differences are $0.00471\%$ and $0.00867\%$. The aerodynamic
coefficients vary across mesh levels, as expected for different spatial
discretizations, but the uniform-flow and ADEx-FNO-based runs converge to
consistent values on each fixed mesh. Accordingly, this experiment assesses
transfer of the initializer across meshes rather than formal grid convergence
of the CFD solution.

\begin{table}[H]
\centering
\caption{Terminal residuals and aerodynamic coefficients for the mesh-refined 2D RANS case.}
\label{tab:rans_2d_mesh_refinement_terminal}
\resizebox{\textwidth}{!}{%
\begin{tabular}{@{}lrrrrrrrr@{}}
    \toprule
    Level
    & $r_\rho^{\mathrm{UF}}$ & $r_\rho^{\mathrm{ADEx-FNO}}$
    & $C_L^{\mathrm{UF}}$ & $C_L^{\mathrm{ADEx-FNO}}$ & $|\Delta C_L|$
    & $C_D^{\mathrm{UF}}$ & $C_D^{\mathrm{ADEx-FNO}}$ & $|\Delta C_D|$ \\
    \midrule
    \texttt{L1}
    & $-12.0$ & $-12.0$
    & $1.147$ & $1.147$ & $\num{1.558e-5}$
    & $0.0181$ & $0.0181$ & $\num{4.782e-7}$ \\
    \texttt{L2}
    & $-12.0$ & $-12.0$
    & $1.151$ & $1.151$ & $\num{2.070e-5}$
    & $0.0172$ & $0.0172$ & $\num{6.177e-7}$ \\
    \texttt{L4}
    & $-12.0$ & $-12.0$
    & $1.149$ & $1.149$ & $\num{5.406e-5}$
    & $0.0171$ & $0.0171$ & $\num{1.484e-6}$ \\
    \bottomrule
\end{tabular}}
\end{table}


\section{Details of the 3D RANS initialization study}
\label{app:rans_3d_appendix}

\subsection{Counting convention, convergence logic, and meshes}

The executed pseudo-time count follows the same zero-based convention as in
the 2D RANS study: a terminal solver index $k$ represents $k+1$ iterations.
The cumulative mean-flow and turbulence linear counts are denoted by
$L_{\mathrm{flow}}$ and $L_{\mathrm{turb}}$, and their unweighted sum is
$L_\Sigma$.

For the 3D RANS cases, convergence monitoring begins after iteration $10$ and
requires the simultaneous conditions
\begin{equation*}
    r_\rho<-12,
    \qquad
    C_{\mathrm{C},D}<10^{-7},
    \qquad
    C_{\mathrm{C},L}<10^{-7},
\end{equation*}
where $C_{\mathrm{C},D}$ and $C_{\mathrm{C},L}$ are the drag- and lift-Cauchy
diagnostics. All $30$ calculations satisfy the complete criterion, with the
density residual controlling termination in every case.

The incidence is represented through the orientation of the case-specific
mesh, while the configured freestream angle-of-attack parameter remains zero.
The reference length is $0.64607\,\mathrm{m}$, the moment origin is
$(0.25,0,0)$, and the reference area is evaluated from the projected geometry.
Table~\ref{tab:rans_3d_meshes} reports the case-specific mesh sizes and
reference areas. Both initializations in each pair use identical geometric and
reference quantities.

\begin{table}[H]
\centering
\caption{Case-specific mesh sizes and reference areas for the 3D RANS wing
cases. Both initializations within each pair use the same mesh.}
\label{tab:rans_3d_meshes}
\begin{tabular}{@{}ccrrr@{}}
\toprule
Wing & $\alpha$ ($^\circ$) & Grid points & Volume elements
& $A_{\mathrm{ref}}$ ($\mathrm{m}^2$) \\
\midrule
1 & -7.0 & 4,047,165 & 5,975,312 & 0.750549 \\
  & -2.5 & 4,206,230 & 6,140,298 & 0.755422 \\
  & 0.0 & 4,228,596 & 6,169,027 & 0.756123 \\
  & 2.5 & 4,208,629 & 6,144,493 & 0.755388 \\
  & 7.0 & 4,221,303 & 6,163,161 & 0.750451 \\
\midrule
2 & -7.0 & 4,214,982 & 6,085,885 & 0.750451 \\
  & -2.5 & 4,467,152 & 6,413,515 & 0.755327 \\
  & 0.0 & 4,866,552 & 6,833,017 & 0.756043 \\
  & 2.5 & 4,883,473 & 6,858,903 & 0.755329 \\
  & 7.0 & 4,908,221 & 6,873,476 & 0.750454 \\
\midrule
3 & -7.0 & 4,052,234 & 5,975,820 & 0.750622 \\
  & -2.5 & 4,206,878 & 6,124,048 & 0.755410 \\
  & 0.0 & 4,217,590 & 6,144,171 & 0.756095 \\
  & 2.5 & 4,482,443 & 6,497,564 & 0.755354 \\
  & 7.0 & 4,971,712 & 7,056,614 & 0.750431 \\
\bottomrule
\end{tabular}
\end{table}

\subsection{Detailed linear-solver effort}
\label{app:rans_3d_linear}

Table~\ref{tab:rans_3d_linear_totals} reports the cumulative mean-flow,
turbulence, and unweighted combined linear iterations for every pair. The
ADEx-FNO-based initialization reduces each cumulative count in all $15$ cases. The
pairwise reductions in $L_\Sigma$ range from $21.10\%$ to $64.15\%$, with a
mean of $43.04\%$.

\begin{table}[H]
\centering
\caption{Per-case cumulative linear-solver iterations for the 3D RANS wing
cases. Superscripts UF and ADEx-FNO denote uniform-flow and ADEx-FNO-based initialization.}
\label{tab:rans_3d_linear_totals}
\resizebox{\textwidth}{!}{%
\begin{tabular}{@{}ccrrcrrcrrc@{}}
\toprule
& & \multicolumn{3}{c}{Mean-flow system} & \multicolumn{3}{c}{Turbulence system} & \multicolumn{3}{c}{Unweighted combined count} \\
\cmidrule(lr){3-5}\cmidrule(lr){6-8}\cmidrule(lr){9-11}
Wing & $\alpha$ & $L_{\mathrm{flow}}^{\mathrm{UF}}$ & $L_{\mathrm{flow}}^{\mathrm{ADEx-FNO}}$ & $G_{\mathrm{flow}}$ (\%) & $L_{\mathrm{turb}}^{\mathrm{UF}}$ & $L_{\mathrm{turb}}^{\mathrm{ADEx-FNO}}$ & $G_{\mathrm{turb}}$ (\%) & $L_\Sigma^{\mathrm{UF}}$ & $L_\Sigma^{\mathrm{ADEx-FNO}}$ & $G_\Sigma$ (\%) \\
\midrule
1 & -7.0 & 87,380 & 50,260 & 42.48 & 71,894 & 41,101 & 42.83 & 159,274 & 91,361 & 42.64 \\
1 & -2.5 & 88,580 & 64,340 & 27.37 & 73,327 & 50,020 & 31.79 & 161,907 & 114,360 & 29.37 \\
1 & 0.0 & 93,760 & 63,720 & 32.04 & 72,741 & 49,455 & 32.01 & 166,501 & 113,175 & 32.03 \\
1 & 2.5 & 91,540 & 61,200 & 33.14 & 71,515 & 48,916 & 31.60 & 163,055 & 110,116 & 32.47 \\
1 & 7.0 & 74,340 & 58,200 & 21.71 & 62,447 & 47,755 & 23.53 & 136,787 & 105,955 & 22.54 \\
2 & -7.0 & 87,660 & 33,180 & 62.15 & 71,825 & 27,031 & 62.37 & 159,485 & 60,211 & 62.25 \\
2 & -2.5 & 87,880 & 45,860 & 47.82 & 69,953 & 36,253 & 48.18 & 157,833 & 82,113 & 47.97 \\
2 & 0.0 & 87,820 & 32,100 & 63.45 & 69,944 & 26,612 & 61.95 & 157,764 & 58,712 & 62.78 \\
2 & 2.5 & 92,080 & 48,680 & 47.13 & 73,045 & 39,384 & 46.08 & 165,125 & 88,064 & 46.67 \\
2 & 7.0 & 75,620 & 59,380 & 21.48 & 62,946 & 49,946 & 20.65 & 138,566 & 109,326 & 21.10 \\
3 & -7.0 & 83,260 & 28,600 & 65.65 & 66,341 & 25,029 & 62.27 & 149,601 & 53,629 & 64.15 \\
3 & -2.5 & 91,160 & 33,040 & 63.76 & 73,173 & 27,536 & 62.37 & 164,333 & 60,576 & 63.14 \\
3 & 0.0 & 99,080 & 55,020 & 44.47 & 79,399 & 44,654 & 43.76 & 178,479 & 99,674 & 44.15 \\
3 & 2.5 & 106,100 & 54,480 & 48.65 & 82,290 & 43,453 & 47.20 & 188,390 & 97,933 & 48.02 \\
3 & 7.0 & 77,900 & 59,120 & 24.11 & 64,247 & 45,712 & 28.85 & 142,147 & 104,832 & 26.25 \\
\bottomrule
\end{tabular}}
\end{table}

Table~\ref{tab:rans_3d_linear_per_outer} reports the mean and maximum linear
iterations per pseudo-time step. The mean-flow solve uses the configured
maximum of $20$ iterations at every step in every run. The turbulence solve
also reaches $20$ at least once in every run, while its per-run mean lies
between $15.46$ and $17.50$. Aggregated over all cases, the mean combined count
per pseudo-time step changes only from $36.09$ for uniform-flow initialization
to $36.14$ for ADEx-FNO-based initialization. The reduction in cumulative linear
work therefore follows from the reduction in pseudo-time iterations rather than
from making an individual pseudo-time step less expensive in terms of Krylov
iterations.

\begin{table}[H]
\centering
\caption{Mean and maximum linear iterations per executed pseudo-time step for the 3D RANS wing cases.}
\label{tab:rans_3d_linear_per_outer}
\resizebox{\textwidth}{!}{%
\begin{tabular}{@{}ccrrrrrrrr@{}}
\toprule
& & \multicolumn{4}{c}{Mean-flow system} & \multicolumn{4}{c}{Turbulence system} \\
\cmidrule(lr){3-6}\cmidrule(lr){7-10}
Wing & $\alpha$ & $\overline{l}_{\mathrm{flow}}^{\mathrm{UF}}$ & $\overline{l}_{\mathrm{flow}}^{\mathrm{ADEx-FNO}}$ & $l_{\mathrm{flow,max}}^{\mathrm{UF}}$ & $l_{\mathrm{flow,max}}^{\mathrm{ADEx-FNO}}$ & $\overline{l}_{\mathrm{turb}}^{\mathrm{UF}}$ & $\overline{l}_{\mathrm{turb}}^{\mathrm{ADEx-FNO}}$ & $l_{\mathrm{turb,max}}^{\mathrm{UF}}$ & $l_{\mathrm{turb,max}}^{\mathrm{ADEx-FNO}}$ \\
\midrule
1 & -7.0 & 20.00 & 20.00 & 20 & 20 & 16.46 & 16.36 & 20 & 20 \\
1 & -2.5 & 20.00 & 20.00 & 20 & 20 & 16.56 & 15.55 & 20 & 20 \\
1 & 0.0 & 20.00 & 20.00 & 20 & 20 & 15.52 & 15.52 & 20 & 20 \\
1 & 2.5 & 20.00 & 20.00 & 20 & 20 & 15.62 & 15.99 & 20 & 20 \\
1 & 7.0 & 20.00 & 20.00 & 20 & 20 & 16.80 & 16.41 & 20 & 20 \\
2 & -7.0 & 20.00 & 20.00 & 20 & 20 & 16.39 & 16.29 & 20 & 20 \\
2 & -2.5 & 20.00 & 20.00 & 20 & 20 & 15.92 & 15.81 & 20 & 20 \\
2 & 0.0 & 20.00 & 20.00 & 20 & 20 & 15.93 & 16.58 & 20 & 20 \\
2 & 2.5 & 20.00 & 20.00 & 20 & 20 & 15.87 & 16.18 & 20 & 20 \\
2 & 7.0 & 20.00 & 20.00 & 20 & 20 & 16.65 & 16.82 & 20 & 20 \\
3 & -7.0 & 20.00 & 20.00 & 20 & 20 & 15.94 & 17.50 & 20 & 20 \\
3 & -2.5 & 20.00 & 20.00 & 20 & 20 & 16.05 & 16.67 & 20 & 20 \\
3 & 0.0 & 20.00 & 20.00 & 20 & 20 & 16.03 & 16.23 & 20 & 20 \\
3 & 2.5 & 20.00 & 20.00 & 20 & 20 & 15.51 & 15.95 & 20 & 20 \\
3 & 7.0 & 20.00 & 20.00 & 20 & 20 & 16.49 & 15.46 & 20 & 20 \\
\bottomrule
\end{tabular}}
\end{table}

\begin{figure}[H]
    \centering
    \includegraphics[width=0.75\linewidth]
    {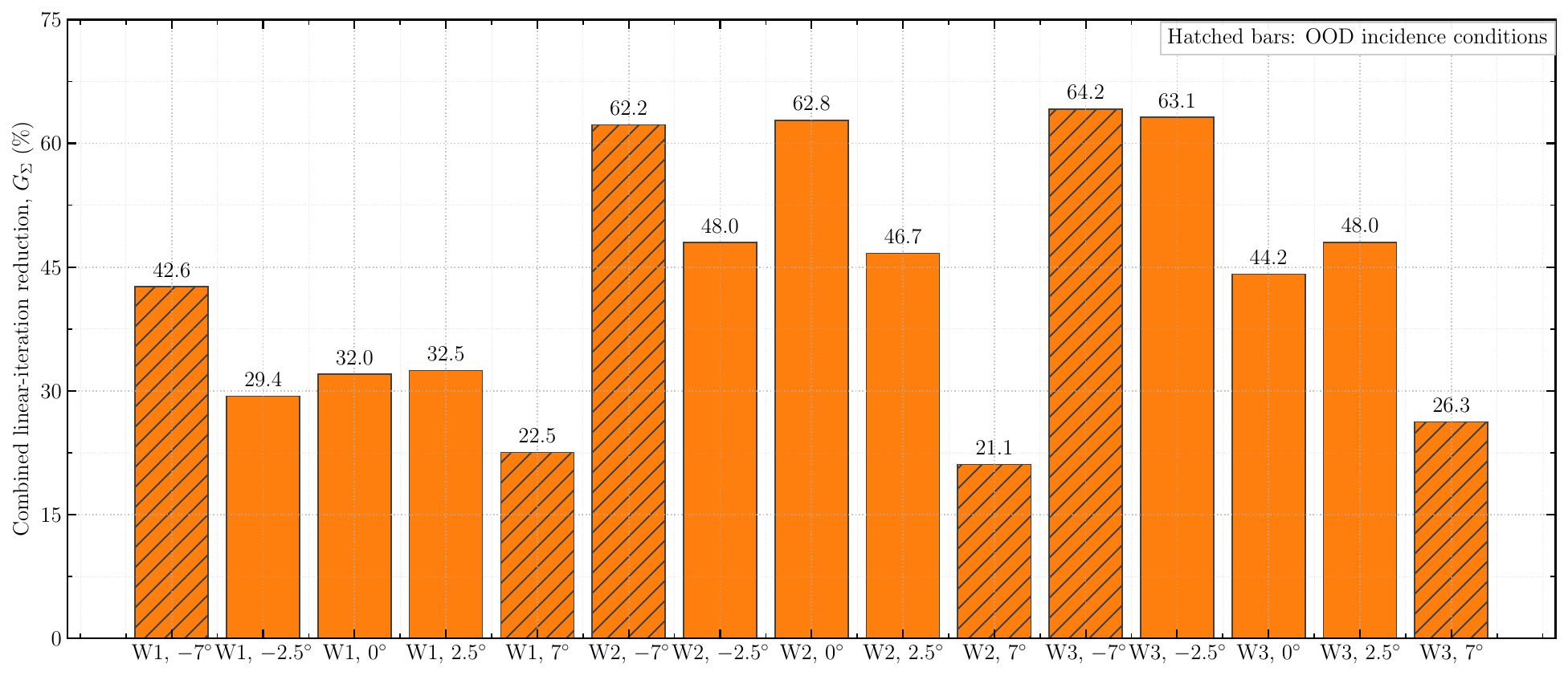}
    \caption{Per-case reduction in the unweighted cumulative linear-solver count
    $L_\Sigma$ for the 3D RANS wing cases. Hatched bars denote OOD
    incidence conditions.}
    \label{fig:rans_3d_linear_reduction}
\end{figure}

\begin{figure}[H]
    \centering
    \includegraphics[width=0.65\linewidth]
    {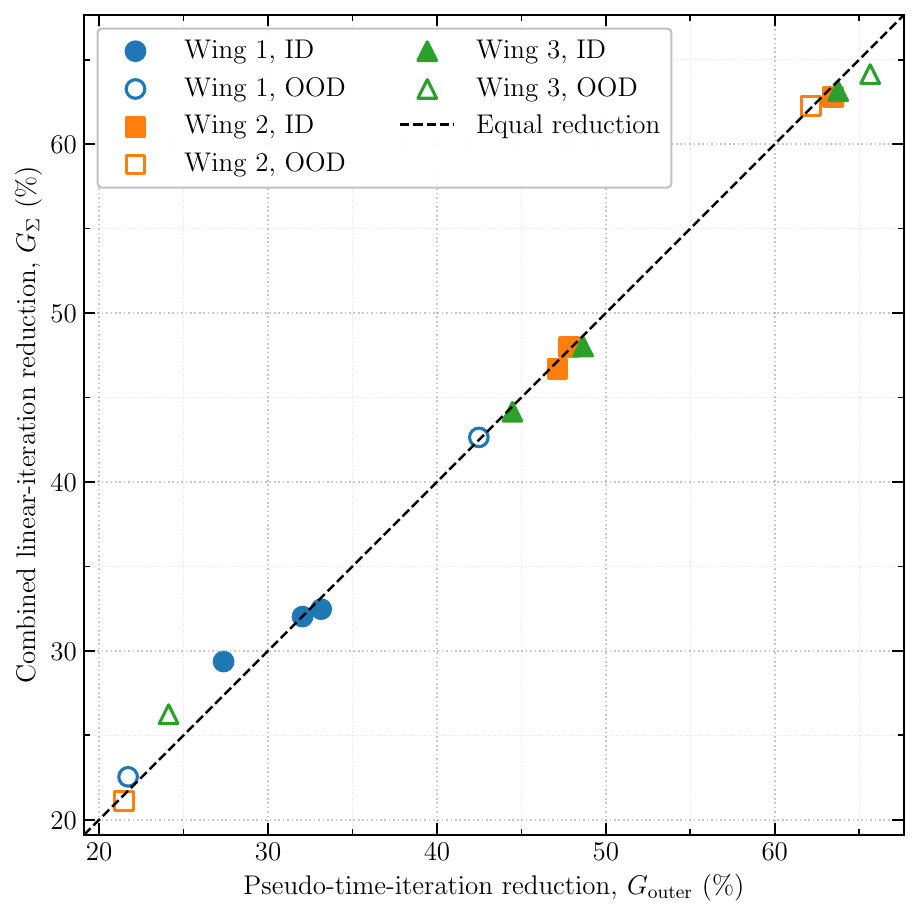}
    \caption{Pairwise comparison of the pseudo-time and unweighted
    cumulative linear-solver reductions for the 3D RANS wing cases. The dotted
    line is the identity relation; filled and open symbols denote ID and OOD
    incidence conditions, respectively.}
    \label{fig:rans_3d_outer_vs_linear}
\end{figure}

Figure~\ref{fig:rans_3d_linear_reduction} shows the per-case reduction in
$L_\Sigma$. Its close correspondence with
Figure~\ref{fig:rans_3d_outer_reduction} is quantified by
Figure~\ref{fig:rans_3d_outer_vs_linear}: for every case, the two reductions
lie near the identity line. The pairwise difference
$G_{\mathrm{outer}}-G_\Sigma$ ranges from $-2.14$ to $1.50$ percentage points
and has a mean of $-0.01$ percentage points.

\subsection{Grouped statistics and dependence on the test parameters}
\label{app:rans_3d_grouped}

Table~\ref{tab:rans_3d_group_statistics} reports descriptive statistics for
$G_{\mathrm{outer}}$ and $G_\Sigma$. The ID cases have mean reductions of
$45.31\%$ and $45.18\%$, respectively, whereas the OOD cases have means of
$39.60\%$ and $39.82\%$. This difference does not establish that interpolation
is intrinsically easier to accelerate: distribution status is confounded with
incidence magnitude, sign, geometry, and the corresponding flow state. All six
OOD pairs retain positive reductions.

\begin{table}[H]
\centering
\caption{Descriptive statistics for the converged 3D RANS pairs.}
\label{tab:rans_3d_group_statistics}
\begin{tabular}{@{}lcrrrr@{}}
\toprule
Group & Metric & $n$ & Mean & Minimum & Maximum \\
\midrule
\multirow{2}{*}{All converged pairs}
    & $G_{\mathrm{outer}}$ & 15 & 43.03 & 21.48 & 65.65 \\
    & $G_\Sigma$          & 15 & 43.04 & 21.10 & 64.15 \\
\addlinespace[2pt]
\multirow{2}{*}{ID angles}
    & $G_{\mathrm{outer}}$ & 9 & 45.31 & 27.37 & 63.76 \\
    & $G_\Sigma$          & 9 & 45.18 & 29.37 & 63.14 \\
\addlinespace[2pt]
\multirow{2}{*}{OOD angles}
    & $G_{\mathrm{outer}}$ & 6 & 39.60 & 21.48 & 65.65 \\
    & $G_\Sigma$          & 6 & 39.82 & 21.10 & 64.15 \\
\addlinespace[2pt]
\multirow{2}{*}{Wing 1}
    & $G_{\mathrm{outer}}$ & 5 & 31.35 & 21.71 & 42.48 \\
    & $G_\Sigma$          & 5 & 31.81 & 22.54 & 42.64 \\
\addlinespace[2pt]
\multirow{2}{*}{Wing 2}
    & $G_{\mathrm{outer}}$ & 5 & 48.40 & 21.48 & 63.45 \\
    & $G_\Sigma$          & 5 & 48.16 & 21.10 & 62.78 \\
\addlinespace[2pt]
\multirow{2}{*}{Wing 3}
    & $G_{\mathrm{outer}}$ & 5 & 49.33 & 24.11 & 65.65 \\
    & $G_\Sigma$          & 5 & 49.14 & 26.25 & 64.15 \\
\bottomrule
\end{tabular}
\end{table}


For Wings~1, 2, and 3, the mean outer-iteration reductions are $31.35\%$,
$48.40\%$, and $49.33\%$, respectively; the corresponding means of
$G_\Sigma$ are $31.81\%$, $48.16\%$, and $49.14\%$. The benefit is therefore
configuration dependent, but three held-out wings are insufficient to infer a
general geometry dependence. At
$|\alpha|=7^\circ$, the negative-incidence case has the larger reduction for
all three wings. At $|\alpha|=2.5^\circ$, however, the ordering changes with
geometry, so the sign of the incidence does not produce a uniform trend across
the complete matrix.

\begin{figure}[H]
    \centering
    \includegraphics[width=0.75\linewidth]
    {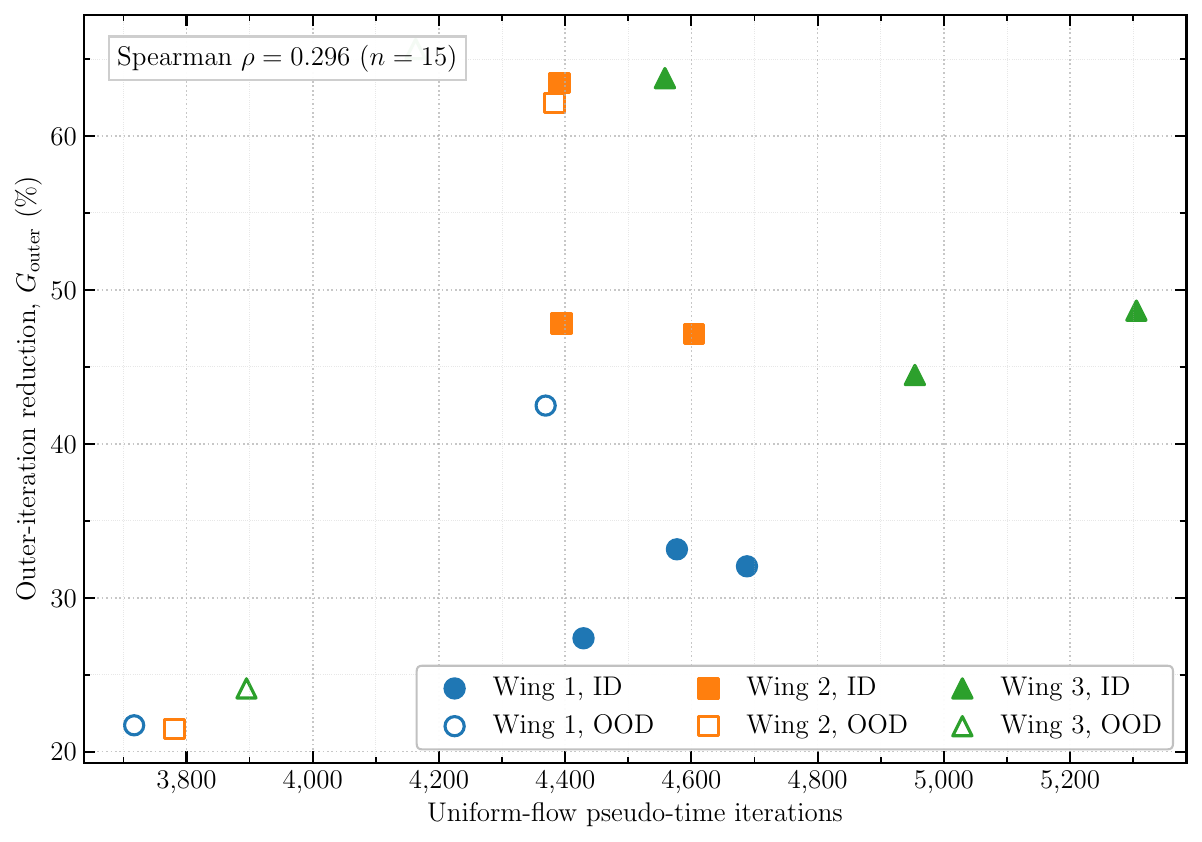}
    \caption{Outer-iteration reduction versus uniform-flow convergence
    effort for the $15$ converged 3D RANS wing pairs. Filled and open symbols
    denote ID and OOD incidence conditions, respectively; color and marker
    shape identify the wing geometry.}
    \label{fig:rans_3d_reduction_vs_difficulty}
\end{figure}

Figure~\ref{fig:rans_3d_reduction_vs_difficulty} plots
$G_{\mathrm{outer}}$ against the uniform-flow convergence effort. Filled symbols denote ID cases inside the training interval
$[-5^\circ,5^\circ]$, and open symbols denote OOD cases at
$\alpha=\pm7^\circ$. Across the $15$ pairs, the Spearman correlation is
$0.296$ ($p=0.283$), while the Pearson correlation is $0.388$ ($p=0.154$).
The present sample therefore does not support a systematic relation between the
uniform-flow iteration count and the percentage reduction obtained from the ADEx-FNO
initialization.

\subsection{Residual histories and convergence phases}
\label{app:rans_3d_convergence}

\begin{figure}[H]
    \centering
    \includegraphics[width=0.99\linewidth]
    {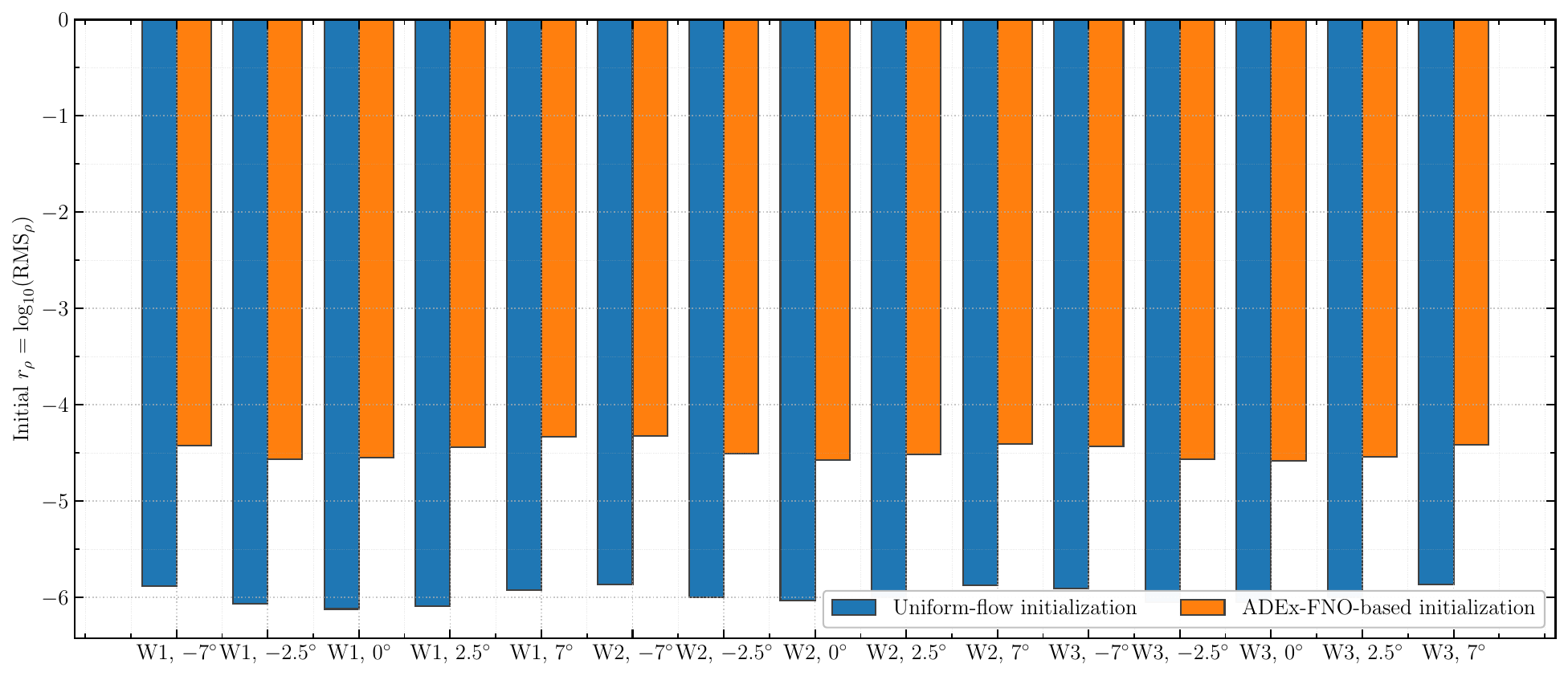}
    \caption{Initial logarithmic RMS-density residuals for the complete 3D
    RANS test matrix. The ADEx-FNO-based histories begin with a higher discrete
    residual in every pair, despite converging in fewer pseudo-time iterations.}
    \label{fig:rans_3d_initial_residual}
\end{figure}

Figure~\ref{fig:rans_3d_initial_residual} compares the initial density
residuals. The ADEx-FNO-based values are higher by $1.45$--$1.66$ logarithmic units
in all $15$ pairs. Thus, the reduction in total solver work is not a consequence
of starting from a smaller discrete residual. The ADEx-FNO-based initial state instead
changes the subsequent trajectory through the coupled nonlinear iteration.

\begin{figure}[H]
    \centering
    \begin{subfigure}[t]{0.485\textwidth}
        \centering
        \includegraphics[width=\linewidth]{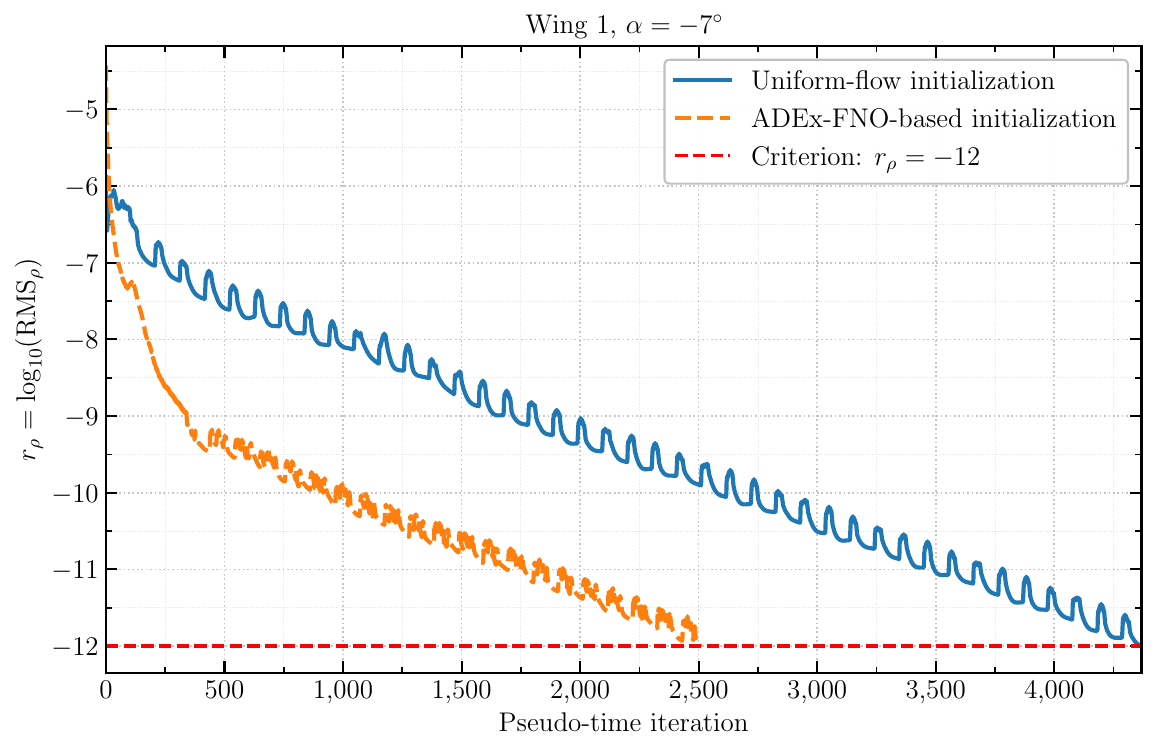}
    \end{subfigure}
    \hfill
    \begin{subfigure}[t]{0.485\textwidth}
        \centering
        \includegraphics[width=\linewidth]{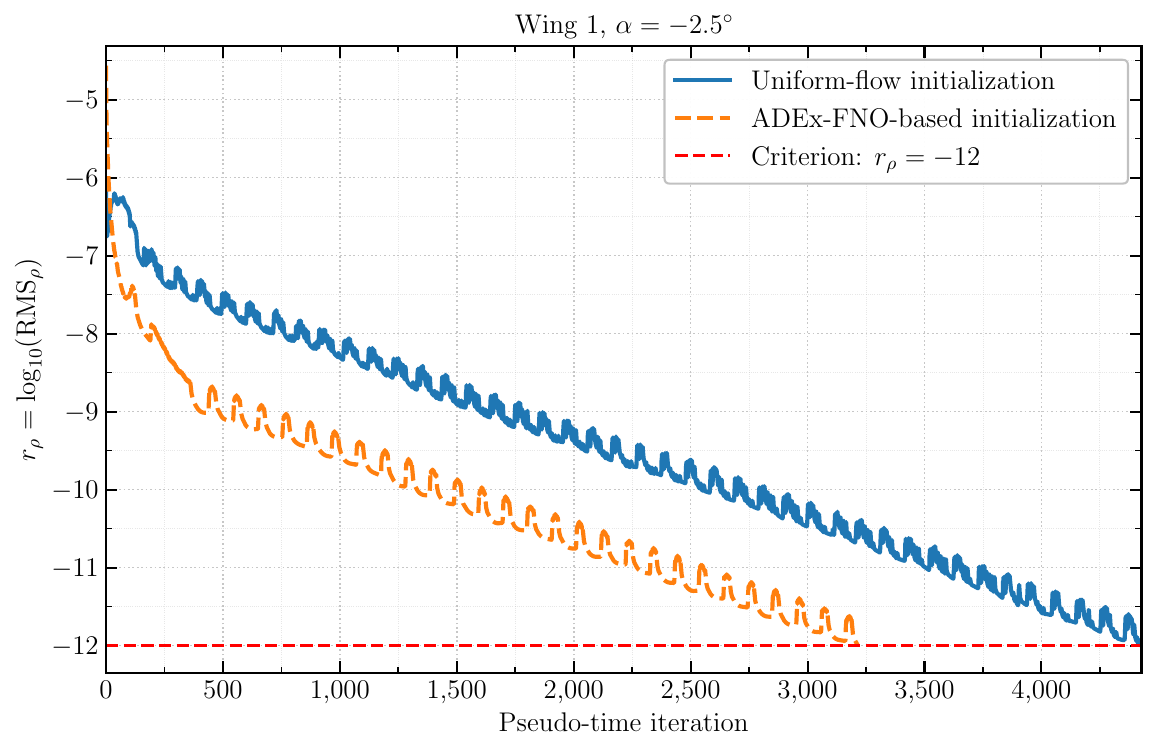}
    \end{subfigure}
    \par\medskip
    \begin{subfigure}[t]{0.485\textwidth}
        \centering
        \includegraphics[width=\linewidth]{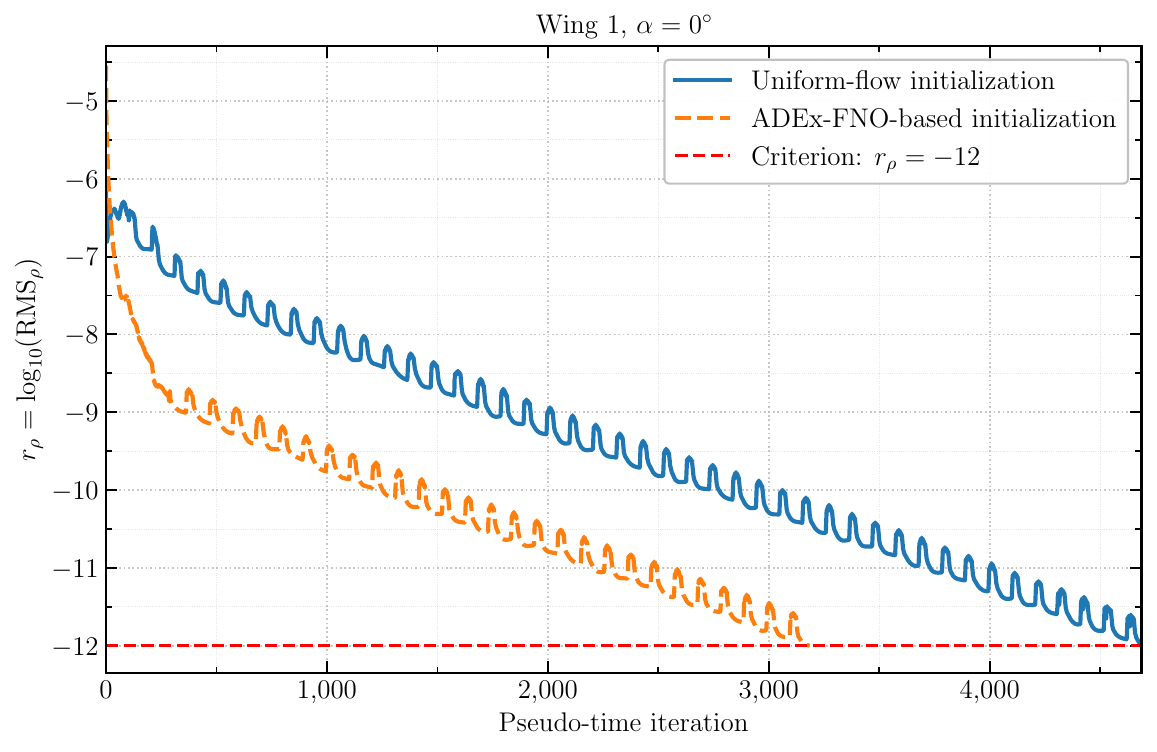}
    \end{subfigure}
    \hfill
    \begin{subfigure}[t]{0.485\textwidth}
        \centering
        \includegraphics[width=\linewidth]{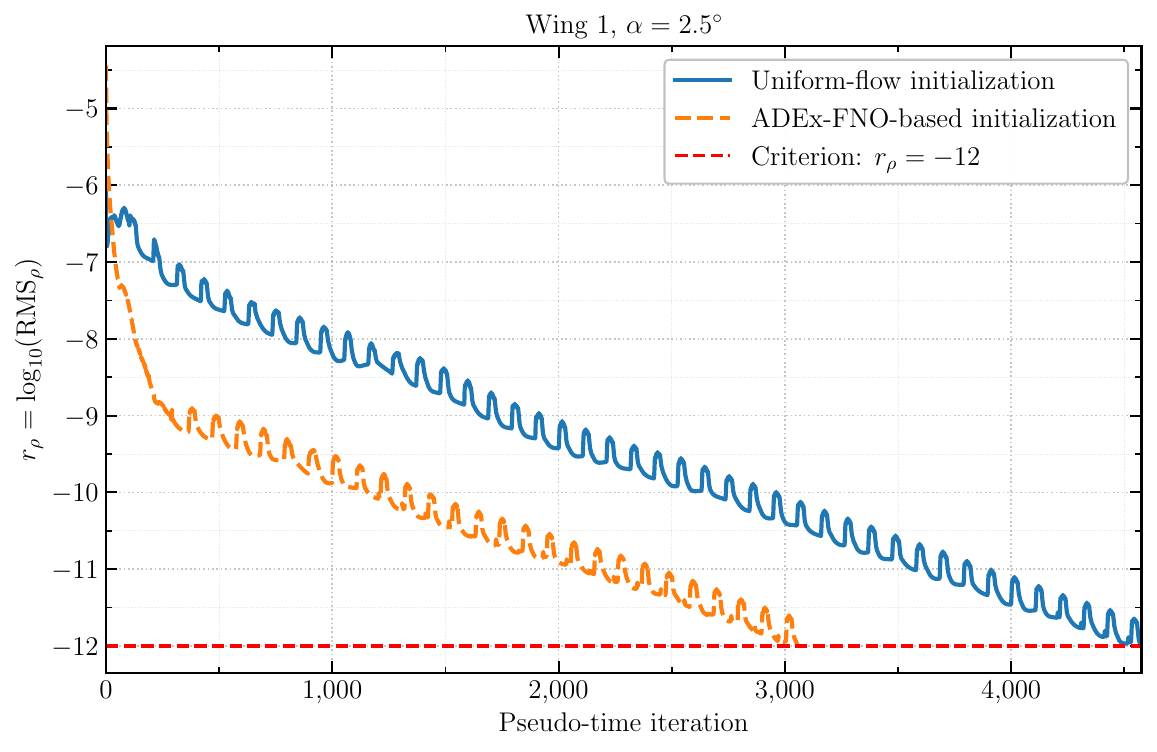}
    \end{subfigure}
    \par\medskip
    \begin{minipage}{\textwidth}
        \centering
        \begin{subfigure}[t]{0.485\textwidth}
            \centering
            \includegraphics[width=\linewidth]{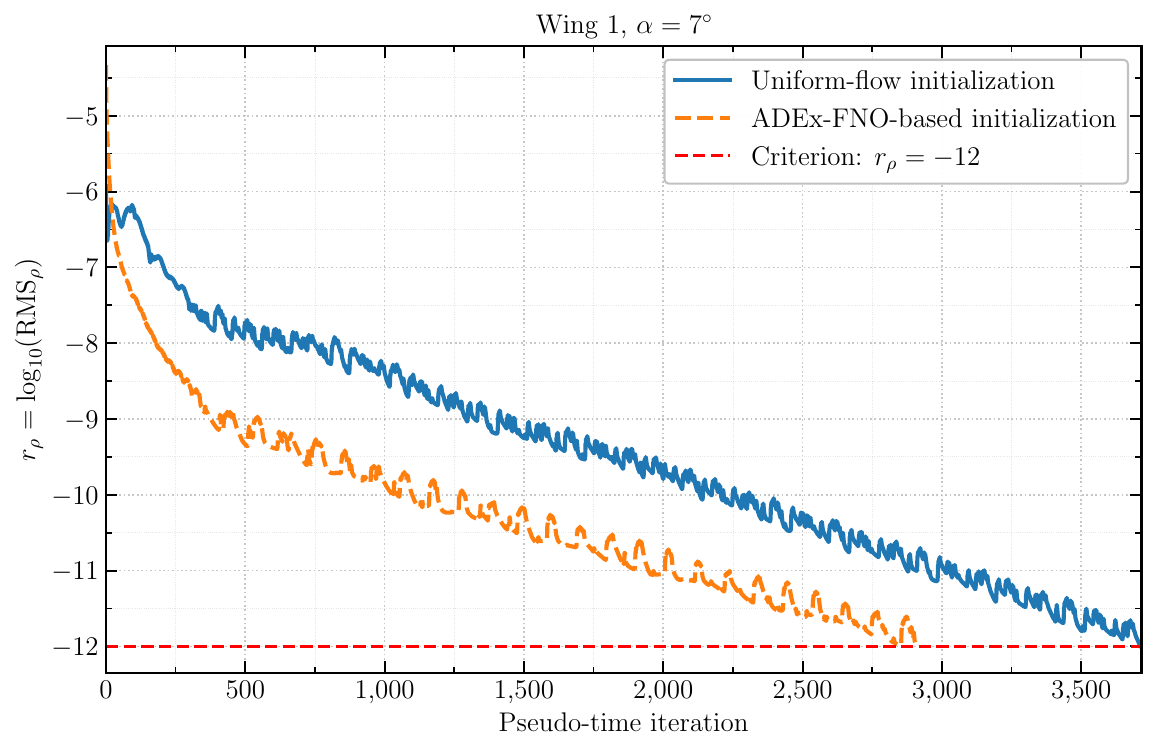}
        \end{subfigure}
    \end{minipage}
    \caption{Paired RMS-density residual histories for the 3D RANS
    Wing~1 cases. Blue solid curves denote uniform-flow initialization, orange
    dashed curves denote ADEx-FNO-based initialization, and the red dashed line is
    the prescribed $r_\rho=-12$ threshold. The panels are ordered by
    increasing incidence.}
    \label{fig:rans_3d_residual_wing1}
\end{figure}

\begin{figure}[H]
    \centering
    \begin{subfigure}[t]{0.485\textwidth}
        \centering
        \includegraphics[width=\linewidth]{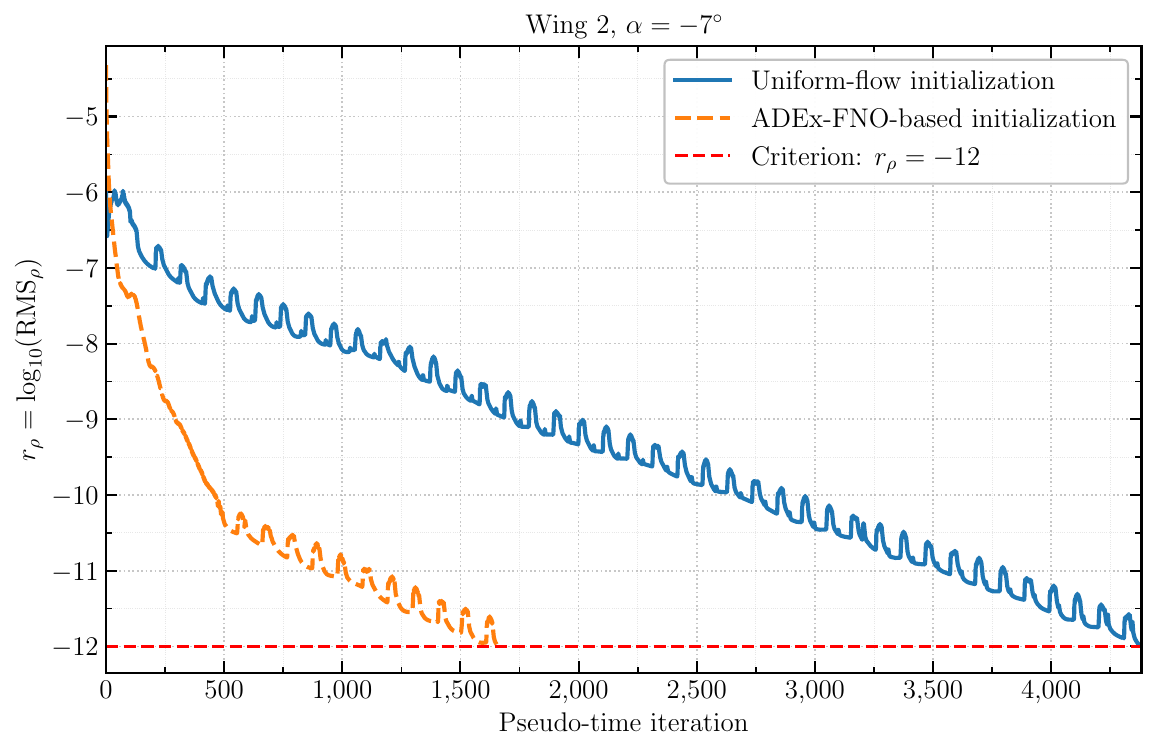}
    \end{subfigure}
    \hfill
    \begin{subfigure}[t]{0.485\textwidth}
        \centering
        \includegraphics[width=\linewidth]{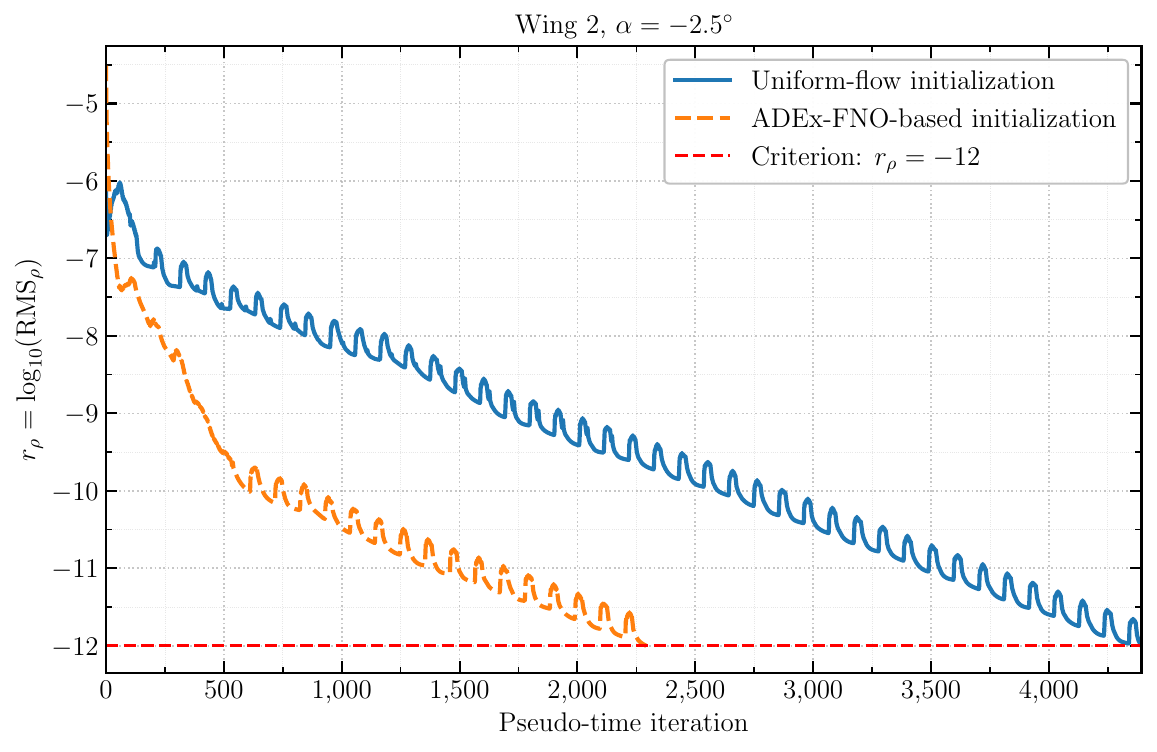}
    \end{subfigure}
    \par\medskip
    \begin{subfigure}[t]{0.485\textwidth}
        \centering
        \includegraphics[width=\linewidth]{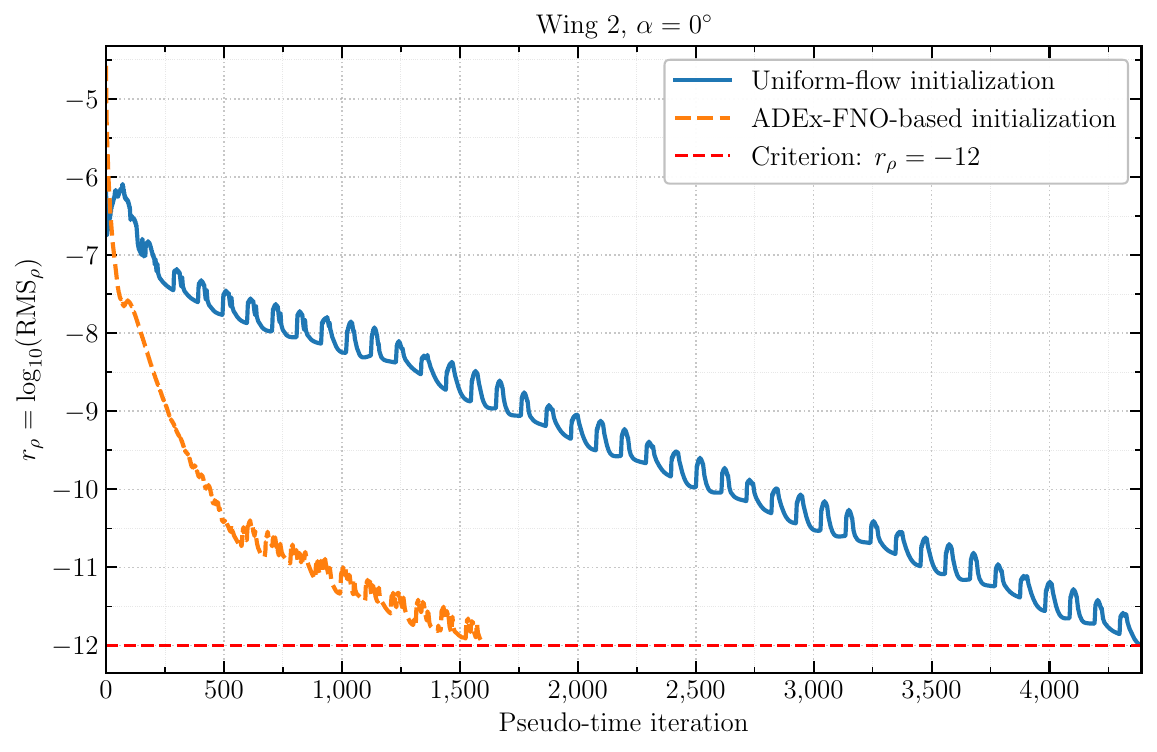}
    \end{subfigure}
    \hfill
    \begin{subfigure}[t]{0.485\textwidth}
        \centering
        \includegraphics[width=\linewidth]{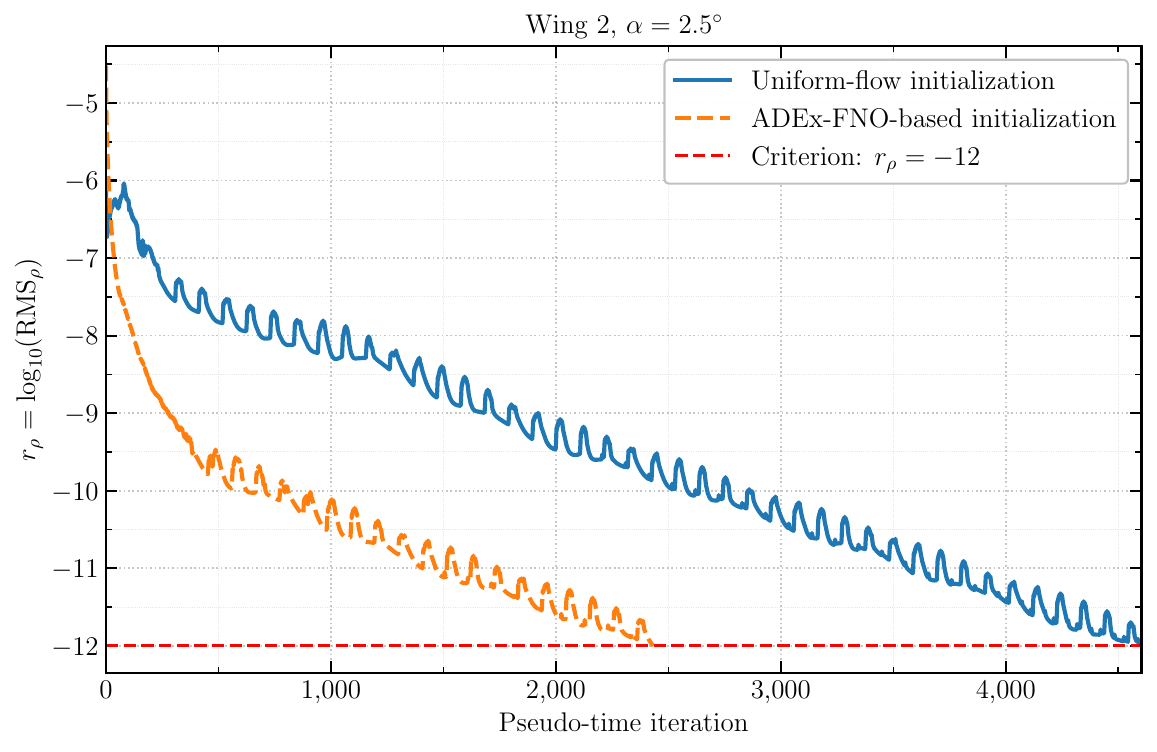}
    \end{subfigure}
    \par\medskip
    \begin{minipage}{\textwidth}
        \centering
        \begin{subfigure}[t]{0.485\textwidth}
            \centering
            \includegraphics[width=\linewidth]{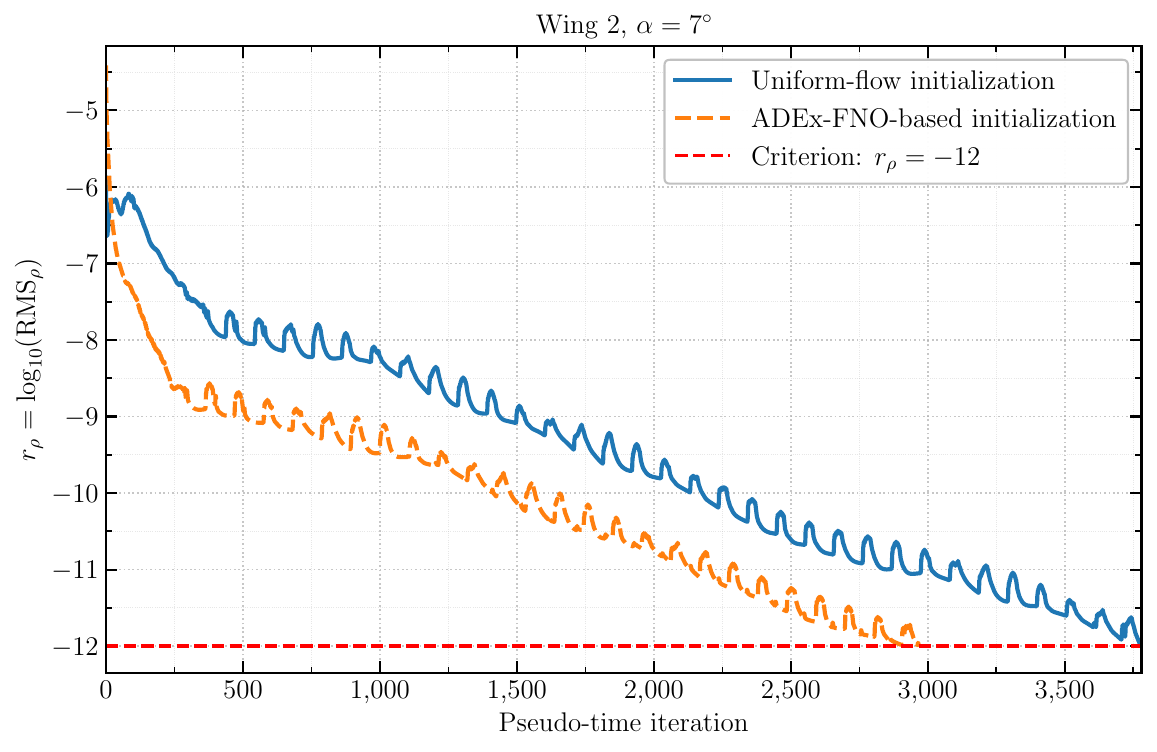}
        \end{subfigure}
    \end{minipage}
    \caption{Paired RMS-density residual histories for the 3D RANS
    Wing~2 cases. The line styles and threshold have the same meaning as in
    Figure~\ref{fig:rans_3d_residual_wing1}.}
    \label{fig:rans_3d_residual_wing2}
\end{figure}

\begin{figure}[H]
    \centering
    \begin{subfigure}[t]{0.485\textwidth}
        \centering
        \includegraphics[width=\linewidth]{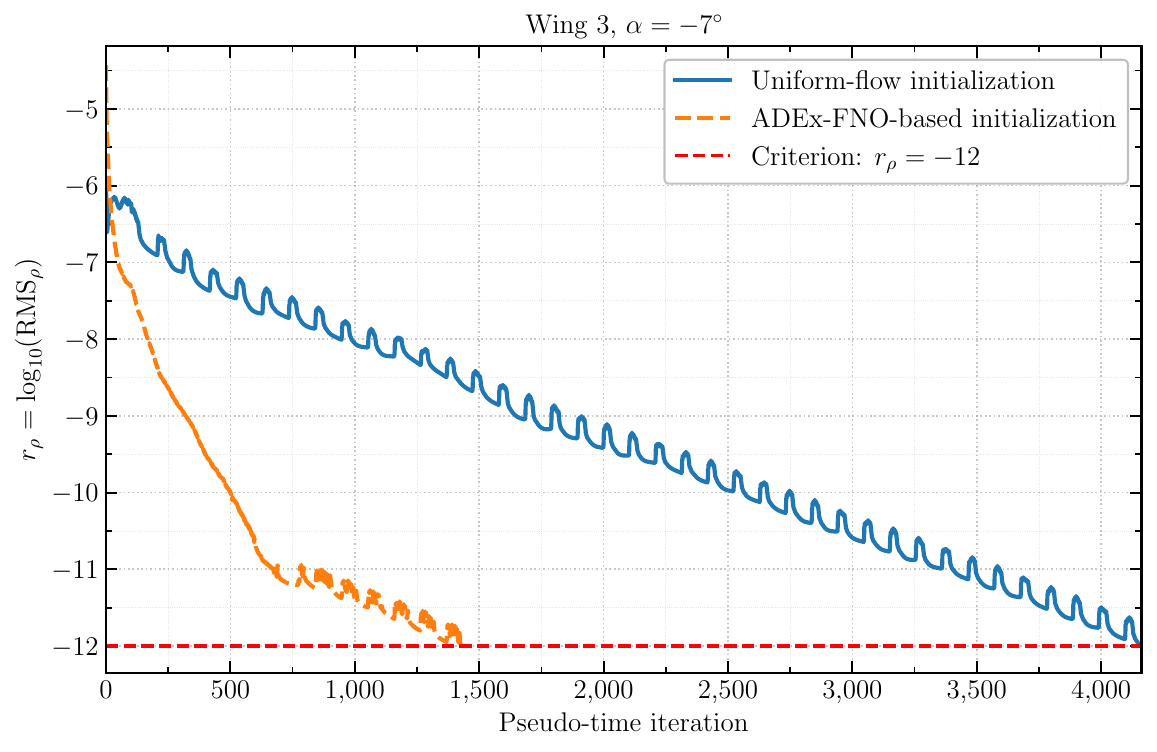}
    \end{subfigure}
    \hfill
    \begin{subfigure}[t]{0.485\textwidth}
        \centering
        \includegraphics[width=\linewidth]{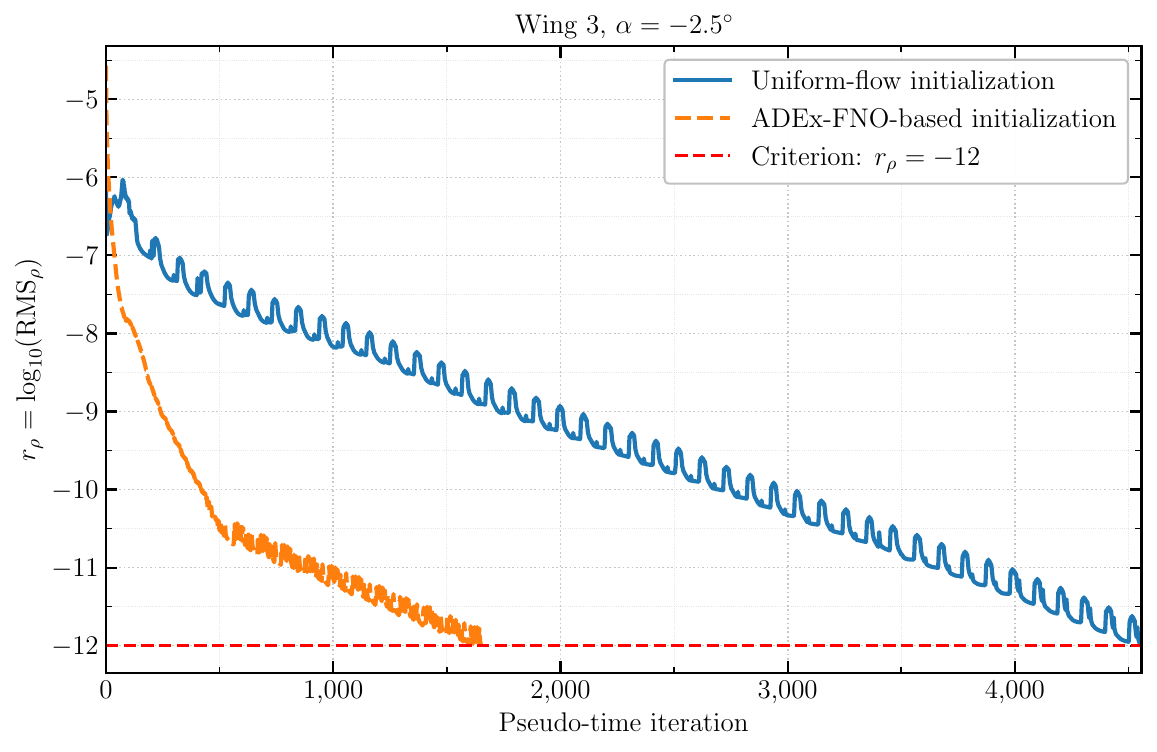}
    \end{subfigure}
    \par\medskip
    \begin{subfigure}[t]{0.485\textwidth}
        \centering
        \includegraphics[width=\linewidth]{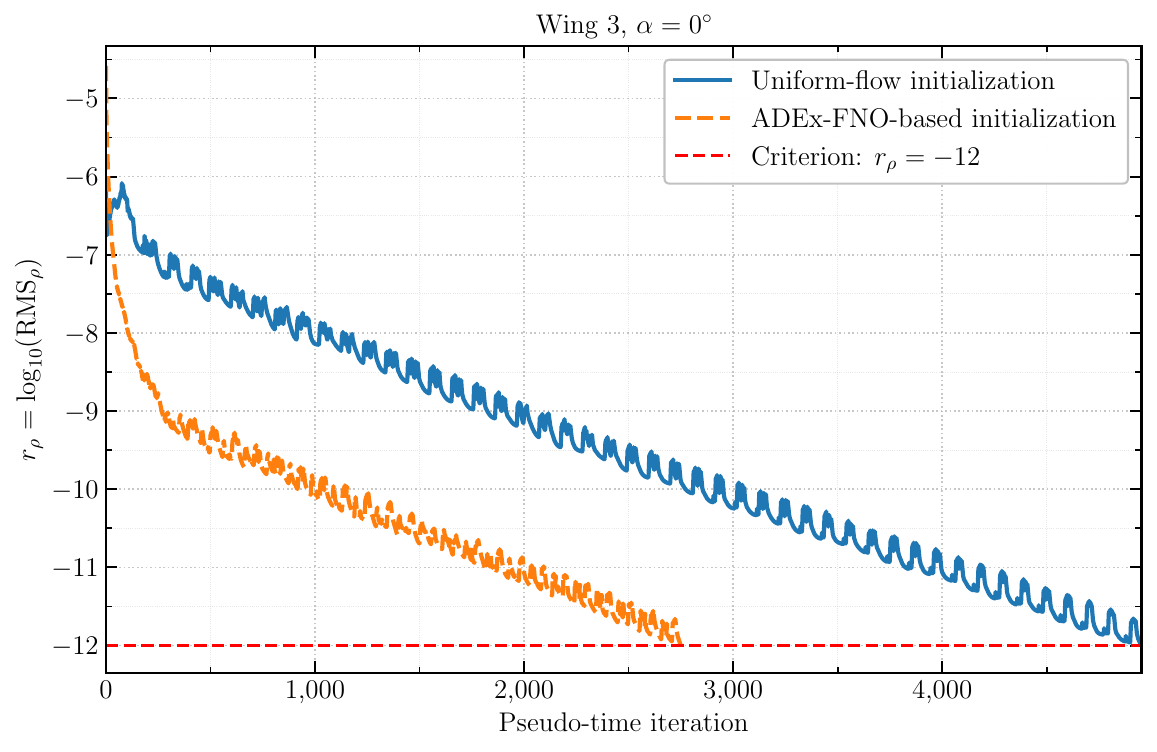}
    \end{subfigure}
    \hfill
    \begin{subfigure}[t]{0.485\textwidth}
        \centering
        \includegraphics[width=\linewidth]{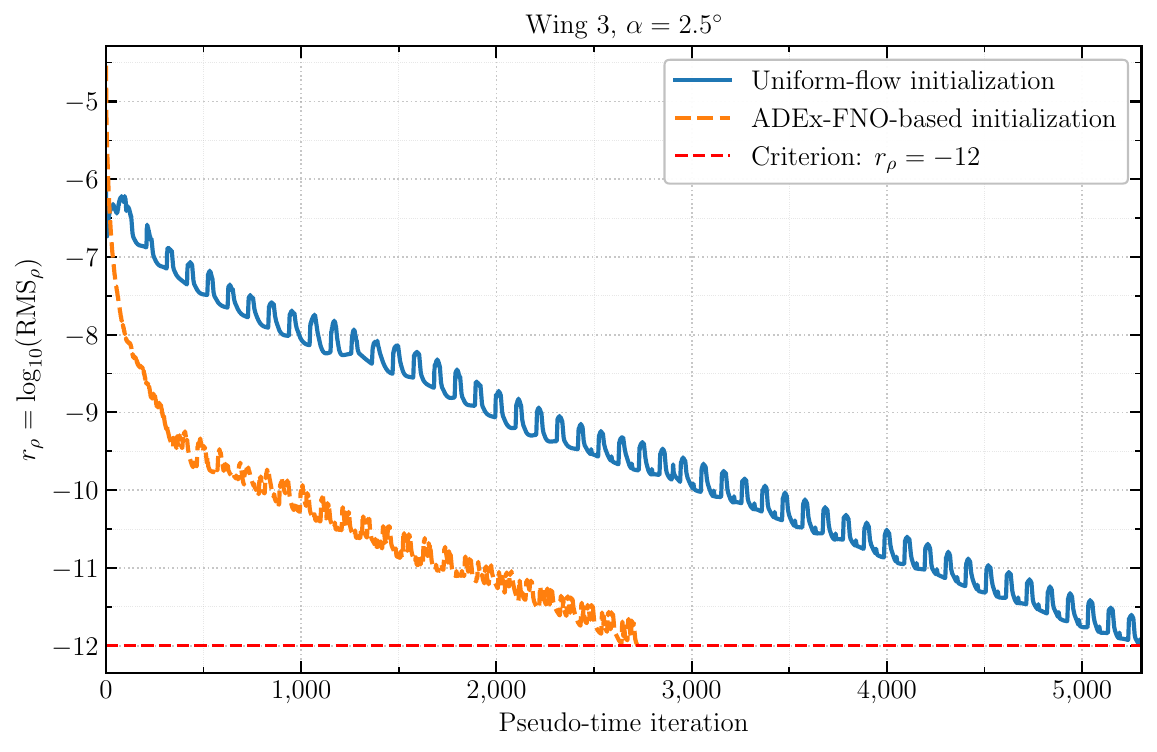}
    \end{subfigure}
    \par\medskip
    \begin{minipage}{\textwidth}
        \centering
        \begin{subfigure}[t]{0.485\textwidth}
            \centering
            \includegraphics[width=\linewidth]{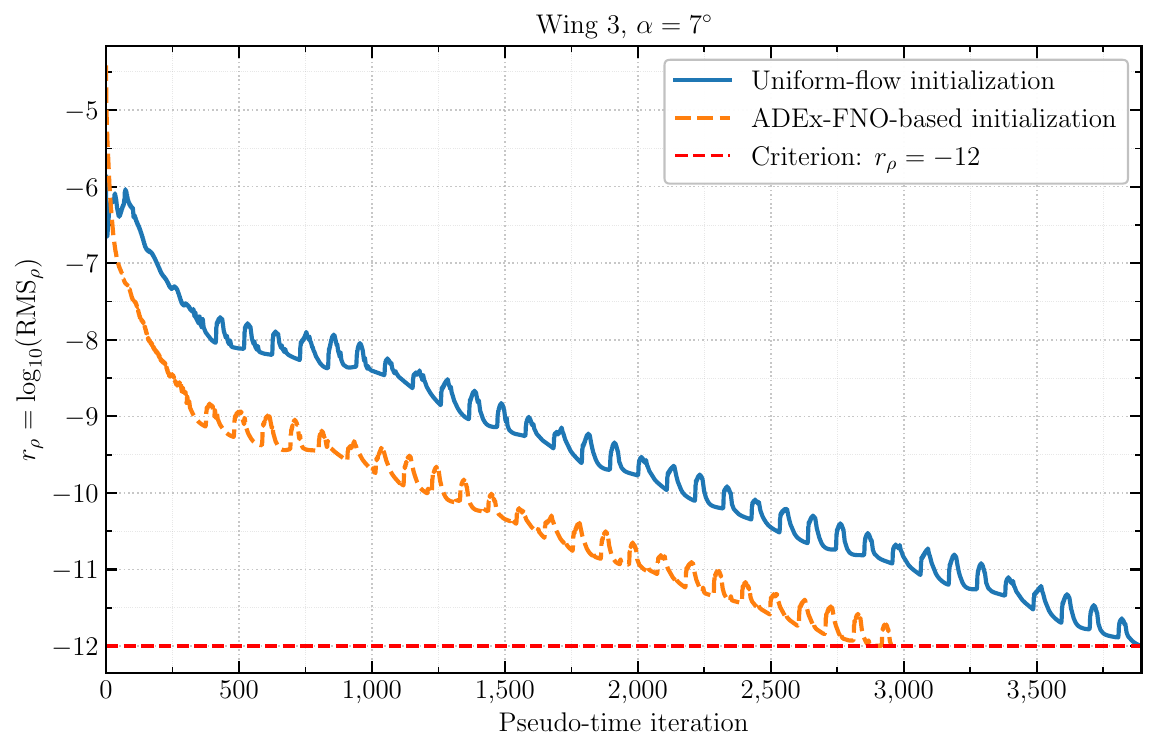}
        \end{subfigure}
    \end{minipage}
    \caption{Paired RMS-density residual histories for the 3D RANS
    Wing~3 cases. The line styles and threshold have the same meaning as in
    Figure~\ref{fig:rans_3d_residual_wing1}.}
    \label{fig:rans_3d_residual_wing3}
\end{figure}

Figures~\ref{fig:rans_3d_residual_wing1},
\ref{fig:rans_3d_residual_wing2}, and
\ref{fig:rans_3d_residual_wing3} show the complete paired residual histories. The ADEx-FNO-based runs first reach
$r_\rho=-10$ using between $34.99\%$ and $84.39\%$ fewer iterations, with a
mean reduction of $65.05\%$. Table~\ref{tab:rans_3d_late_stage} then isolates
the interval between the first $r_\rho<-10$ crossing and final convergence.
This interval is shortened in $10$ pairs and lengthened in $5$ pairs. The
late-stage reduction has a mean of $12.38\%$ and ranges from $-17.68\%$ to
$41.83\%$. The gain is therefore concentrated mainly in the initial and
intermediate convergence transient. The data do not establish a uniformly
improved asymptotic residual-decay rate.

\begin{table}[H]
\centering
\caption{Pseudo-time iterations required after first reaching $r_\rho=-10$ in the 3D RANS wing cases. A negative reduction means that the ADEx-FNO-based run required more iterations than the uniform-flow run between the first $r_\rho<-10$ crossing and final convergence.}
\label{tab:rans_3d_late_stage}
\begin{tabular}{@{}ccrrrrr@{}}
\toprule
Wing & $\alpha$ & \makecell{UF index at\\$r_\rho=-10$} & \makecell{ADEx-FNO \\ idx at $r_\rho=-10$} & \makecell{UF iter\\$-10\rightarrow-12$} & \makecell{ADEx-FNO iter\\$-10\rightarrow-12$} & \makecell{Reduction\\ (\%)} \\
\midrule
1 & -7.0 & 2,577 & 933 & 1,791 & 1,579 & 11.84 \\
1 & -2.5 & 2,548 & 1,330 & 1,880 & 1,886 & -0.32 \\
1 & 0.0 & 2,774 & 1,252 & 1,913 & 1,933 & -1.05 \\
1 & 2.5 & 2,678 & 1,158 & 1,898 & 1,901 & -0.16 \\
1 & 7.0 & 2,132 & 1,045 & 1,584 & 1,864 & -17.68 \\
2 & -7.0 & 2,675 & 460 & 1,707 & 1,198 & 29.82 \\
2 & -2.5 & 2,597 & 604 & 1,796 & 1,688 & 6.01 \\
2 & 0.0 & 2,555 & 442 & 1,835 & 1,162 & 36.68 \\
2 & 2.5 & 2,585 & 626 & 2,018 & 1,807 & 10.46 \\
2 & 7.0 & 2,177 & 1,415 & 1,603 & 1,553 & 3.12 \\
3 & -7.0 & 2,565 & 500 & 1,597 & 929 & 41.83 \\
3 & -2.5 & 2,696 & 420 & 1,861 & 1,231 & 33.85 \\
3 & 0.0 & 2,772 & 964 & 2,181 & 1,786 & 18.11 \\
3 & 2.5 & 3,021 & 774 & 2,283 & 1,949 & 14.63 \\
3 & 7.0 & 2,170 & 1,206 & 1,724 & 1,749 & -1.45 \\
\bottomrule
\end{tabular}
\end{table}

\subsection{Consistency of the terminal aerodynamic coefficients}
\label{app:rans_3d_coefficients}

Table~\ref{tab:rans_3d_terminal_diagnostics} lists the terminal residual and
Cauchy diagnostics. All residuals are below $-12$, and all drag and lift Cauchy
values are below $10^{-7}$, confirming that every coefficient comparison uses
runs satisfying the same complete stopping criterion.

\begin{table}[H]
\centering
\caption{Terminal residual and aerodynamic-coefficient Cauchy diagnostics for the 3D RANS wing cases. Every value satisfies the corresponding threshold.}
\label{tab:rans_3d_terminal_diagnostics}
\resizebox{\textwidth}{!}{%
\begin{tabular}{@{}ccrrrrrr@{}}
\toprule
Wing & $\alpha$ & $r_\rho^{\mathrm{UF}}$ & $r_\rho^{\mathrm{ADEx-FNO}}$ & $C_{\mathrm{C},D}^{\mathrm{UF}}$ & $C_{\mathrm{C},D}^{\mathrm{ADEx-FNO}}$ & $C_{\mathrm{C},L}^{\mathrm{UF}}$ & $C_{\mathrm{C},L}^{\mathrm{ADEx-FNO}}$ \\
\midrule
1 & -7.0 & $-12.0$ & $-12.0$ & $\num{1.688e-8}$ & $\num{1.306e-8}$ & $\num{1.160e-8}$ & $\num{1.184e-8}$ \\
1 & -2.5 & $-12.0$ & $-12.0$ & $\num{2.249e-8}$ & $\num{1.690e-8}$ & $\num{2.701e-8}$ & $\num{2.700e-8}$ \\
1 &  0.0 & $-12.0$ & $-12.0$ & $\num{1.842e-8}$ & $\num{1.909e-8}$ & $\num{2.180e-8}$ & $\num{2.241e-8}$ \\
1 &  2.5 & $-12.0$ & $-12.0$ & $\num{1.943e-8}$ & $\num{1.894e-8}$ & $\num{1.384e-8}$ & $\num{1.191e-8}$ \\
1 &  7.0 & $-12.0$ & $-12.0$ & $\num{2.011e-8}$ & $\num{2.280e-8}$ & $\num{8.727e-9}$ & $\num{9.940e-9}$ \\
2 & -7.0 & $-12.0$ & $-12.0$ & $\num{1.571e-8}$ & $\num{1.487e-8}$ & $\num{3.223e-8}$ & $\num{4.033e-8}$ \\
2 & -2.5 & $-12.0$ & $-12.0$ & $\num{1.670e-8}$ & $\num{1.615e-8}$ & $\num{3.913e-8}$ & $\num{3.593e-8}$ \\
2 &  0.0 & $-12.0$ & $-12.0$ & $\num{2.851e-8}$ & $\num{2.126e-8}$ & $\num{3.593e-8}$ & $\num{2.811e-8}$ \\
2 &  2.5 & $-12.0$ & $-12.0$ & $\num{2.182e-8}$ & $\num{2.271e-8}$ & $\num{1.293e-8}$ & $\num{1.515e-8}$ \\
2 &  7.0 & $-12.0$ & $-12.0$ & $\num{2.668e-8}$ & $\num{2.028e-8}$ & $\num{1.741e-8}$ & $\num{1.294e-8}$ \\
3 & -7.0 & $-12.0$ & $-12.0$ & $\num{1.366e-8}$ & $\num{1.044e-8}$ & $\num{4.813e-9}$ & $\num{3.268e-9}$ \\
3 & -2.5 & $-12.0$ & $-12.0$ & $\num{2.170e-8}$ & $\num{2.030e-8}$ & $\num{2.076e-8}$ & $\num{1.911e-8}$ \\
3 &  0.0 & $-12.0$ & $-12.0$ & $\num{1.875e-8}$ & $\num{2.680e-8}$ & $\num{1.336e-8}$ & $\num{6.742e-9}$ \\
3 &  2.5 & $-12.0$ & $-12.0$ & $\num{1.184e-8}$ & $\num{1.091e-8}$ & $\num{4.595e-9}$ & $\num{4.213e-9}$ \\
3 &  7.0 & $-12.0$ & $-12.0$ & $\num{1.915e-8}$ & $\num{2.022e-8}$ & $\num{1.146e-8}$ & $\num{1.236e-8}$ \\
\bottomrule
\end{tabular}}
\end{table}

\begin{figure}[H]
    \centering
    \begin{subfigure}[t]{0.48\textwidth}
        \centering
        \includegraphics[width=\linewidth]{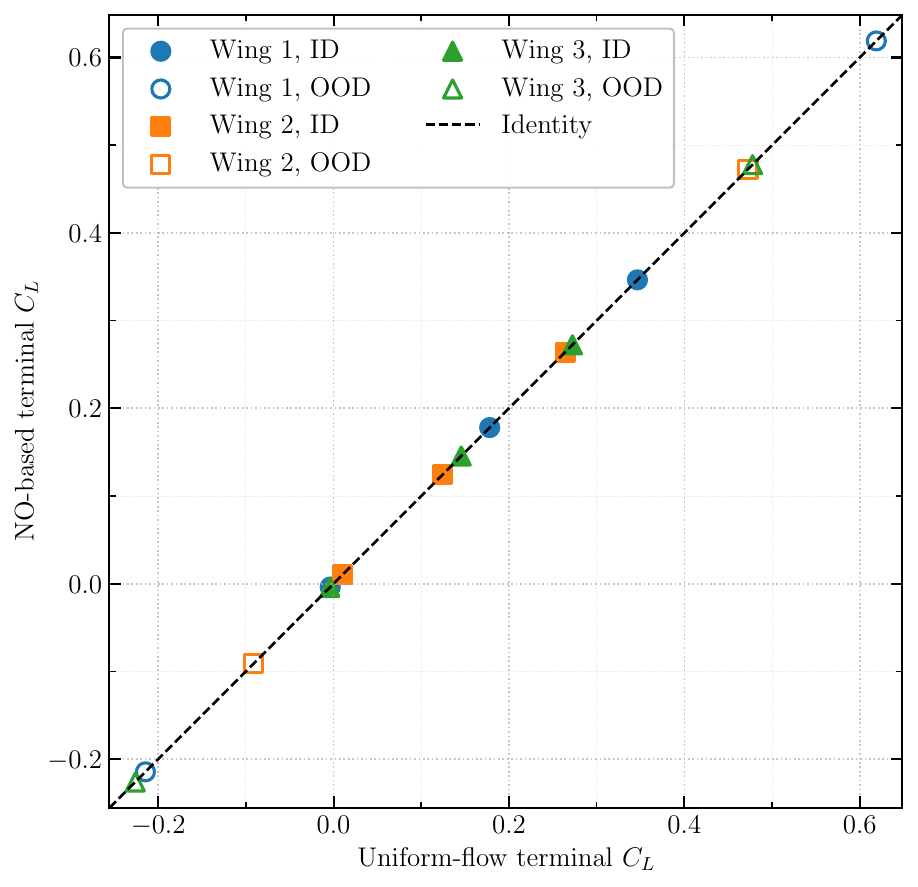}
        \caption{Lift coefficient.}
    \end{subfigure}
    \hfill
    \begin{subfigure}[t]{0.48\textwidth}
        \centering
        \includegraphics[width=\linewidth]{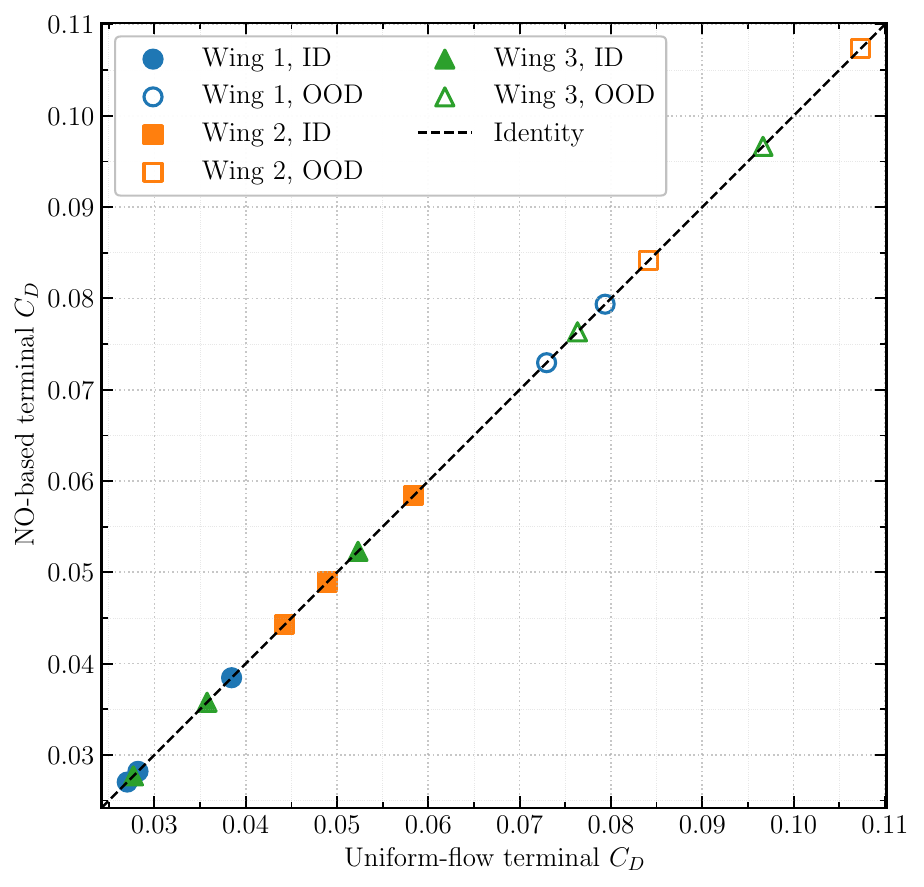}
        \caption{Drag coefficient.}
    \end{subfigure}
    \caption{Consistency of the terminal aerodynamic coefficients for the
    converged 3D RANS wing pairs obtained with uniform-flow and ADEx-FNO-based
    initialization.}
    \label{fig:rans_3d_coefficient_parity}
\end{figure}

Figure~\ref{fig:rans_3d_coefficient_parity} compares the terminal lift and drag
coefficients. The points lie close to the identity line for all three wings.
Table~\ref{tab:rans_3d_terminal_coefficients} gives the complete numerical
values. The maximum absolute differences are $4.1349\times10^{-5}$ for $C_L$
and $4.2101\times10^{-6}$ for $C_D$. The maximum relative differences are
$0.0382\%$ and $0.00577\%$, respectively. The differences are therefore small
relative to the corresponding coefficients and show that the reduction in
solver work does not alter the converged aerodynamic state at the precision
established by the stopping criteria.

\begin{table}[H]
\centering
\caption{Terminal aerodynamic coefficients for the 3D RANS wing pairs. All values correspond to runs satisfying the same residual and Cauchy criteria.}
\label{tab:rans_3d_terminal_coefficients}
\resizebox{\textwidth}{!}{%
\begin{tabular}{@{}ccrrrrrr@{}}
\toprule
Wing & $\alpha$ & $C_L^{\mathrm{UF}}$ & $C_L^{\mathrm{ADEx-FNO}}$ & $|\Delta C_L|$ & $C_D^{\mathrm{UF}}$ & $C_D^{\mathrm{ADEx-FNO}}$ & $|\Delta C_D|$ \\
\midrule
1 & -7.0 & $-0.214$ & $-0.215$ & $\num{4.135e-5}$ & $0.0730$ & $0.0730$ & $\num{4.210e-6}$ \\
1 & -2.5 & $-0.004$ & $-0.004$ & $\num{1.567e-7}$ & $0.0282$ & $0.0282$ & $\num{8.724e-8}$ \\
1 &  0.0 & $ 0.178$ & $ 0.178$ & $\num{4.630e-8}$ & $0.0270$ & $0.0270$ & $\num{1.449e-7}$ \\
1 &  2.5 & $ 0.347$ & $ 0.347$ & $\num{2.434e-7}$ & $0.0385$ & $0.0385$ & $\num{1.176e-8}$ \\
1 &  7.0 & $ 0.619$ & $ 0.619$ & $\num{5.270e-8}$ & $0.0794$ & $0.0794$ & $\num{3.355e-8}$ \\
2 & -7.0 & $-0.091$ & $-0.091$ & $\num{2.390e-5}$ & $0.0841$ & $0.0841$ & $\num{3.183e-6}$ \\
2 & -2.5 & $ 0.010$ & $ 0.010$ & $\num{1.837e-6}$ & $0.0489$ & $0.0489$ & $\num{2.195e-7}$ \\
2 &  0.0 & $ 0.124$ & $ 0.124$ & $\num{3.466e-7}$ & $0.0443$ & $0.0443$ & $\num{1.423e-7}$ \\
2 &  2.5 & $ 0.264$ & $ 0.264$ & $\num{2.900e-9}$ & $0.0584$ & $0.0584$ & $\num{1.322e-8}$ \\
2 &  7.0 & $ 0.472$ & $ 0.472$ & $\num{2.328e-6}$ & $0.1073$ & $0.1073$ & $\num{7.700e-8}$ \\
3 & -7.0 & $-0.226$ & $-0.226$ & $\num{4.581e-7}$ & $0.0763$ & $0.0763$ & $\num{2.437e-6}$ \\
3 & -2.5 & $-0.004$ & $-0.004$ & $\num{1.615e-6}$ & $0.0277$ & $0.0277$ & $\num{1.536e-6}$ \\
3 &  0.0 & $ 0.146$ & $ 0.146$ & $\num{5.260e-8}$ & $0.0358$ & $0.0358$ & $\num{9.979e-8}$ \\
3 &  2.5 & $ 0.272$ & $ 0.272$ & $\num{2.730e-8}$ & $0.0523$ & $0.0523$ & $\num{6.058e-8}$ \\
3 &  7.0 & $ 0.478$ & $ 0.478$ & $\num{1.692e-6}$ & $0.0966$ & $0.0966$ & $\num{1.182e-7}$ \\
\bottomrule
\end{tabular}}
\end{table}

\section{Physical and numerical parameters of the URANS test cases}
\label{app:urans_configurations}

Table~\ref{tab:urans_su2_configurations} summarizes the physical and numerical
settings used in the 2D and 3D URANS studies. The numerical-method acronyms
follow the definitions given in~\ref{app:rans_configurations}.

\begin{table}[H]
\centering
\caption{Principal physical and numerical parameters of the 2D and 3D URANS
cases.}
\label{tab:urans_su2_configurations}
\renewcommand{\arraystretch}{1.06}
\setlength{\tabcolsep}{4pt}
\begin{tabularx}{\linewidth}{@{}
    >{\raggedright\arraybackslash}p{0.26\linewidth}
    >{\raggedright\arraybackslash}X
    >{\raggedright\arraybackslash}X
@{}}
\toprule
Parameter & \textbf{2D URANS} & \textbf{3D URANS} \\
\midrule

Reference configuration
& NACA~0012 airfoil
& Held-out transonic wing \\

Flow regime
& Low-Mach-number external flow
& Transonic external flow \\

Mach number
& $M_\infty=0.15$
& $M_\infty=0.8395$ \\

Reynolds number
& $Re=6.0\times10^{6}$
& $Re=1.5\times10^{5}$ \\

Reference length
& $L_{\mathrm{ref}}=1$
& $L_{\mathrm{ref}}=0.64607$ \\

Freestream temperature
& $T_\infty=300\,\mathrm{K}$
& $T_\infty=288.15\,\mathrm{K}$ \\

Boundary conditions
& Adiabatic wall and far field
& Adiabatic wall, far field, and symmetry plane \\

\addlinespace[2pt]

Turbulence closure
& Spalart--Allmaras-\texttt{noft2}
& Spalart--Allmaras-\texttt{noft2} \\

Mean-flow discretization
& Second-order JST
& Second-order JST \\

Turbulence discretization
& Scalar-upwind
& Scalar-upwind \\

Physical-time integration
& Second-order dual-time stepping
& Second-order dual-time stepping \\

Physical time step
& $\Delta t=5\times10^{-4}\,\mathrm{s}$
& $\Delta t=5\times10^{-4}\,\mathrm{s}$ \\

Maximum inner iterations
& $80$
& $80$ \\

Linear solver and preconditioner
& ILU(0)-preconditioned FGMRES
& ILU(0)-preconditioned FGMRES \\

Maximum linear iterations
& $20$
& $20$ \\

\addlinespace[2pt]

Inner convergence threshold
& $r_\rho<-11$
& $r_\rho<-11$ \\

Window function
& Squared Hann
& Squared Hann \\

Monitored time averages
& $C_D$ and $C_L$
& $C_D$ and $C_L$ \\

Window start
& After $300$ physical-time iterations
& After $150$ physical-time iterations \\

Cauchy tolerance
& $10^{-4}$
& $10^{-5}$ \\

Number of Cauchy entries
& $50$
& $50$ \\

\bottomrule
\end{tabularx}
\end{table}

\section{Details of the 2D URANS initialization study}
\label{app:urans_2d_appendix}

\subsection{Physical-time, pseudo-time, and linear-solver effort}
\label{app:urans_2d_effort}

The URANS solver has a three-level work hierarchy: physical-time advancement,
inner pseudo-time iterations at each physical step, and linear iterations for
the mean-flow and turbulence systems within each inner solve. Table~
\ref{tab:urans_2d_work_detail} reports the cumulative work after activation
of the windowed averaging procedure and over the complete simulations.

\begin{table}[H]
\centering
\caption{Detailed solver work for the 2D NACA~0012 URANS case. The combined count
$L_\Sigma=L_\mathrm{flow}+L_\mathrm{turb}$ is an unweighted
algebraic-work diagnostic.}
\label{tab:urans_2d_work_detail}
\begin{tabularx}{0.98\textwidth}{@{}Xrrr@{}}
    \toprule
    Quantity
    & \makecell{Uniform-flow\\initialization}
    & \makecell{ADEx-FNO-based\\initialization}
    & $G_Q$ (\%) \\
    \midrule
    \multicolumn{4}{@{}l}{\textit{After window activation}} \\
    Physical-time advances
    & 189
    & 137
    & 27.51 \\
    Inner pseudo-time iterations
    & 1,649
    & 1,223
    & 25.83 \\
    Mean-flow linear iterations
    & 32,980
    & 23,524
    & 28.67 \\
    Turbulence linear iterations
    & 23,471
    & 6,179
    & 73.67 \\
    Unweighted combined linear iterations
    & 56,451
    & 29,703
    & 47.38 \\
    \addlinespace[3pt]
    \multicolumn{4}{@{}l}{\textit{Complete simulation}} \\
    Executed physical-time iterations
    & 490
    & 439
    & 10.41 \\
    Inner pseudo-time iterations
    & 10,055
    & 7,158
    & 28.81 \\
    Mean-flow linear iterations
    & 201,100
    & 141,844
    & 29.47 \\
    Turbulence linear iterations
    & 138,175
    & 38,164
    & 72.38 \\
    Unweighted combined linear iterations
    & 339,275
    & 180,008
    & 46.94 \\
    \bottomrule
\end{tabularx}
\end{table}

The post-window mean number of inner iterations per physical-time advance is
8.72 for uniform-flow initialization and
8.93 for ADEx-FNO-based initialization. Thus, the
ADEx-FNO-based run does not reduce the pseudo-time count of an individual remaining
physical step. Its reduction in cumulative inner iterations follows from
reaching the time-convergence criterion after
52
fewer physical-time advances.

The linear-solver behavior provides an additional distinction. After window
activation, the mean-flow linear count is reduced by
$28.67\%$, while the turbulence linear count
is reduced by $73.67\%$. The latter
reduction is substantially larger because the average turbulence-solver count
per inner iteration decreases from 14.23
to 5.05. Consequently, the combined
linear work decreases by $47.38\%$,
exceeding the reduction in the number of inner iterations. The same pattern
persists over the complete simulations, for which $L_\Sigma$ decreases by
$46.94\%$.

\begin{figure}[H]
    \centering
    \includegraphics[width=0.75\linewidth]
    {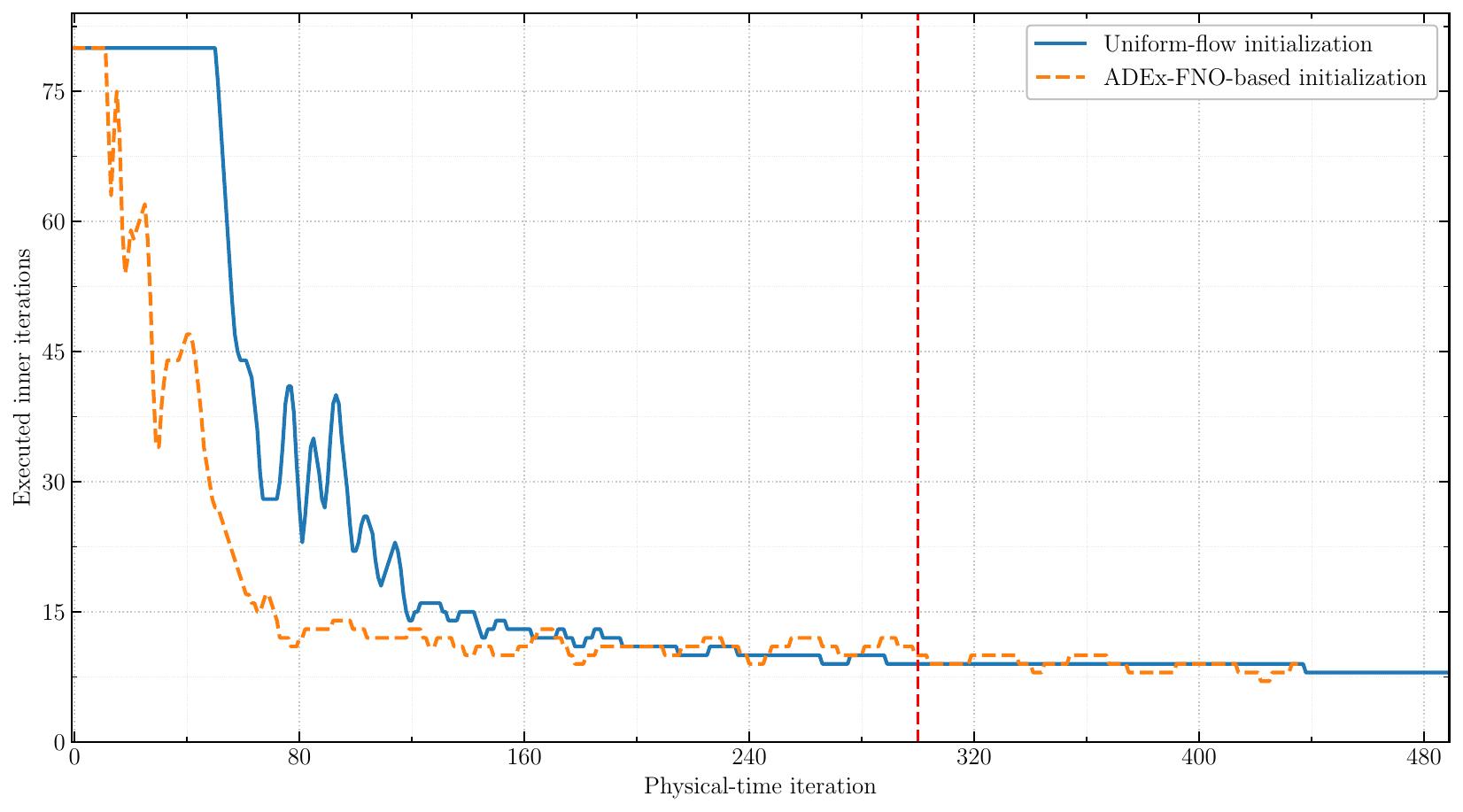}
    \caption{Executed inner pseudo-time iterations at each physical-time
    iteration for the 2D NACA~0012 URANS case. Blue solid and orange dashed
    curves denote uniform-flow and ADEx-FNO-based initialization, respectively. The
    vertical red dashed line marks activation of the averaging window after
    $300$ physical-time iterations.}
    \label{fig:urans_2d_inner_iterations}
\end{figure}

Figure~\ref{fig:urans_2d_inner_iterations} shows that both calculations
initially reach the imposed limit of $80$ inner iterations. The limit is reached
at $51$ physical-time iterations for the uniform-flow
calculation and at $13$ physical-time iterations for the ADEx-FNO-based
calculation; it is not reached after window activation in either run. The long
uniform-flow transient therefore accounts for a large fraction of the complete
pseudo-time work.

\begin{figure}[H]
    \centering
    \includegraphics[width=0.75\linewidth]
    {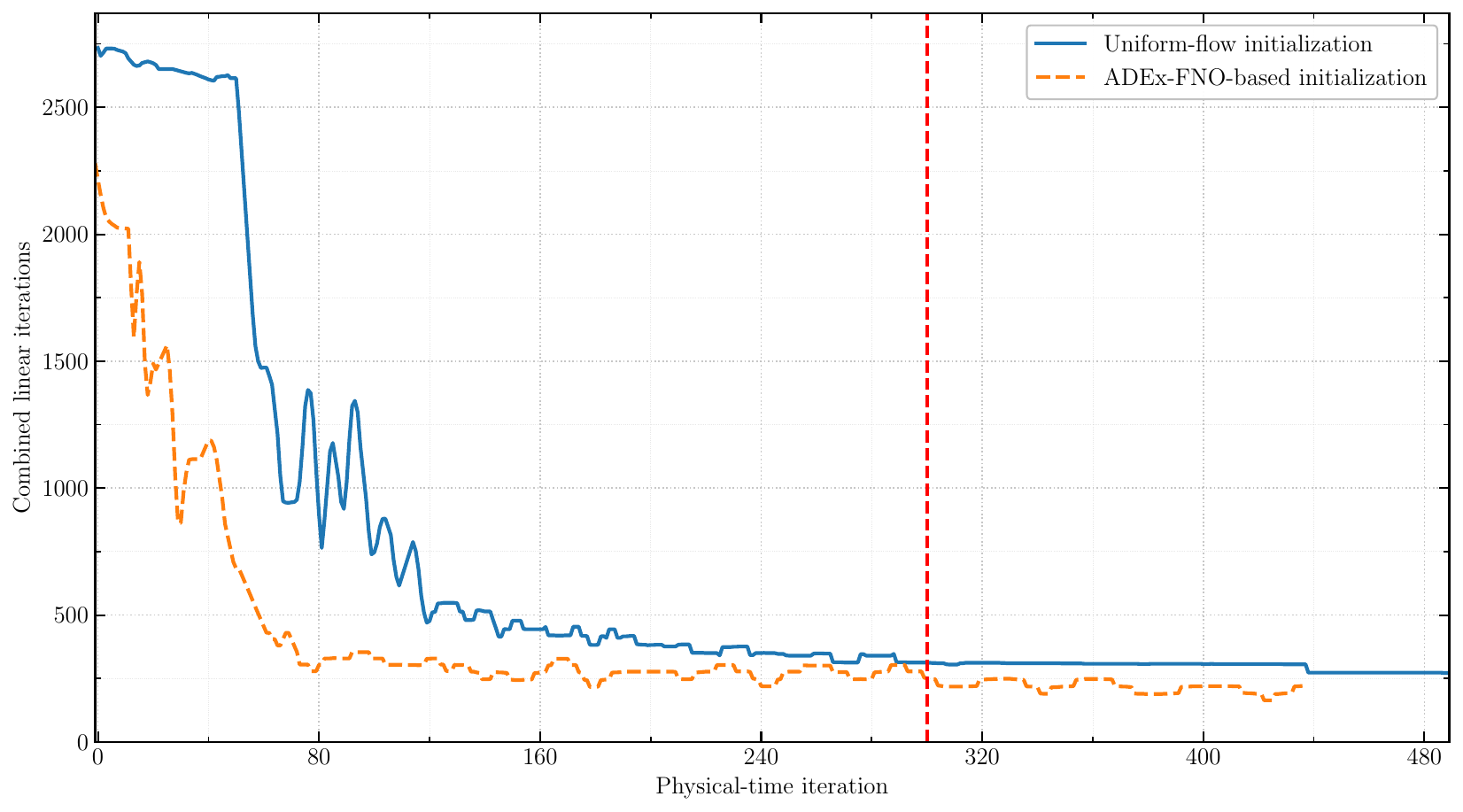}
    \caption{Unweighted combined linear iterations accumulated at each
    physical-time iteration for the 2D NACA~0012 URANS case. The line styles
    and averaging-window marker have the same meaning as in
    Figure~\ref{fig:urans_2d_inner_iterations}.}
    \label{fig:urans_2d_linear_per_step}
\end{figure}

Figure~\ref{fig:urans_2d_linear_per_step} confirms that the ADEx-FNO-based
calculation reduces not only the number of physical and inner iterations but
also the linear-algebra work accumulated within the remaining pseudo-time
solves. The reduction is particularly pronounced for the turbulence system.

\subsection{Windowed time-convergence diagnostics}
\label{app:urans_2d_window}

Table~\ref{tab:urans_2d_cauchy_terminal} reports the configured window
starts, the first iterations at which both Cauchy diagnostics are below
$10^{-4}$, and the actual solver-termination iterations. In each run, SU2
exits one physical-time iteration after the first simultaneous
crossing. The lift diagnostic controls termination in both calculations.

\begin{table}[H]
\centering
\caption{Window-start, Cauchy-crossing, and termination data for the 2D
NACA~0012 URANS calculations.}
\label{tab:urans_2d_cauchy_terminal}
\begin{tabular}{@{}lrrrr@{}}
    \toprule
    Initialization
    & \makecell{Window\\start}
    & \makecell{First simultaneous\\Cauchy crossing}
    & \makecell{Solver\\termination}
    & \makecell{Terminal diagnostics\\$(C_D,C_L)$} \\
    \midrule
    Uniform-flow
    & 300
    & 488
    & 489
    & $(5.9590\times10^{-6},\,9.9431\times10^{-5})$ \\
    ADEx-FNO-based
    & 800
    & 936
    & 937
    & $(4.0894\times10^{-6},\,9.8166\times10^{-5})$ \\
    \bottomrule
\end{tabular}
\end{table}

\begin{figure}[H]
    \centering
    \begin{subfigure}[t]{0.48\textwidth}
        \centering
        \includegraphics[width=\linewidth]
        {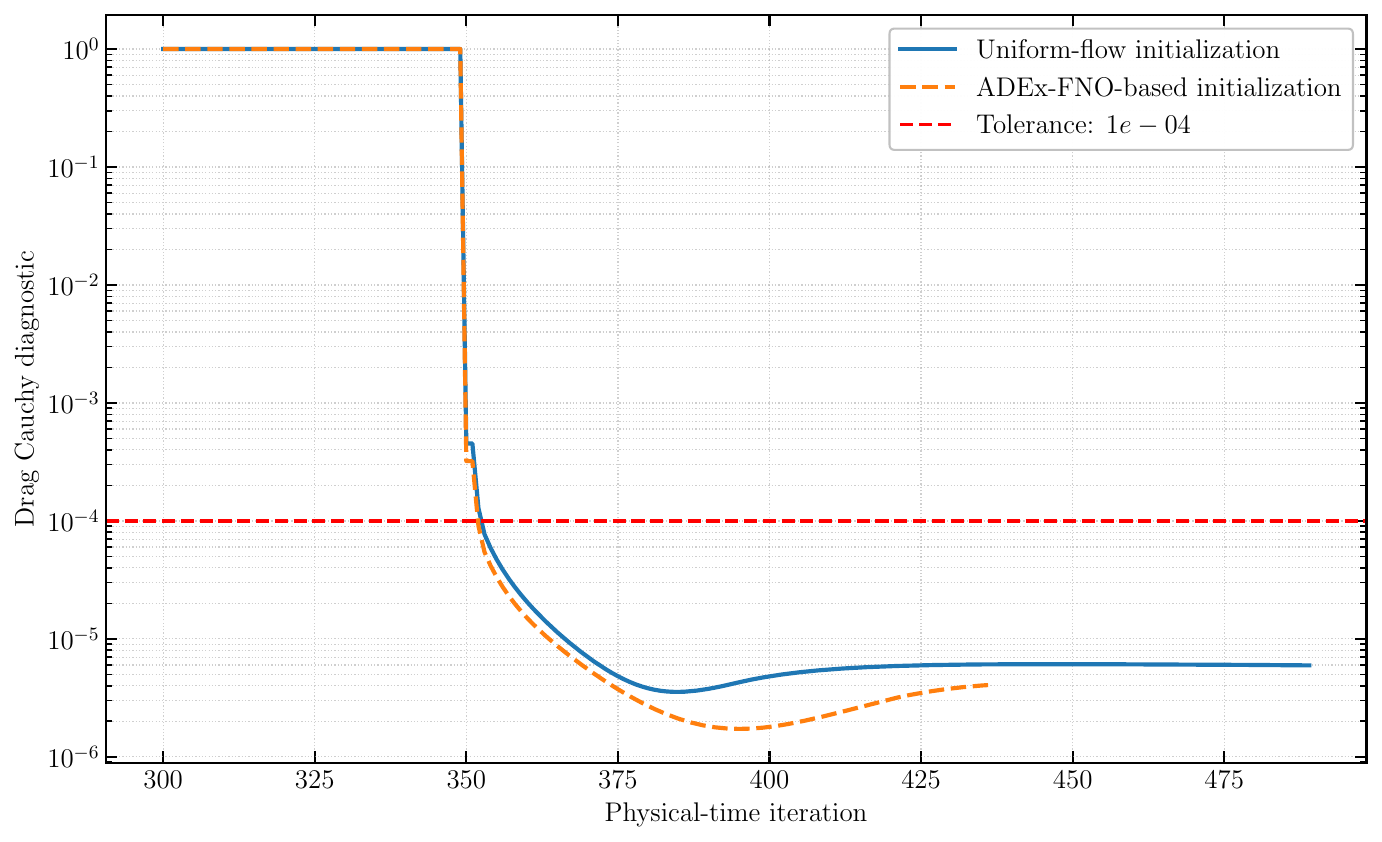}
        \caption{Drag-coefficient diagnostic.}
    \end{subfigure}
    \hfill
    \begin{subfigure}[t]{0.48\textwidth}
        \centering
        \includegraphics[width=\linewidth]
        {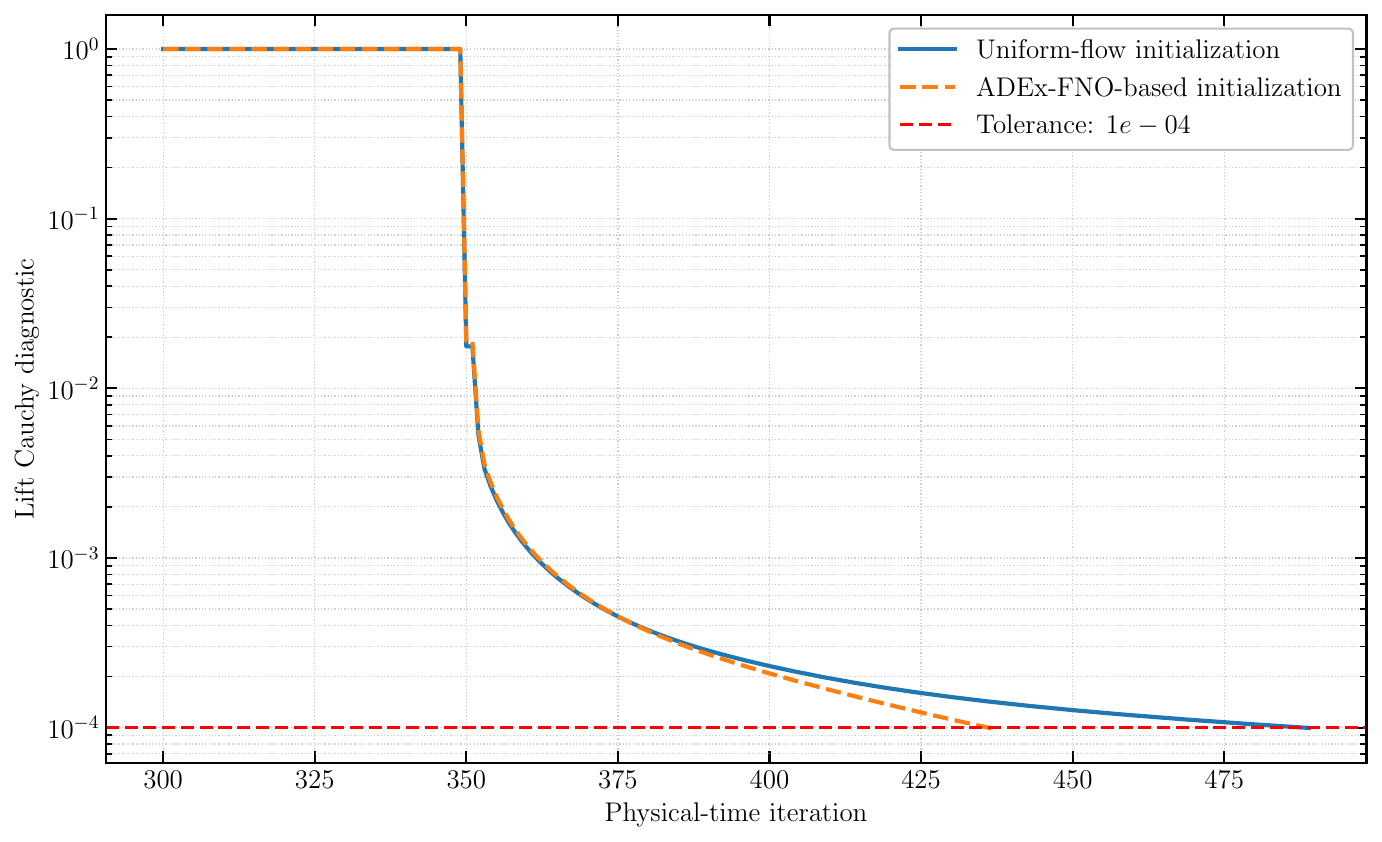}
        \caption{Lift-coefficient diagnostic.}
    \end{subfigure}
    \caption{Cauchy diagnostics for the squared-Hann-windowed aerodynamic
    coefficients in the 2D NACA~0012 URANS case. The red dashed line marks the
    prescribed tolerance $10^{-4}$.}
    \label{fig:urans_2d_cauchy}
\end{figure}

Figure~\ref{fig:urans_2d_cauchy} shows that the drag diagnostic falls below the
tolerance earlier, whereas the lift diagnostic determines the final stopping
point. The terminal RMS-density residuals are
-11.1411 and -11.0827,
and the corresponding minimum values are
-11.3162 and -11.2280.
Thus, the final inner solves in both calculations also satisfy the common
$r_\rho<-11$ threshold.

\subsection{Consistency of the time-averaged aerodynamic coefficients}
\label{app:urans_2d_statistics}

Table~\ref{tab:urans_2d_terminal_statistics} compares the terminal
instantaneous coefficients and the terminal squared-Hann averages. Although
both calculations independently satisfy the same windowed Cauchy tolerance,
the terminal averages are not mutually consistent to the level observed for
the steady RANS cases.

\begin{table}[H]
\centering
\caption{Terminal instantaneous and squared-Hann-windowed aerodynamic coefficients for the 2D NACA~0012 URANS case. Relative differences use the magnitude of the uniform-flow value as reference.}
\label{tab:urans_2d_terminal_statistics}
\begin{tabular}{@{}lrrrr@{}}
    \toprule
    Quantity
    & \makecell{Uniform-flow\\initialization}
    & \makecell{ADEx-FNO-based\\initialization}
    & Absolute difference
    & Relative difference (\%) \\
    \midrule
    Terminal $C_L$
    & $-0.726$
    & $-0.762$
    & $\num{3.664e-2}$
    & 5.05 \\
    Terminal $C_D$
    & $0.0158$
    & $0.0120$
    & $\num{3.828e-3}$
    & 24.16 \\
    Terminal $\overline{C_L}$
    & $-0.711$
    & $-0.753$
    & $\num{4.246e-2}$
    & 5.97 \\
    Terminal $\overline{C_D}$
    & $0.0168$
    & $0.0122$
    & $\num{4.550e-3}$
    & 27.12 \\
    \bottomrule
\end{tabular}
\end{table}

\begin{figure}[H]
    \centering
    \begin{subfigure}[t]{0.48\textwidth}
        \centering
        \includegraphics[width=\linewidth]
        {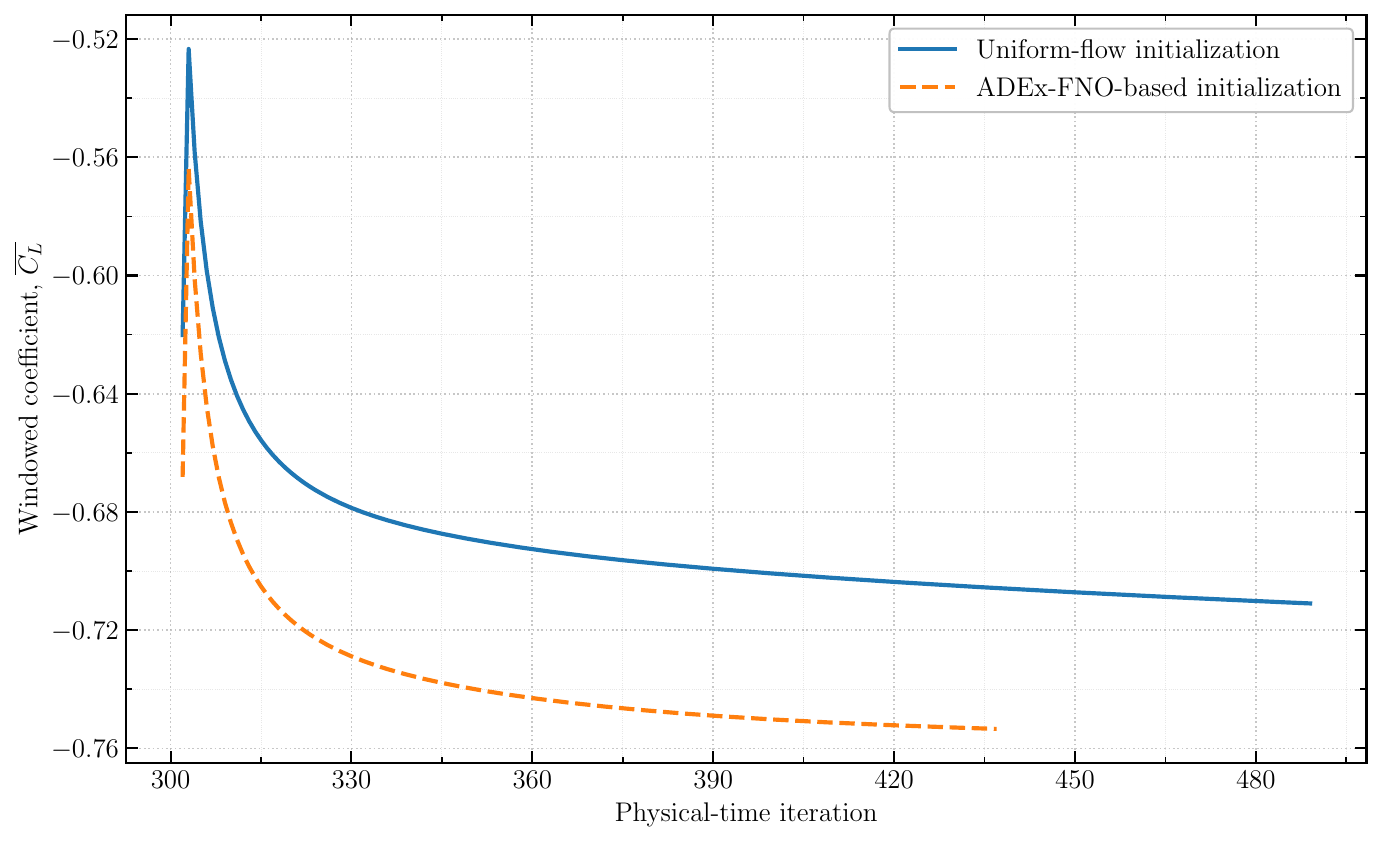}
        \caption{Windowed lift coefficient.}
    \end{subfigure}
    \hfill
    \begin{subfigure}[t]{0.48\textwidth}
        \centering
        \includegraphics[width=\linewidth]
        {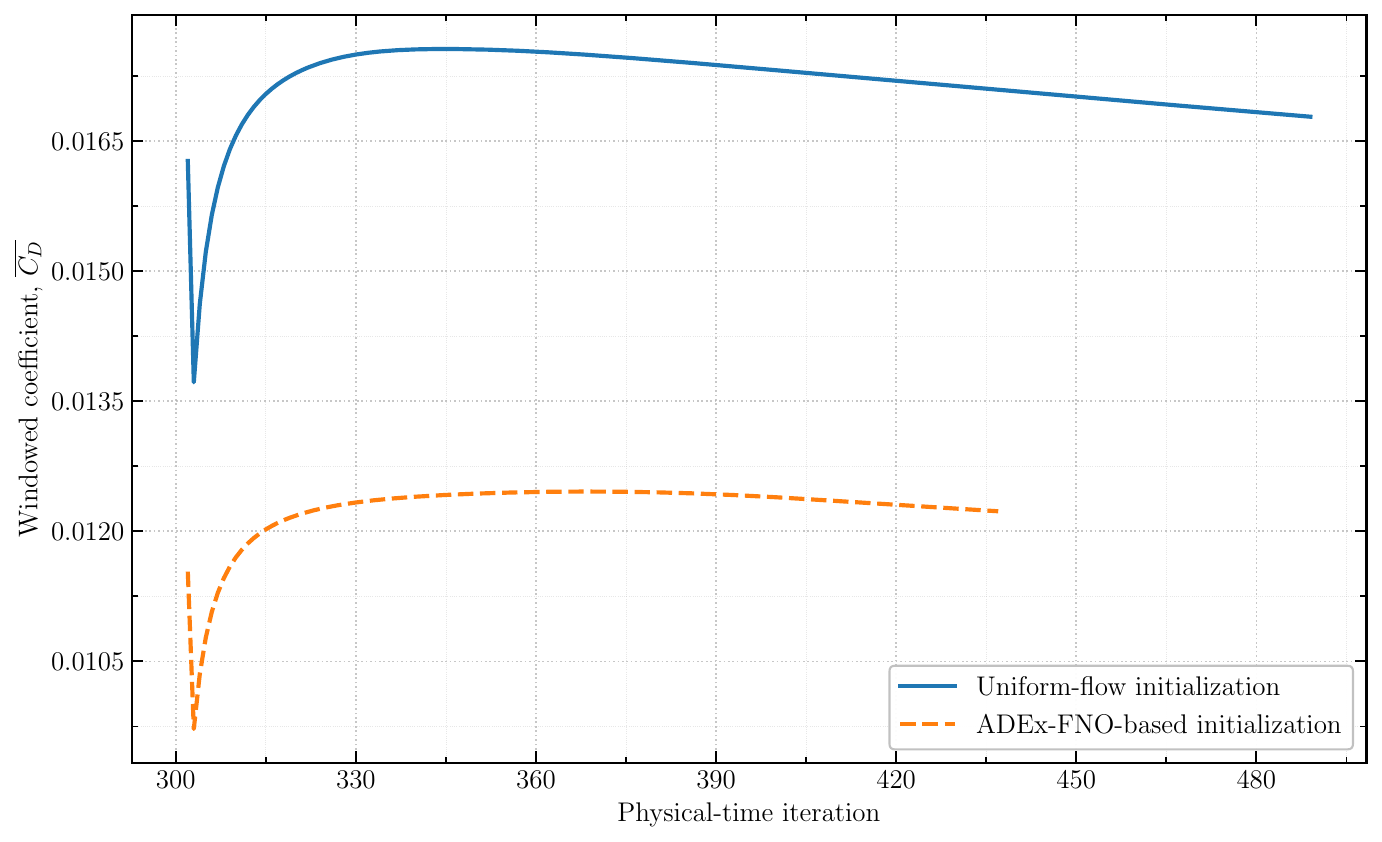}
        \caption{Windowed drag coefficient.}
    \end{subfigure}
    \caption{Squared-Hann-windowed aerodynamic coefficients for the 2D
    NACA~0012 URANS case after activation of the averaging window.}
    \label{fig:urans_2d_tavg}
\end{figure}

Figure~\ref{fig:urans_2d_tavg} makes the discrepancy explicit. The ADEx-FNO-based
averages meet the self-convergence criterion after a shorter physical-time
interval, but their terminal values differ from those of the uniform-flow run
by $5.97\%$ in lift and
$27.12\%$ in drag. Consequently, the current
evidence supports a reduction in solver work to the configured stopping
condition, not yet acceleration to an independently verified common
long-time-averaged state. A longer common integration interval or a stricter
cross-initialization consistency requirement is needed to complete that
validation.

\section{Details of the 3D URANS initialization study}
\label{app:urans_3d_appendix}

\subsection{Physical-time, pseudo-time, and linear-solver effort}
\label{app:urans_3d_effort}

The 3D URANS calculation has a three-level work hierarchy: physical-time
advancement, inner pseudo-time iterations at each physical step, and linear
iterations for the mean-flow and turbulence systems within each inner solve.
An inner solve ending at zero-based index $k$ contains $k+1$ executed
pseudo-time iterations. The cumulative mean-flow and turbulence linear counts
are obtained by summing the linear iterations over all inner solves.

Table~\ref{tab:urans_3d_work_detail} reports the work after activation of the
windowed averaging procedure and over the complete simulations.

\begin{table}[H]
\centering
\caption{Detailed solver work for the 3D transonic-wing URANS case. The combined count
$L_\Sigma=L_{\mathrm{flow}}+L_{\mathrm{turb}}$ is an unweighted
algebraic-work diagnostic.}
\label{tab:urans_3d_work_detail}
\begin{tabularx}{0.98\textwidth}{@{}Xrrr@{}}
    \toprule
    Quantity
    & \makecell{Uniform-flow\\initialization}
    & \makecell{ADEx-FNO-based\\initialization}
    & $G_Q$ (\%) \\
    \midrule
    \multicolumn{4}{@{}l}{\textit{After window activation}} \\
    Physical-time advances
    & 162 & 132 & 18.52 \\
    Inner pseudo-time iterations
    & 2,382 & 1,710 & 28.21 \\
    Mean-flow linear iterations
    & 46,169 & 32,245 & 30.16 \\
    Turbulence linear iterations
    & 46,771 & 28,231 & 39.64 \\
    Unweighted combined linear iterations
    & 92,940 & 60,476 & 34.93 \\
    \addlinespace[3pt]
    \multicolumn{4}{@{}l}{\textit{Complete simulation}} \\
    Executed physical-time iterations
    & 313 & 284 & 9.27 \\
    Inner pseudo-time iterations
    & 7,860 & 5,409 & 31.18 \\
    Mean-flow linear iterations
    & 155,612 & 104,729 & 32.70 \\
    Turbulence linear iterations
    & 155,795 & 94,453 & 39.37 \\
    Unweighted combined linear iterations
    & 311,407 & 199,182 & 36.04 \\
    \bottomrule
\end{tabularx}
\end{table}

After window activation, the mean number of inner iterations per physical-time
advance decreases from $14.70$ to $12.95$, a reduction of $11.90\%$. The mean
combined linear count per inner iteration also decreases from $39.02$ to
$35.37$, a reduction of $9.36\%$. The larger decrease in the turbulence linear
count is associated with a reduction of its mean count per inner solve from
$19.64$ to $16.51$. Hence, the reduction in cumulative work has two components:
fewer physical-time advances and lower pseudo-time and linear work within the
remaining advances.

\begin{figure}[H]
    \centering
    \includegraphics[width=0.75\linewidth]
    {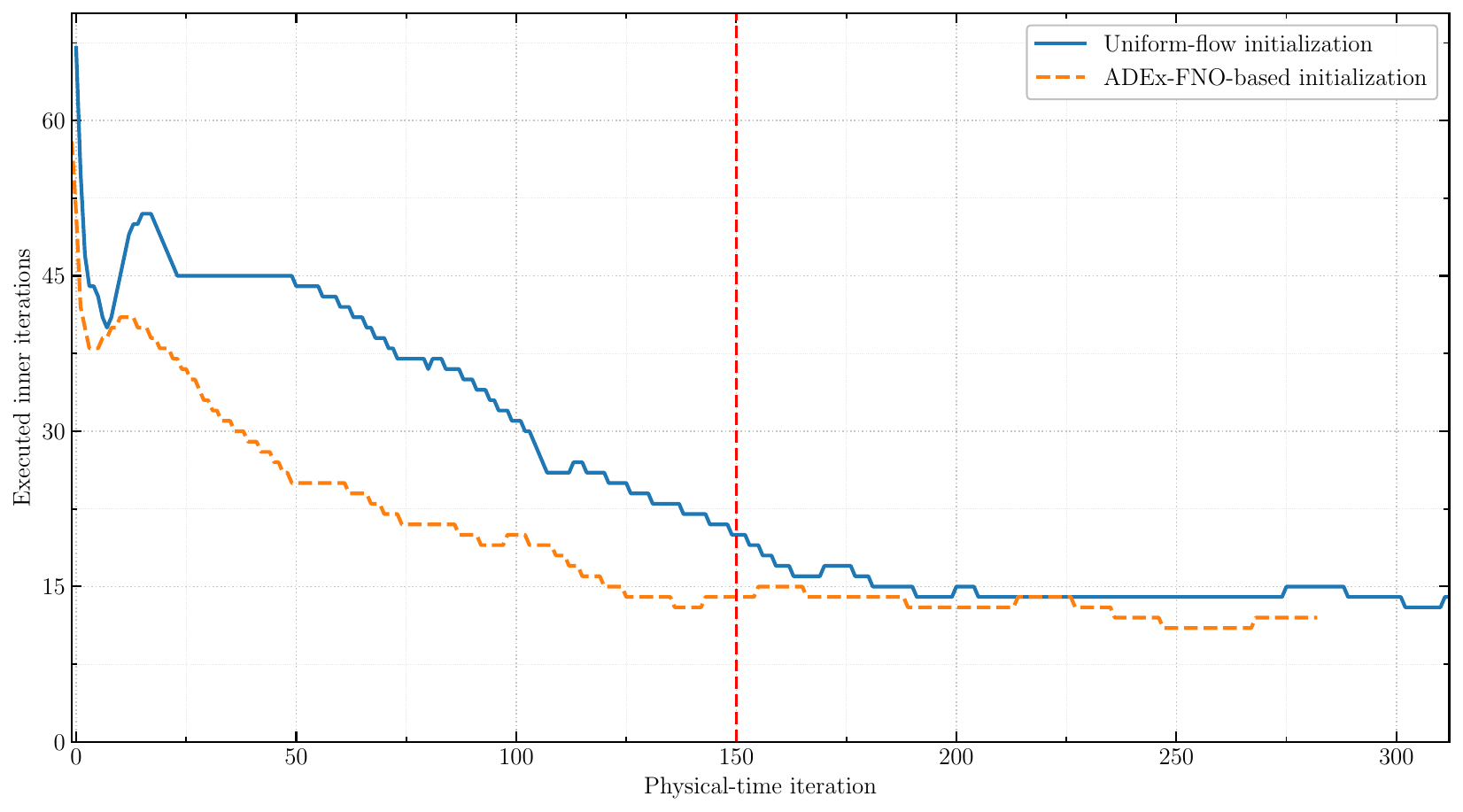}
    \caption{Executed inner pseudo-time iterations at each physical-time
    iteration for the 3D transonic-wing URANS case. Blue solid and orange
    dashed curves denote uniform-flow and ADEx-FNO-based initialization,
    respectively. The vertical red dashed line marks activation of the
    averaging window after $150$ physical-time iterations.}
    \label{fig:urans_3d_inner_iterations}
\end{figure}

Figure~\ref{fig:urans_3d_inner_iterations} shows that neither calculation
reaches the imposed limit of $80$ inner iterations. The maximum observed counts
are $67$ and $58$ for the uniform-flow and ADEx-FNO-based calculations,
respectively. After window activation, the corresponding ranges are
$13$--$20$ and $11$--$15$ inner iterations per physical-time advance.

\begin{figure}[H]
    \centering
    \includegraphics[width=0.75\linewidth]
    {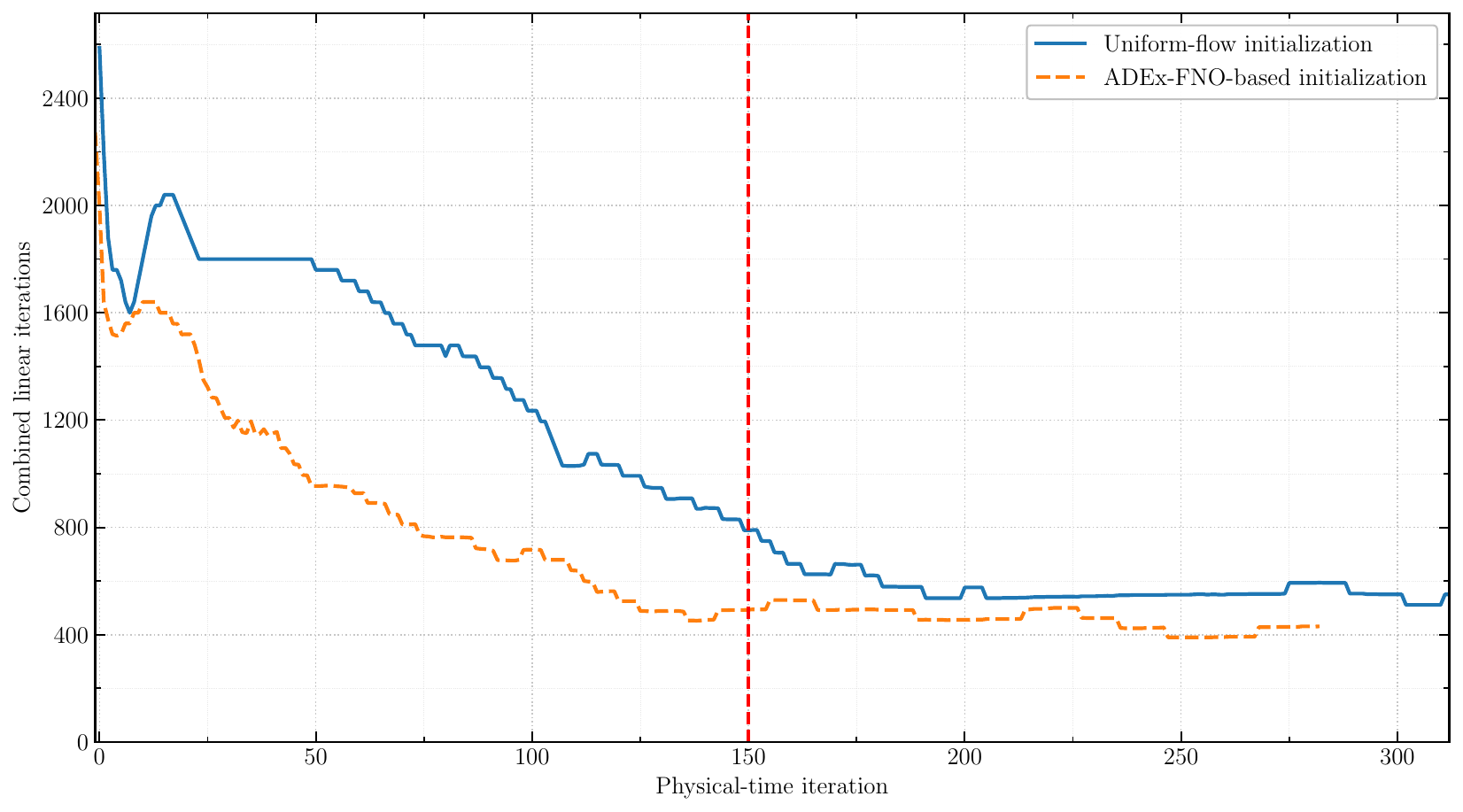}
    \caption{Unweighted combined linear iterations accumulated at each
    physical-time iteration for the 3D transonic-wing URANS case. The line
    styles and averaging-window marker have the same meaning as in
    Figure~\ref{fig:urans_3d_inner_iterations}.}
    \label{fig:urans_3d_linear_per_step}
\end{figure}

Figure~\ref{fig:urans_3d_linear_per_step} confirms that the ADEx-FNO-based
initialization reduces the algebraic work accumulated within the physical-time
steps as well as the number of such steps required before termination. The
initial ADEx-FNO-based residual is not smaller: the initial values of
$r_\rho$ are $-5.9071$ and $-4.4336$ for the uniform-flow and ADEx-FNO-based
calculations, respectively. The reduction therefore cannot be attributed
simply to a lower initial discrete residual.

\subsection{Windowed time-convergence diagnostics}
\label{app:urans_3d_window}

Table~\ref{tab:urans_3d_cauchy_terminal} reports the configured window starts,
the first iterations at which both Cauchy diagnostics fall below $10^{-5}$,
and the actual solver-termination iterations. In each run, SU2 terminates one
physical-time iteration after the first simultaneous crossing. The
lift diagnostic controls termination in both calculations.

\begin{table}[H]
\centering
\caption{Window-start, Cauchy-crossing, and termination data for the 3D transonic-wing URANS calculations.}
\label{tab:urans_3d_cauchy_terminal}
\begin{tabular}{@{}lrrrr@{}}
    \toprule
    Initialization
    & \makecell{Window\\start}
    & \makecell{First simultaneous\\Cauchy crossing}
    & \makecell{Solver\\termination}
    & \makecell{Terminal diagnostics\\$(C_D,C_L)$} \\
    \midrule
    Uniform-flow
    & 150 & 311 & 312
    & $(\num{4.019e-6},\,\num{9.675e-6})$ \\
    ADEx-FNO-based
    & 650 & 781 & 782
    & $(\num{6.164e-6},\,\num{9.197e-6})$ \\
    \bottomrule
\end{tabular}
\end{table}

\begin{figure}[H]
    \centering
    \begin{subfigure}[t]{0.48\textwidth}
        \centering
        \includegraphics[width=\linewidth]
        {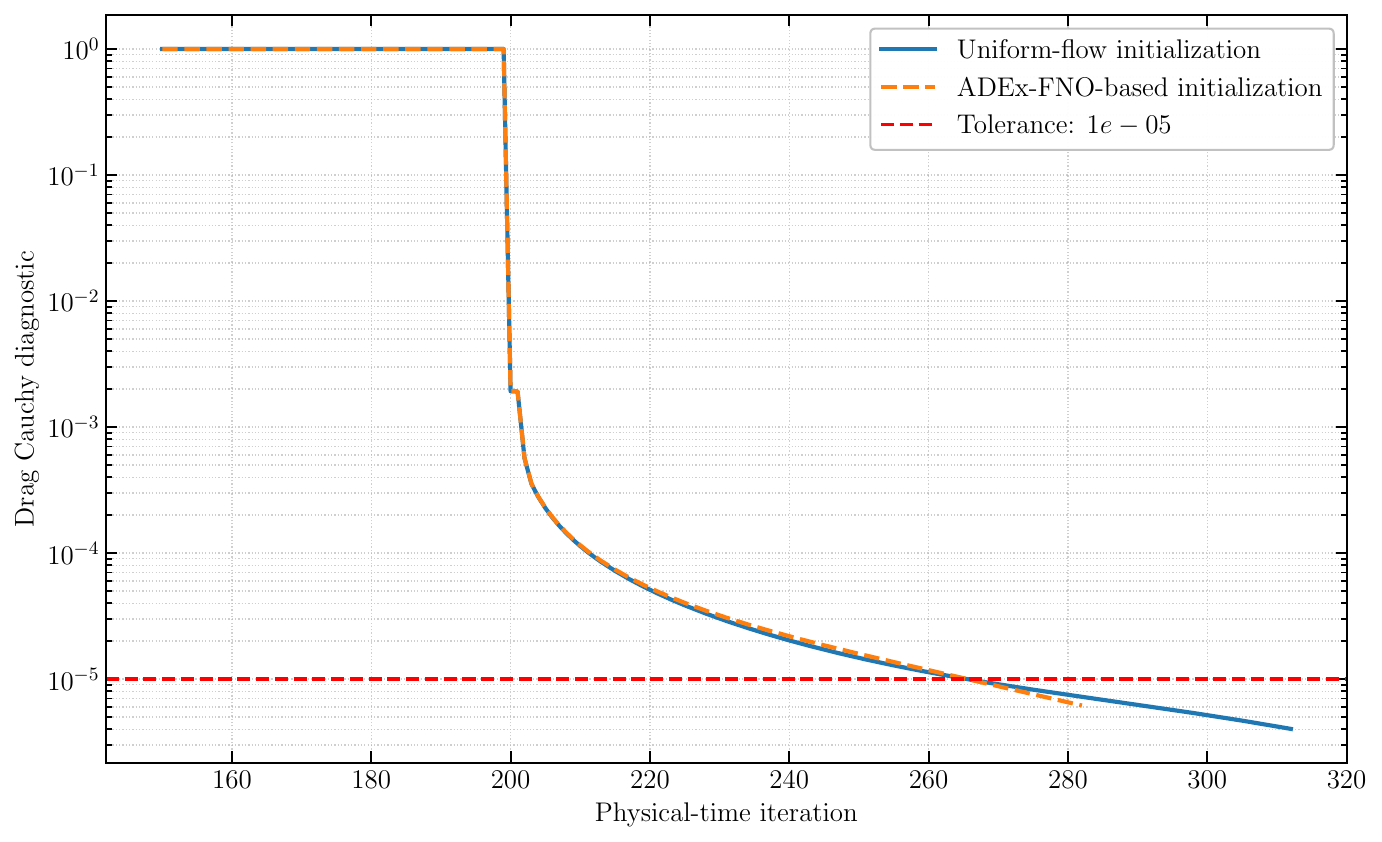}
        \caption{Drag-coefficient diagnostic.}
    \end{subfigure}
    \hfill
    \begin{subfigure}[t]{0.48\textwidth}
        \centering
        \includegraphics[width=\linewidth]
        {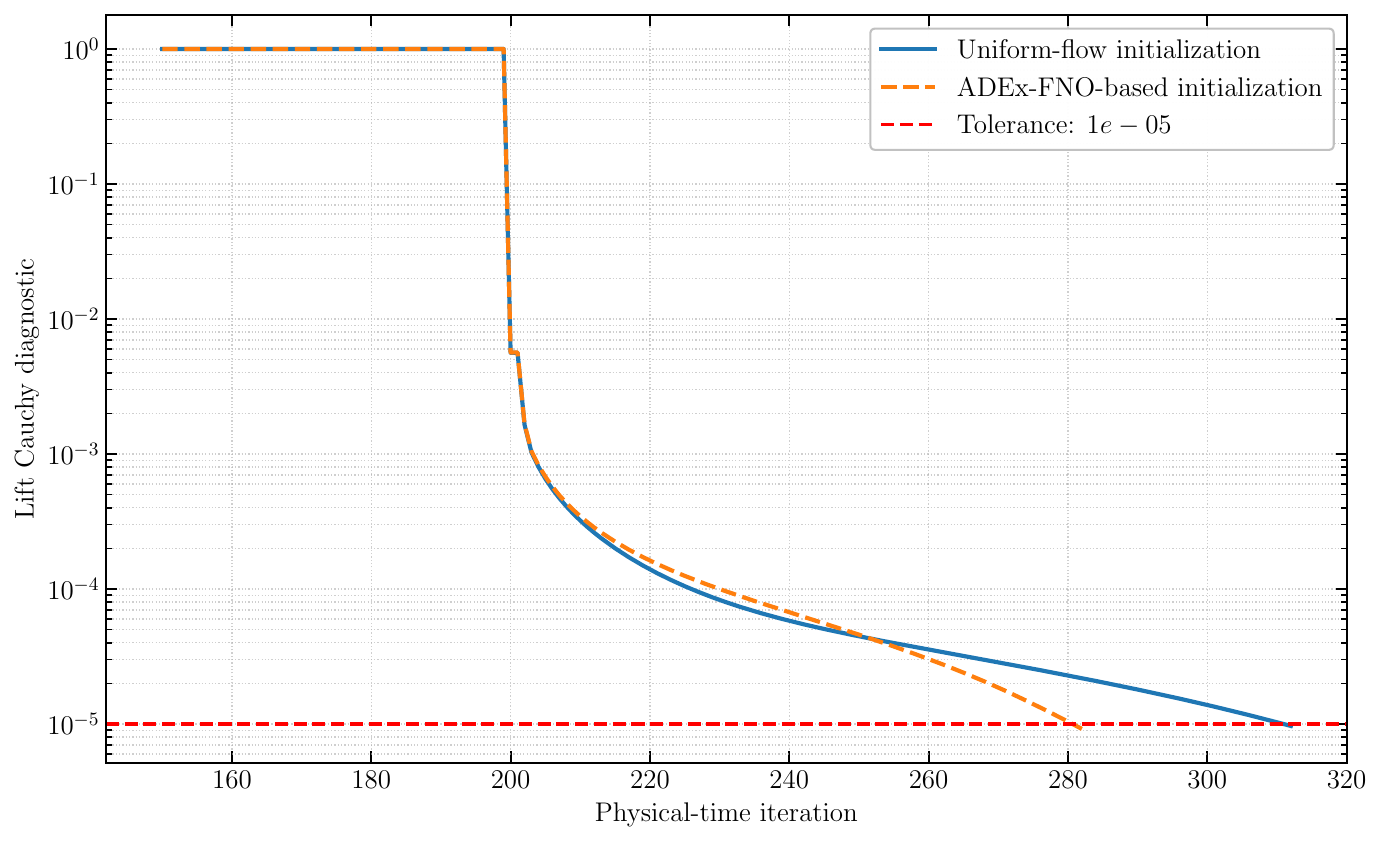}
        \caption{Lift-coefficient diagnostic.}
    \end{subfigure}
    \caption{Cauchy diagnostics for the squared-Hann-windowed aerodynamic
    coefficients in the 3D transonic-wing URANS case. The red dashed line marks
    the prescribed tolerance $10^{-5}$.}
    \label{fig:urans_3d_cauchy}
\end{figure}

Figure~\ref{fig:urans_3d_cauchy} shows that the drag diagnostic crosses the
tolerance first, whereas the lift diagnostic determines the final stopping
point. The terminal RMS-density residuals are $-11.0955$ and $-11.0341$, and
the minimum final-inner residuals attained during the calculations are $-11.1049$
and $-11.1492$. Thus, the final inner solves in both calculations satisfy the
common $r_\rho<-11$ threshold.

\subsection{Consistency of the time-averaged aerodynamic coefficients}
\label{app:urans_3d_statistics}

Table~\ref{tab:urans_3d_terminal_statistics} compares the terminal
instantaneous coefficients and the terminal squared-Hann averages. Both
calculations independently satisfy the same coefficient-based tolerance, and
the terminal averages are mutually consistent at the sub-percent level.

\begin{table}[H]
\centering
\caption{Terminal instantaneous and squared-Hann-windowed aerodynamic
coefficients for the 3D transonic-wing URANS case. Relative differences use
magnitude of the uniform-flow value as reference.}
\label{tab:urans_3d_terminal_statistics}
\begin{tabular}{@{}lrrrr@{}}
    \toprule
    Quantity
    & \makecell{Uniform-flow\\initialization}
    & \makecell{ADEx-FNO-based\\initialization}
    & Absolute difference
    & Relative difference (\%) \\
    \midrule
    Terminal $C_L$
    & $-0.224$ & $-0.225$
    & $\num{5.719e-4}$ & 0.255 \\
    Terminal $C_D$
    & $0.0762$ & $0.0763$
    & $\num{1.121e-4}$ & 0.147 \\
    Terminal $\overline{C_L}$
    & $-0.222$ & $-0.223$
    & $\num{7.988e-4}$ & 0.360 \\
    Terminal $\overline{C_D}$
    & $0.0758$ & $0.0758$
    & $\num{5.493e-5}$ & 0.072 \\
    \bottomrule
\end{tabular}
\end{table}

\begin{figure}[H]
    \centering
    \begin{subfigure}[t]{0.48\textwidth}
        \centering
        \includegraphics[width=\linewidth]
        {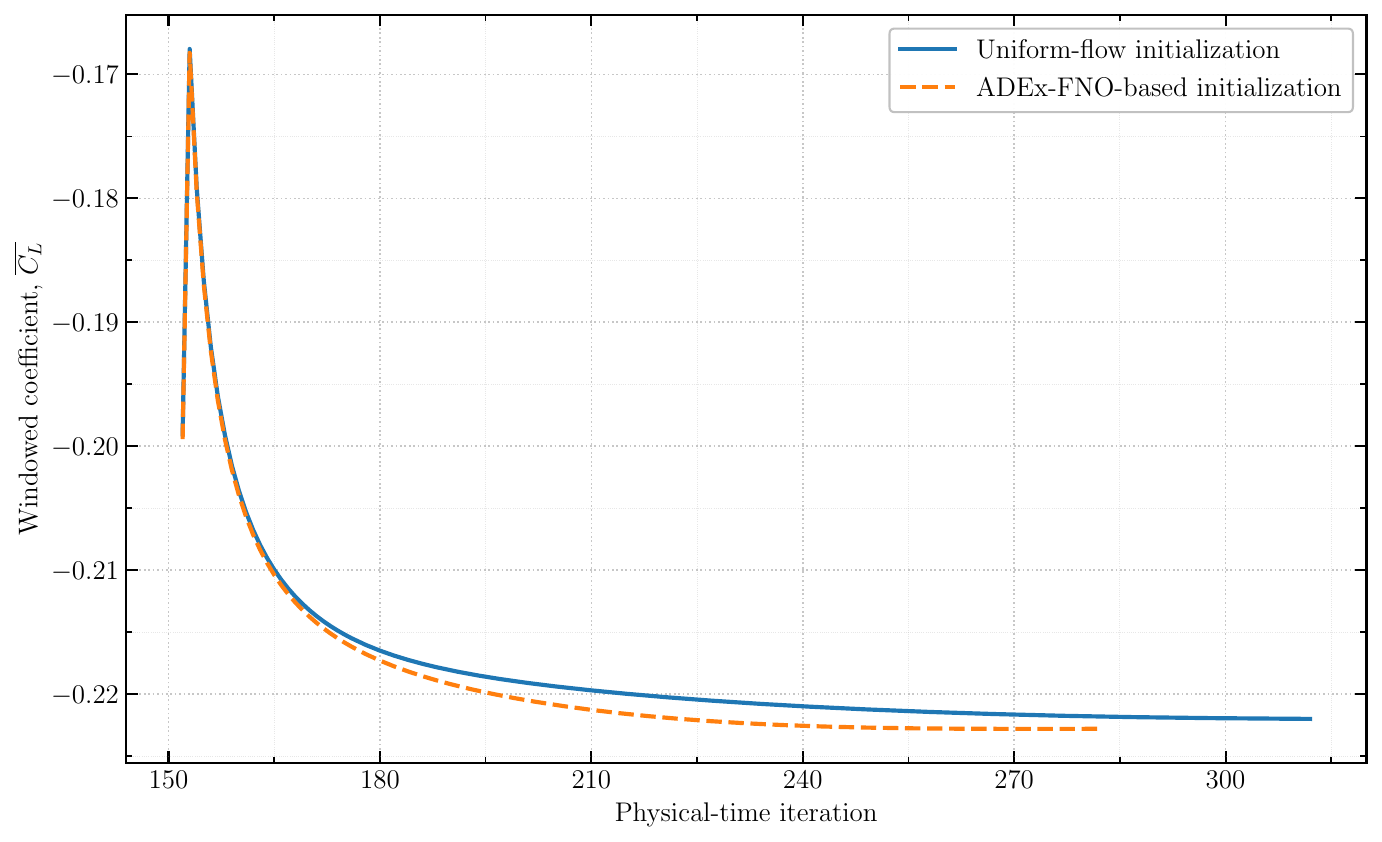}
        \caption{Windowed lift coefficient.}
    \end{subfigure}
    \hfill
    \begin{subfigure}[t]{0.48\textwidth}
        \centering
        \includegraphics[width=\linewidth]
        {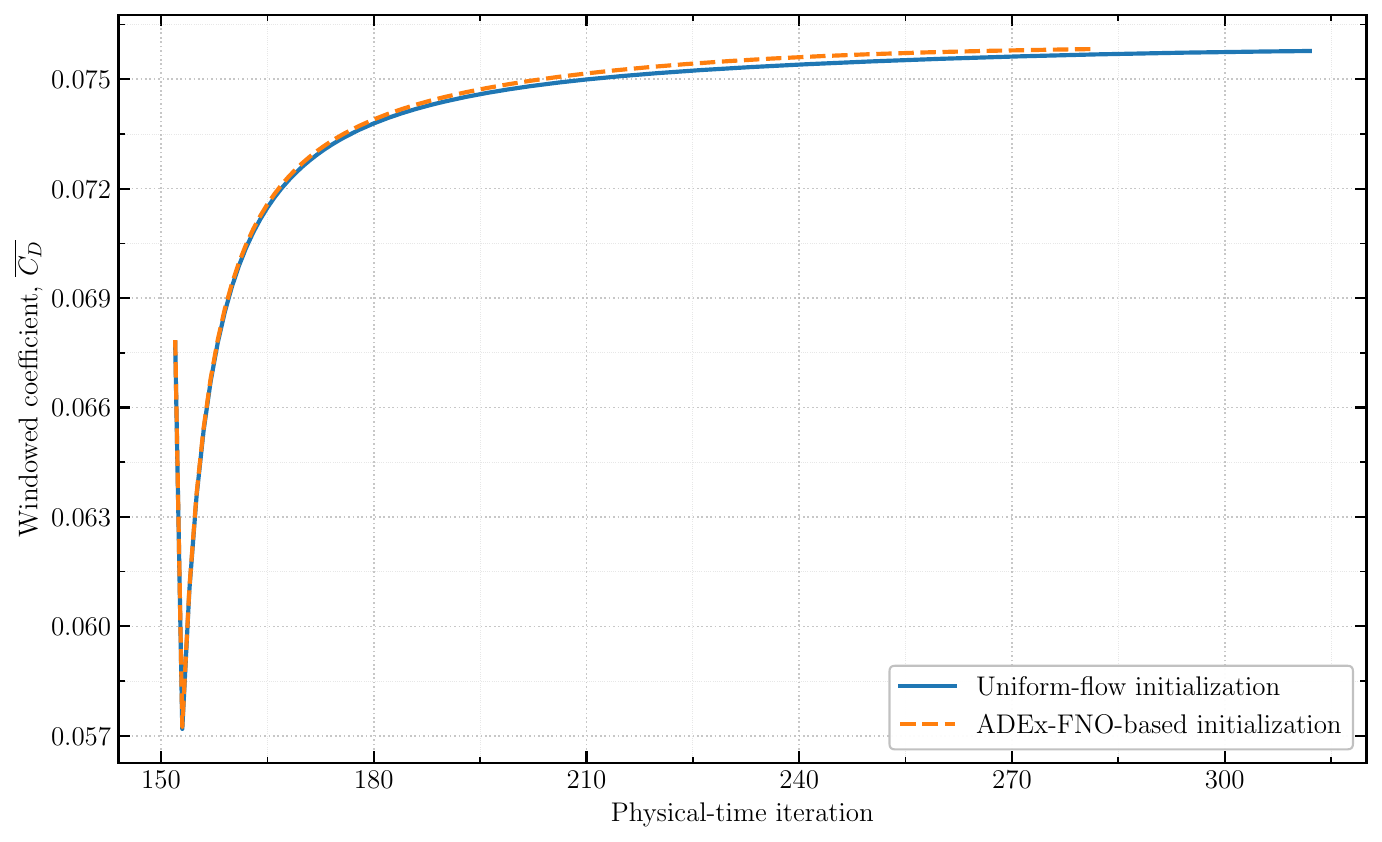}
        \caption{Windowed drag coefficient.}
    \end{subfigure}
    \caption{Squared-Hann-windowed aerodynamic coefficients for the 3D
    transonic-wing URANS case after activation of the averaging window.}
    \label{fig:urans_3d_tavg}
\end{figure}

Figure~\ref{fig:urans_3d_tavg} shows that the two running averages approach closely
comparable terminal values despite their different startup transients. The
terminal differences are $0.360\%$ in lift and $0.072\%$ in drag. These values
do not imply exact equality, but they support the conclusion that the ADEx-FNO-based
initialization reduces the solver work without materially changing the
terminal time-averaged aerodynamic state for this case.

\section{Details of the DNS bootstrap analysis}
\label{app:dns_bootstrap}

\subsection{Force-based truncation and sensitivity}

The force-based bootstrap estimate uses the same MSER procedure for both initializations. For each prescribed block duration $B$, the $C_L$ and $C_D$ histories are replaced by common block-averaged signals, and MSER determines a candidate truncation time for each coefficient. The bootstrap time is the later of the two coefficient-based truncation times. Six values in the range $0.02\leq B\leq0.50$ are evaluated, and the median of the resulting bootstrap times is used in the main text. The sensitivity shown in Figure~\ref{fig:dns_mser_sensitivity} gives reductions from $34.15\%$ to $38.10\%$, which supports the central value of $35.96\%$ without selecting a block duration after observing the outcome.

Because the DNS uses adaptive explicit time stepping, equal convective times do not correspond to equal solver-step counts. The central MSER times $t^*=5.075$ and $3.250$ correspond to step indices $1{,}243{,}957$ and $829{,}132$, respectively. Thus, the CTU reduction is $35.96\%$, while the step-count reduction is $33.35\%$.

\subsection{Semi-discrete right-hand-side histories}

\begin{figure}[H]
    \centering
    \includegraphics[width=0.95\linewidth]{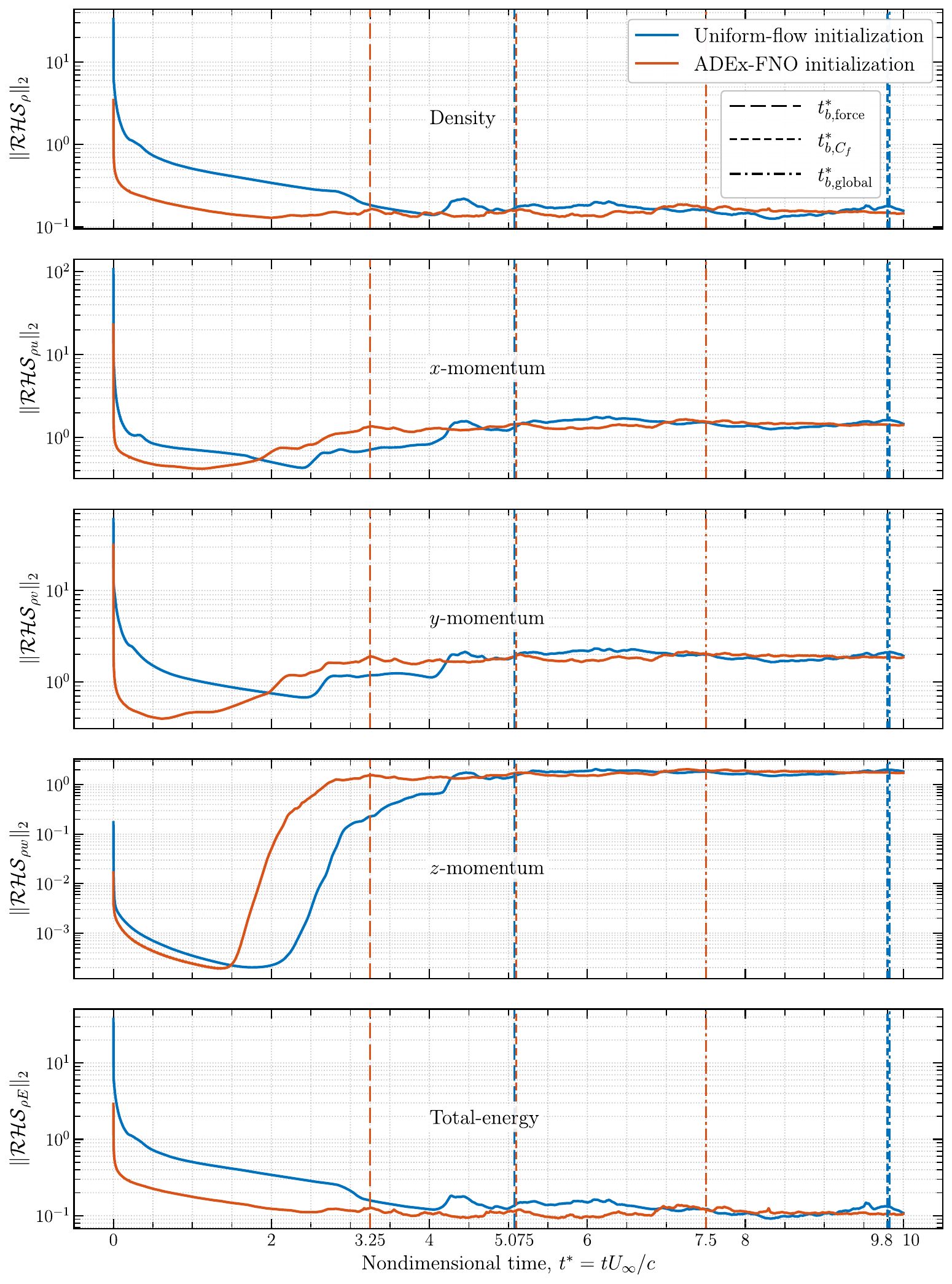}
    \caption{Component-wise $L^2$ norms of the semi-discrete DNS right-hand side for uniform-flow initialization (blue curve) and ADEx-FNO-based initialization (orange curve). From top to bottom, the panels show density, $x$-momentum, $y$-momentum, $z$-momentum, and total energy components. The dashed orange and blue vertical lines mark the different bootstrap estimates at $t^*=3.250$, $5.075$, $7.5$, and $9.8$. The $\rho w$ component indicates the development of spanwise activity but does not identify a particular instability mechanism.}
    \label{fig:dns_rhs_components}
\end{figure}

Figure~\ref{fig:dns_rhs_components} reports the five recorded componentwise $L^2$ norms of the semi-discrete right-hand side,
\begin{equation*}
    \left\lVert\mathcal{R}_{\rho}\right\rVert_2,
    \quad
    \left\lVert\mathcal{R}_{\rho u}\right\rVert_2,
    \quad
    \left\lVert\mathcal{R}_{\rho v}\right\rVert_2,
    \quad
    \left\lVert\mathcal{R}_{\rho w}\right\rVert_2,
    \quad
    \left\lVert\mathcal{R}_{\rho E}\right\rVert_2.
\end{equation*}
The five equations represent different conserved quantities and are therefore shown separately rather than being combined into a single unnormalized norm.
All components exhibit a substantially smaller initial imbalance under ADEx-FNO-based initialization. They do not decay to zero after the flow develops; instead, they approach bounded fluctuating levels, as expected for a physically unsteady DNS. The $\rho w$ component is particularly informative: it remains very small during the initially quasi-2D stage and then grows as spanwise motion develops. Its growth is evidence of emerging 3D activity, but it is not by itself a diagnostic of a specific instability mechanism.

\subsection{3D wake development}

\begin{figure}[H]
    \centering
    \begin{tabular}{
        @{}
        >{\centering\arraybackslash}m{0.095\textwidth}
        @{\hspace{0.01\textwidth}}
        >{\centering\arraybackslash}m{0.435\textwidth}
        @{\hspace{0.01\textwidth}}
        >{\centering\arraybackslash}m{0.435\textwidth}
        @{}}

        & \text{Uniform-flow initialization} & \text{ADEx-FNO-based initialization} \\

        $t^*=1.5$ &
        \includegraphics[width=\linewidth]{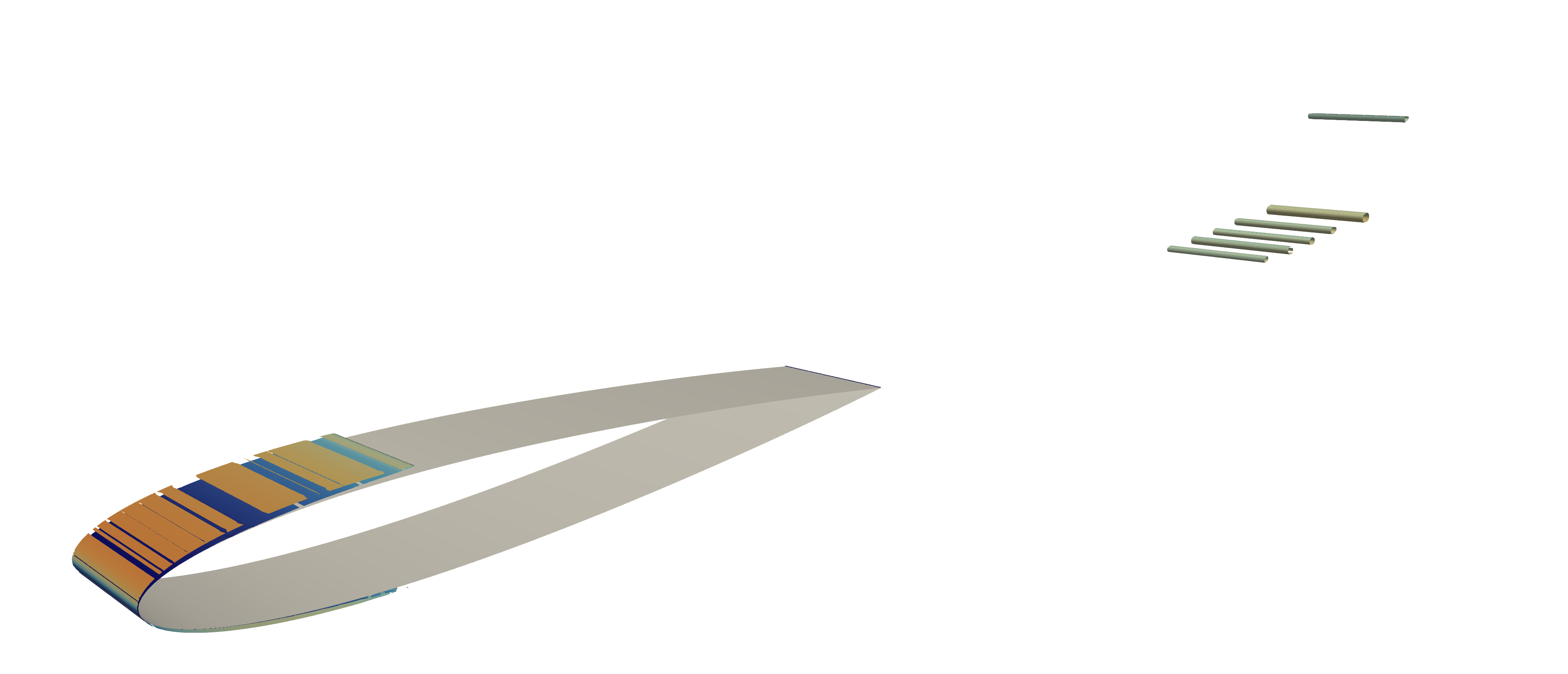} &
        \includegraphics[width=\linewidth]{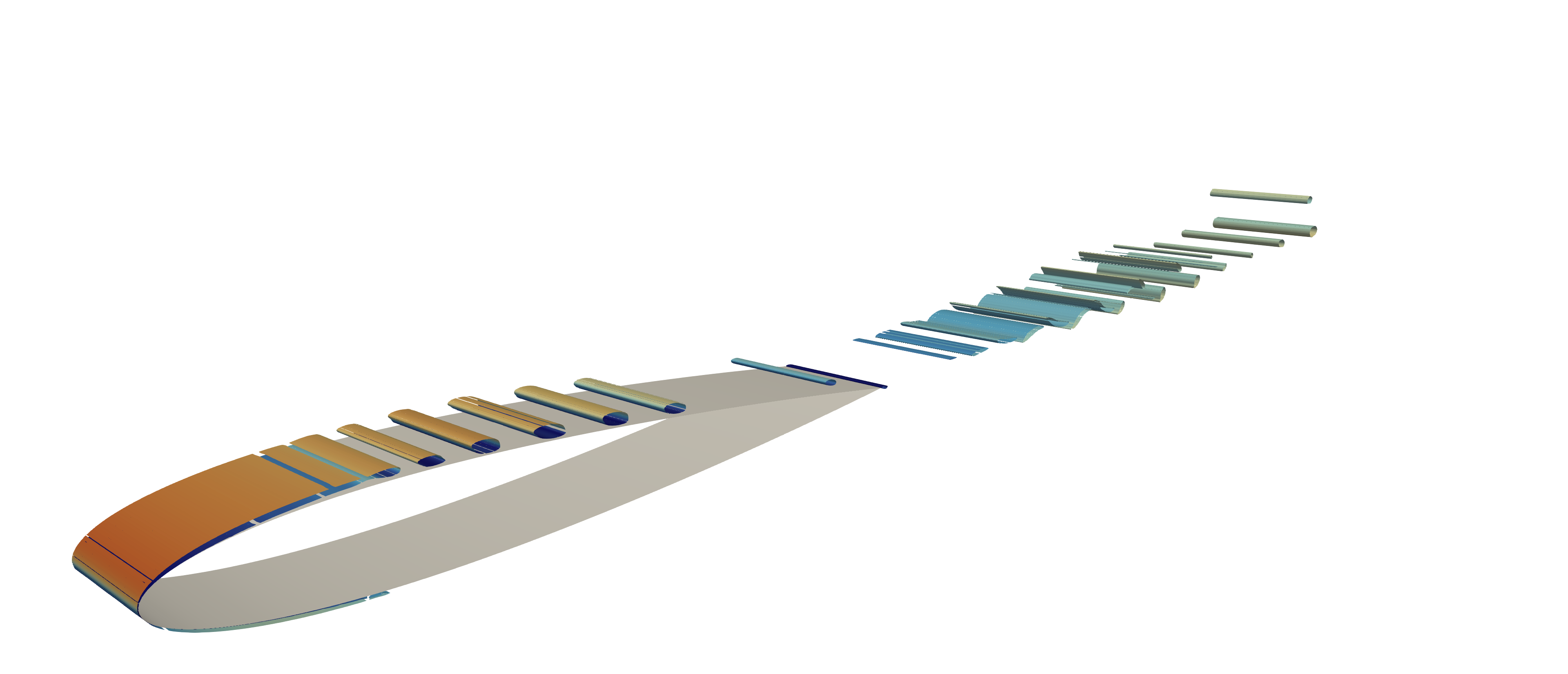} \\[0.8em]

        $t^*=3.25$ &
        \includegraphics[width=\linewidth]{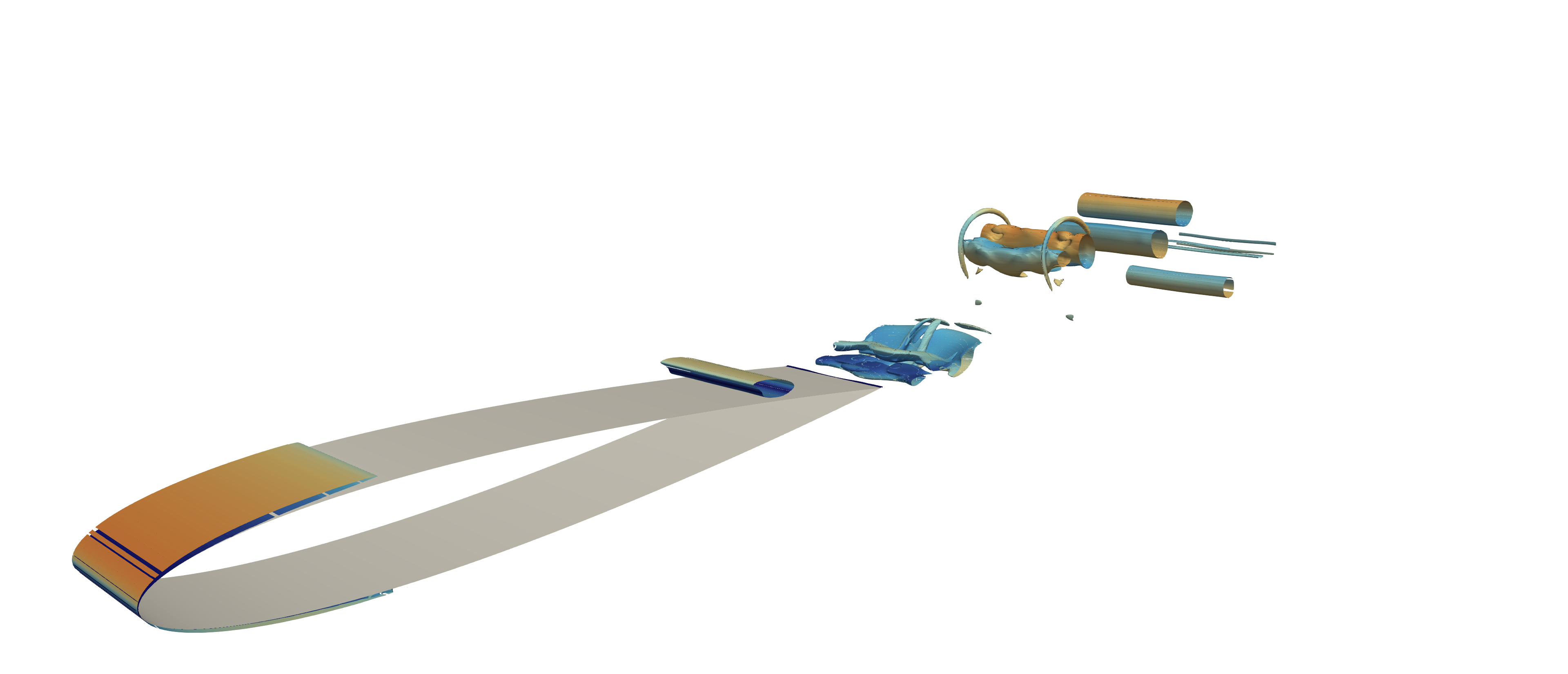} &
        \includegraphics[width=\linewidth]{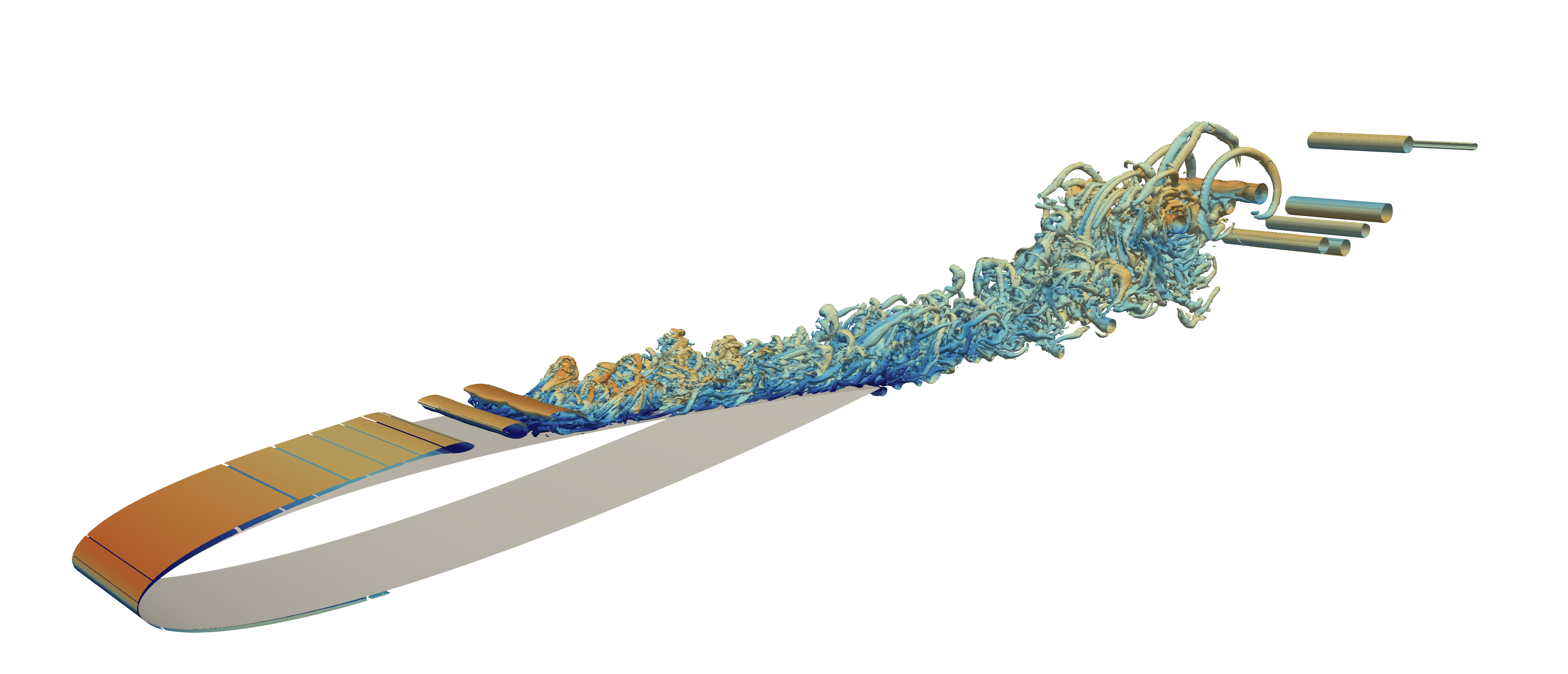} \\[0.8em]

        $t^*=5.1$ &
        \includegraphics[width=\linewidth]{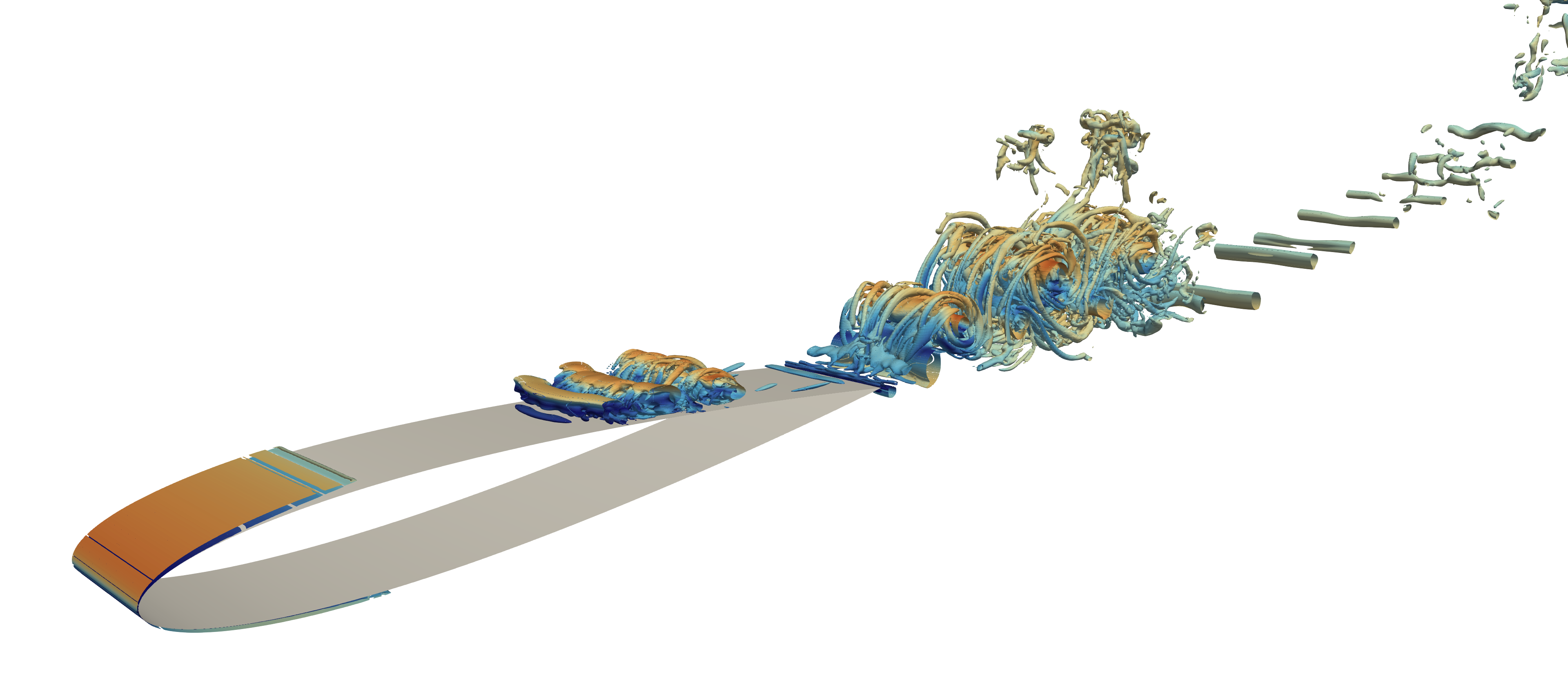} &
        \includegraphics[width=\linewidth]{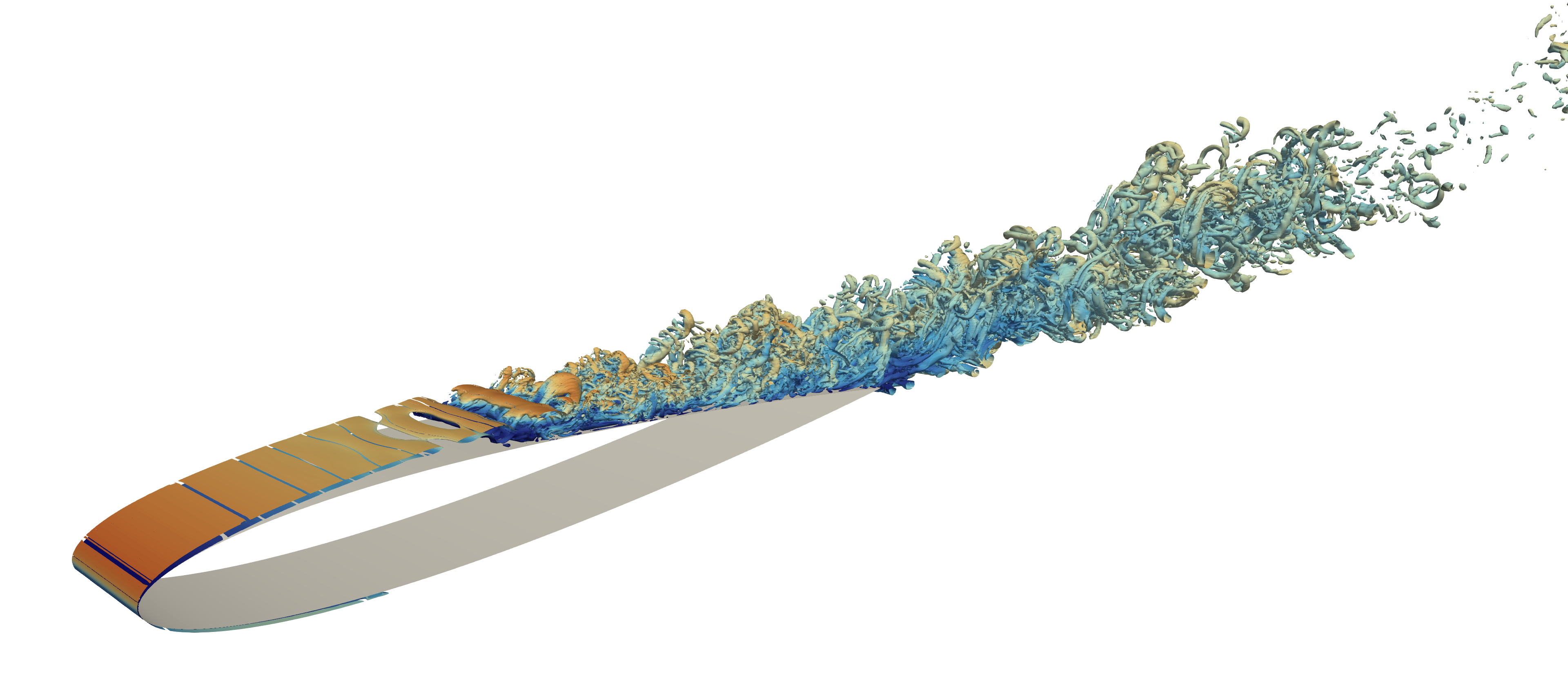} \\[0.8em]

        $t^*=7.5$ &
        \includegraphics[width=\linewidth]{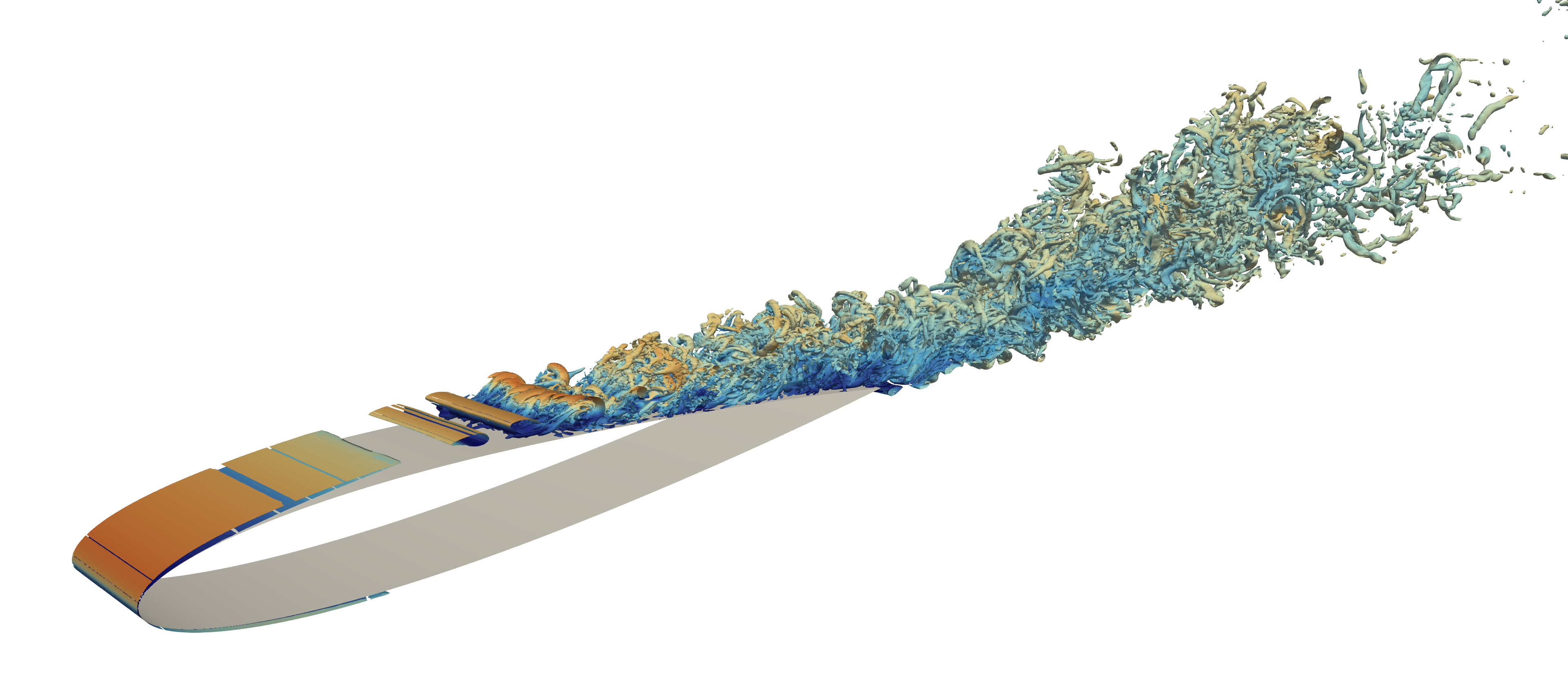} &
        \includegraphics[width=\linewidth]{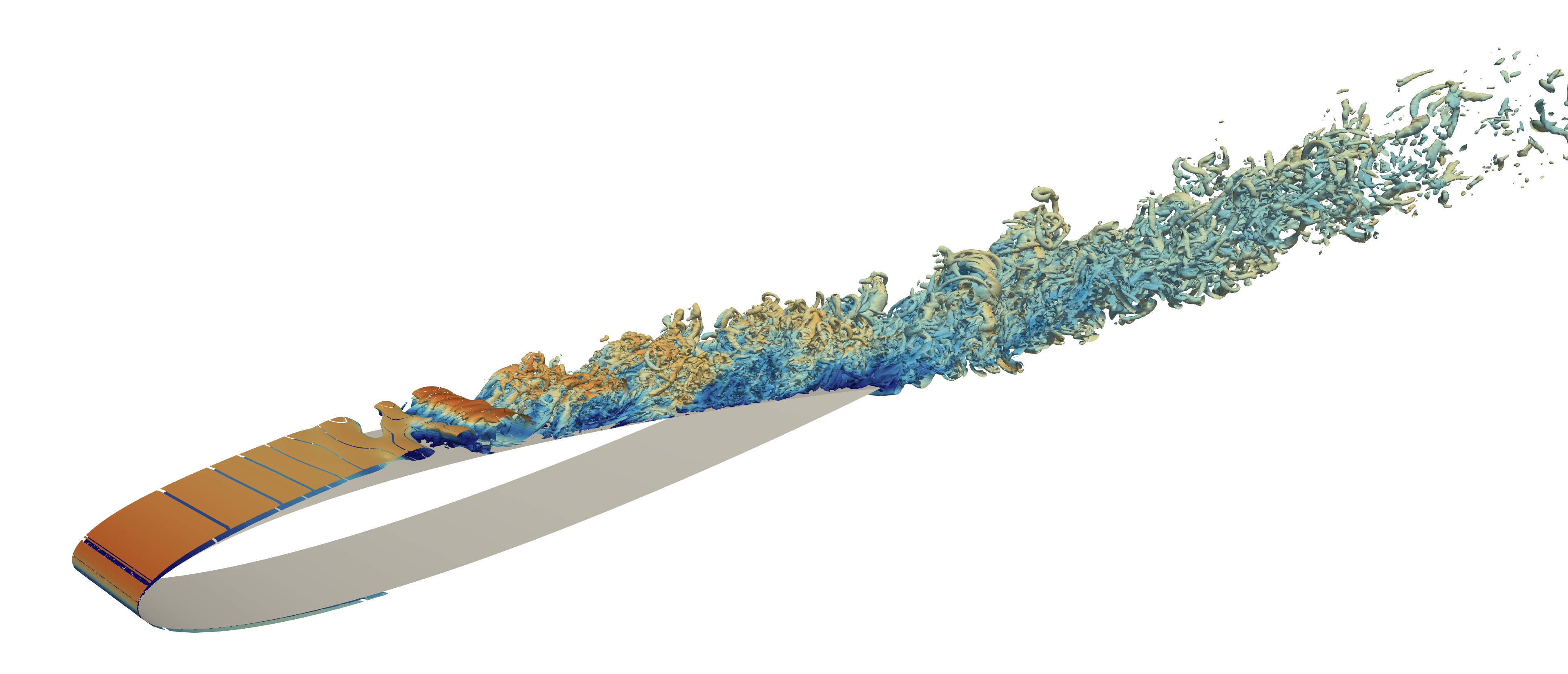} \\[0.8em]

        \multicolumn{3}{@{}c@{}}{
        \includegraphics[width=0.6\textwidth]{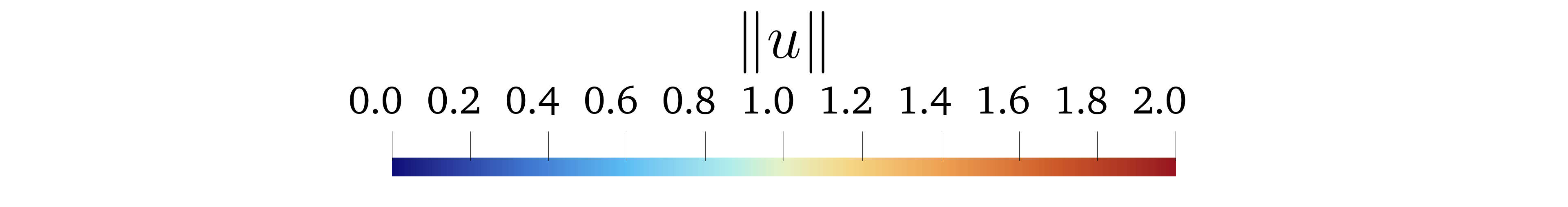}}
    \end{tabular}

    \caption{DNS $Q$-criterion iso-surfaces for uniform-flow and ADEx-FNO initialization at common convective times. The iso-surfaces correspond to $Q=50$ and are colored by velocity magnitude, $\lVert u \rVert$. The sequence provides qualitative evidence of the formation and downstream development of 3D wake structures.}
    \label{fig:QCriterion_DNS}
\end{figure}

The $Q$-criterion sequence in Figure~\ref{fig:QCriterion_DNS} provides a qualitative view of the emergence of 3D structures. At $t^*=1.5$, both calculations remain dominated by organized structures, but the ADEx-FNO-based solution already shows more developed near-body and wake activity. At $t^*=3.25$, a sustained 3D wake is visible for ADEx-FNO-based initialization, whereas the uniform-flow calculation retains a shorter and more weakly developed wake. The uniform-flow calculation develops comparable complexity only at later times. These images support the force- and skin-friction-based results, but no exact bootstrap time is assigned from the $Q$-criterion rendering alone.

\subsection{Full-domain wake flushing}

\begin{figure}[H]
    \centering
    \begin{tabular}{
        @{}
        >{\centering\arraybackslash}m{0.095\textwidth}
        @{\hspace{0.01\textwidth}}
        >{\centering\arraybackslash}m{0.435\textwidth}
        @{\hspace{0.01\textwidth}}
        >{\centering\arraybackslash}m{0.435\textwidth}
        @{}}

        & \text{Uniform-flow initialization} & \text{ADEx-FNO-based initialization} \\ \\

        $t^*=3.25$ &
        \includegraphics[width=\linewidth]{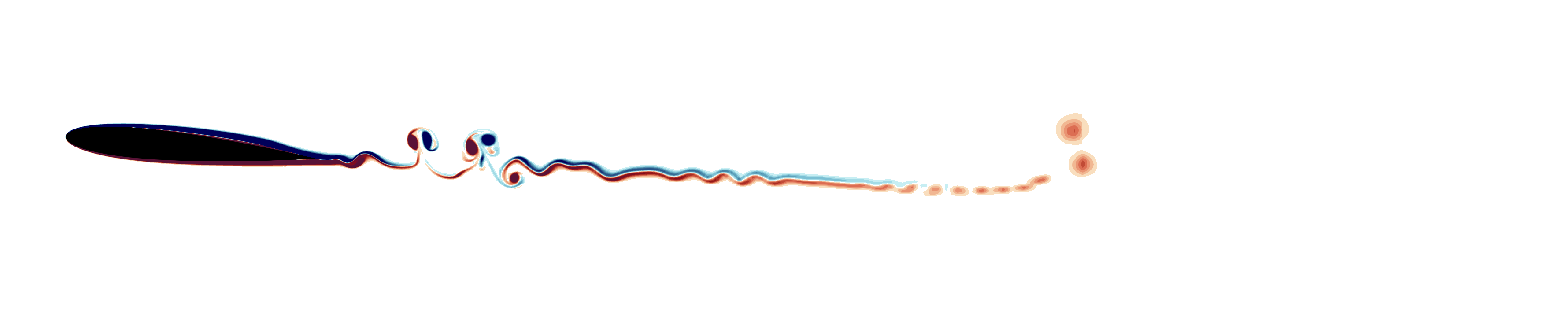} &
        \includegraphics[width=\linewidth]{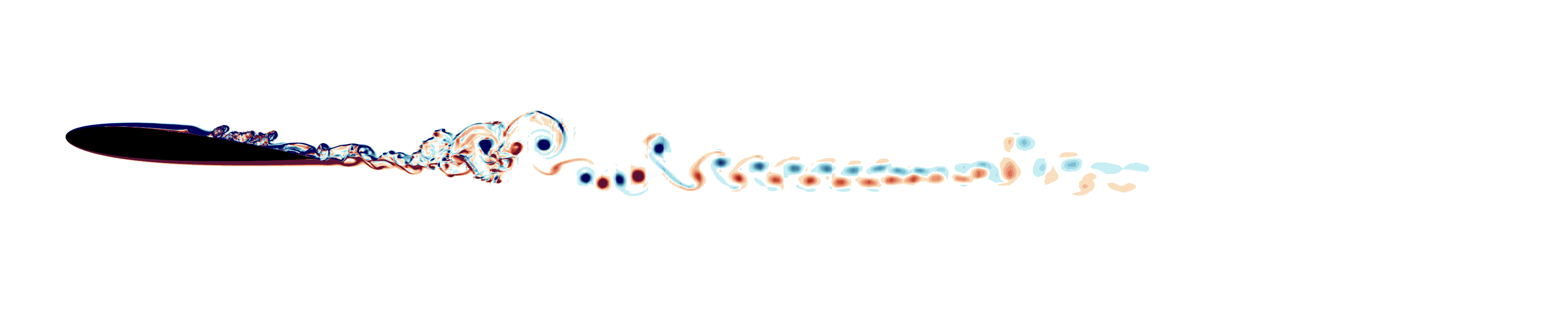} \\[0.8em]

        $t^*=5.1$ &
        \includegraphics[width=\linewidth]{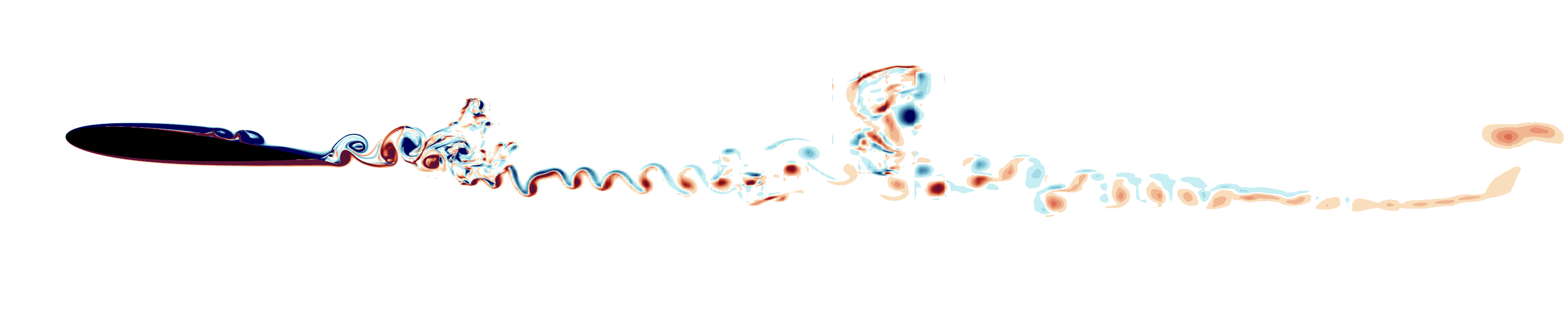} &
        \includegraphics[width=\linewidth]{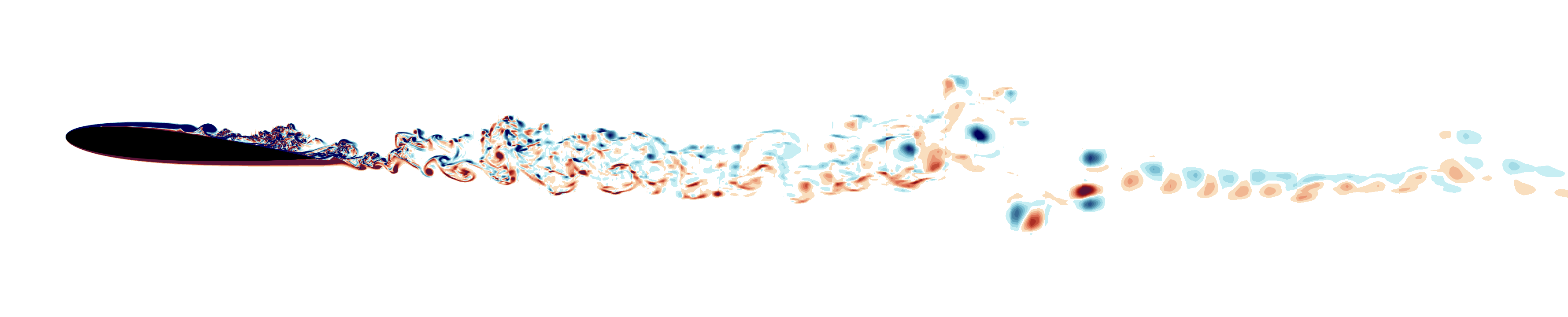} \\[0.8em]

        $t^*=7.5$ &
        \includegraphics[width=\linewidth]{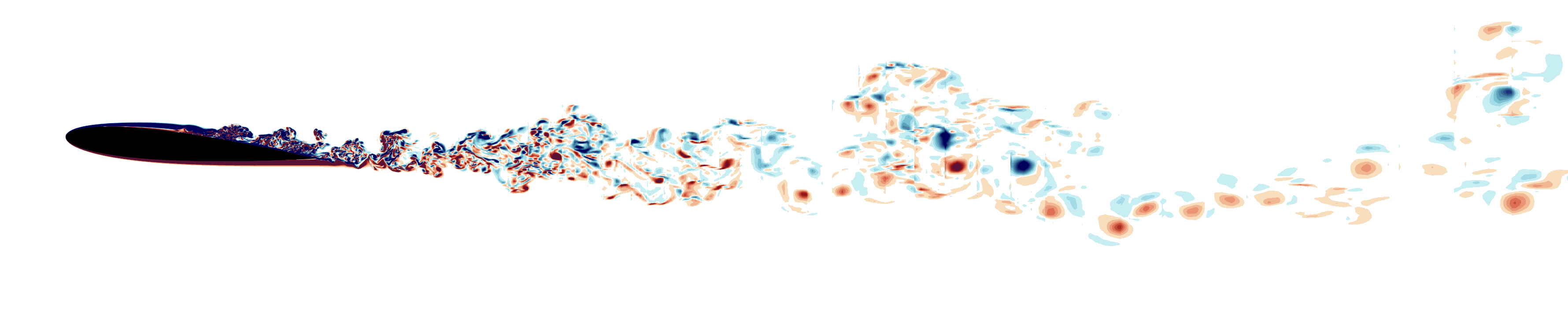} &
        \includegraphics[width=\linewidth]{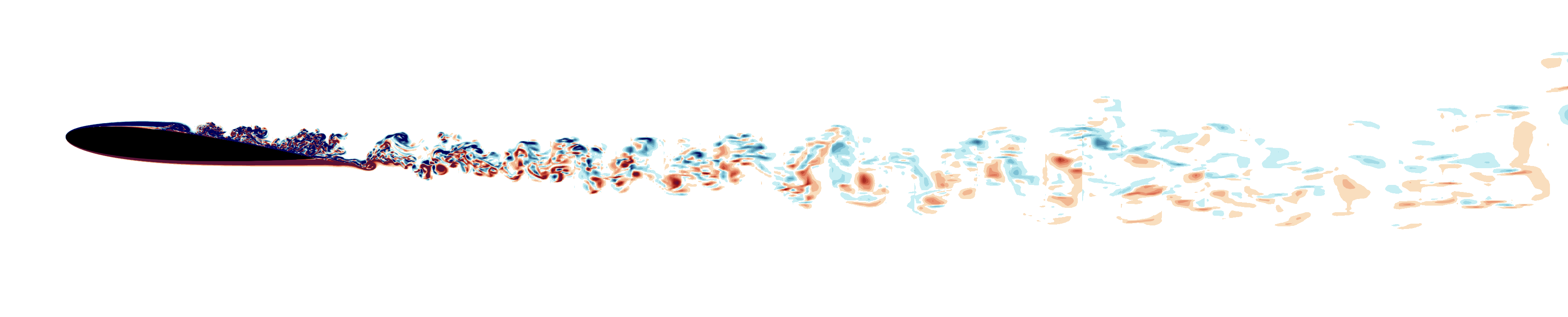} \\[0.8em]

        $t^*=9.8$ &
        \includegraphics[width=\linewidth]{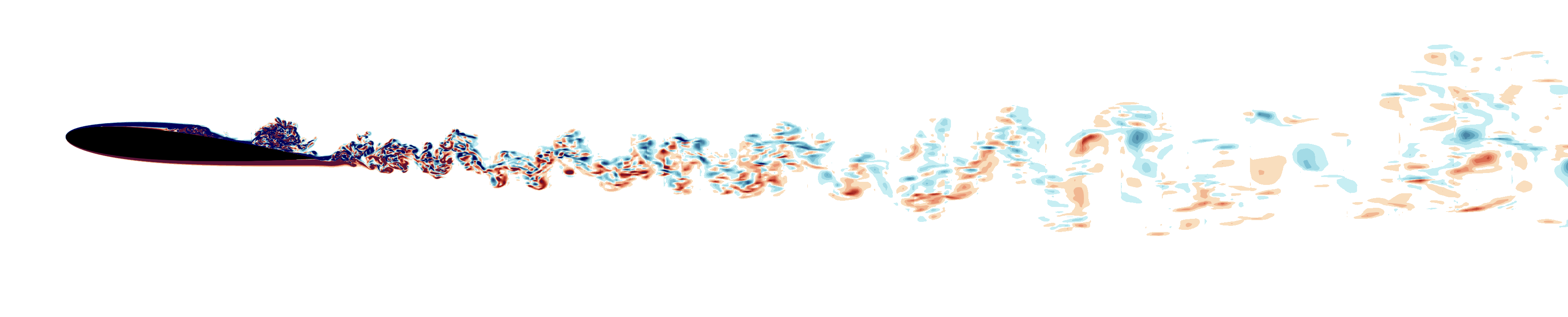} &
        \includegraphics[width=\linewidth]{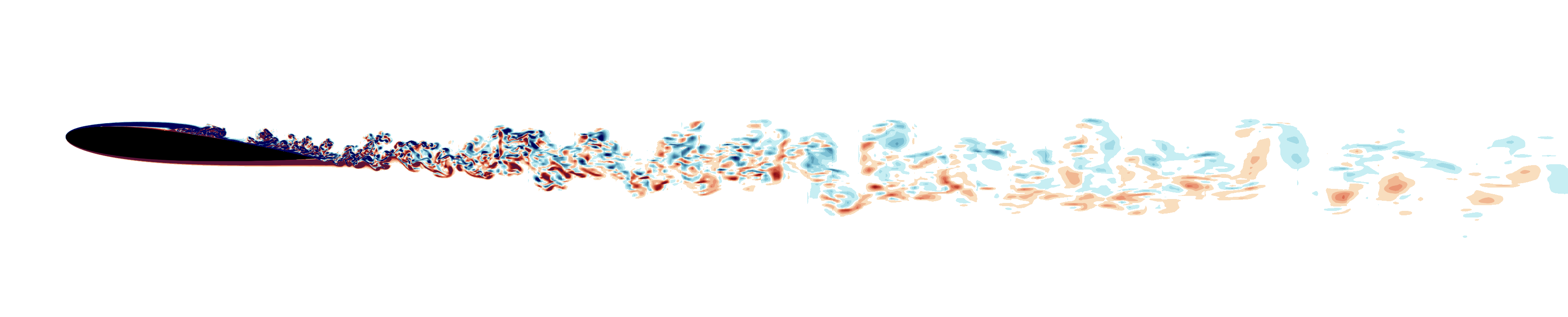} \\[0.8em]

        \multicolumn{3}{@{}c@{}}{
        \includegraphics[width=0.6\textwidth]{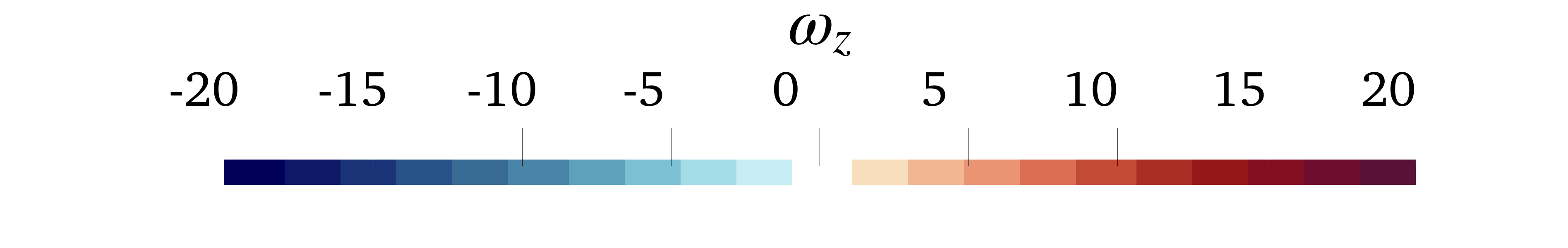}}
    \end{tabular}

    \caption{DNS spanwise-vorticity contours, $\omega_z$, on the section $z=0.1$ for uniform-flow and ADEx-FNO initialization at common convective times. The qualitative sequence distinguishes near-wake development from flushing of initialization-generated structures through the downstream domain.}
    \label{fig:Vorticity_DNS}
\end{figure}

Figure~\ref{fig:Vorticity_DNS} shows the spanwise vorticity field on the section $z=0.1$. The right edge of each panel coincides with the downstream extent of the displayed far-field domain. At $t^*=5.1$, the ADEx-FNO near wake is developed, but coherent structures generated during initialization are still visible downstream. These structures have left the domain by $t^*=7.5$. For uniform-flow initialization, comparable flushing is not observed until $t^*=9.8$. This distinction motivates the global bootstrap times reported in the main text.

\subsection{Density disturbances}

\begin{figure}[H]
    \centering
    \begin{tabular}{
        @{}
        >{\centering\arraybackslash}m{0.095\textwidth}
        @{\hspace{0.01\textwidth}}
        >{\centering\arraybackslash}m{0.435\textwidth}
        @{\hspace{0.01\textwidth}}
        >{\centering\arraybackslash}m{0.435\textwidth}
        @{}}

        & \text{Uniform-flow initialization} & \text{ADEx-FNO-based initialization} \\ \\

        $t^*=0.4$ &
        \includegraphics[width=0.6\linewidth]{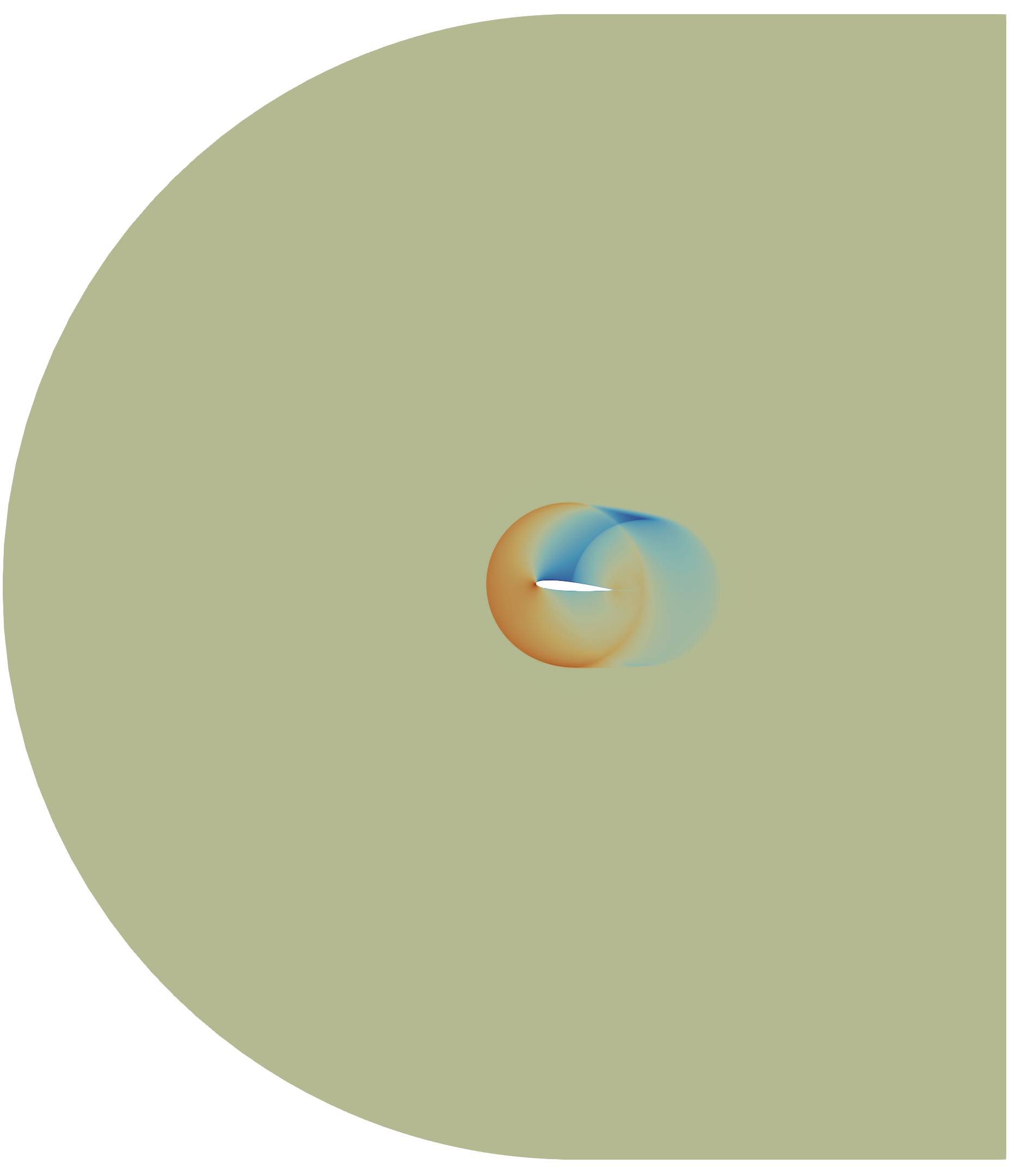} &
        \includegraphics[width=0.6\linewidth]{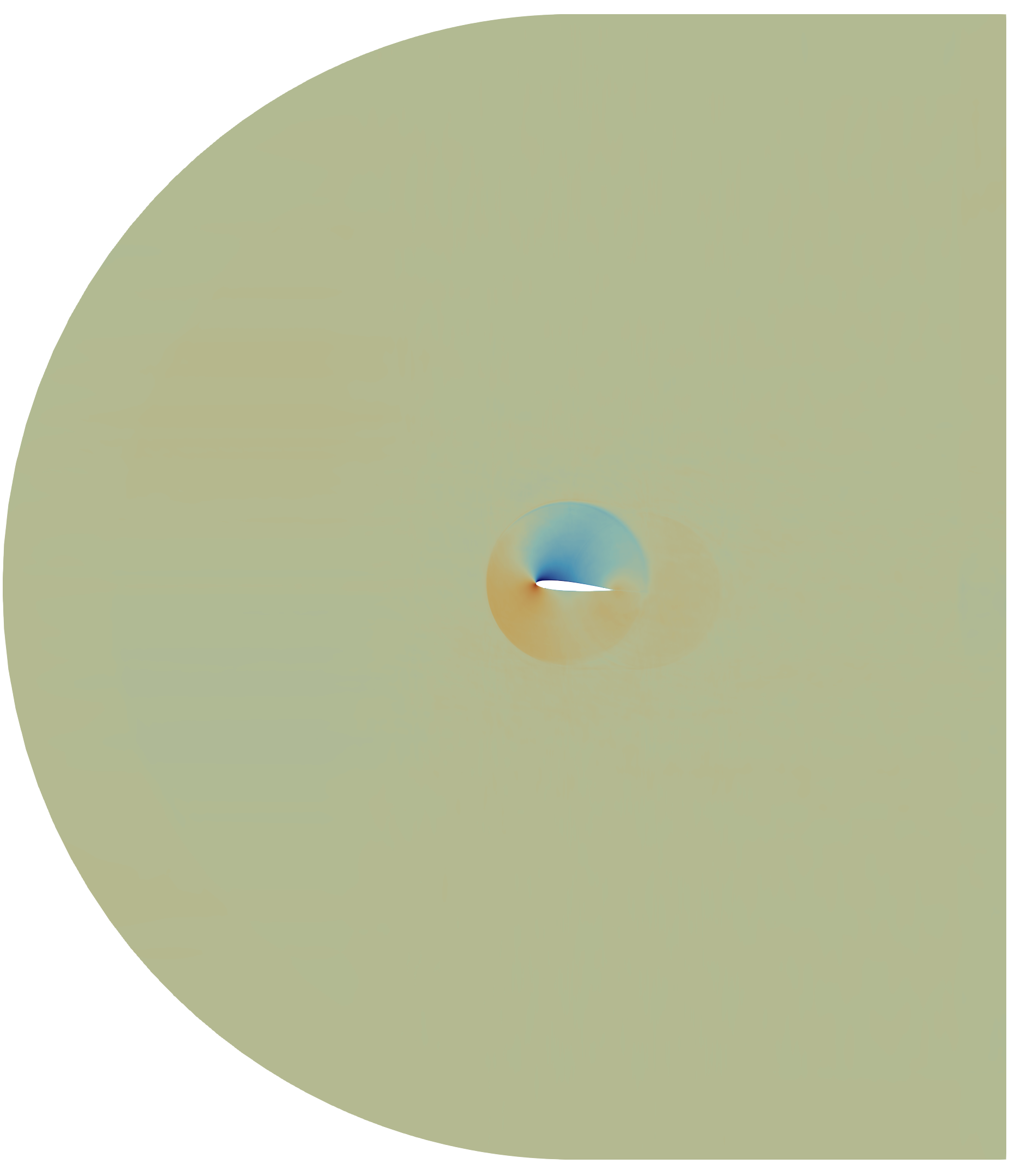} \\[0.8em]

        $t^*=3.25$ &
        \includegraphics[width=0.6\linewidth]{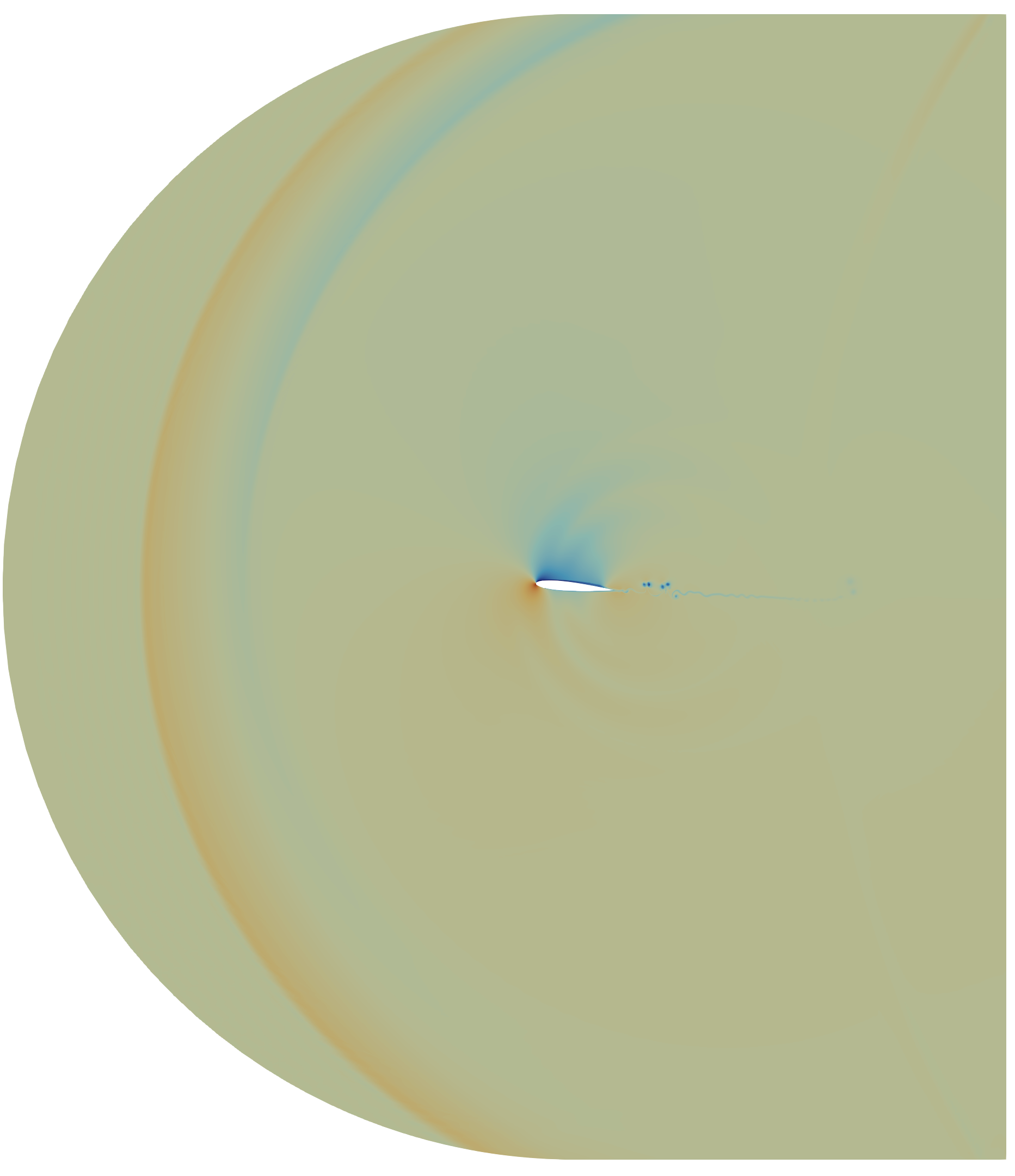} &
        \includegraphics[width=0.6\linewidth]{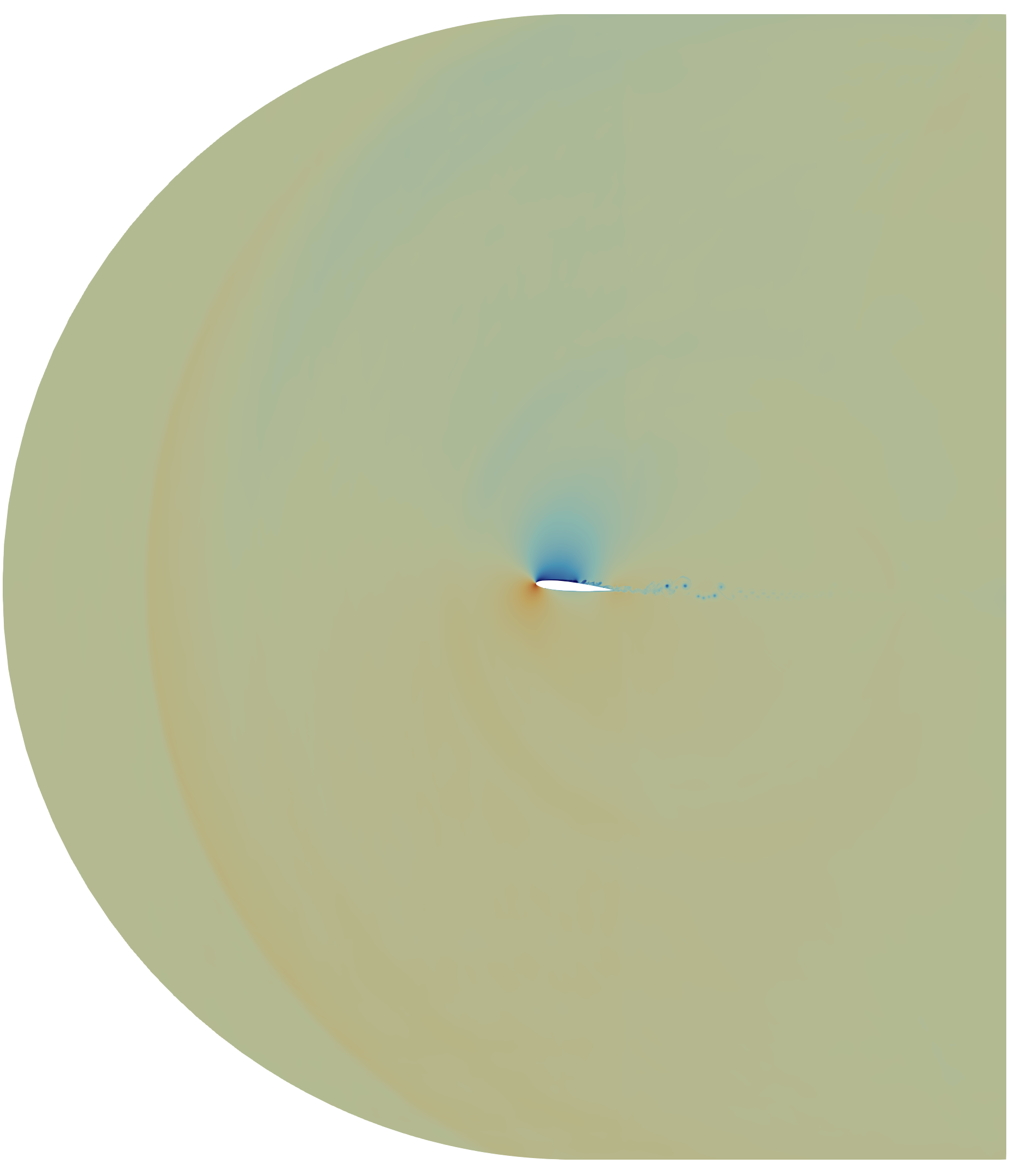} \\[0.8em]

        $t^*=5.1$ &
        \includegraphics[width=0.6\linewidth]{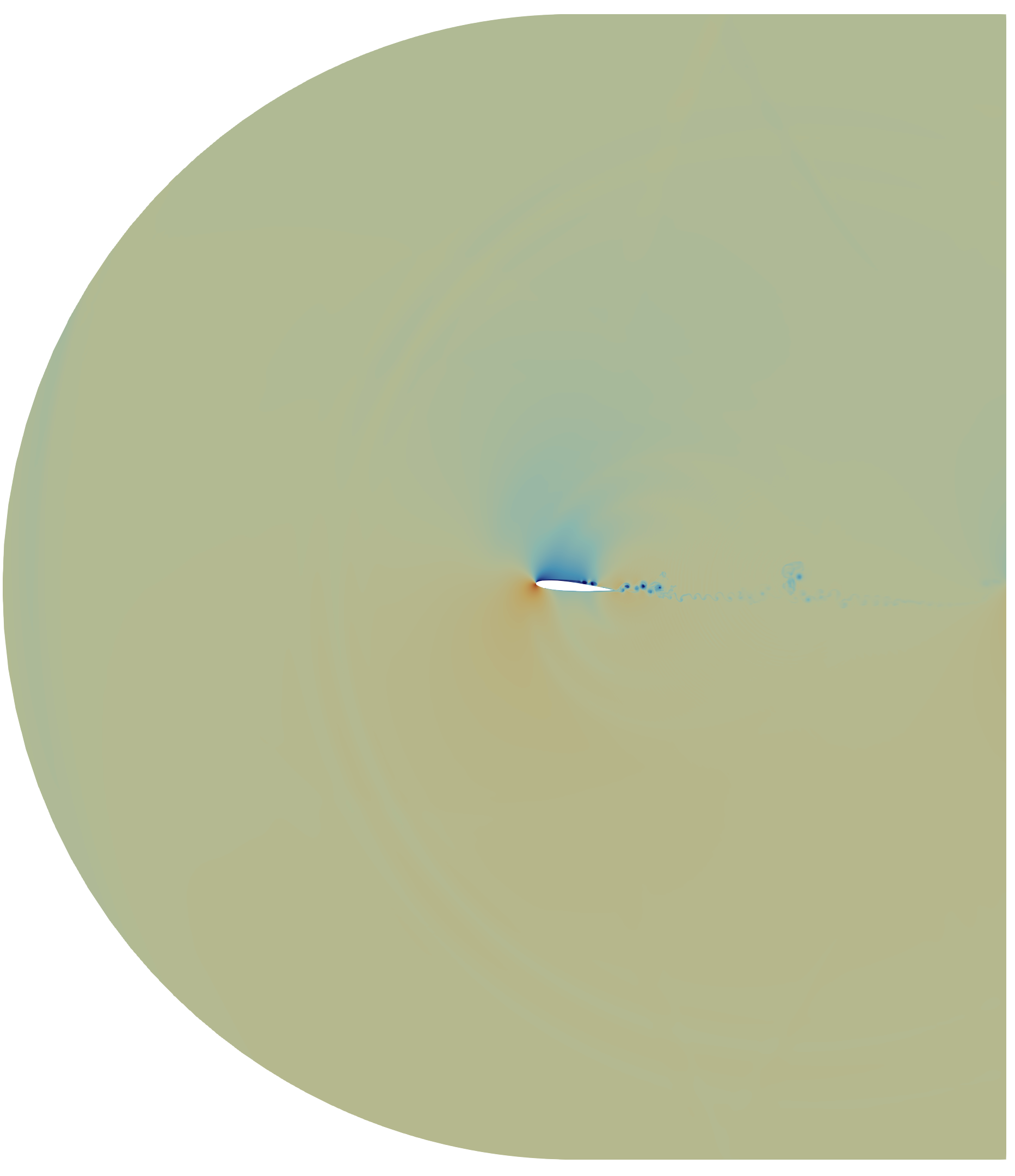} &
        \includegraphics[width=0.6\linewidth]{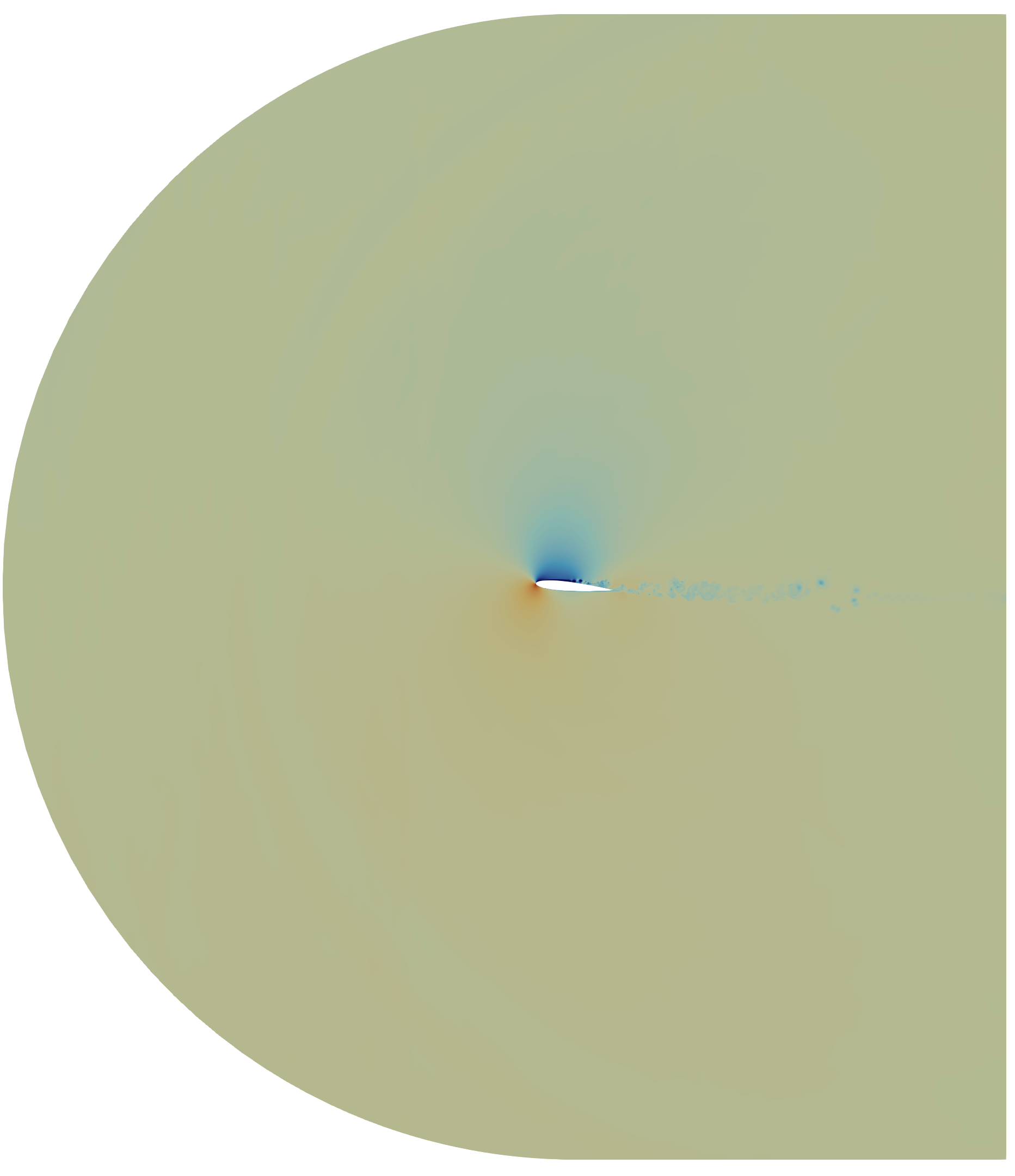} \\[0.8em]

        $t^*=7.5$ &
        \includegraphics[width=0.6\linewidth]{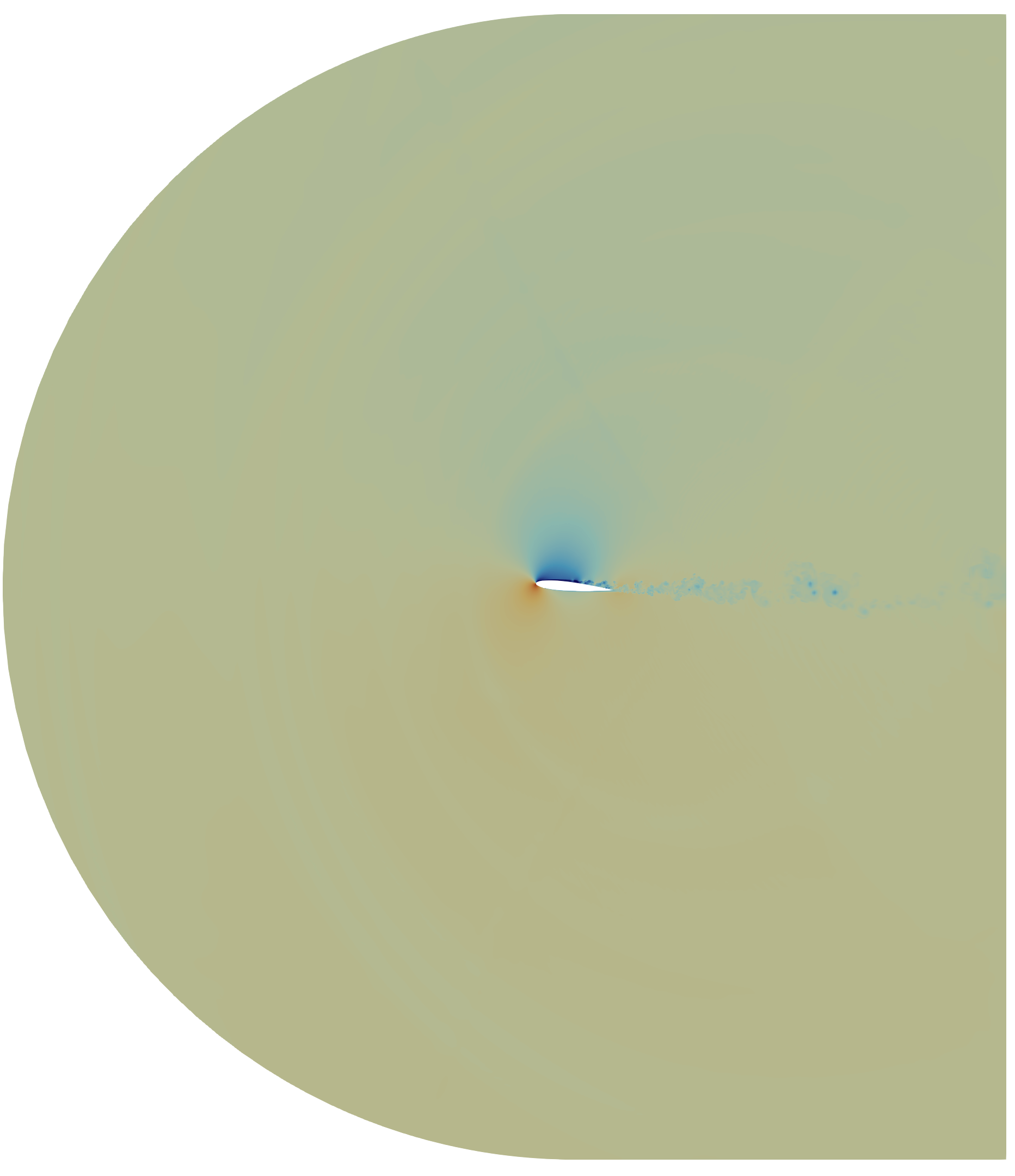} &
        \includegraphics[width=0.6\linewidth]{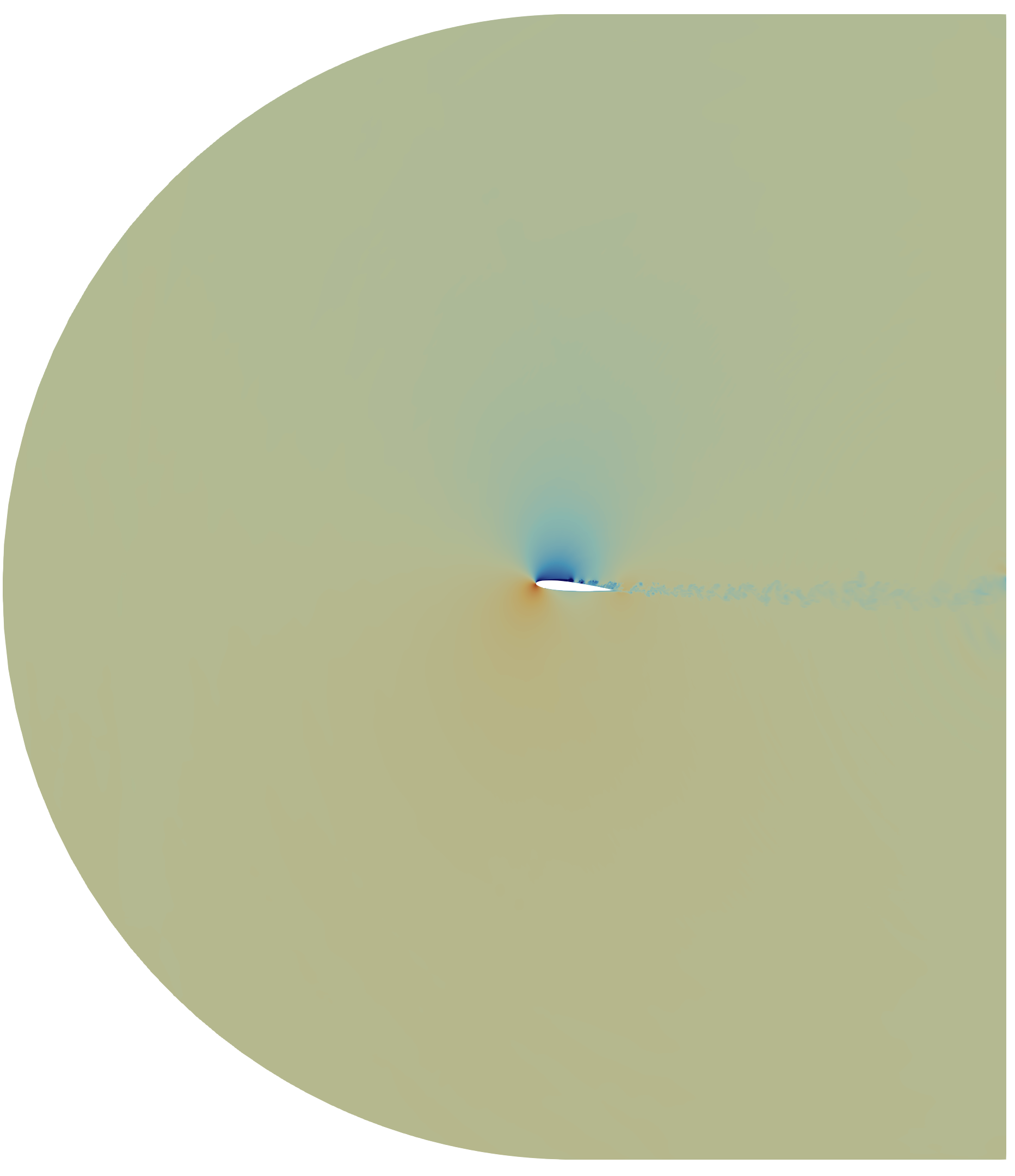} \\[0.8em]

        \multicolumn{3}{@{}c@{}}{
        \includegraphics[width=0.6\textwidth]{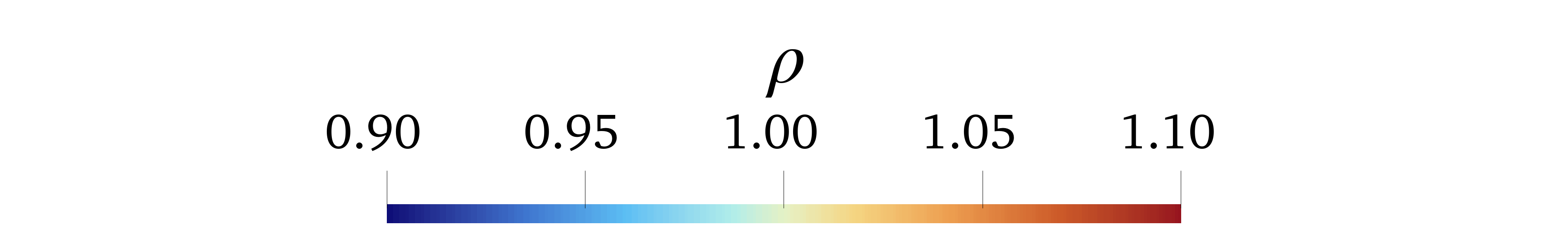}}
    \end{tabular}
    \caption{DNS density fields, $\rho$, for uniform-flow and ADEx-FNO-based initialization at $t^*=0.4$, $3.25$, $5.10$, and $7.5$. The columns compare the two initialization strategies using the displayed common density scale. The sequence illustrates the outward propagation of the startup disturbance and is used as qualitative supporting evidence only.}
    \label{fig:dns_density_frames}
\end{figure}

The density fields in Figure~\ref{fig:dns_density_frames} show the outward propagation of the startup disturbance in both calculations. The images are used only to document this transient and its propagation away from the wing. They do not establish convergence of acoustic spectra or long-time acoustic statistics.

\subsection{Late-time force comparison and limitations}

The common interval $7\leq t^*\leq10$ is selected independently of the MSER bootstrap times and is used to compare the late force histories. The mean values are summarized in Table~\ref{tab:dns_late_force_statistics}. The standard deviations differ by $4.19\%$ for $C_L$ and $5.29\%$ for $C_D$.
The correlation-aware $95\%$ confidence intervals for the differences between
the ADEx-FNO-based and uniform-flow mean forces include zero. Thus, within the
uncertainty of the available records, no statistically distinguishable
difference is detected between their late-time mean force levels.


\begin{table}[H]
\centering 
\caption{Mean DNS force coefficients over the common late-time interval
$7\leq t^*\leq10$.}
\label{tab:dns_late_force_statistics}
\begin{tabular}{lcrrrr}
\toprule
Setup & $t_{avg} (C/U_\infty)$ & $C_L$ & $C_D$ & $C_{D_P}$ & $C_{D_{sf}}$ \\
\midrule
Ref. \cite{karp2026EffectsOfLowerFloatingPointPrecision} & 12 & 0.610 & 0.0355 & 0.0264 & 0.0090 \\
Ref. \cite{jones2008naca0012} & 7.7 & 0.621 & 0.0358 & 0.0220 & 0.0087 \\
\midrule
Uniform-flow initialization & 3 & 0.591 & 0.0359 & 0.0271 & 0.0088 \\
ADEx-FNO initialization & 3 & 0.595 & 0.0346 & 0.0253 & 0.0093 \\
\bottomrule
\end{tabular}
\end{table}


\end{document}